\documentclass[fleqn,10pt]{article}

\usepackage{setspace}
\usepackage[utf8]{inputenc}
\usepackage{enumitem}
\setlist{nolistsep}

\usepackage{amsmath,amssymb,amsthm}
\usepackage[hidelinks]{hyperref}
\usepackage{graphicx}
\usepackage{graphics}
\usepackage{color}
\usepackage{multirow}
\usepackage{subfig}
\usepackage{bm} 
\usepackage{enumerate}
\usepackage{array}

\usepackage{natbib}
\newcommand{\dx}{\ensuremath{{\,\hspace{-0.06em}\rm{d}}\bm{x}}}
\newcommand{\ds}{\ensuremath{{\,\hspace{-0.06em}\rm{d}}s(\bm{x})}}
\newcommand{\dt}{\ensuremath{{\,\hspace{-0.06em}\rm{d}}t}}
\newcommand{\diff}[1]{\ensuremath{{\,\hspace{-0.06em}\rm{d}}{#1}}}

\theoremstyle{plain}

\newtheorem{Theorem}{Theorem}[section]
\newtheorem{Lemma}[Theorem]{Lemma}
\newtheorem{Def}[Theorem]{Definition}
\newtheorem{Remark}[Theorem]{Remark}
{\normalfont}

\title{
  Error Analysis of Physics-Informed Neural Networks for Approximating Dynamic PDEs
  of Second Order in Time
} 
\author{
  Yanxia Qian$^a$, Yongchao Zhang$^b$, Yunqing Huang$^{a*}$,
  Suchuan Dong$^c$\thanks{Authors of correspondence.
    Emails: yxqian0520@xtu.edu.cn (Y.~Qian), yoczhang@nwu.edu.cn (Y.~Zhang),  huangyq@xtu.edu.cn (Y.~Huang), sdong@purdue.edu (S.~Dong)
  }
   \\[0.1in]
   $^a$School of Mathematics and Computational Science, Xiangtan University,\\
   Xiangtan, Hunan, 411105, P.R.~China\\
   $^b$School of Mathematics, Northwest University,\\
   Xi'an, Shaanxi 710069, P.R.~China\\
   $^c$Center for Computational and Applied Mathematics,\\ 
  Department of Mathematics,
  Purdue University, USA 
 } 

\date{(March 21, 2023)}

\begin{document}
\maketitle


\begin{abstract}

We consider the approximation of a class of dynamic partial differential equations (PDE) of second order in time by the physics-informed neural network (PINN) approach, and provide an error analysis of PINN for the wave equation, the Sine-Gordon equation and the linear elastodynamic equation. Our analyses show that, with feed-forward neural networks having two hidden layers and the $\tanh$ activation function, the PINN approximation errors for the solution field, its time derivative and its gradient field can be effectively bounded by the training loss and the number of training data points (quadrature points). Our analyses further suggest new forms for the training loss function, which contain certain residuals that are crucial to the error estimate but would be absent from the canonical PINN loss formulation. Adopting these new forms for the loss function leads to a variant PINN algorithm. We present ample numerical experiments with the new PINN algorithm for the wave equation, the Sine-Gordon equation and the linear elastodynamic equation, which show that the method can capture the solution well.

\end{abstract}


\vspace{0.05cm}
Keywords: {\em
  physics informed neural network; neural network; error estimate; PDE; scientific machine learning
}



	\section{Introduction}\label{sec:intro}


Deep neural networks (DNN) have achieved a great success in a number of fields in science and engineering \cite{LeCun2015DP} such as natural language processing, robotics, computer vision, speech and image recognition, to name but a few. 
This has inspired a great deal of research efforts in the past few years to adapt such techniques to scientific computing.
DNN-based techniques seem particularly promising for problems in higher dimensions, e.g.~high-dimensional partial differential equations (PDE), since
traditional numerical methods for high-dimensional problems can quickly 
become infeasible due to the exponential increase in the computational effort (so-called  curse of dimensionality). 
Under these circumstances deep-learning algorithms can be helpful. In particular, the neural networks approach for PDE problems provide implicit regularization and can alleviate and perhaps overcome the curse of high dimensions  \cite{Beck2019Machine,Berner2020Analysis}. This approach also provides a natural framework for estimating the unknown parameters \cite{Fang2020NN,Raissi2019pinn,Raissi2018Hidden,Thuerey2020Deep,Wang2017pinn}.

As deep neural networks are universal function approximators, it is natural to employ them as ansatz spaces for solutions of (ordinary or partial) differential equations. This paves the way for their  use in physical modeling and scientific computing and gives rise to the field of scientific machine learning \cite{Karniadakisetal2021,SirignanoS2018,Raissi2019pinn,EY2018,Lu2021DeepXDE}. 
The physics-informed neural network (PINN) approach was introduced in \cite{Raissi2019pinn}. 
It has been successfully applied to a variety of forward and inverse PDE problems and has become one of the most commonly-used  methods in scientific machine learning (see e.g.~\cite{Raissi2019pinn,HeX2019,CyrGPPT2020,JagtapKK2020,WangL2020,JagtapK2020,CaiCLL2020,Tartakovskyetal2020,DongN2021,TangWL2021,DongL2021,CalabroFS2021,WanW2022,FabianiCRS2021,KrishnapriyanGZKM2021,DongY2022,DongY2022rm,WangYP2022,Pateletal2022,DongW2022,Siegeletal2022,HuLWX2022,Penwardenetal2023}, among others). 
The references \cite{Karniadakisetal2021,Cuomo2022Scientific} provide a comprehensive review of the literature on PINN and about the  benefits and drawbacks of this approach.

The mathematical foundation for PINN aiming at the approximation of PDE solution is currently an active area of research. It is important to account for different components of the neural-network error: optimization error, approximation error, and estimation error~\cite{Niyogi1999Generalization,Shin2020On}. 
Approximation error refers to the discrepancy between the exact functional map and the neural network mapping function on a given network architecture~\cite{Calin2020Deep,Elbrachter2021deep}. 
Estimation error arises when the network is trained on a finite data set to get a mapping on the target domain. The generalization error is the combination of approximation and estimation errors and defines the accuracy of the neural-network predicted solution trained on the given set of data. 

Theoretical understanding of PINN has been advanced by a number of recent works. 
In~\cite{Shin2020On} Shin et al. rigorously justify why PINN works and shows its consistency for linear elliptic and parabolic PDEs under certain assumptions. 
These results are extended in \cite{Shin2010.08019} to a general abstract framework for analyzing PINN for linear problems with the loss function formulated in terms of the strong or weak forms of the equations. In \cite{Mishra2022Estimates} Mishra and Molinaro provide an abstract framework on PINN for forward PDE problems, and estimate the generalization error by means of the training error and the number of training data points. This framework is extended in~\cite{Mishra2022inverse} to study several inverse PDE problems, including the Poisson, heat, wave and Stokes equations. Bai and Koley \cite{Bai2021PINN} investigate the PINN approximation of nonlinear dispersive PDEs such as the KdV-Kawahara, Camassa-Holm and Benjamin-Ono equations. In~\cite{Biswas2022Error} Biswa et al.~provide explicit error estimates (in suitable norms) and stability analyses for the incompressible Navier–Stokes equations. Zerbinati~\cite{Zerbinati2022pinns} presents PINN as an under-determined point matching collocation method, reveals its connection with Galerkin Least Squares (GALS) method, and establishes an a priori error estimate for elliptic problems. 

An important theoretical result on the approximation errors from the recent work~\cite{DeRyck2021On} establishes  that a feed-forward neural network $\hat{u}_{\theta }$ with a $\tanh$ activation function and two hidden layers may approximate a function $u$ with a bound in a Sobolev space,
\[
\|\hat{u}_{\theta N}-u\|_{w^{k,\infty}}\leq C{\rm ln}(cN)^k/N^{s-k}. 
\]
Here $u \in w^{s,\infty}([0,1]^d)$, $d$ is the dimension of the problem, $N$ is the number of training points, and $c, C > 0$ are explicitly known constants independent of $N$. 
Based on this result,
De Ryck et al.~\cite{2023_IMA_Mishra_NS} have studied the PINN for the Navier–Stokes equations and shown that a small training error implies a small generalization error.
In particular, Hu et al.~\cite{Ruimeng2209.11929} provide the higher-order (spatial Sobolev norm) error estimates for the primitive equations, which improve the existing results in the PINN literature that only involve $L^2$ errors.
In~\cite{DeRyck2022Estimates} it has been shown that,
with a sufficient number of  randomly chosen training points,
the total $L^2$ error can be bounded by the generalization error for Kolmogorov-type PDEs, which in turn is bounded by the training error. It is proved that the size of the PINN and the number of training samples only increase polynomially with the problem dimension, thus enabling PINN to overcome the curse of dimensionality in this case. 
In~\cite{Mishra2021pinn} the authors
investigate the high-dimensional radiative transfer equation and prove that the generalization error is bounded by the training error and the number of training points, where the upper bound depends on the dimension only through 
a logarithmic factor. 
Hence PINN does not suffer from the curse of dimensionality, provided that the training errors do not depend on the underlying dimension. 

Although PINN has been widely used for approximating PDEs, theoretical investigations on its convergence and errors are still quite limited and are largely confined to elliptic and parabolic PDEs. 
There seems to be less (or little) theoretical analysis on the convergence of PINN for hyperbolic type PDEs. 
In this paper, we consider a class of 
dynamic PDEs of second order in time, which are hyperbolic in nature, and provide an analysis of the convergence and errors of the PINN algorithm applied to such problems. 
We have focused on the wave equation, the Sine-Gordon equation and the linear elastodynamic equation in our analyses. 
Building upon the result of~\cite{DeRyck2021On,2023_IMA_Mishra_NS} on $\tanh$ neural networks with two hidden layers, we have shown that for these three kinds of PDEs: 
\begin{itemize}
\item
The underlying PDE residuals in PINN can be made arbitrarily small with $\tanh$ neural networks having two hidden layers. 
\item
The total error of the PINN approximation is bounded by the generalization error of PINN.
\item
The total error of PINN approximations for the solution field, its time derivative and its gradient is bounded by the training error (training loss) of PINN and the number of quadrature points (training data points).
\end{itemize}


Furthermore, our theoretical analyses have suggested PINN training loss functions for these PDEs that are somewhat different in form than from the canonical PINN formulation. 
These lie in two aspects: (i) Our analyses require certain residual terms (such as the gradient of the initial condition, the time derivative of the boundary condition, or in the case of linear elastodynamic equation the strain and divergence of the initial condition) in the training loss, which would be absent from the canonical PINN formulation of the loss function. (ii) Our analyses may require, depending on the type of boundary conditions, a norm other than the $L^2$ norm for certain boundary residuals in the training loss, which is different from the commonly-used $L^2$ norm in the canonical PINN formulation of the loss function.

These new forms for the training loss function suggested by the theoretical analyses lead to a variant PINN algorithm. We have implemented the PINN algorithm based on these new forms of the training loss function for the wave equation, the Sine-Gordon equation and the linear elastodynamic equation. Ample numerical experiments based on this algorithm have been presented. The simulation results indicate that the method has captured the solution field reasonably well for these PDEs. The numerical results also to some extent corroborate the theoretical relation between the approximation error and the PINN training loss obtained from the error analysis.

The rest of this paper is organized as follows. In Section \ref{PINN} we present an overview of PINN for dynamic PDEs of second order in time. In Sections~\ref{Wave}, \ref{Sine-Gordon} and \ref{Elasto-dynamics}, we present an error analysis of the PINN algorithm for approximating the wave equation, Sine-Gordon equation, and the linear elastodynamic equation, respectively. 
Section~\ref{numerical examples} summarizes a set of
numerical experiments with these three PDEs to supplement and support our theoretical analyses. Section~\ref{Conclusion} concludes the presentation with some closing remarks.
Finally, the appendix (Section~\ref{Appendix}) recalls some auxiliary results for our analysis and provides the proofs of the main theorems in Sections~\ref{Sine-Gordon} and~\ref{Elasto-dynamics}.

	\section{Physics Informed Neural Networks (PINN) for Approximating PDEs}\label{PINN}

\subsection{Generic PDE of Second Order in Time}\label{PINN0}

Consider a compact domain $D\subset \mathbb{R}^d$ ($d>0$ being an integer), and let $\mathcal{D}$ and $\mathcal{B}$ denote the differential and boundary operators. We consider the following general form of an initial boundary value problem with a generic PDE of second order in time. For any $\bm{x} \in D$, $\bm{y} \in \partial D$ and $t\in[0,T]$,
\begin{subequations}\label{general}
	\begin{align}\label{general_eq1}
		&\frac{\partial^2u}{\partial t^2}(\bm{x},t)+
  \mathcal{D}[u](\bm{x},t)=0,\\
		\label{general_eq2}
		&\mathcal{B}u(\bm{y},t)=u_{d}(\bm{y},t),\\
		\label{general_eq3}
		&u(\bm{x},0)=u_{in}(\bm{x}),
  \quad \frac{\partial u}{\partial t}(\bm{x},0) = v_{in}(\bm{x}).
	\end{align}
\end{subequations}
Here, $u(\bm x,t)$ is the unknown field solution, $u_{d}$ denotes the boundary data, and $u_{in}$ and $v_{in}$ are the initial distributions for $u$ and 
$\frac{\partial u}{\partial t}$.
We assume that in $\mathcal{D}$ the highest derivative with respect to the time variable $t$, if any,
is of first order.

\subsection{Neural Network Representation of a Function}\label{PINN1}

Let $\sigma: \mathbb{R} \rightarrow \mathbb{R}$ denote an activation function that is at least twice continuously differentiable. 
For any $n\in \mathbb{N}$ and $z\in \mathbb{R}^n$, we define $\sigma(z): = (\sigma(z_1),\cdots,\sigma(z_n))$, where $z_i$ ($1\leq i\leq n$) are the components of $z$. We adopt the following formal definition for a feedforward neural network as given in~\cite{2023_IMA_Mishra_NS}.

\begin{Def}[\cite{2023_IMA_Mishra_NS}]\label{pre_lem1}
Let $R\in (0,\infty]$, $L, W \in \mathbb{N}$ and $l_0, \cdots,l_L \in \mathbb{N}$. Let $\sigma: \mathbb{R} \rightarrow \mathbb{R}$ be a twice differentiable function and define
	\begin{equation}\label{pre_eq1}
		\Theta = \Theta_{L,W,R}:= \bigcup_{L' \in \mathbb{N},L' \leq L} \, \bigcup_{l_0, \cdots,l_L \in \{1,\cdots,W\}} \, {\rlap{$ \diagdown $}\diagup}_{k=1}^{L'}([-R,R]^{l_k\times l_{k-1}}\times[-R,R]^{l_k}).
	\end{equation}
\end{Def}

For $\theta\in \Theta$, we define $\theta_k:=(W_k,b_k)$ and $\mathcal{A}_k^{\theta}: \mathbb{R}^{l_{k-1}}\rightarrow \mathbb{R}^{l_k}$ by $ z\mapsto W_kz+b_k$ for $1\leq k\leq L$, and we define $f_k^{\theta}: \mathbb{R}^{l_{k-1}}\rightarrow \mathbb{R}^{l_k}$ by
\begin{align}\label{pre_eq2}
	f_k^{\theta}=\left\{\begin {array}{lll}
	\mathcal{A}_L^{\theta}(z) & k=L,\\
	(\sigma\circ\mathcal{A}_k^{\theta})(z)& 1\leq k <L.
\end{array}\right.
\end{align}
Denote $u_{\theta}: \mathbb{R}^{l_0}\rightarrow \mathbb{R}^{l_L}$ the function that satisfies for all $z\in \mathbb{R}^{l_0}$ that
\begin{equation}\label{pre_eq3}
u_{\theta}(z)=(f_{L}^{\theta}\circ f_{L-1}^{\theta}\circ\cdots\circ f_{1}^{\theta})(z)\qquad z\in \mathbb{R}^{l_0}.
\end{equation}
We set $z=(\bm{x},t)$ and $l_0=d+1$ for approximating  the PDE problem~\eqref{general}. 

$u_{\theta}$ as defined above is the neural-network representation of a parameterized function  associated with the parameter $\theta$. This neural network contains $(L+1)$ layers ($L\geq 2$), with widths $(l_0,l_1, \cdots,l_L)$ for each layer. The input layer has a width $l_0$, and the output layer has a width $l_L$. The $(L-1)$ layers between the input/output layers are the hidden layers, with widths $l_k$ ($1\leq k\leq L-1$). $W_k$ and $b_k$ are the weight/bias coefficients corresponding to layer $k$ for $1\leq k\leq L$. 
From layer to layer the network logic represents an affine transform, followed by a function composition with the activation function $\sigma$.
Note that no activation function is applied to the output layer.
We refer to $u_{\theta}$ with $L=2$ (i.e.~single hidden layer) as a shallow neural network, and  $u_{\theta}$ with $L\geq 3$ (i.e.~multiple hidden layers) as a deeper or deep neural network.

\subsection{Physics Informed Neural Network for Initial/Boundary Value Problem}\label{PINN2}

Let $\Omega = D\times [0,T]$ and $\Omega_* = \partial D\times [0,T]$ be the spatial-temporal domain
domain. We approximate the solution $u$ to the problem \eqref{general} by a neural network $u_{\theta}: \Omega \rightarrow \mathbb{R}^n$. 
With PINN we consider the residual function of the initial/boundary value problem~\eqref{general}, defined for any sufficiently smooth function $u: \Omega \rightarrow \mathbb{R}^n$ as,
for any $\bm{x} \in D$, $\bm{y} \in \partial D$ and $t\in[0,T]$,
\begin{subequations}\label{general_pinn}
\begin{align}
	\label{general_pinn_eq1}
	&\mathcal{R}_{int}[u](\bm{x},t)=
	\frac{\partial^2u}{\partial t^2}(\bm{x},t)+
	\mathcal{D}[u](\bm{x},t),\\
	\label{general_pinn_eq2}
	&\mathcal{R}_{sb}[u](\bm{y},t)=\mathcal{B}u(\bm{y},t)-u_{d}(\bm{y},t),\\
	\label{general_pinn_eq3}
	&\mathcal{R}_{tb1}[u](\bm{x},0)=u(\bm{x},0)-u_{in}(\bm{x}),\\
	\label{general_pinn_eq4}
	& \mathcal{R}_{tb2}[u](\bm{x},0)=\frac{\partial u}{\partial t}(\bm{x},0)-v_{in}(\bm{x}).
\end{align}
\end{subequations}
These residuals chacracterize how well a given function $u$ satisfies the initial/boundary value problem~\eqref{general}. If $u$ is the exact solution, $\mathcal{R}_{int}[u]=\mathcal{R}_{sb}[u]=\mathcal{R}_{tb1}[u]
=\mathcal{R}_{tb2}[u]=0$.

To facilitate the subsequent analyses,
we introduce an auxiliary function $v=\frac{\partial u}{\partial t}$ and rewrite  $\mathcal{R}_{tb2}$ as 
\begin{equation}\label{general_pinn_eq4.1}
	\mathcal{R}_{tb2}[v](\bm{x},0)=v(\bm{x},0)-v_{in}(\bm{x}).
\end{equation}
We reformulate~\eqref{general_eq1} into two equations, thus separating the interior residual into the following two components:
\begin{align}
		\label{general_pinn_eq1.1}
		&\mathcal{R}_{int1}[u,v](\bm{x},t)=
		\frac{\partial u}{\partial t}(\bm{x},t)-v(\bm{x},t),\\
		\label{general_pinn_eq1.2}
		&\mathcal{R}_{int2}[u,v](\bm{x},t)=
		\frac{\partial v}{\partial t}(\bm{x},t)+
		\mathcal{D}[u](\bm{x},t).
\end{align}
With PINN, we seek a neural network $(u_{\theta},v_{\theta})$ to minimize the following quantity,
\begin{equation}\label{general_0G}
	\begin{split}
		\mathcal{E}_G(\theta)^2=
		&\int_{\Omega}|R_{int1}[u_{\theta}, v_{\theta}](\bm{x},t)|^2\dx
		+\int_{\Omega}|R_{int2}[u_{\theta}, v_{\theta}](\bm{x},t)|^2\dx
		+\int_{D}|R_{tb1}[u_{\theta}](\bm{x})|^2\dx \\
		&+ \int_{D}|R_{tb2}[v_{\theta}](\bm{x})|^2\dx 
		+ \int_{\Omega_*}|R_{sb}[u_{\theta}](\bm{x},t)|^2\ds\dt.
	\end{split}
\end{equation}
The different terms of \eqref{general_0G} may be rescaled by  different weights (penalty coefficients). For simplicity, we set all these weights to one in the analysis. $\mathcal{E}_G$ as defined above is often referred to as the generalization error. 
Because of the integrals involved therein, $\mathcal{E}_G$ can be hard to minimize. In practice, one will approximate  \eqref{general_0G} by an appropriate numerical quadrature rule, as follows
\begin{align}\label{general_0T}
\mathcal{E}_T(\theta,\mathcal{S})^2=&\mathcal{E}_T^{int1}(\theta,\mathcal{S}_{int})^2
	+\mathcal{E}_T^{int2}(\theta,\mathcal{S}_{int})^2
	+\mathcal{E}_T^{tb1}(\theta,\mathcal{S}_{tb})^2 
	+ \mathcal{E}_T^{tb2}(\theta,\mathcal{S}_{tb})^2
	+ \mathcal{E}_T^{sb}(\theta,\mathcal{S}_{sb})^2,
\end{align}
where
\begin{subequations}\label{general_TT}
	\begin{align}
		\label{general_T1}
		\mathcal{E}_T^{int1}(\theta,\mathcal{S}_{int})^2 &= \sum_{n=1}^{N_{int}}\omega_{int}^n|R_{int1}[u_{\theta}, v_{\theta}](\bm{x}_{int}^n,t_{int}^n)|^2,\\
		\label{general_T1.1}
		\mathcal{E}_T^{int2}(\theta,\mathcal{S}_{int})^2 &= \sum_{n=1}^{N_{int}}\omega_{int}^n|R_{int2}[u_{\theta}, v_{\theta}](\bm{x}_{int}^n,t_{int}^n)|^2,\\
		\label{general_T2}
		\mathcal{E}_T^{tb1}(\theta,\mathcal{S}_{tb})^2 &= \sum_{n=1}^{N_{tb}}\omega_{tb}^n|R_{tb1}[u_{\theta}](\bm{x}_{tb}^n)|^2,\\
		\label{general_T2.1}
		\mathcal{E}_T^{tb2}(\theta,\mathcal{S}_{tb})^2 &= \sum_{n=1}^{N_{tb}}\omega_{tb}^n|R_{tb2}[v_{\theta}](\bm{x}_{tb}^n)|^2,\\
		\label{general_T3}
		\mathcal{E}_T^{sb}(\theta,\mathcal{S}_{sb})^2 &= \sum_{n=1}^{N_{sb}}\omega_{sb}^n|R_{sb}[u_{\theta}](\bm{x}_{sb}^n,t_{sb}^n)|^2.
	\end{align}
\end{subequations}
The quadrature points in the spatial-temporal domain and on the spatial and temporal boundaries,  $\mathcal{S}_{int}=\{(\bm{x}_{int}^n,t_{int}^n)\}_{n=1}^{N_{int}}$, $\mathcal{S}_{sb}=\{(\bm{x}_{sb}^n,t_{sb}^n)\}_{n=1}^{N_{sb}}$ and $\mathcal{S}_{tb}=\{(\bm{x}_{tb}^n,t_{tb}^n=0)\}_{n=1}^{N_{tb}}$, constitute the input data sets to the neural network. 
In the above equations $\mathcal{E}_T(\theta,\mathcal{S})^2$ is referred to as the training error (or training loss), and $\omega_{\star}^n$ are suitable quadrature weights for $\star=int$, $sb$ and $tb$. Therefore, PINN attempts to minimize the training error $\mathcal{E}_T(\theta,\mathcal{S})^2$ over the network parameters $\theta$, and upon  convergence of optimization the trained $u_{\theta}$ contains the  approximation of the solution $u$ to the problem \eqref{general}. 

\begin{Remark}
The generalization error \eqref{general_0G} (with the corresponding training error \eqref{general_0T}) is the standard (canonical) PINN form if one introduces $v=\frac{\partial u}{\partial t}$ and reformulates~\eqref{general_eq1} into two equations. 
%
We would like to emphasize that our analyses below suggest  alternative forms for the generalization error, e.g.
\begin{equation}\label{general_G}
	\begin{split}
		\mathcal{E}_G(\theta)^2=
		&\int_{\Omega}|R_{int1}[u_{\theta}, v_{\theta}](\bm{x},t)|^2\dx
		+\int_{\Omega}|R_{int2}[u_{\theta}, v_{\theta}](\bm{x},t)|^2\dx
		+\int_{\Omega}|\nabla R_{int1}[u_{\theta}, v_{\theta}](\bm{x},t)|^2\dx\\
		&+\int_{D}|R_{tb1}[u_{\theta}](\bm{x})|^2\dx
		+ \int_{D}|R_{tb2}[v_{\theta}](\bm{x})|^2\dx
		+\int_{D}|\nabla R_{tb1}[u_{\theta}](\bm{x})|^2\dx\\
		&+ \left(
		\int_{\Omega_*}|R_{sb}[u_{\theta}](\bm{x},t)|^2\ds\dt
		\right)^{\frac{1}{2}},
	\end{split}
\end{equation}
which differs from~\eqref{general_0G} in the terms $\nabla R_{int1}$,  $\nabla R_{tb1}$ and the last term. The corresponding training error is,
\begin{align}\label{general_T}
\mathcal{E}_T(\theta,\mathcal{S})^2=&\mathcal{E}_T^{int1}(\theta,\mathcal{S}_{int})^2
	+\mathcal{E}_T^{int2}(\theta,\mathcal{S}_{int})^2
	+\mathcal{E}_T^{int3}(\theta,\mathcal{S}_{int})^2
	+\mathcal{E}_T^{tb1}(\theta,\mathcal{S}_{tb})^2 
	\nonumber\\
	&+ \mathcal{E}_T^{tb2}(\theta,\mathcal{S}_{tb})^2
	+ \mathcal{E}_T^{tb3}(\theta,\mathcal{S}_{tb})^2
	+ \mathcal{E}_T^{sb}(\theta,\mathcal{S}_{sb}),
\end{align}
where
\begin{equation}\left\{
\begin{split}
&\mathcal{E}_T^{int3}(\theta,\mathcal{S}_{int})^2 = \sum_{n=1}^{N_{int}}\omega_{int}^n|\nabla R_{int1}[u_{\theta}, v_{\theta}](\bm{x}_{int}^n,t_{int}^n)|^2, \\
&
\mathcal{E}_T^{tb3}(\theta,\mathcal{S}_{tb})^2 = \sum_{n=1}^{N_{tb}}\omega_{tb}^n|\nabla R_{tb1}[u_{\theta}](\bm{x}_{tb}^n)|^2.
\end{split}
\right.
\end{equation}
The error analyses also suggest additional terms in the generalization error for different equations.

\end{Remark}

\subsection{Numerical Quadrature Rules}\label{PINN3}

As discussed above, we need to approximate the integrals of functions. The analysis in the subsequent sections requires well-known results on numerical quadrature rules as reviewed below. 

Given $\Lambda \subset \mathbb{R}^d$ and a function $f\in L^1(\Lambda)$, we would like to approximate $\int_{\Lambda}f(z)\diff{z}$. 
A  quadrature rule provides an approximation by
\begin{equation}\label{int}
\int_{\Lambda}f(z)\diff{z}\approx
\frac{1}{M}\sum_{n=1}^M \omega_n f(z_n),
\end{equation}
where $z_n\in \Lambda$ ($1\leq n\leq M$) are the quadrature points and $\omega_n$ ($1\leq n\leq M$) denote the appropriate quadrature weights. 
The approximation accuracy is influenced by the type of quadrature rule, the number of quadrature points ($M$), and the regularity of $f$. For the mid-point rule, which is assumed in the analysis in the current work, the approximation accuracy is given by
\begin{equation}\label{int1}
\left|\int_{\Lambda}f(z)\diff{z}-
\frac{1}{M}\sum_{n=1}^M \omega_n f(z_n)
\right|\leq C_fM^{-2/d},
\end{equation}
where $C_f\lesssim \|f\|_{C^2(D)}$ ($a\lesssim b$ denotes $a\leq Cb$) and $D$ has been partitioned into $M\sim N^d$ cubes 
and $z_n$ ($1\leq n\leq M$) denote the midpoints of these cubes~\cite{DavisR2007}.
In this paper, we use $C$ to denote a universal constant, which may depend on $k, d, T, u$ and $v$ but not on $N$. And we use the subscript to emphasize its dependence when necessary, e.g.~$C_{d}$ is a constant depending only on $d$.

We focus on PDE problems in relatively low dimensions ($d\leq 3$) in this paper and employ the standard quadrature rules. 
We note that in higher dimensions the standard quadrature rules may not be favorable. In this case the random training points or low-discrepancy training points 
\cite{Mishra2021Enhancing} may be preferred. 


\vspace{0.2in} 
In subsequent sections we focus on three representative dynamic equations of second order in time (the wave equation, Sine-Gordon equation, and the linear elastodynamic equation), and provide the error estimate for approximating these equations by PINN. We note that these  analyses suggest alternative forms for the training loss function that are somewhat different from the standard PINN forms~\cite{Raissi2019pinn}. The PINN numerical results based on the standard  form for the loss function, and based on the alternative forms as suggested by
the error estimate, will be provided after the presentation of the theoretical analysis. In what follows,  for brevity we adopt the notation of $\mathcal{F}_\Xi=\frac{\partial \mathcal{F}}{\partial \Xi}$, $\mathcal{F}_{\Xi\Upsilon}=\frac{\partial^2 \mathcal{F}}{\partial\Xi \partial \Upsilon}$ ($\Xi,\Upsilon\in\{t, x\}$),
for any sufficiently smooth function $\mathcal{F}: \Omega\rightarrow \mathbb R^n$.

	\section{Physics Informed Neural Networks for Approximating Wave Equation}\label{Wave}

\subsection{Wave Equation}

Consider the following wave equations on the torus $D=[0,1)^d \subset \mathbb{R}^d$ with periodic boundary conditions:
\begin{subequations}\label{wave}
	\begin{align}
		\label{wave_eq0}
		&u_{t} - v = 0  \ \, \qquad\qquad\qquad \text{in}\ D\times[0,T],\\
		\label{wave_eq1}
		&v_{t} -\Delta u = f   \  \ \quad\qquad\qquad \text{in}\ D\times[0,T],\\
		\label{wave_eq2}
		&u(\bm{x},0)=\psi_{1}(\bm{x})\qquad\qquad \ \text{in}\ D,\\
		\label{wave_eq3}
		&v(\bm{x},0)=\psi_{2}(\bm{x}) \,\qquad\qquad \ \text{in}\ D,\\
		\label{wave_eq4}
		&u(\bm{x},t)=u(\bm{x}+1,t)  \qquad\ \ \text{in}\ \partial D\times[0,T],\\
		\label{wave_eq5}
		&\nabla u(\bm{x},t)=\nabla u(\bm{x}+1,t)  \quad \text{in}\ \partial D\times[0,T].
	\end{align}
\end{subequations}

The regularity results for linear evolution equations of the second order in time have been studied in the Book \cite{Temam1997Infinite}. 
When the self-adjoint operator $\mathcal{A}$ takes $\Delta$, the linear evolution equations of second order in time become the classical wave equations, and then we can also obtain the following regularity results.

\begin{Lemma}\label{sec5_Lemma1} Let $r\geq1$, $\psi_{1}\in H^{r}(D)$, $\psi_{2}\in H^{r-1}(D)$ and $f\in L^{2}([0,T];H^{r-1}(D))$, then there exists a unique solution $u$ to the classical wave equations such that $u \in C([0,T];H^{r}(D))$ and $u_t \in C([0,T];H^{r-1}(D))$.
\end{Lemma}

\begin{Lemma}\label{sec5_Lemma2} Let $k\in\mathbb N$, $\psi_{1}\in H^{r}(D)$, $\psi_{2}\in H^{r-1}(D)$ and $f\in C^{k-1}([0,T];H^{r-1}(D))$ with $r>\frac{d}{2}+k$, then there exists $T>0$ and a classical solution $u$ to the wave equations such that $u(t = 0) = \psi_1$, $u_t(t = 0) = \psi_2$, $u \in C^k(D\times[0,T])$ and $v \in C^{k-1}(D\times[0,T])$.
\end{Lemma}
\begin{proof}  By Lemma \ref{sec5_Lemma1}, there exists $T > 0$ and the solution $(u, v)$ to the wave equations such that $u(t = 0) = \psi_1$, $v(t = 0) = \psi_2$, $u \in C([0,T];H^{r}(D))$ and $v \in C([0,T];H^{r-1}(D))$. As $r>\frac{d}{2}+k$, $H^{r-k}(D)$ is a Banach algebra.
	
	For $k = 1$, since $u \in C([0,T];H^{r}(D))$, $v \in C([0,T];H^{r-1}(D))$ and $f\in C([0,T];H^{r-1}(D))$, we have $u_t = v\in C([0,T];H^{r-1}(D))$ and $v_t = \Delta u +f\in C([0,T];H^{r-2}(D))$. Then, it implies that $u\in C^1([0,T];H^{r-1}(D))$ and $v\in C^1([0,T];H^{r-2}(D))$.
	
	For $k = 2$, by $f\in C^{1}([0,T];H^{r-1}(D))$, we have $u_{tt} = v_t \in C([0,T];H^{r-2}(D))$ and $v_{tt} = \Delta u_t +f_t \in C([0,T];H^{r-3}(D))$. Then, it implies that $u\in C^2([0,T];H^{r-2}(D))$ and $v\in C^2([0,T];H^{r-3}(D))$.
	
	Repeating the same argument, we have $u \in \cap_{l=0}^k C^l([0,T];H^{r-l}(D))$ and $v \subset \cap_{l=0}^k C^l([0,T];H^{r-l-1}(D))$. Then, applying the Sobolev embedding theorem and $r>\frac{d}{2}+k$, it holds $H^{r-l}(D) \subset C^{r-l}(D)$ and $H^{r-l-1}(D) \subset C^{r-l-1}(D)$ for $0\leq l\leq k$. Therefore, $u \in C^k(D\times[0,T])$ and $v \in C^{k-1}(D\times[0,T])$.
\end{proof} 

\subsection{Physics Informed Neural Networks}

We would like to approximate the solutions to  the problem~\eqref{wave} with PINN.
%
We seek deep neural networks $u_{\theta} : D\times [0,T] \rightarrow \mathbb{R}$ and $v_{\theta} : D\times [0,T] \rightarrow \mathbb{R}$, parameterized by $\theta \in \Theta$, 
that approximate the solution $u$ and $v$ of \eqref{wave}. Define residuals,
\begin{subequations}\label{wave_pinn}
	\begin{align}
		\label{wave_pinn_eq1}
		&R_{int1}[u_{\theta},v_{\theta}](\bm{x},t) =u_{\theta t}-v_{\theta},\\
		\label{wave_pinn_eq2}
		&R_{int2}[u_{\theta},v_{\theta}](\bm{x},t) =v_{\theta t}-\Delta u_{\theta} - f,\\
		\label{wave_pinn_eq3}
		&R_{tb1}[u_{\theta}](\bm{x}) =u_{\theta}(\bm{x},0)-\psi_{1}(\bm{x}),\\
		\label{wave_pinn_eq4}
		&R_{tb2}[v_{\theta}](\bm{x}) =v_{\theta}(\bm{x},0)-\psi_{2}(\bm{x}),\\
		\label{wave_pinn_eq5}
		&R_{sb1}[v_{\theta}](\bm{x},t) =v_{\theta}(\bm{x},t)-v_{\theta}(\bm{x}+1,t),\\
		\label{wave_pinn_eq6}
		&R_{sb2}[u_{\theta}](\bm{x},t) =\nabla u_{\theta}(\bm{x},t)-\nabla u_{\theta}(\bm{x}+1,t).
	\end{align}
\end{subequations}
Note that for the exact solution $R_{int1}[u,v]=R_{int2}[u,v]=R_{tb1}[u]=R_{tb2}[v]=R_{sb1}[v]=R_{sb2}[u]=0$. 
Let $\Omega = D\times [0,T]$ and $\Omega_* = \partial D\times [0,T]$ be the space-time domain. 
With PINN, we minimize the the following generalization error,
\begin{align}\label{wave_G}
\mathcal{E}_G(\theta)^2&=\int_{\Omega}|R_{int1}[u_{\theta},v_{\theta}](\bm{x},t)|^2\dx+\int_{\Omega}|R_{int2}[u_{\theta},v_{\theta}](\bm{x},t)|^2\dx+\int_{\Omega}|\nabla R_{int1}[u_{\theta},v_{\theta}](\bm{x},t)|^2\dx
	\nonumber\\
	&+\int_{D}|R_{tb1}[u_{\theta}](\bm{x})|^2\dx+\int_{D}|R_{tb2}[v_{\theta}](\bm{x})|^2\dx
	+\int_{D}|\nabla R_{tb1}[u_{\theta}](\bm{x})|^2\dx
	\nonumber\\
	&
	+\int_{\Omega_*}|R_{sb1}[v_{\theta}](\bm{x},t)|^2\ds\dt +\int_{\Omega_*}|R_{sb2}[u_{\theta}](\bm{x},t)|^2\ds\dt.
\end{align}
The form of different terms in this expression will become clearer below.

To complete the PINN formulation, we will choose the training set $\mathcal{S} \subset \overline{D}\times [0,T]$ based on suitable quadrature points. We divide the full training set $\mathcal{S} = \mathcal{S}_{int} \cup \mathcal{S}_{sb} \cup \mathcal{S}_{tb}$ into the following three components: 
\begin{itemize}
    \item Interior training points $\mathcal{S}_{int}=\{{z}_n\}$ for $1\leq n \leq N_{int}$, with each ${z}_n= (\bm{x},t)_n \in D \times(0,T)$.
    
    \item Spatial boundary training points $\mathcal{S}_{sb}=\{{z}_n\}$ for $1\leq n \leq N_{sb}$, with each ${z}_n= (\bm{x},t)_n \in \partial D\times (0,T)$.
    
    \item Temporal boundary training points $\mathcal{S}_{tb}=\{\bm{x}_n\}$ for $1\leq n \leq N_{tb}$ with  each $\bm{x}_n \in D$.
\end{itemize}
%
We define the PINN training loss, $\theta \mapsto \mathcal{E}_T(\theta,\mathcal{S})^2$, as follows,
\begin{align}
	\label{wave_T}
	\mathcal{E}_T(\theta,\mathcal{S})^2&
	=\mathcal{E}_T^{int1}(\theta,\mathcal{S}_{int})^2+\mathcal{E}_T^{int2}(\theta,\mathcal{S}_{int})^2+\mathcal{E}_T^{int3}(\theta,\mathcal{S}_{int})^2	+\mathcal{E}_T^{tb1}(\theta,\mathcal{S}_{tb})^2\nonumber\\
	& +\mathcal{E}_T^{tb2}(\theta,\mathcal{S}_{tb})^2 +\mathcal{E}_T^{tb3}(\theta,\mathcal{S}_{tb})^2 
	+\mathcal{E}_T^{sb1}(\theta,\mathcal{S}_{sb})^2+\mathcal{E}_T^{sb2}(\theta,\mathcal{S}_{sb})^2,
\end{align}
where
\begin{subequations}\label{wave_G_add1}
	\begin{align}
		\label{wave_T1}
		\mathcal{E}_T^{int1}(\theta,\mathcal{S}_{int})^2 &= \sum_{n=1}^{N_{int}}\omega_{int}^n|R_{int1}[u_{\theta},v_{\theta}](\bm{x}_{int}^n,t_{int}^n)|^2,\\
		\label{wave_T01}
		\mathcal{E}_T^{int2}(\theta,\mathcal{S}_{int})^2 &= \sum_{n=1}^{N_{int}}\omega_{int}^n|R_{int2}[u_{\theta},v_{\theta}]](\bm{x}_{int}^n,t_{int}^n)|^2,\\
		\label{wave_T001}
		\mathcal{E}_T^{int3}(\theta,\mathcal{S}_{int})^2 &= \sum_{n=1}^{N_{int}}\omega_{int}^n|\nabla R_{int1}[u_{\theta},v_{\theta}](\bm{x}_{int}^n,t_{int}^n)|^2,\\
		\label{wave_T2}
		\mathcal{E}_T^{tb1}(\theta,\mathcal{S}_{tb})^2 &= \sum_{n=1}^{N_{tb}}\omega_{tb}^n|R_{tb1}[u_{\theta}](\bm{x}_{tb}^n)|^2,\\
		\label{wave_T02}
		\mathcal{E}_T^{tb2}(\theta,\mathcal{S}_{tb})^2 &= \sum_{n=1}^{N_{tb}}\omega_{tb}^n|R_{tb2}[v_{\theta}](\bm{x}_{tb}^n)|^2,\\
		\label{wave_T002}
		\mathcal{E}_T^{tb3}(\theta,\mathcal{S}_{tb})^2 &= \sum_{n=1}^{N_{tb}}\omega_{tb}^n|\nabla R_{tb1}[u_{\theta}](\bm{x}_{tb}^n)|^2,\\
		\label{wave_T3}
		\mathcal{E}_T^{sb1}(\theta,\mathcal{S}_{sb})^2 &= \sum_{n=1}^{N_{sb}}\omega_{sb}^n|R_{sb1}[v_{\theta}](\bm{x}_{sb}^n,t_{sb}^n)|^2,\\
		\label{wave_T03}
		\mathcal{E}_T^{sb2}(\theta,\mathcal{S}_{sb})^2 &= \sum_{n=1}^{N_{sb}}\omega_{sb}^n|R_{sb2}[u_{\theta}](\bm{x}_{sb}^n,t_{sb}^n)|^2.
	\end{align}
\end{subequations}
Here the quadrature points in space-time constitute the data sets $\mathcal{S}_{int} = \{(\bm{x}_{int}^n,t_{int}^n)\}_{n=1}^{N_{int}}$, $\mathcal{S}_{tb} = \{\bm{x}_{tb}^n)\}_{n=1}^{N_{tb}}$ and $\mathcal{S}_{sb} = \{(\bm{x}_{sb}^n,t_{sb}^n)\}_{n=1}^{N_{sb}}$, and $\omega_{\star}^n$ are suitable quadrature weights with $\star$ denoting $int$, $tb$ or $sb$. 

Let 
\[
\hat{u} = u_{\theta}-u, \qquad \hat{v} = v_{\theta}-v,	
\]
denote the difference between the solution to the wave equations and the PINN approximation of the solution.
We define the total error of the PINN approximation by
\begin{equation}\label{wave_total}
	\mathcal{E}(\theta)^2=\int_0^{T}\int_{D}(|\hat{u}(\bm{x},t)|^2+|\nabla\hat{u}(\bm{x},t)|^2+|\hat{v}(\bm{x},t)|^2)\dx\dt.
\end{equation}

\subsection{Error Analysis}  

In light of the wave equations \eqref{wave} and the definitions for  different residuals \eqref{wave_pinn}, we have
\begin{subequations}\label{wave_error}
	\begin{align}
		\label{wave_error_eq1}
		&R_{int1}=\hat{u}_t-\hat{v},\\
		\label{wave_error_eq2}
		&R_{int2}=\hat{v}_t-\Delta \hat{u}\\
		\label{wave_error_eq3}
		&R_{tb1}=\hat{u}(\bm{x},0),\\
		\label{wave_error_eq4}
		&R_{tb2}=\hat{v}(\bm{x},0),\\
		\label{wave_error_eq5}
		&R_{sb1}=\hat{v}(\bm{x},t)-\hat{v}(\bm{x}+1,t),\\
		\label{wave_error_eq6}
		&R_{sb2}=\nabla\hat{u}(\bm{x},t)-\nabla\hat{u}(\bm{x}+1,t).
	\end{align}
\end{subequations}

\subsubsection{Bound on the Residuals}  
\begin{Theorem}\label{sec5_Theorem1} 
	Let $d$, $r$, $k \in \mathbb{N}$ with $k\geq 3$. Let $\psi_1 \in H^{r}(D)$, $\psi_2 \in H^{r-1}(D)$ and $f\in C^{k-1}([0,T];H^{r-1}(D))$ with $r>\frac{d}{2}+k$. For every integer $N>5$, there exist $\tanh$ neural networks $u_{\theta}$ and $v_{\theta}$, each with two hidden layers, of widths at most $3\lceil\frac{k}{2}\rceil|P_{k-1,d+2}| + \lceil NT\rceil+ d(N-1)$ and $3\lceil\frac{d+3}{2}\rceil|P_{d+2,d+2}| \lceil NT\rceil N^d$, such that
	\begin{subequations}
		\begin{align}
			\label{lem5.1}
			&\|R_{int1}\|_{L^2(\Omega)},\|R_{tb1}\|_{L^2(D)}\lesssim {\rm ln}NN^{-k+1},\\
			\label{lem5.2}
			&\|R_{int2}\|_{L^2(\Omega)},\|\nabla R_{int1}\|_{L^2(\Omega)}, \|\nabla R_{tb1}\|_{L^2(D)}, \|R_{sb2}\|_{L^2(\partial D\times [0,t])}\lesssim {\rm ln}^2NN^{-k+2},\\
			\label{lem5.3}
			&\|R_{tb2}\|_{L^2(D)},\|R_{sb1}\|_{L^2(\partial D\times [0,t])}\lesssim {\rm ln}NN^{-k+2}.
		\end{align}
	\end{subequations}
\end{Theorem}
\begin{proof} Based on Lemma \ref{sec5_Lemma2}, it holds that $u\in H^k(\Omega)$ and $v\in H^{k-1}(\Omega)$. In light of Lemma \ref{Ar_4}, there  exists neural networks $u_{\theta}$ and $v_{\theta}$, with the same two hidden layers and widths $3\lceil\frac{k}{2}\rceil|P_{k-1,d+2}| + \lceil NT\rceil+ d(N-1)$ and $3\lceil\frac{d+3}{2}\rceil|P_{d+2,d+2}| \lceil NT\rceil N^d$, such that for every $0 \leq l \leq 2$ and $0 \leq s \leq 2$,
\begin{align}
\label{lem5.1_eq1}
&\|u_{\theta}-u\|_{H^l(\Omega)}\leq C_{l,k,d+1,u}\lambda_{l,u}(N)N^{-k+l},\\
\label{lem5.1_eq2}
&\|v_{\theta}-v\|_{H^s(\Omega)}\leq C_{s,k-1,d+1,v}\lambda_{s,v}(N)N^{-k+1+s},
\end{align}
where $\lambda_{l,u} = 2^l3^{d+1}(1+\sigma){\rm ln}^l\left(\beta_{l,\sigma,d+1,u}N^{d+k+3}\right)$, $\sigma = \frac{1}{100}$, $\lambda_{s,v}=2^s3^{d+1}(1+\sigma){\rm ln}^s\left(\beta_{s,\sigma,d+1,v}N^{d+k+2}\right)$, and the definition for the other constants can be found in Lemma \ref{Ar_4}. 

In light of Lemma \ref{Ar_2}, we can bound the PINN residual terms,
	\begin{align*}
&\|\hat{u}_t\|_{L^2(\Omega)}\leq\|\hat{u}\|_{H^1(\Omega)},\qquad \|\hat{v}_t\|_{L^2(\Omega)}\leq\|\hat{v}\|_{H^1(\Omega)},\\
		&\|\Delta\hat{u}\|_{L^2(\Omega)}\leq\|\hat{u}\|_{H^2(\Omega)},\qquad \|\nabla\hat{u}_t\|_{L^2(\Omega)}\leq\|\hat{u}\|_{H^2(\Omega)},\\
		&\|\nabla\hat{v}\|_{L^2(\Omega)}\leq\|\hat{v}\|_{H^1(\Omega)},\\
		&\|\hat{u}\|_{L^2(D)}\leq \|\hat{u}\|_{L^2(\partial\Omega)}\leq C_{h_{\Omega},d+1,\rho_{\Omega}}\|\hat{u}\|_{H^1(\Omega)},\\
		&\|\hat{v}\|_{L^2(D)}\leq \|\hat{v}\|_{L^2(\partial\Omega)}\leq C_{h_{\Omega},d+1,\rho_{\Omega}}\|\hat{v}\|_{H^1(\Omega)},\\
		&\|\nabla\hat{u}\|_{L^2(D)}\leq \|\nabla\hat{u}\|_{L^2(\partial\Omega)}\leq C_{h_{\Omega},d+1,\rho_{\Omega}}\|\hat{u}\|_{H^2(\Omega)},\\
		&\|\hat{v}\|_{L^2(\partial D\times [0,t])}\leq \|\hat{v}\|_{L^2(\partial\Omega)}\leq C_{h_{\Omega},d+1,\rho_{\Omega}}\|\hat{v}\|_{H^1(\Omega)},\\
		&\|\nabla\hat{u}\|_{L^2(\partial D\times [0,t])}\leq \|\nabla\hat{u}\|_{L^2(\partial\Omega)}\leq C_{h_{\Omega},d+1,\rho_{\Omega}}\|\hat{u}\|_{H^2(\Omega)}.
	\end{align*}
	By combining these relations  with \eqref{lem5.1_eq1} and \eqref{lem5.1_eq2}, we can obtain
	\begin{align*}
&\|R_{int1}\|_{L^2(\Omega)}=\|\hat{u}_t-\hat{v}\|_{L^2(\Omega)}\leq
		\|\hat{u}\|_{H^1(\Omega)}+\|\hat{v}\|_{L^2(\Omega)}\\
		&\qquad\leq C_{1,k,d+1,u}\lambda_{1,u}(N)N^{-k+1} + C_{0,k-1,d+1,v}\lambda_{0,v}(N)N^{-k+1}\lesssim {\rm ln}NN^{-k+1},\\
		&\|R_{int2}\|_{L^2(\Omega)}=\|\hat{v}_t-\Delta\hat{u}\|_{L^2(\Omega)}\leq
		\|\hat{v}\|_{H^1(\Omega)}+\|\hat{u}\|_{H^2(\Omega)}\\
		&\qquad\leq C_{2,k,d+1,u}\lambda_{2,u}(N)N^{-k+2} + C_{1,k-1,d+1,v}\lambda_{1,v}(N)N^{-k+2}\lesssim {\rm ln}^2NN^{-k+2},\\		
		&\|\nabla R_{int1}\|_{L^2(\Omega)}=\|\nabla(\hat{u}_t-\hat{v})\|_{L^2(\Omega)}\leq
		\|\hat{u}\|_{H^2(\Omega)}+\|\hat{v}\|_{H^1(\Omega)}\\
		&\qquad\leq C_{2,k,d+1,u}\lambda_{2,u}(N)N^{-k+2} + C_{1,k-1,d+1,v}\lambda_{1,v}(N)N^{-k+2}\lesssim {\rm ln}^2NN^{-k+2},\\
		&\|R_{tb1}\|_{L^2(D)}\leq C_{h_{\Omega},d+1,\rho_{\Omega}}\|\hat{u}\|_{H^1(\Omega)}\lesssim {\rm ln}NN^{-k+1},\\	
		&\|R_{tb2}\|_{L^2(D)}, \|R_{sb1}\|_{L^2(\partial D\times [0,t])}\leq C_{h_{\Omega},d+1,\rho_{\Omega}}\|\hat{v}\|_{H^1(\Omega)}\lesssim {\rm ln}NN^{-k+2},\\	
		&\|\nabla R_{tb1}\|_{L^2(D)}, \|R_{sb2}\|_{L^2(\partial D\times [0,t])}\leq C_{h_{\Omega},d+1,\rho_{\Omega}}\|\hat{u}\|_{H^2(\Omega)}\lesssim {\rm ln}^2NN^{-k+2}.
	\end{align*}
\end{proof}

Theorem \ref{sec5_Theorem1} implies that one can make the PINN residuals~\eqref{wave_pinn} arbitrarily small by choosing $N$ to be sufficiently large. It follows that the generalization error $\mathcal{E}_G(\theta)^2$ in~\eqref{wave_G} can be made arbitrarily small.  

\subsubsection{Bounds on the Total Approximation Error} 

We next show that the total error $\mathcal{E}(\theta)^2$ is also small when the generalization error $\mathcal{E}_G(\theta)^2$ is small with the PINN approximation $(u_{\theta},v_{\theta})$. 
Then we prove that the total error $\mathcal{E}(\theta)^2$ can be arbitrarily small, provided that the training error $\mathcal{E}_T(\theta,\mathcal{S})^2$ is sufficiently small  and the sample set is sufficiently large. 

\begin{Theorem}\label{sec5_Theorem2} Let $d\in \mathbb{N}$, $u\in C^1(\Omega)$ and $v\in C^0(\Omega)$ be the classical solution to the wave equations \eqref{wave}. Let $u_{\theta}$ and $v_{\theta}$ denote the PINN approximation with parameter $\theta$. Then the following relation holds,
	\begin{equation}\label{lem5.4}
			\mathcal{E}(\theta)^2=\int_0^{T}\int_{D}(|\hat{u}(\bm{x},t)|^2+|\nabla\hat{u}(\bm{x},\tau)|^2+|\hat{v}(\bm{x},t)|^2)\dx\dt\leq C_GT\exp(2T),
	\end{equation}
	where 
	\begin{align*}
		C_G&=\int_{D}(|R_{tb1}|^2+|R_{tb2}|^2+|\nabla R_{tb1}|^2)\dx +  \int_{0}^{T}\int_{D}(|R_{int1}|^2+|R_{int2}|^2+|\nabla R_{int1}|^2)\dx\dt
		\\
		&\quad+
		\int_{0}^{T}\int_{\partial D}(|R_{sb1}|^2+|R_{sb2}|^2)\ds\dt.
	\end{align*}
\end{Theorem}
\begin{proof} By taking the inner product of \eqref{wave_error_eq1} and \eqref{wave_error_eq2} with $\hat{u}$ and $\hat{v}$ and integrating over $D$, respectively, we have
	\begin{align}
		\label{sec5_eq0}
		\frac{d}{2dt}\int_{D} |\hat{u}|^2\dx &= \int_{D}\hat{u}\hat{v}\dx+\int_{D} R_{int1}\hat{u}\dx\leq \int_{D} |\hat{u}|^2\dx+\frac{1}{2}\int_{D} |R_{int1}|^2\dx+\frac{1}{2}\int_{D} |\hat{v}|^2\dx,\\
		\label{sec5_eq1}
		\frac{d}{2dt}\int_{D} |\hat{v}|^2\dx &=- \int_{D}\nabla\hat{u}\cdot\nabla\hat{v}\dx+\int_{\partial D} \hat{v}\nabla\hat{u}\cdot\bm{n}\ds+\int_{D} R_{int2}\hat{v}\dx
		\nonumber\\
		&=- \int_{D}\nabla\hat{u}\cdot\nabla\hat{u}_t\dx
		+\int_{D}\nabla\hat{u}\cdot\nabla R_{int1}\dx
		+\int_{\partial D} \hat{v}\nabla\hat{u}\cdot\bm{n}\ds
		+\int_{D} R_{int2}\hat{v}\dx
		\nonumber\\
		&=-\frac{d}{2dt}\int_{D} |\nabla\hat{u}|^2\dx 
		+\int_{D}\nabla\hat{u}\cdot\nabla R_{int1}\dx
		+\int_{\partial D} R_{sb1}R_{sb2}\cdot\bm{n}\ds
		+\int_{D} R_{int2}\hat{v}\dx
		\nonumber\\
		&\leq-\frac{d}{2dt}\int_{D} |\nabla\hat{u}|^2\dx +\frac{1}{2}\int_{D} |\nabla\hat{u}|^2\dx+\frac{1}{2}\int_{D} |\nabla R_{int1}|^2\dx
		\nonumber\\
		&\qquad+\frac{1}{2}\int_{\partial D}(|R_{sb1}|^2+|R_{sb2}|^2)\ds+\frac{1}{2}\int_{D} |\hat{v}|^2\dx+\frac{1}{2}\int_{D} |R_{int2}|^2\dx.
	\end{align}
	Here, we have used $\hat{v}=\hat{u}_t-R_{int1}$. 
 
	By adding \eqref{sec5_eq0} to \eqref{sec5_eq1}, we have
	\begin{align}\label{sec5_eq2}
		&\frac{d}{2dt}\int_{D} |\hat{u}|^2\dx+\frac{d}{2dt}\int_{D} |\nabla\hat{u}|^2\dx +\frac{d}{2dt}\int_{D} |\hat{v}|^2\dx
		\nonumber\\
		&\qquad \leq \int_{D} |\hat{u}|^2\dx+\frac{1}{2}\int_{D} |\nabla\hat{u}|^2\dx+\int_{D} |\hat{v}|^2\dx+\frac{1}{2}\int_{D} |R_{int1}|^2\dx
		\nonumber\\
		&\qquad+\frac{1}{2}\int_{D} |R_{int2}|^2\dx+\frac{1}{2}\int_{D} |\nabla R_{int1}|^2\dx+\frac{1}{2}\int_{\partial D}(|R_{sb1}|^2+|R_{sb2}|^2)\ds.
	\end{align}
	Integrating \eqref{sec5_eq2} over $[0,
	\tau]$ for any $\tau \leq T$ and applying the Cauchy–Schwarz inequality, we obtain
	\begin{align*}
		&\int_{D} |\hat{u}(\bm{x},\tau)|^2\dx 
		+\int_{D} |\nabla\hat{u}(\bm{x},\tau)|^2\dx
		+\int_{D} |\hat{v}(\bm{x},\tau)|^2\dx\\
		&\qquad\leq\int_{D}|R_{tb1}|^2\dx +\int_{D}|R_{tb2}|^2\dx+\int_{D}|\nabla R_{tb1}|^2\dx
		+2\int_{0}^{\tau}\int_{D} \left(|\hat{u}|^2+ |\nabla\hat{u}|^2+ |\hat{v}|^2 \right)\dx\dt
		\\
		&\qquad+ \int_{0}^{T}\int_{D}\left(|R_{int1}|^2+|R_{int2}|^2+|\nabla R_{int1}|^2\right)\dx\dt+\int_{0}^{T}\int_{\partial D}(|R_{sb1}|^2+|R_{sb2}|^2)\ds\dt.
	\end{align*}

 We apply the integral form of the Gr${\rm\ddot{o}}$nwall inequality to the above inequality to get
	\[
	\int_{D} \left(|\hat{u}(\bm{x},\tau)|^2+|\nabla\hat{u}(\bm{x},\tau)|^2+|\hat{v}(\bm{x},\tau)|^2\right)\dx\leq C_G\exp(2T),
	\]
	where
	\begin{align*}
		C_G&=\int_{D}(|R_{tb1}|^2+|R_{tb2}|^2+|\nabla R_{tb1}|^2)\dx +  \int_{0}^{T}\int_{D}(|R_{int1}|^2+|R_{int2}|^2+|\nabla R_{int1}|^2)\dx\dt
		\\
		&\qquad+\int_{0}^{T}\int_{\partial D}(|R_{sb1}|^2+|R_{sb2}|^2)\ds\dt.
	\end{align*}
	Then, we integrate the above inequality over $[0,T]$ to yield \eqref{lem5.4}.
\end{proof} 

\begin{Remark}\label{sec5_Remark1} 
For the wave equations~\eqref{wave} with periodic boundary, we would like to mention below two other forms for the generalization error (and the related training loss). Compared with \eqref{wave_G}, they  differ only on the spatial boundary $\Omega_{*}$, i.e.,
	\begin{align}\label{wave_G1}
		\mathcal{E}_G(\theta)^2&=\int_{\Omega}|R_{int1}[u_{\theta},v_{\theta}](\bm{x},t)|^2\dx\dt+\int_{\Omega}|R_{int2}[u_{\theta},v_{\theta}](\bm{x},t)|^2\dx\dt+\int_{\Omega}|\nabla R_{int1}[u_{\theta},v_{\theta}](\bm{x},t)|^2\dx\dt
		\nonumber\\
		&+\int_{D}|R_{tb1}[u_{\theta}](\bm{x})|^2\dx+\int_{D}|R_{tb2}[v_{\theta}](\bm{x})|^2\dx
		+\int_{D}|\nabla R_{tb1}[u_{\theta}](\bm{x})|^2\dx
		\nonumber\\
		&
		+\left(\int_{\Omega_*}|R_{sb1}[v_{\theta}](\bm{x},t)|^2\ds\dt\right)^{\frac{1}{2}},
	\end{align}
	and
	\begin{align}\label{wave_G2}
		\mathcal{E}_G(\theta)^2&=\int_{\Omega}|R_{int1}[u_{\theta},v_{\theta}](\bm{x},t)|^2\dx\dt+\int_{\Omega}|R_{int2}[u_{\theta},v_{\theta}](\bm{x},t)|^2\dx\dt+\int_{\Omega}|\nabla R_{int1}[u_{\theta},v_{\theta}](\bm{x},t)|^2\dx\dt
		\nonumber\\
		&+\int_{D}|R_{tb1}[u_{\theta}](\bm{x})|^2\dx+\int_{D}|R_{tb2}[v_{\theta}](\bm{x})|^2\dx
		+\int_{D}|\nabla R_{tb1}[u_{\theta}](\bm{x})|^2\dx
		\nonumber\\
		&
		+\left(\int_{\Omega_*}|R_{sb2}[u_{\theta}](\bm{x},t)|^2\ds\dt\right)^{\frac{1}{2}}.
	\end{align}
	The related training loss functions are given by
	\begin{align}\label{wave_TT1}
		\mathcal{E}_T(\theta,\mathcal{S})^2&
		=\mathcal{E}_T^{int1}(\theta,\mathcal{S}_{int})^2+\mathcal{E}_T^{int2}(\theta,\mathcal{S}_{int})^2+\mathcal{E}_T^{int3}(\theta,\mathcal{S}_{int})^2	+\mathcal{E}_T^{tb1}(\theta,\mathcal{S}_{tb})^2
		\nonumber\\
		& +\mathcal{E}_T^{tb2}(\theta,\mathcal{S}_{tb})^2 +\mathcal{E}_T^{tb3}(\theta,\mathcal{S}_{tb})^2 
		+\mathcal{E}_T^{sb1}(\theta,\mathcal{S}_{sb}),
	\end{align}
	or
	\begin{align}\label{wave_TT2}
		\mathcal{E}_T(\theta,\mathcal{S})^2&
		=\mathcal{E}_T^{int1}(\theta,\mathcal{S}_{int})^2+\mathcal{E}_T^{int2}(\theta,\mathcal{S}_{int})^2+\mathcal{E}_T^{int3}(\theta,\mathcal{S}_{int})^2	+\mathcal{E}_T^{tb1}(\theta,\mathcal{S}_{tb})^2
		\nonumber\\
		& +\mathcal{E}_T^{tb2}(\theta,\mathcal{S}_{tb})^2 +\mathcal{E}_T^{tb3}(\theta,\mathcal{S}_{tb})^2+\mathcal{E}_T^{sb2}(\theta,\mathcal{S}_{sb}).
	\end{align}
	
 These three forms for the generalization error result from different treatments of the boundary term $\int_{\partial D}\hat{v}\nabla\hat{u}\cdot\bm{n}$ in the proof of Theorem~\ref{sec5_Theorem2}:
	\begin{align*}
		&\int_{\partial D} \hat{v}\nabla\hat{u}\cdot\bm{n}\ds=\int_{\partial D} R_{sb1}\nabla\hat{u}\cdot\bm{n}\ds\leq |\partial D|^{\frac{1}{2}}(\|u\|_{C^1(\partial D\times[0,t])}+||u_{\theta}||_{C^1(\partial D\times[0,t])})\left(\int_{\partial D}|R_{sb1}|^2\ds\right)^{\frac{1}{2}},\\
		&\int_{\partial D} \hat{v}\nabla\hat{u}\cdot\bm{n}\ds=\int_{\partial D} \hat{v}R_{sb2}\cdot\bm{n}\ds\leq |\partial D|^{\frac{1}{2}}(\|v\|_{C^0(\partial D\times[0,t])}+||v_{\theta}||_{C^0(\partial D\times[0,t])})\left(\int_{\partial D}|R_{sb2}|^2\ds\right)^{\frac{1}{2}},\\
		&\int_{\partial D} \hat{v}\nabla\hat{u}\cdot\bm{n}\ds=\int_{\partial D} R_{sb1}R_{sb2}\cdot\bm{n}\ds\leq \frac{1}{2}\left(\int_{\partial D} |R_{sb1}|^2\ds+\int_{\partial D}|R_{sb2}|^2\ds\right).
	\end{align*}
	Our numerical experiments indicate that 
 adopting the training loss \eqref{wave_TT1} or \eqref{wave_TT2} seems to lead to poorer simulation results. 
 For the periodic boundary, both terms $R_{sb1}$ and $R_{sb2}$ may be needed for the periodicity information.
 We suspect that this may be why only a single boundary term ($R_{sb1}$ or $R_{sb2}$), as given by~\eqref{wave_TT1} and~\eqref{wave_TT2}, leads to poorer numerical results.
\end{Remark}

\begin{Theorem}\label{sec5_Theorem3} Let $d\in \mathbb{N}$ and $T>0$. Let $u\in C^4(\Omega)$ and $v\in C^3(\Omega)$ be the classical solution of
	the wave equations \eqref{wave},  and let $(u_{\theta},v_{\theta})$ denote the PINN approximation with parameter $\theta \in \Theta$. Then the total error satisfies
	\begin{align}\label{lem5.5}
		&\int_0^{T}\int_{D}(|\hat{u}(\bm{x},t)|^2+|\nabla \hat{u}(\bm{x},t)|^2+|\hat{v}(\bm{x},t)|^2)\dx\dt\leq C_TT\exp(2T)
		\nonumber\\
		&\qquad=\mathcal{O}(\mathcal{E}_T(\theta,\mathcal{S})^2 + M_{int}^{-\frac{2}{d+1}} +M_{tb}^{-\frac{2}{d}}+M_{sb}^{-\frac{2}{d}}).
	\end{align}
	The constant $C_T$ is defined as
	\begin{align*}
		C_T =&C_{({R_{tb1}^2})}M_{tb}^{-\frac{2}{d}}+\mathcal{Q}_{M_{tb}}^{D}(R_{tb1}^2)+C_{({R_{tb2}^2})}M_{tb}^{-\frac{2}{d}}+\mathcal{Q}_{M_{tb}}^{D}(R_{tb2}^2)+C_{(|\nabla R_{tb1}|^2)}M_{tb}^{-\frac{2}{d}}+\mathcal{Q}_{M_{tb}}^{D}(|\nabla R_{tb1}|^2)\\
		&+C_{({R_{int1}^2})}M_{int}^{-\frac{2}{d+1}}+\mathcal{Q}_{M_{int}}^{\Omega}(R_{int1}^2)+C_{(R_{int2}^2)}M_{int}^{-\frac{2}{d+1}}+\mathcal{Q}_{M_{int}}^{\Omega}(R_{int2}^2)+C_{(|\nabla R_{int1}|^2)}M_{int}^{-\frac{2}{d+1}}\\
		&+\mathcal{Q}_{M_{int}}^{\Omega}(|\nabla R_{int1}|^2)+C_{({R_{sb1}^2})}M_{sb}^{-\frac{2}{d}}+\mathcal{Q}_{M_{sb}}^{\Omega_*}(R_{sb1}^2)+C_{({R_{sb2}^2})}M_{sb}^{-\frac{2}{d}}+\mathcal{Q}_{M_{sb}}^{\Omega_*}(R_{sb2}^2),
	\end{align*}
where
\begin{align*}
	&C_{({R_{tb1}^2})}\lesssim\|\hat{u}\|_{C^2}^2, \quad C_{({R_{tb2}^2})}\lesssim \|\hat{v}\|_{C^2}^2, \quad C_{(|\nabla R_{tb1}|^2)}\lesssim \|\hat{u}\|_{C^3}^2, \quad C_{({R_{int1}^2})}\lesssim \|\hat{u}\|_{C^3}^2+\|\hat{u}\|_{C^2}^2,\\
	&\qquad C_{(R_{int2}^2)}, C_{(|\nabla R_{int1}|^2)}\lesssim \|\hat{u}\|_{C^4}^2+\|\hat{v}\|_{C^3}^2, \quad C_{({R_{sb1}^2})}\lesssim \|\hat{v}\|_{C^3}^2,\quad C_{({R_{sb2}^2})}\lesssim\|\hat{u}\|_{C^4}^2,
\end{align*}
and the bounds $\|u_{\theta}\|_{C^n}$ and $\|v_{\theta}\|_{C^n}$ ($n \in \mathbb{N}$) are given by Lemma \ref{Ar_3}. 
\end{Theorem}
\begin{proof} By combining Theorem \ref{sec5_Theorem2} with the quadrature error formula \eqref{int1}, we have
	\begin{align*}
		\int_{D}|R_{tb1}|^2\dx&=\int_{D}|R_{tb1}|^2\dx-\mathcal{Q}_{M_{tb}}^{D}(R_{tb1}^2)+\mathcal{Q}_{M_{tb}}^{D}(R_{tb1}^2)\\
		&\leq C_{({R_{tb1}^2})}M_{tb}^{-\frac{2}{d}}+\mathcal{Q}_{M_{tb}}^{D}(R_{tb1}^2),\\
		\int_{D}|R_{tb2}|^2\dx&=\int_{D}|R_{tb2}|^2\dx-\mathcal{Q}_{M_{tb}}^{D}(R_{tb2}^2)+\mathcal{Q}_{M_{tb}}^{D}(R_{tb2}^2)\\
		&\leq C_{({R_{tb2}^2})}M_{tb}^{-\frac{2}{d}}+\mathcal{Q}_{M_{tb}}^{D}(R_{tb2}^2),\\
		\int_{D}|\nabla R_{tb1}|^2\dx&=\int_{D}|\nabla R_{tb1}|^2\dx-\mathcal{Q}_{M_{tb}}^{D}(|\nabla R_{tb1}|^2)+\mathcal{Q}_{M_{tb}}^{D}(|\nabla R_{tb1}|^2)\\
		&\leq C_{(|\nabla R_{tb1}|^2)}M_{tb}^{-\frac{2}{d}}+\mathcal{Q}_{M_{tb}}^{D}(|\nabla R_{tb1}|^2),\\
		\int_{\Omega}|R_{int1}|^2\dx\dt&=\int_{\Omega}|R_{int1}|^2\dx\dt-\mathcal{Q}_{M_{int}}^{\Omega}(R_{int1}^2)+\mathcal{Q}_{M_{int}}^{\Omega}(R_{int1}^2)\\
		&\leq C_{({R_{int1}^2})}M_{int}^{-\frac{2}{d+1}}+\mathcal{Q}_{M_{int}}^{\Omega}(R_{int1}^2),\\
		\int_{\Omega}|R_{int2}|^2\dx\dt&=\int_{\Omega}|R_{int2}|^2\dx\dt-\mathcal{Q}_{M_{int}}^{\Omega}(R_{int2}^2)+\mathcal{Q}_{M_{int}}^{\Omega}(R_{int2}^2)\\
		&\leq C_{(R_{int2}^2)}M_{int}^{-\frac{2}{d+1}}+\mathcal{Q}_{M_{int}}^{\Omega}(R_{int2}^2),\\
		\int_{\Omega}|\nabla R_{int1}|^2\dx\dt&=\int_{\Omega}|\nabla R_{int1}|^2\dx\dt-\mathcal{Q}_{M_{int}}^{\Omega}(|\nabla R_{int1}|^2)+\mathcal{Q}_{M_{int}}^{\Omega}(|\nabla R_{int1}|^2)\\
		&\leq C_{(|\nabla R_{int1}|^2)}M_{int}^{-\frac{2}{d+1}}+\mathcal{Q}_{M_{int}}^{\Omega}(|\nabla R_{int1}|^2),\\
		\int_{\Omega_*}|R_{sb1}|^2\ds\dt&=\int_{\Omega_*}|R_{sb1}|^2\ds\dt-\mathcal{Q}_{M_{sb}}^{\Omega_*}(R_{sb1}^2)+\mathcal{Q}_{M_{sb}}^{\Omega_*}(R_{sb1}^2)\\
		&\leq C_{({R_{sb1}^2})}M_{sb}^{-\frac{2}{d}}+\mathcal{Q}_{M_{sb}}^{\Omega_*}(R_{sb1}^2),\\
		\int_{\Omega_*}|R_{sb2}|^2\ds\dt&=\int_{\Omega_*}|R_{sb2}|^2\ds\dt-\mathcal{Q}_{M_{sb}}^{\Omega_*}(R_{sb2}^2)+\mathcal{Q}_{M_{sb}}^{\Omega_*}(R_{sb2}^2)\\
		&\leq C_{({R_{sb2}^2})}M_{sb}^{-\frac{2}{d}}+\mathcal{Q}_{M_{sb}}^{\Omega_*}(R_{sb2}^2).
	\end{align*}
	By the above inequalities and \eqref{lem5.4}, it holds that
	\begin{equation*}
		\int_0^{T}\int_{D}(|\hat{u}(\bm{x},t)|^2+|\nabla \hat{u}(\bm{x},t)|^2+|\hat{v}(\bm{x},t)|^2)\dx\dt\leq  C_TT\exp(2T),
	\end{equation*}
	where 
\begin{align*}
	C_T =&C_{({R_{tb1}^2})}M_{tb}^{-\frac{2}{d}}+\mathcal{Q}_{M_{tb}}^{D}(R_{tb1}^2)+C_{({R_{tb2}^2})}M_{tb}^{-\frac{2}{d}}+\mathcal{Q}_{M_{tb}}^{D}(R_{tb2}^2)+C_{(|\nabla R_{tb1}|^2)}M_{tb}^{-\frac{2}{d}}+\mathcal{Q}_{M_{tb}}^{D}(|\nabla R_{tb1}|^2)\\
	&+C_{({R_{int1}^2})}M_{int}^{-\frac{2}{d+1}}+\mathcal{Q}_{M_{int}}^{\Omega}(R_{int1}^2)+C_{(R_{int2}^2)}M_{int}^{-\frac{2}{d+1}}+\mathcal{Q}_{M_{int}}^{\Omega}(R_{int2}^2)+C_{(|\nabla R_{int1}|^2)}M_{int}^{-\frac{2}{d+1}}\\
	&+\mathcal{Q}_{M_{int}}^{\Omega}(|\nabla R_{int1}|^2)+C_{({R_{sb1}^2})}M_{sb}^{-\frac{2}{d}}+\mathcal{Q}_{M_{sb}}^{\Omega_*}(R_{sb1}^2)+C_{({R_{sb2}^2})}M_{sb}^{-\frac{2}{d}}+\mathcal{Q}_{M_{sb}}^{\Omega_*}(R_{sb2}^2).
\end{align*}
The complexities of the constants $C_{({R_{q}^2})}$ are given by Lemma \ref{Ar_3}, and we observe that for every residual $R_q$, it holds that $\|R_q^2\|_{C^n}\leq 2^n\|R_q\|_{C^n}^2$ ($n \in \mathbb{N}$) for $R_q=R_{tb1}$, $R_{tb2}$, $\nabla R_{tb1}$, $R_{int1}$, $R_{int2}$, $\nabla R_{int1}$ and $R_{sb2}$. 
\end{proof} 

	\section{Physics Informed Neural Networks for Approximating the Sine-Gordon Equation}\label{Sine-Gordon}
\subsection{Sine-Gordon Equation}
Let $D\subset \mathbb{R}^d$ be  an open connected bounded set with a boundary $\partial D$. We consider the following Sine-Gordon equation:
\begin{subequations}\label{SG}
             \begin{align}
                          \label{SG_eq0}
                          &u_{t} - v = 0  \ \qquad\qquad\qquad\qquad\qquad\qquad\quad\ \ \, \text{in}\ D\times[0,T],\\
                          \label{SG_eq1}
                          &\varepsilon^2v_{t} = a^2\Delta u - \varepsilon_1^2u-g(u)+f   \  \ \quad\qquad\qquad \text{in}\ D\times[0,T],\\
                          \label{SG_eq2}
                          &u(\bm{x},0)=\psi_{1}(\bm{x})\qquad\qquad\qquad\qquad\qquad\qquad \text{in}\ D,\\
                          \label{SG_eq3}
                          &v(\bm{x},0)=\psi_{2}(\bm{x}) \qquad\qquad\qquad\qquad\qquad\qquad \, \text{in}\ D,\\
                          \label{SG_eq4}
                          &u(\bm{x},t)|_{\partial D}=u_{d}(t)  \qquad\qquad\qquad\qquad\qquad\ \ \ \, \text{in}\ \partial D\times[0,T],
             \end{align}
\end{subequations}
where $u$ and $v$ are the field functions to be solved for, $f$ is a source term, and $u_d$, $\psi_1$ and $\psi_2$ denote the boundary/initial conditions. $\varepsilon>0$, $a>0$ and $\varepsilon_1\geq 0$ are  constants.
$g(u)$ is a nonlinear term. We assume that the nonlinearity is globally Lipschitz, i.e., there exists a constant $L$ (independent of $v$ and $w$) such that
\begin{equation}\label{non2}
             |g(v) - g(w)|\leq L|v-w|, \qquad \forall v, \, w \in \mathbb{R}.
\end{equation}

\begin{Remark}\label{sec6_Remark1} 
The existence and regularity of the solution to the Sine-Gordon equation with different nonlinear terms have been the subject of several studies in the literature; see~\cite{Baoxiang1997Classical,Kubota2001Global,Shatah1982Global,Shatah1985Normal,Temam1997Infinite}.
             
             The book \cite{Temam1997Infinite} provides the existence and regularity result of the following Sine-Gordon equation, 
             \[
             u_{tt} + \alpha u_t - \Delta u + g(u)=f.
             \]
             Let $\alpha \in \mathbb{R}$, $g(u)$ be a $C^2$ function from $\mathbb{R}$ to $\mathbb{R}$ and satisfy certain assumptions. If $f\in C([0,T]; L^{2}(D))$, $\psi_{1}\in H^{1}(D)$ and $\psi_{2}\in L^{2}(D)$, then there exists a unique solution $u$ to this Sine-Gordon equation such that $u \in C([0,T];H^{1}(D))$ and $u_t \in C([0,T];L^{2}(D))$. Furthermore, 
             $f'\in C([0,T]; L^{2}(D))$, $\psi_{1}\in H^{2}(D)$ and $\psi_{2}\in H^{1}(D)$, it holds $u \in C([0,T];H^{2}(D))$ and $u_t \in C([0,T];H^{1}(D))$. 
             
             Let $g$ be a smooth function of degree 2. The following equation is studied in \cite{Shatah1985Normal}, 
             \[
             u_{tt} - \Delta u  + u+ g(u,u_t,u_{tt})=0,
             \]
             where it is reformulated as
             \[
             \bm{u}_{t} = A\bm{u} + G(\bm{u}),
             \]
             in which $\bm{u}=\begin{pmatrix}
             	u\\	
             	u_t
             \end{pmatrix}$, $A =\begin{pmatrix}
             0 & 1 \\
             \Delta-1  & 0
         \end{pmatrix}$ 
     and $G=\begin{pmatrix}
         0, \\
         -g(u,u_t,u_{tt})
     \end{pmatrix}$. Set $X=H^k(\mathbb{R}^n) \bigoplus H^{k-1}(\mathbb{R}^n)$, $k>n+2+2a$ with $a>1$. Given $\bm{u}_0 =\begin{pmatrix}
     \psi_1\\	
     \psi_2
 \end{pmatrix} \in X$ and $\|\bm{u}_0\|_{X} = \sigma$, there exists a $T_0 = T_0(\sigma)$ depending on the size of the initial data $\sigma$ and a unique solution $\bm{u} \in C([0,T_0],X)$.

The reference~\cite{Baoxiang1997Classical} provides the following result.
Under certain conditions for the nonlinear term $g(u)$, with $f = 0$, $d\leq 5$, $k \geq \frac{d}{2} + 1$, $\psi_{1}\in H^{k}(D)$ and $\psi_{2}\in H^{k-1}(D)$, there exists a unique solution $u \in C((0,\infty);H^{k}(D))$ of nonlinear Klein–Gordon equation.  

The following result is due to~\cite{Kubota2001Global}.
Under certain conditions for the nonlinear term $g(u)$, with $f = 0$, $\psi_{1}\in H^{k}(D)$ and $\psi_{2}\in H^{k-1}(D)$ with a positive constant $k\geq 4$, there exists a positive constant $T_k$ and a unique solution $u \in C([0,T_k];H^{k}(D))\cap C^1([0,T_k];H^{k-1}(D))\cap C^2([0,T_k];H^{k-2}(D))$ to the nonlinear wave equations with different speeds of propagation.
\end{Remark}

A survey of literature indicates that,
while several works have touched on the regularity of the solution to the Sine-Gordon equations, none of them is comprehensive. 
To facilitate the subsequent analyses, we make the following  assumption in light of Remark \ref{sec6_Remark1}. Let $k\geq1$, $g(u)$ and $f$ be sufficiently smooth and bounded. Given $\psi_{1}\in H^{r}(D)$ and $\psi_{2}\in H^{r-1}(D)$ with $r \geq \frac{d}{2} + k$,  we assume that there exists $T>0$ and a classical solution $u$ and $v$ to the Sine-Gordon equations \eqref{SG} such that $u \in C([0,T];H^{k}(D))$ and $v \in C([0,T];H^{k-1}(D))$. Then, it follows that $u \in C^k(D\times[0,T])$ and $v \in C^{k-1}(D\times[0,T])$ based on the Sobolev embedding theorem.

\subsection{Physics Informed Neural Networks}

Let $\Omega = D\times [0,T]$ 
and $\Omega_* = \partial D\times [0,T]$ be the space-time domain. We define the following residuals for the PINN approximation, $u_{\theta}: \Omega \rightarrow \mathbb{R}$ and $v_{\theta}: \Omega \rightarrow \mathbb{R}$, for the Sine-Gordon equations \eqref{SG}:
\begin{subequations}\label{SG_pinn}
	\begin{align}
		\label{SG_pinn_eq1}
		&R_{int1}[u_{\theta},v_{\theta}](\bm{x},t) =u_{\theta t}-v_{\theta},\\
		\label{SG_pinn_eq2}
		&R_{int2}[u_{\theta},v_{\theta}](\bm{x},t) =\varepsilon^2v_{\theta t}-a^2\Delta u_{\theta} +\varepsilon_1^2u_{\theta}+g(u_{\theta})-f,\\
		\label{SG_pinn_eq3}
		&R_{tb1}[u_{\theta}](\bm{x}) =u_{\theta}(\bm{x},0)-\psi_{1}(\bm{x}),\\
		\label{SG_pinn_eq4}
		&R_{tb2}[v_{\theta}](\bm{x}) =v_{\theta}(\bm{x},0)-\psi_{2}(\bm{x}),\\
		\label{SG_pinn_eq5}
		&R_{sb}[v_{\theta}](\bm{x},t) =v_{\theta}(\bm{x},t)|_{\partial D}-u_{dt}(t),
	\end{align}
\end{subequations}
where $u_{dt}=\frac{\partial u_d}{\partial t}$.
Note that for the exact solution $(u,v)$, $R_{int1}[u,v]=R_{int2}[u,v]=R_{tb1}[u]=R_{tb2}[v]=R_{sb}[v]=0$. 
With PINN we minimize the following generalization error,
\begin{align}\label{SG_G}
	\mathcal{E}_G(\theta)^2&=\int_{\Omega}|R_{int1}[u_{\theta},v_{\theta}](\bm{x},t)|^2\dx\dt+\int_{\Omega}|R_{int2}[u_{\theta},v_{\theta}](\bm{x},t)|^2\dx\dt+\int_{\Omega}|\nabla R_{int1}[u_{\theta},v_{\theta}](\bm{x},t)|^2\dx\dt
	\nonumber\\
	&+\int_{D}|R_{tb1}[u_{\theta}](\bm{x})|^2\dx+\int_{D}|R_{tb2}[v_{\theta}](\bm{x})|^2\dx
	+\int_{D}|\nabla R_{tb1}[u_{\theta}](\bm{x})|^2\dx
	\nonumber\\
	&
	+\left(\int_{\Omega_*}|R_{sb}[v_{\theta}](\bm{x},t)|^2\ds\dt\right)^{\frac{1}{2}}.
\end{align}

Let 
\begin{equation*}
\hat{u} = u_{\theta}-u, \quad \hat{v} = v_{\theta}-v,
\end{equation*}
where $(u,v)$ denotes the exact solution.
We define the total error of the PINN approximation of the Sine-Gordon equations \eqref{SG} as, 
\begin{equation}\label{SG_total}
	\mathcal{E}(\theta)^2=\int_{\Omega}(|\hat{u}(\bm{x},t)|^2+a^2|\nabla \hat{u}(\bm{x},t)|^2+\varepsilon^2|\hat{v}(\bm{x},t)|^2)\dx\dt.
\end{equation}
Then we choose the training set $\mathcal{S} \subset \overline{D}\times [0,T]$ with $\mathcal{S} = \mathcal{S}_{int} \cup \mathcal{S}_{sb} \cup \mathcal{S}_{tb}$, based on suitable quadrature points:
\begin{itemize}
    \item Interior training points $\mathcal{S}_{int}=\{{z}_n\}$ for $1\leq n \leq N_{int}$, with each ${z}_n= (\bm{x},t)_n \in D \times(0,T)$.
    
    \item Spatial boundary training points $\mathcal{S}_{sb}=\{{z}_n\}$ for $1\leq n \leq N_{sb}$, with each ${z}_n= (\bm{x},t)_n \in \partial D\times (0,T)$.
 
    \item Temporal boundary training points $\mathcal{S}_{tb}=\{\bm{x}_n\}$ for $1\leq n \leq N_{tb}$ with  each $\bm{x}_n \in D$.
\end{itemize}
The integrals in \eqref{SG_G} are approximated by a  numerical quadrature rule, resulting in the training loss,
\begin{align}\label{SG_T}
	\mathcal{E}_T(\theta,\mathcal{S})^2&
	=\mathcal{E}_T^{int1}(\theta,\mathcal{S}_{int})^2+\mathcal{E}_T^{int2}(\theta,\mathcal{S}_{int})^2+\mathcal{E}_T^{int3}(\theta,\mathcal{S}_{int})^2
	\nonumber\\
	&
	+\mathcal{E}_T^{tb1}(\theta,\mathcal{S}_{tb})^2 +\mathcal{E}_T^{tb2}(\theta,\mathcal{S}_{tb})^2 +\mathcal{E}_T^{tb3}(\theta,\mathcal{S}_{tb})^2 
	+\mathcal{E}_T^{sb}(\theta,\mathcal{S}_{sb}),
\end{align}
where
\begin{subequations}
	\begin{align}
		\label{SG_T1}
		\mathcal{E}_T^{int1}(\theta,\mathcal{S}_{int})^2 &= \sum_{n=1}^{N_{int}}\omega_{int}^n|R_{int1}[u_{\theta},v_{\theta}](\bm{x}_{int}^n,t_{int}^n)|^2,\\
		\label{SG_T01}
		\mathcal{E}_T^{int2}(\theta,\mathcal{S}_{int})^2 &= \sum_{n=1}^{N_{int}}\omega_{int}^n|R_{int2}[u_{\theta},v_{\theta}](\bm{x}_{int}^n,t_{int}^n)|^2,\\
		\label{SG_T001}
		\mathcal{E}_T^{int3}(\theta,\mathcal{S}_{int})^2 &= \sum_{n=1}^{N_{int}}\omega_{int}^n|\nabla R_{int1}[u_{\theta},v_{\theta}](\bm{x}_{int}^n,t_{int}^n)|^2,\\
		\label{SG_T2}
		\mathcal{E}_T^{tb1}(\theta,\mathcal{S}_{tb})^2 &= \sum_{n=1}^{N_{tb}}\omega_{tb}^n|R_{tb1}[u_{\theta}](\bm{x}_{tb}^n)|^2,\\
		\label{SG_T02}
		\mathcal{E}_T^{tb2}(\theta,\mathcal{S}_{tb})^2 &= \sum_{n=1}^{N_{tb}}\omega_{tb}^n|R_{tb2}[v_{\theta}](\bm{x}_{tb}^n)|^2,\\
		\label{SG_T002}
		\mathcal{E}_T^{tb3}(\theta,\mathcal{S}_{tb})^2 &= \sum_{n=1}^{N_{tb}}\omega_{tb}^n|\nabla R_{tb1}[u_{\theta}](\bm{x}_{tb}^n)|^2,\\
		\label{SG_T3}
		\mathcal{E}_T^{sb}(\theta,\mathcal{S}_{sb})^2&= \sum_{n=1}^{N_{sb}}\omega_{sb}^n|R_{sb}[v_{\theta}](\bm{x}_{sb}^n,t_{sb}^n)|^2.
	\end{align}
\end{subequations}
Here the quadrature points in space-time constitute the data sets $\mathcal{S}_{int} = \{(\bm{x}_{int}^n,t_{int}^n)\}_{n=1}^{N_{int}}$, $\mathcal{S}_{tb} = \{\bm{x}_{tb}^n)\}_{n=1}^{N_{tb}}$ and $\mathcal{S}_{sb} = \{(\bm{x}_{sb}^n,t_{sb}^n)\}_{n=1}^{N_{sb}}$, and $\omega_{\star}^n$ are the quadrature weights with $\star$ being $int$, $tb$ or $sb$.

\subsection{Error Analysis}  


By substracting the Sine-Gordon equations \eqref{SG} from the residual equations~\eqref{SG_pinn}, we get,
\begin{subequations}\label{SG_error}
	\begin{align}
		\label{SG_error_eq1}
		&R_{int1}=\hat{u}_t-\hat{v},\\
		\label{SG_error_eq2}
		&R_{int2}=\varepsilon^2\hat{v}_{t}-a^2\Delta \hat{u} +\varepsilon_1^2\hat{u}+g(u_{\theta})-g(u),\\
		\label{SG_error_eq3}
		&R_{tb1}=\hat{u}(\bm{x},0),\\
		\label{SG_error_eq4}
		&R_{tb2}=\hat{v}(\bm{x},0),\\
		\label{SG_error_eq5}
		&R_{sb}=\hat{v}(\bm{x},t)|_{\partial D}.
	\end{align}
\end{subequations}
The results on the PINN approximations to the Sine-Gordon equations are summarized in the following theorems.

\begin{Theorem}\label{sec6_Theorem1} 
	Let $d$, $r$, $k \in \mathbb{N}$ with $k\geq 3$. Assume that $g(u)$ is Lipschitz continuous, $u \in C^k(D\times[0,T])$ and $v \in C^{k-1}(D\times[0,T])$. Then for every integer $N>5$, there exist $\tanh$ neural networks $u_{\theta}$ and $v_{\theta}$, each with two hidden layers, of widths at most $3\lceil\frac{k}{2}\rceil|P_{k-1,d+2}| + \lceil NT\rceil+ d(N-1)$ and $3\lceil\frac{d+3}{2}\rceil|P_{d+2,d+2}| \lceil NT\rceil N^d$, such that
	\begin{subequations}
		\begin{align}
			\label{lem6.1}
			&\|R_{int1}\|_{L^2(\Omega)},\|R_{tb1}\|_{L^2(D)}\lesssim {\rm ln}NN^{-k+1},\\
			\label{lem6.2}
			&\|R_{int2}\|_{L^2(\Omega)},\|\nabla R_{int1}\|_{L^2(\Omega)}, \|\nabla R_{tb1}\|_{L^2(D)}\lesssim {\rm ln}^2NN^{-k+2},\\
			\label{lem6.3}
			&\|R_{tb2}\|_{L^2(D)},\|R_{sb}\|_{L^2(\partial D\times [0,t])}\lesssim {\rm ln}NN^{-k+2}.
		\end{align}
	\end{subequations}
\end{Theorem}  
\noindent The proof of this theorem is provided in the Appendix~\ref{Proof_sinegordon}.

Theorem \ref{sec6_Theorem1} implies that the PINN residuals in \eqref{SG_pinn} can be made arbitrarily small by choosing a  sufficiently large $N$. Therefore, the generalization error $\mathcal{E}_G(\theta)^2$ can be made arbitrarily small. 

We next show that the PINN total approximation error $\mathcal{E}(\theta)^2$ can be controlled by the generalization error $\mathcal{E}_G(\theta)^2$ (Theorem \ref{sec6_Theorem2} below), and by the training error~$\mathcal{E}_T(\theta,\mathcal{S})^2$ (Theorem \ref{sec6_Theorem3} below). The proofs for Theorem \ref{sec6_Theorem2} and Theorem \ref{sec6_Theorem3} are provided in the Appendix \ref{Proof_sinegordon}.

\begin{Theorem}\label{sec6_Theorem2} Let $d\in \mathbb{N}$, $u\in C^1(\Omega)$ and $v\in C^0(\Omega)$ be the classical solution of
	the Sine-Gordon equation \eqref{SG}. Let $(u_{\theta},v_{\theta})$ denote the  PINN approximation with parameter $\theta$. Then the following relation holds, 
	\begin{equation}\label{lem6.09}
		\mathcal{E}(\theta)^2 = \int_0^{T}\int_{D}(|\hat{u}(\bm{x},t)|^2+a^2|\nabla \hat{u}(\bm{x},t)|^2+\varepsilon^2|\hat{v}(\bm{x},t)|^2)\dx\dt
		\leq C_GT\exp\left((2+\varepsilon_1^2+L+a^2)T\right),
	\end{equation}
	where
	\begin{align*}
	&C_G=\int_{D}(|R_{tb1}|^2+a^2|\nabla R_{tb1}|^2+\varepsilon^2|R_{tb2}|^2)\dx + \int_{0}^{T}\int_{D}(|R_{int1}|^2+|R_{int2}|^2+a^2|\nabla R_{int1}|^2)\dx\dt\\
	&\qquad +  2C_{\partial D}|T|^{\frac{1}{2}}\left(\int_{0}^{T}\int_{\partial D}|R_{sb}|^2\ds\dt\right)^{\frac{1}{2}},
	\end{align*}
and $C_{\partial D}=a^2|\partial D|^{\frac{1}{2}}(\|u\|_{C^1(\partial D\times[0,t])}+||u_{\theta}||_{C^1(\partial D\times[0,t])})$.
\end{Theorem}

\begin{Theorem}\label{sec6_Theorem3} Let $d\in \mathbb{N}$ and $T>0$, and let $u\in C^4(\Omega)$ and $v\in C^3(\Omega)$ be the classical solution to the Sine-Gordon equation \eqref{SG}. Let $(u_{\theta},v_{\theta})$ denote the PINN approximation with parameter $\theta \in \Theta$. Then the following relation holds, 
             \begin{align}\label{lem6.9}
                          &\int_0^{T}\int_{D}(|\hat{u}(\bm{x},t)|^2+a^2|\nabla \hat{u}(\bm{x},t)|^2+\varepsilon^2|\hat{v}(\bm{x},t)|^2)\dx\dt\leq C_TT\exp\left((2+\varepsilon_1^2+L+a^2)T\right) 
                          	\nonumber\\
                          &\qquad=\mathcal{O}(\mathcal{E}_T(\theta,\mathcal{S})^2 + M_{int}^{-\frac{2}{d+1}} +M_{tb}^{-\frac{2}{d}}+M_{sb}^{-\frac{1}{d}}),                       
             \end{align}
             where the constant $C_T$ is defined by 
             \begin{align*}
             	C_T=&C_{({R_{tb1}^2})}M_{tb}^{-\frac{2}{d}}+\mathcal{Q}_{M_{tb}}^{D}(R_{tb1}^2)+\varepsilon^2\left(C_{({R_{tb2}^2})}M_{tb}^{-\frac{2}{d}}+\mathcal{Q}_{M_{tb}}^{D}(R_{tb2}^2) \right)\\
             	&+a^2\left( C_{(|\nabla R_{tb1}|^2)}M_{tb}^{-\frac{2}{d}}+\mathcal{Q}_{M_{tb}}^{D}(|\nabla R_{tb1}|^2) \right)+C_{({R_{int1}^2})}M_{int}^{-\frac{2}{d+1}}+\mathcal{Q}_{M_{int}}^{\Omega}(R_{int1}^2)\\
             	&+C_{({R_{int2}^2})}M_{int}^{-\frac{2}{d+1}}+\mathcal{Q}_{M_{int}}^{\Omega}(R_{int2}^2)
             	+a^2\left(C_{(|\nabla R_{int1}|^2)}M_{int}^{-\frac{2}{d+1}}+\mathcal{Q}_{M_{int}}^{\Omega}(|\nabla R_{int1}|^2)\right),\\
             	&+2C_{\partial D}|T|^{\frac{1}{2}}\left(C_{({R_{sb}^2})}M_{sb}^{-\frac{2}{d}}+\mathcal{Q}_{M_{sb}}^{\Omega_*}(R_{sb}^2)\right)^{\frac{1}{2}}.
             \end{align*}
\end{Theorem}

It follows from Theorem \ref{sec6_Theorem3} that the PINN approximation error $\mathcal{E}(\theta)^2$ can be arbitrarily small, provided that the training error $\mathcal{E}_T(\theta,\mathcal{S})^2$ is sufficiently small  and the sample set is sufficiently large.

	\section{Physics Informed Neural Networks for Approximating  Linear Elastodynamic Equation}\label{Elasto-dynamics}
\subsection{Linear Elastodynamic Equation}

Consider an elastic body occupying an open, bounded convex polyhedral domain $D\subset \mathbb{R}^d$. The boundary $\partial D = \Gamma_D\cup \Gamma_N$, with the outward  unit normal vector $\bm{n}$, is assumed to be composed of two disjoint portions $\Gamma_D\neq \emptyset$ and $\Gamma_N$, with $\Gamma_D\cap \Gamma_N=\emptyset$. Given a suitable external load $\bm{f} \in L^2((0,T];\bm{L}^2(D))$, and suitable initial/boundary data $\bm{g} \in C^1((0,T];\bm{H}^{\frac{1}{2}}(\Gamma_N))$, $\bm{\psi}_{1}\in \bm{H}_{0,\Gamma_D}^{\frac{1}{2}}(D)$ and $\bm{\psi}_{2}\in \bm{L}^2(D)$, we consider the linear elastodynamic equations, 
\begin{subequations}\label{elast}
	\begin{align}
		\label{elast_eq0}
		&\bm{u}_{t} - \bm{v} = 0  \ \quad\quad\qquad\qquad\qquad\qquad\qquad\quad\ \ \, \text{in}\ D\times [0,T],\\
		\label{elast_eq1}
		&\rho\bm{v}_{t} - 2\mu\nabla\cdot(\underline{\bm{\varepsilon}}(\bm{u})) -\lambda\nabla(\nabla\cdot\bm{u})= \bm{f}  \, \quad\qquad\text{in}\ D\times [0,T],\\
		\label{elast_eq2}
		&\bm{u}=\bm{u}_{d} \ \ \, \quad\qquad\qquad\qquad\qquad\qquad\qquad\qquad \text{in}\ \Gamma_D\times [0,T],\\
		\label{elast_eq3}
		&2\mu\underline{\bm{\varepsilon}}(\bm{u})\bm{n} +\lambda(\nabla\cdot\bm{u})\bm{n}=\bm{g}\ \ \quad\qquad\qquad\qquad \text{in}\ \Gamma_N\times [0,T],\\
		\label{elast_eq4}
		&\bm{u}=\bm{\psi}_{1} \ \  \quad\qquad\qquad\qquad\qquad\qquad\qquad\qquad \, \text{in}\ D\times\{0\},\\
		\label{elast_eq5}
		&\bm{v}=\bm{\psi}_{2} \qquad\qquad\qquad\qquad\qquad\qquad\qquad\qquad \text{in}\ D\times\{0\}.
	\end{align}
\end{subequations}
In the above system, $\bm{u}=(u_1,u_2,\cdots,u_d)$ and $\bm{v}=(v_1,v_2,\cdots,v_d)$ denote the displacement and the velocity, respectively, and $[0,T]$ (with $T>0$) denotes the time domain. $\underline{\bm{\varepsilon}}(\bm{u})$ is the strain tensor, $\underline{\bm{\varepsilon}}(\bm{u})=\frac{1}{2}(\nabla\bm{u}+\nabla\bm{u}^T)$. The constants $\lambda$ and $\mu$ are the first and the second Lam${\rm \acute{e}}$ parameters, respectively.  

Combining \eqref{elast_eq0} and \eqref{elast_eq1}, we can recover the classical linear elastodynamics equation:
\begin{equation}\label{elast1}
	\rho\bm{u}_{tt} - 2\mu\nabla\cdot(\underline{\bm{\varepsilon}}(\bm{u})) -\lambda\nabla(\nabla\cdot\bm{u})= \bm{f} \qquad\text{in}\ D\times [0,T].
\end{equation}
The well-posedness of this equation is established in~\cite{Hughes1978Classical}.

\begin{Lemma}[\cite{Hughes1978Classical,Yosida1980Functional}]\label{sec9_Lemma1} 
Let $\bm{\psi}_{1}\in H^{r}(D)$, $\bm{\psi}_{2}\in H^{r-1}(D)$ and $\bm{f}\in H^{r-1}(D\times [0,T])$ with $r\geq1$. Then there exists a unique solution $\bm{u}$ to the classical linear elastodynamic equation \eqref{elast1} such that $\bm{u}(t = 0) =\bm{\psi}_{1}$, $\bm{u}_t(t = 0) = \bm{\psi}_{2}$ and $\bm{u} \in C^l([0,T];H^{r-l}(D))$ with $0\leq l \leq r$.
\end{Lemma}

\begin{Lemma}\label{sec9_Lemma2} Let $k\in \mathbb{N}$, $\bm{\psi}_{1}\in H^{r}(D)$, $\bm{\psi}_{2}\in H^{r-1}(D)$ and $\bm{f}\in H^{r-1}(D\times [0,T])$ with $r>\frac{d}{2}+k$, then there exists $T>0$ and a classical solution $(\bm{u},\bm{v})$ to the elastodynamic equations \eqref{elast}  such that $\bm{u}(t = 0) =\bm{\psi}_{1}$, $\bm{u}_t(t = 0) = \bm{\psi}_{2}$, $\bm{u} \in C^k(D\times [0,T])$ and $\bm{v} \in C^{k-1}(D\times [0,T])$.
\end{Lemma}
\begin{proof}  As $r>\frac{d}{2}+k$, $H^{r-k}(D)$ is a Banach algebra. By Lemma \ref{sec9_Lemma1}, there exists $T > 0$ and the solution $(\bm{u}, \bm{v})$ to the linear elastodynamics equations such that $\bm{u}(t = 0) =\bm{\psi}_{1}$, 
	$\bm{v}(t = 0) = \bm{\psi}_{2}$, $\bm{u} \in C^l([0,T];H^{r-l}(D))$ with $0\leq l \leq r$ and $\bm{v} \in C^l([0,T];H^{r-1-l}(D))$ with $0\leq l \leq r-1$. 
	
	Since $\bm{u} \in \cap_{l=0}^k C^l([0,T];H^{r-l}(D))$ and $\bm{v} \subset \cap_{l=0}^{k-1} C^l([0,T];H^{r-l-1}(D))$. By applying the Sobolev embedding theorem and $r>\frac{d}{2}+k$, we obtain $H^{r-l}(D) \subset C^{r-l}(D)$ and $H^{r-l-1}(D) \subset C^{r-l-1}(D)$ for $0\leq l\leq k$. Therefore, $\bm{u} \in C^k(D\times [0,T])$ and $\bm{v} \in C^{k-1}(D\times [0,T])$.
\end{proof} 

\subsection{Physics Informed Neural Networks}

We now consider the PINN approximation of the linear elastodynamic equations \eqref{elast}.
Let $\Omega = D\times [0,T]$, 
$\Omega_{D} = \Gamma_D\times [0,T]$ and $\Omega_{N} = \Gamma_N\times [0,T]$ denote the space-time domain. 
Define the following residuals for the PINN approximation $\bm{u}_{\theta}:\Omega \rightarrow \mathbb{R}$ and $\bm{v}_{\theta}: \Omega \rightarrow \mathbb{R}$ for the elastodynamic equations \eqref{elast}:
\begin{subequations}\label{elast_pinn}
	\begin{align}
		\label{elast_pinn_eq1}
		&\bm{R}_{int1}[\bm{u}_{\theta},\bm{v}_{\theta}](\bm{x},t) =\bm{u}_{\theta t}-\bm{v}_{\theta},\\
		\label{elast_pinn_eq2}
		&\bm{R}_{int2}[\bm{u}_{\theta},\bm{v}_{\theta}](\bm{x},t) =\rho \bm{v}_{\theta t}-2\mu\nabla\cdot(\underline{\bm{\varepsilon}}(\bm{u}_{\theta})) -\lambda\nabla(\nabla\cdot\bm{u}_{\theta})-\bm{f},\\
		\label{elast_pinn_eq3}
		&\bm{R}_{tb1}[\bm{u}_{\theta}](\bm{x}) =\bm{u}_{\theta}(\bm{x}, 0)-\bm{\psi}_{1}(\bm{x}),\\
		\label{elast_pinn_eq4}
		&\bm{R}_{tb2}[\bm{v}_{\theta}](\bm{x}) =\bm{v}_{\theta}(\bm{x}, 0)-\bm{\psi}_{2}(\bm{x}),\\
		\label{elast_pinn_eq5}
		&\bm{R}_{sb1}[\bm{v}_{\theta}](\bm{x},t) =\bm{v}_{\theta}|_{\Gamma_D}-\bm{u}_{dt},\\
		\label{elast_pinn_eq6}
		&\bm{R}_{sb2}[\bm{u}_{\theta}](\bm{x},t) =(2\mu\underline{\bm{\varepsilon}}(\bm{u}_{\theta})\bm{n} +\lambda(\nabla\cdot\bm{u}_{\theta})\bm{n})|_{\Gamma_N}-\bm{g}.
	\end{align}
\end{subequations}
Note that for the exact solution $(\bm{u},\bm{v})$,  we have $\bm{R}_{int1}[\bm{u},\bm{v}]=\bm{R}_{int2}[\bm{u},\bm{v}]=\bm{R}_{tb1}[\bm{u}]=\bm{R}_{tb2}[\bm{v}]=\bm{R}_{sb1}[\bm{v}]=\bm{R}_{sb2}[\bm{u}]=0$. 
With PINN we minimize the the following generalization error,
\begin{align}\label{elast_G}
	\mathcal{E}_G(\theta)^2&=\int_{\Omega}|\bm{R}_{int1}[\bm{u}_{\theta},\bm{v}_{\theta}](\bm{x},t)|^2\dx\dt+\int_{\Omega}|\bm{R}_{int2}[\bm{u}_{\theta},\bm{v}_{\theta}](\bm{x},t)|^2\dx\dt+\int_{\Omega}|\underline{\bm{\varepsilon}}(\bm{R}_{int1}[\bm{u}_{\theta},\bm{v}_{\theta}](\bm{x},t))|^2\dx\dt
	\nonumber\\
	&+\int_{\Omega}|\nabla\cdot(\bm{R}_{int1}[\bm{u}_{\theta},\bm{v}_{\theta}](\bm{x},t))|^2\dx\dt+\int_{D}|\bm{R}_{tb1}[\bm{u}_{\theta}](\bm{x})|^2\dx+\int_{D}|\bm{R}_{tb2}[\bm{v}_{\theta}](\bm{x})|^2\dx
	\nonumber\\
	&
	+\int_{D}|\underline{\bm{\varepsilon}}(\bm{R}_{tb1}[\bm{u}_{\theta}](\bm{x}))|^2\dx+\int_{D}|\nabla\cdot \bm{R}_{tb1}[\bm{u}_{\theta}](\bm{x})|^2\dx
	\nonumber\\
	&+\left(\int_{\Omega_D}|\bm{R}_{sb1}[\bm{v}_{\theta}](\bm{x},t)|^2\ds\dt\right)^{\frac{1}{2}}
	+\left(\int_{\Omega_N}|\bm{R}_{sb2}[\bm{u}_{\theta}](\bm{x},t)|^2\ds\dt\right)^{\frac{1}{2}}.
\end{align}

Let 
\begin{equation*}
\hat{\bm{u}} = \bm{u}_{\theta}-\bm{u}, \quad \hat{\bm{v}} = \bm{v}_{\theta}-\bm{v}
\end{equation*}
denote the difference between the solution to the elastodynamic equations~\eqref{elast} and the PINN approximation with parameter $\theta$. We define the total error of the PINN approximation as,
\begin{equation}\label{elast_total}
	\mathcal{E}(\theta)^2=\int_{\Omega}( |\hat{\bm{u}}(\bm{x},t)|^2+2\mu|\underline{\bm{\varepsilon}}(\hat{\bm{u}}(\bm{x},t))|^2
	+\lambda|\nabla\cdot\hat{\bm{u}}(\bm{x},t)|^2+\rho|\hat{\bm{v}}(\bm{x},t)|^2)\dx\dt.
\end{equation}

We choose the training set $\mathcal{S} \subset \overline{D}\times [0,T]$ based on suitable quadrature points. The full training set is defined by $\mathcal{S} = \mathcal{S}_{int} \cup \mathcal{S}_{sb} \cup \mathcal{S}_{tb}$, and $\mathcal{S}_{sb}=\mathcal{S}_{sb1}\cup \mathcal{S}_{sb2}$:
\begin{itemize}
    \item Interior training points $\mathcal{S}_{int}=\{{z}_n\}$ for $1\leq n \leq N_{int}$, with each ${z}_n= (\bm{x},t)_n \in D \times(0,T)$. 
    
    \item Spatial boundary training points $\mathcal{S}_{sb1}=\{{z}_n\}$ for $1\leq n \leq N_{sb1}$, with each ${z}_n= (\bm{x},t)_n \in \Gamma_D \times (0,T)$, and $\mathcal{S}_{sb2}=\{{z}_n\}$ for $1\leq n \leq N_{sb2}$, with each ${z}_n= (\bm{x},t)_n \in \Gamma_N\times (0,T)$.
 
	\item Temporal boundary training points $\mathcal{S}_{tb}=\{\bm{x}_n\}$ for $1\leq n \leq N_{tb}$ with  each $\bm{x}_n \in D$.
\end{itemize}

Then, the integrals in \eqref{elast_G} can be approximated by a suitable numerical quadrature, resulting in the following training loss,
\begin{align}\label{elast_T}
	\mathcal{E}_T(\theta,\mathcal{S})^2&
	=\mathcal{E}_T^{int1}(\theta,\mathcal{S}_{int})^2+\mathcal{E}_T^{int2}(\theta,\mathcal{S}_{int})^2+\mathcal{E}_T^{int3}(\theta,\mathcal{S}_{int})^2+\mathcal{E}_T^{int4}(\theta,\mathcal{S}_{int})^2+\mathcal{E}_T^{tb1}(\theta,\mathcal{S}_{tb})^2
	\nonumber\\
	&\quad
	+\mathcal{E}_T^{tb2}(\theta,\mathcal{S}_{tb})^2 +\mathcal{E}_T^{tb3}(\theta,\mathcal{S}_{tb})^2
	+\mathcal{E}_T^{tb4}(\theta,\mathcal{S}_{tb})^2
	+\mathcal{E}_T^{sb1}(\theta,\mathcal{S}_{sb1}) 
	+\mathcal{E}_T^{sb2}(\theta,\mathcal{S}_{sb2}),
\end{align}
where, 
\begin{subequations}\label{elast_TT}
	\begin{align}
		\label{elast_T1}
		\mathcal{E}_T^{int1}(\theta,\mathcal{S}_{int})^2 &= \sum_{n=1}^{N_{int}}\omega_{int}^n|\bm{R}_{int1}[\bm{u}_{\theta},\bm{v}_{\theta}](\bm{x}_{int}^n,t_{int}^n)|^2,\\
		\label{elast_T2}
		\mathcal{E}_T^{int2}(\theta,\mathcal{S}_{int})^2 &= \sum_{n=1}^{N_{int}}\omega_{int}^n|\bm{R}_{int2}[\bm{u}_{\theta},\bm{v}_{\theta}](\bm{x}_{int}^n,t_{int}^n)|^2,\\
		\label{elast_T3}
		\mathcal{E}_T^{int3}(\theta,\mathcal{S}_{int})^2 &= \sum_{n=1}^{N_{int}}\omega_{int}^n|\underline{\bm{\varepsilon}}(\bm{R}_{int1}[\bm{u}_{\theta},\bm{v}_{\theta}](\bm{x}_{int}^n,t_{int}^n))|^2,\\
		\label{elast_T33}
		\mathcal{E}_T^{int4}(\theta,\mathcal{S}_{int})^2 &= \sum_{n=1}^{N_{int}}\omega_{int}^n|\nabla\cdot \bm{R}_{int1}[\bm{u}_{\theta},\bm{v}_{\theta}](\bm{x}_{int}^n,t_{int}^n)|^2,\\
		\label{elast_T4}
		\mathcal{E}_T^{tb1}(\theta,\mathcal{S}_{tb})^2 &= \sum_{n=1}^{N_{tb}}\omega_{tb}^n|\bm{R}_{tb1}[\bm{u}_{\theta}](\bm{x}_{tb}^n)|^2,\\
		\label{elast_T5}
		\mathcal{E}_T^{tb2}(\theta,\mathcal{S}_{tb})^2 &= \sum_{n=1}^{N_{tb}}\omega_{tb}^n|\bm{R}_{tb2}[\bm{v}_{\theta}](\bm{x}_{tb}^n)|^2,\\
		\label{elast_T6}
		\mathcal{E}_T^{tb3}(\theta,\mathcal{S}_{tb})^2 &= \sum_{n=1}^{N_{tb}}\omega_{tb}^n|\underline{\bm{\varepsilon}}(\bm{R}_{tb1}[\bm{u}_{\theta}](\bm{x}_{tb}^n))|^2,\\
		\label{elast_T66}
		\mathcal{E}_T^{tb4}(\theta,\mathcal{S}_{tb})^2 &= \sum_{n=1}^{N_{tb}}\omega_{tb}^n|\nabla\cdot \bm{R}_{tb1}[\bm{u}_{\theta}](\bm{x}_{tb}^n)|^2,\\
		\label{elast_T7}
		\mathcal{E}_T^{sb1}(\theta,\mathcal{S}_{sb1})^2&= \sum_{n=1}^{N_{sb1}}\omega_{sb1}^n|\bm{R}_{sb1}[\bm{v}_{\theta}](\bm{x}_{sb}^n,t_{sb}^n)|^2,\\
		\label{elast_T8}
		\mathcal{E}_T^{sb2}(\theta,\mathcal{S}_{sb2})^2&= \sum_{n=1}^{N_{sb2}}\omega_{sb2}^n|\bm{R}_{sb2}[\bm{u}_{\theta}](\bm{x}_{sb}^n,t_{sb}^n)|^2.
	\end{align}
\end{subequations}
Here the quadrature points in space-time constitute the data sets $\mathcal{S}_{int} = \{(\bm{x}_{int}^n,t_{int}^n)\}_{n=1}^{N_{int}}$, $\mathcal{S}_{tb} = \{\bm{x}_{tb}^n)\}_{n=1}^{N_{tb}}$, $\mathcal{S}_{sb1} = \{(\bm{x}_{sb1}^n,t_{sb1}^n)\}_{n=1}^{N_{sb1}}$ and $\mathcal{S}_{sb2} = \{(\bm{x}_{sb2}^n,t_{sb2}^n)\}_{n=1}^{N_{sb2}}$. $\omega_{\star}^n$ denote the suitable quadrature weights with $\star$ being $int$, $tb$, $sb1$ and $sb2$.

\subsection{Error Analysis}  

Subtracting the elastodynamic equations \eqref{elast} from the residual equations~\eqref{elast_pinn}, we obtain
\begin{subequations}\label{elast_error}
	\begin{align}
		\label{elast_error_eq1}
		&\bm{R}_{int1}=\hat{\bm{u}}_t-\hat{\bm{v}},\\
		\label{elast_error_eq2}
		&\bm{R}_{int2}=\rho\hat{\bm{v}}_t-2\mu\nabla\cdot(\underline{\bm{\varepsilon}}(\hat{\bm{u}})) -\lambda\nabla(\nabla\cdot\hat{\bm{u}}),\\
		\label{elast_error_eq3}
		&\bm{R}_{tb1}=\hat{\bm{u}}|_{t=0},\\
		\label{elast_error_eq4}
		&\bm{R}_{tb2}=\hat{\bm{v}}|_{t=0},\\
		\label{elast_error_eq5}
		&\bm{R}_{sb1}=\hat{\bm{v}}|_{\Gamma_D},\\
		\label{elast_error_eq6}
		&\bm{R}_{sb2}=(2\mu\underline{\bm{\varepsilon}}(\hat{\bm{u}})\bm{n} +\lambda(\nabla\cdot\hat{\bm{u}})\bm{n})|_{\Gamma_N}.
	\end{align}
\end{subequations}	
%
The PINN approximation results are summarized in the following three theorems. The proofs of these theorems are provided in the Appendix~\ref{Proof_elasto}.

\begin{Theorem}\label{sec9_Theorem1} 
		Let $d$, $r$, $k \in \mathbb{N}$ with $k\geq 3$. Let $\bm{\psi}_{1}\in H^{r}(D)$, $\bm{\psi}_{2}\in H^{r-1}(D)$ and $\bm{f}\in H^{r-1}(D\times [0,T])$ with $r>\frac{d}{2}+k$. 
		For every integer $N>5$, there exist $\tanh$ neural networks $(\bm{u}_j)_{\theta}$ and $(\bm{v}_j)_{\theta}$, with $j=1,2,\cdots,d$, each with two hidden layers, of widths at most $3\lceil\frac{k}{2}\rceil|P_{k-1,d+2}| + \lceil NT\rceil+ d(N-1)$ and $3\lceil\frac{d+3}{2}\rceil|P_{d+2,d+2}| \lceil NT\rceil N^d$, such that
	\begin{subequations}
		\begin{align}
			\label{lem9.1}
			&\|\bm{R}_{int1}\|_{L^2(\Omega)},\|\bm{R}_{tb1}\|_{L^2(\Omega)}\lesssim {\rm ln}NN^{-k+1},\\
			\label{lem9.2}
			&\|\bm{R}_{int2}\|_{L^2(\Omega)},\|\underline{\bm{\varepsilon}}(\bm{R}_{int1})\|_{L^2(\Omega)},\|\nabla\cdot\bm{R}_{int1}\|_{L^2(\Omega)}\lesssim {\rm ln}^2NN^{-k+2},\\
			\label{lem9.3}
			&\|\underline{\bm{\varepsilon}}(\bm{R}_{tb1})\|_{L^2(D)},\|\nabla\cdot\bm{R}_{tb1}\|_{L^2(D)},\|\bm{R}_{sb2}\|_{L^2(\Gamma_N\times [0,t])}\lesssim {\rm ln}^2NN^{-k+2},\\
			\label{lem9.4}
			&\|\bm{R}_{tb2}\|_{L^2(D)}, \|\bm{R}_{sb1}\|_{L^2(\Gamma_D\times [0,t])}\lesssim {\rm ln}NN^{-k+2}.
		\end{align}
	\end{subequations}
\end{Theorem}

It follows from Theorem \ref{sec9_Theorem1} that,
by choosing a sufficiently large $N$, 
one can make the PINN residuals in~\eqref{elast_pinn}, and thus  the generalization error $\mathcal{E}_G(\theta)^2$ in~\eqref{elast_G}, arbitrarily small. 

\begin{Theorem}\label{sec9_Theorem2} Let $d\in \mathbb{N}$, $\bm{u} \in C^1(\Omega)$ and $\bm{v}\in C(\Omega)$ be the classical solution to the linear elastodynamic equation~\eqref{elast}. Let $(\bm{u}_{\theta},\bm{v}_{\theta})$ denote the PINN approximation with the parameter $\theta$. Then the following relation holds,
	\begin{equation*}
		\int_0^{T}\int_{D}( |\hat{\bm{u}}(\bm{x},t)|^2+2\mu|\underline{\bm{\varepsilon}}(\hat{\bm{u}}(\bm{x},t))|^2
		+\lambda|\nabla\cdot\hat{\bm{u}}(\bm{x},t)|^2+\rho|\hat{\bm{v}}(\bm{x},t)|^2)\dx\dt
		\leq C_GT\exp\left((2+2\mu+\lambda)T\right),
	\end{equation*}
	where
	\begin{align*}
		&C_G=\int_{D}|\bm{R}_{tb1}|^2\dx+\int_{D}2\mu|\underline{\bm{\varepsilon}}(\bm{R}_{tb1})|^2\dx+\int_{D}\lambda|\nabla\cdot \bm{R}_{tb1}|^2\dx+\rho\int_{D}|\bm{R}_{tb2}|^2\dx\\
		&\qquad+\int_{0}^{T}\int_{D}\left(|\bm{R}_{int1}|^2+2\mu|\underline{\bm{\varepsilon}}(\bm{R}_{int1})|^2+\lambda|\nabla\cdot \bm{R}_{int1}|^2+|\bm{R}_{int2}|^2\right)\dx\dt\\
		&\qquad +2|T|^{\frac{1}{2}}C_{\Gamma_D}\left(\int_{0}^{T}\int_{\Gamma_D}|\bm{R}_{sb1}|^2\ds\dt\right)^{\frac{1}{2}}+2|T|^{\frac{1}{2}}C_{\Gamma_N}\left(\int_{0}^{T}\int_{\Gamma_N}|\bm{R}_{sb2}|^2\ds\dt\right)^{\frac{1}{2}},
	\end{align*}
	with $C_{\Gamma_D}=(2\mu+\lambda)|\Gamma_D|^{\frac{1}{2}}\|\bm{u}\|_{C^1(\Gamma_D\times [0,T])}+(2\mu+\lambda)|\Gamma_D|^{\frac{1}{2}}||\bm{u}_{\theta}||_{C^1(\Gamma_D\times [0,T])}$ and $C_{\Gamma_N}=|\Gamma_N|^{\frac{1}{2}}(\|\bm{v}\|_{C(\Gamma_N\times [0,T])}+||\bm{v}_{\theta}||_{C(\Gamma_N\times [0,T])})$.
\end{Theorem}
 Theorem \ref{sec9_Theorem2} shows that the total error of the PINN approximation $\mathcal{E}(\theta)^2$ can be controlled by the generalization error $\mathcal{E}_G(\theta)^2$.

\begin{Theorem}\label{sec9_Theorem3} Let $d\in \mathbb{N}$, $\bm{u}\in C^4(\Omega)$ and $\bm{v}\in C^3(\Omega)$ be the classical solution to
	the linear elastodynamic equation~\eqref{elast}. Let $(\bm{u}_{\theta},\bm{v}_{\theta})$ denote the PINN approximation with the parameter $\theta$. Then the following relation holds,
  \begin{align}\label{lem9.9}
	&\int_0^{T}\int_{D}( |\hat{\bm{u}}(\bm{x},t)|^2+2\mu|\underline{\bm{\varepsilon}}(\hat{\bm{u}}(\bm{x},t))|^2
	+\lambda|\nabla\cdot\hat{\bm{u}}(\bm{x},t)|^2+\rho|\hat{\bm{v}}(\bm{x},t)|^2)\dx\dt
	\leq C_TT\exp\left((2+2\mu+\lambda)T\right)
		\nonumber\\
	&\qquad=\mathcal{O}(\mathcal{E}_T(\theta)^2 + M_{int}^{-\frac{2}{d+1}} +M_{tb}^{-\frac{2}{d}}+M_{sb}^{-\frac{1}{d}}),        
  \end{align}
	where 
	\begin{align*}
		C_T=&C_{({\bm{R}_{tb1}^2})}M_{tb}^{-\frac{2}{d}}+\mathcal{Q}_{M_{tb}}^{D}(\bm{R}_{tb1}^2)
		+\rho\left(C_{({\bm{R}_{tb2}^2})}M_{tb}^{-\frac{2}{d}}+\mathcal{Q}_{M_{tb}}^{D}(\bm{R}_{tb2}^2)\right)
		+ 2\mu\left(C_{(|\underline{\bm{\varepsilon}}(\bm{R}_{tb1})|^2)}M_{tb}^{-\frac{2}{d}}+\mathcal{Q}_{M_{tb}}^{D}(|\underline{\bm{\varepsilon}}(\bm{R}_{tb1})|^2)\right)\\
		&+\lambda\left(C_{(|\nabla\cdot \bm{R}_{tb1}|^2)}M_{tb}^{-\frac{2}{d}}+\mathcal{Q}_{M_{tb}}^{D}(|\nabla\cdot \bm{R}_{tb1}|^2)\right)+C_{({\bm{R}_{int1}^2})}M_{int}^{-\frac{2}{d+1}}+\mathcal{Q}_{M_{int}}^{\Omega}(\bm{R}_{int1}^2)\\
		&+ C_{({\bm{R}_{int2}^2})}M_{int}^{-\frac{2}{d+1}}+\mathcal{Q}_{M_{int}}^{\Omega}(\bm{R}_{int2}^2)
		+2\mu\left(C_{(|\underline{\bm{\varepsilon}}(\bm{R}_{int1})|^2)}M_{int}^{-\frac{2}{d+1}}+\mathcal{Q}_{M_{int}}^{\Omega}(|\underline{\bm{\varepsilon}}(\bm{R}_{int1})|^2)\right)\\
		&+\lambda\left(C_{(|\nabla\cdot\bm{R}_{int1}|^2)}M_{int}^{-\frac{2}{d+1}}+\mathcal{Q}_{M_{int}}^{\Omega}(|\nabla\cdot\bm{R}_{int1}|^2)\right)
		+ 2|T|^{\frac{1}{2}}C_{\Gamma_D}\left(C_{({\bm{R}_{sb1}^2})}M_{sb1}^{-\frac{2}{d}}+\mathcal{Q}_{M_{sb1}}^{\Omega_D}(\bm{R}_{sb1}^2)\right)^{\frac{1}{2}}\\
		&+2|T|^{\frac{1}{2}}C_{\Gamma_N}\left(C_{({\bm{R}_{sb2}^2})}M_{sb2}^{-\frac{2}{d}}+\mathcal{Q}_{M_{sb2}}^{\Omega_N}(\bm{R}_{sb2}^2)\right)^{\frac{1}{2}}.
	\end{align*}
\end{Theorem}
Theorem \ref{sec9_Theorem3} shows that the PINN approximation error $\mathcal{E}(\theta)^2$ can be controlled by the training error $\mathcal{E}_T(\theta,\mathcal{S})^2$ with a large enough sample set $\mathcal{S}$.

	\section{Numerical Examples}\label{numerical examples}

The theoretical analyses from Sections~\ref{Wave} to \ref{Elasto-dynamics} suggest several forms for the PINN loss function with the wave, Sine-Gordon and the linear elastodynamic equations. These forms contain certain non-standard terms, such as the square root of the residuals or the gradient terms on some boundaries, which would generally be  absent from the canonical PINN formulation of the loss function. The presence of such non-standard terms is crucial to bounding the PINN approximation errors, as shown in the error analyses.

These non-standard forms of the loss function lead to a variant PINN algorithm. In this section we illustrate the performance of the variant PINN algorithm as suggested by the theoretical analysis, as well as the more standard PINN algorithm, using several numerical examples in one spatial dimension (1D) plus time for the wave equation and the Sine-Gordon equation, and in two spatial dimensions (2D) plus time for the linear elastodynamic equation.

The following settings are common to all the numerical simulations in this section. Let $(\bm x, t )\in D\times[0,T] $ denote the spatial and temporal coordinates in the spatial-temporal domain, where $ \bm x = x $ and $ \bm x = (x, y) $ for one and two spatial dimensions, respectively. For the wave equation and the Sine-Gordon equation, the neural networks contain two nodes in the input layer (representing $x$ and $t$), two hidden layers with the number of nodes to be specified later, and two nodes in the output layer (representing the solution $u$ and its time detivative $v=\frac{\partial u}{\partial t}$). For the linear elastodynamic equaton, three input nodes and four output nodes are employed in the neural network, as will be explained in more detail later.
We employ the $\tanh$ (hyperbolic tangent) activation function for all the hidden nodes, and no activation function is applied to the output nodes (i.e.~linear).
For training the neural networks, we employ $ N $ collocation points within the spatial-temporal domain drawn from a uniform random distribution, and also $N$ uniform random points on  each spatial boundary and on the initial boundary. In the simulations the value of $ N $ is varied systematically among $ 1000 $, $ 1500 $, $ 2000 $, $ 2500 $ and $ 3000 $. After the neural networks are trained, for the wave equation and the Sine-Gordon equation, we compare the PINN solution and the exact solution on a set of $N_{ev} =3000 \times 3000 $ uniform spatial-temporal grid points (evaluation points) $(x, t)_n\in D\times[0,T] $ ($n=1,\cdots, N_{ev} $) that covers the problem domain and the boundaries. For the elastodynamic equation, we compare the PINN solution and the exact solution at different time instants, and at each time instant the corresponding solutions are evaluated at a uniform set of $ N_{ev} =1500 \times 1500 $  grid points in the spatial domain, $ \bm x_n = (x, y)_n\in D $ ($n=1,\cdots, N_{ev} $).

The PINN errors reported below are computed as follows. Let $ z_n = (\bm x, t)_n $ ($(\bm x, t)_n\in D\times[0,T], n=1,\cdots, N_{ev} $) denote the set of uniform grid  points, where $ N_{ev} $ denote the number of evaluation points. The errors of PINN are defined by,
\begin{subequations}\label{PINN_num_eq1}
	\begin{align}\label{PINN_num_eq1_1}
		&l_2\text{-error} = \frac{\sqrt{\sum_{n=1}^{N_{ev}} |u(z_n) - u_{\theta}(z_n)|^2}}{\sqrt{\sum_{n=1}^{N_{ev}} u(z_n)^2}} =  \frac{\sqrt{\left(\sum_{n=1}^{N_{ev}} |u(z_n) - u_{\theta}(z_n)|^2\right)/N_{ev}}}{\sqrt{\left(\sum_{n=1}^{N_{ev}} u(z_n)^2\right)/N_{ev}}}, \\ 
		\label{PINN_num_eq1_2}
		&l_\infty\text{-error}= \frac{\max\{|u(z_n) - u_{\theta}(z_n)|\}_{n=1}^{N_{ev}} }{\sqrt{\left(\sum_{n=1}^{N_{ev}} u(z_n)^2\right)/N_{ev}}},
	\end{align}
\end{subequations}
where $ u_\theta $ denotes the PINN solution and $u$ denotes the exact solution. 

Our implementation of the PINN algorithm is based on the PyTorch library (\href{https://pytorch.org/}{pytorch.org}). In all the following  numerical examples, we combine the Adam \cite{kingma2014adam} optimizer and the L-BFGS \cite{2006_NumericalOptimization} optimizer (in batch mode) to train the neural network. We first employ the Adam optimizer to train the network for 100 epochs/iterations, and then employ the L-BFGS optimizer to continue the network training for another 30000 iterations. 
We employ the default parameter values in Adam, with the learning rate $0.001$,  $\beta_1=0.9$ and $\beta_2=0.99$. The initial learning rate $1.0$ is adopted in the L-BFGS optimizer.

\subsection{Wave Equation}

\begin{figure}[tb]
	\centerline{
	\subfloat[True solution for $ u $]{\includegraphics[width=0.25\linewidth]{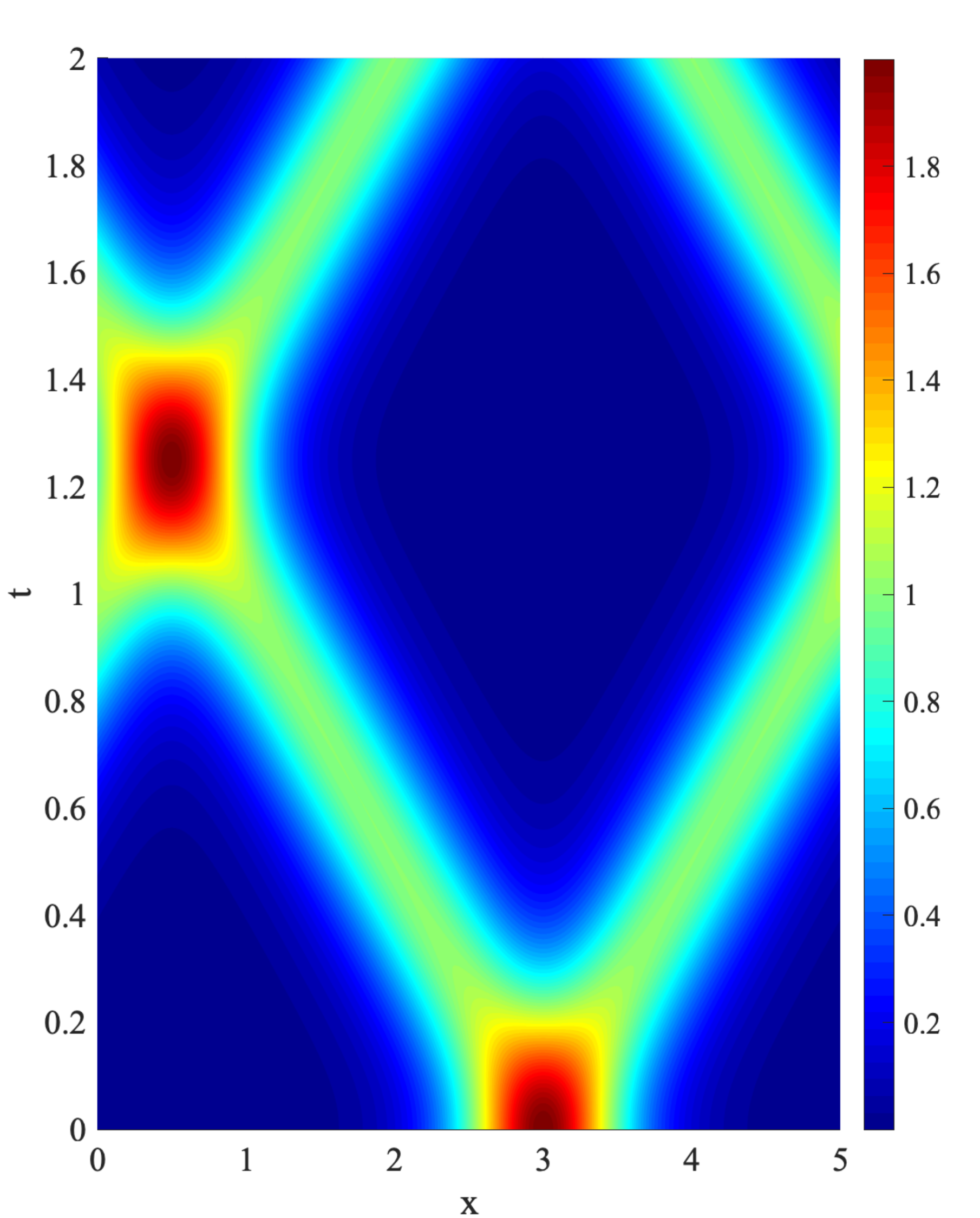}}\hspace{0.2em}
	\subfloat[PINN solution $ u_\theta $]{\includegraphics[width=0.25\linewidth]{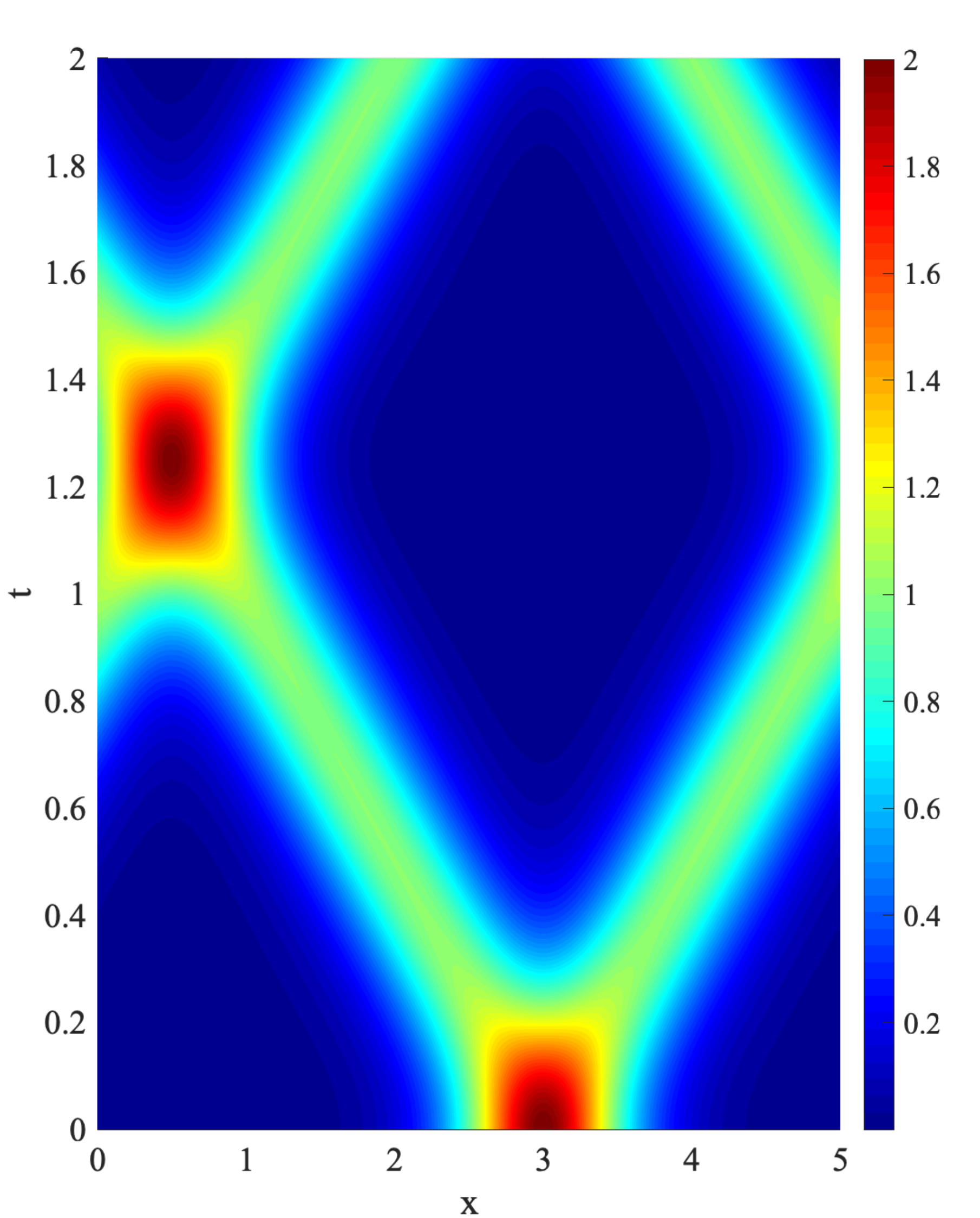}}\hspace{0.2em}
	\subfloat[{PINN} {absolute} error for $u$]{\includegraphics[width=0.25\linewidth]{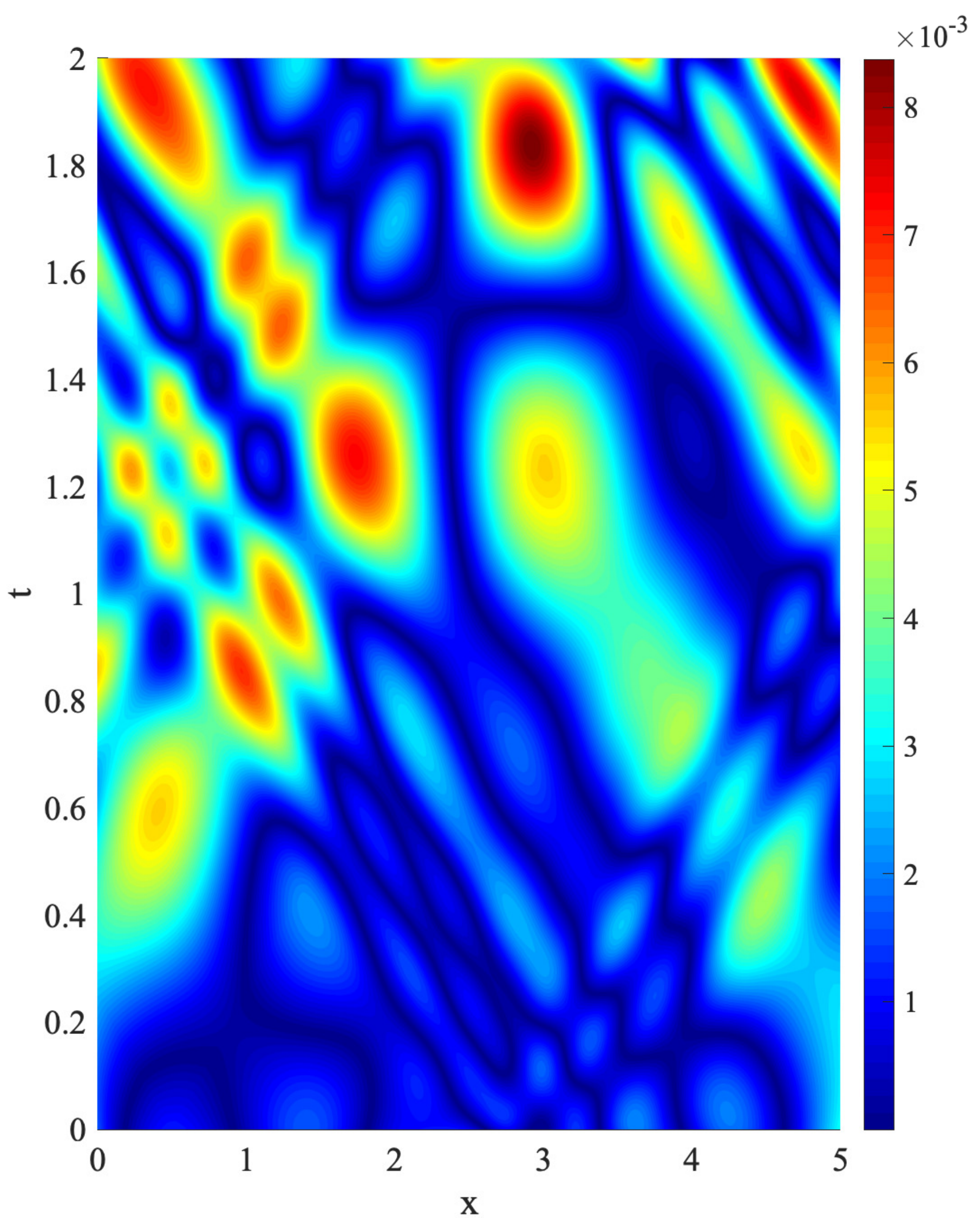}}
 }
	\centerline{
	\subfloat[True solution for $ v $]{\includegraphics[width=0.25\linewidth]{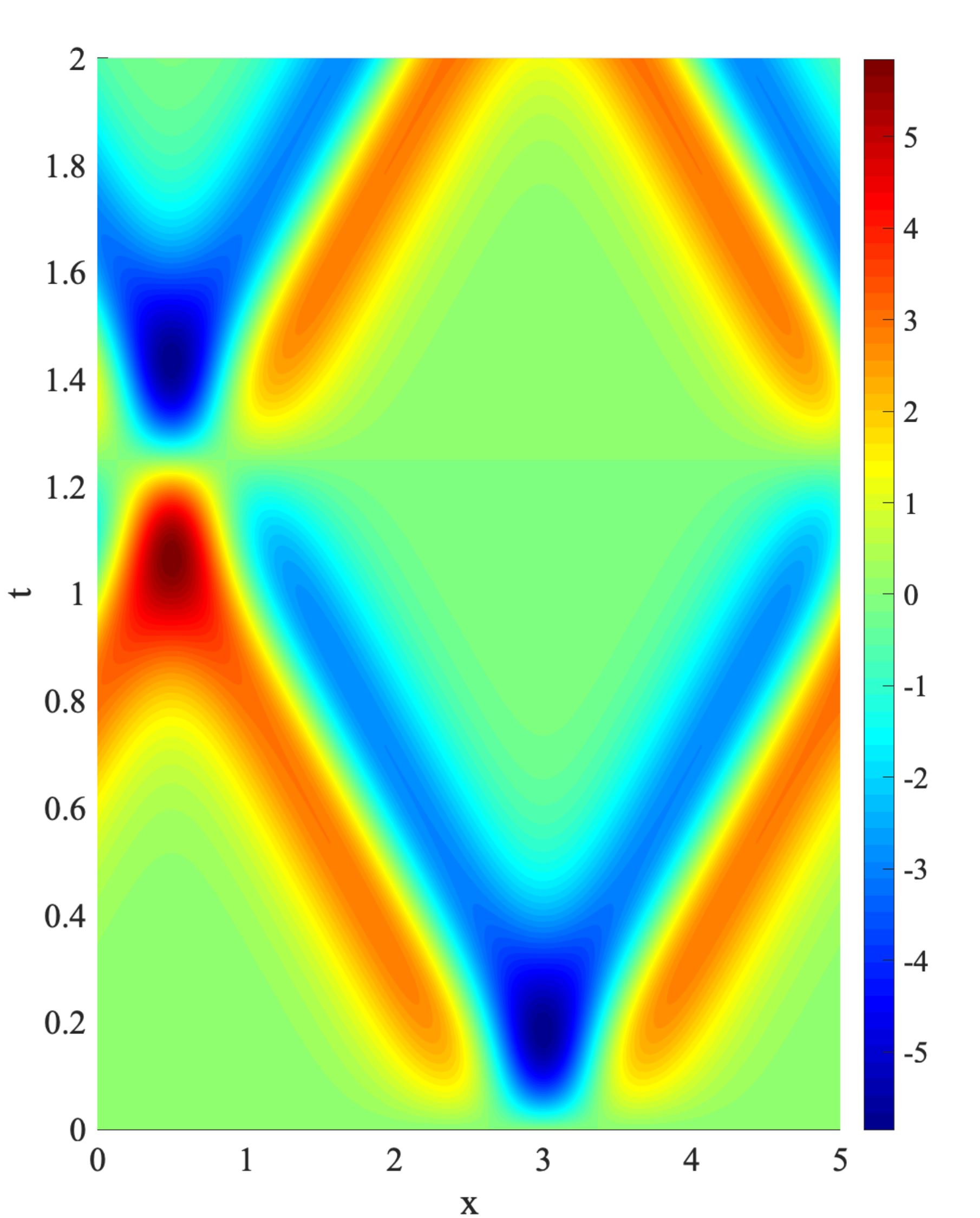}}\hspace{0.2em}
	\subfloat[PINN solution $ v_\theta $]{\includegraphics[width=0.25\linewidth]{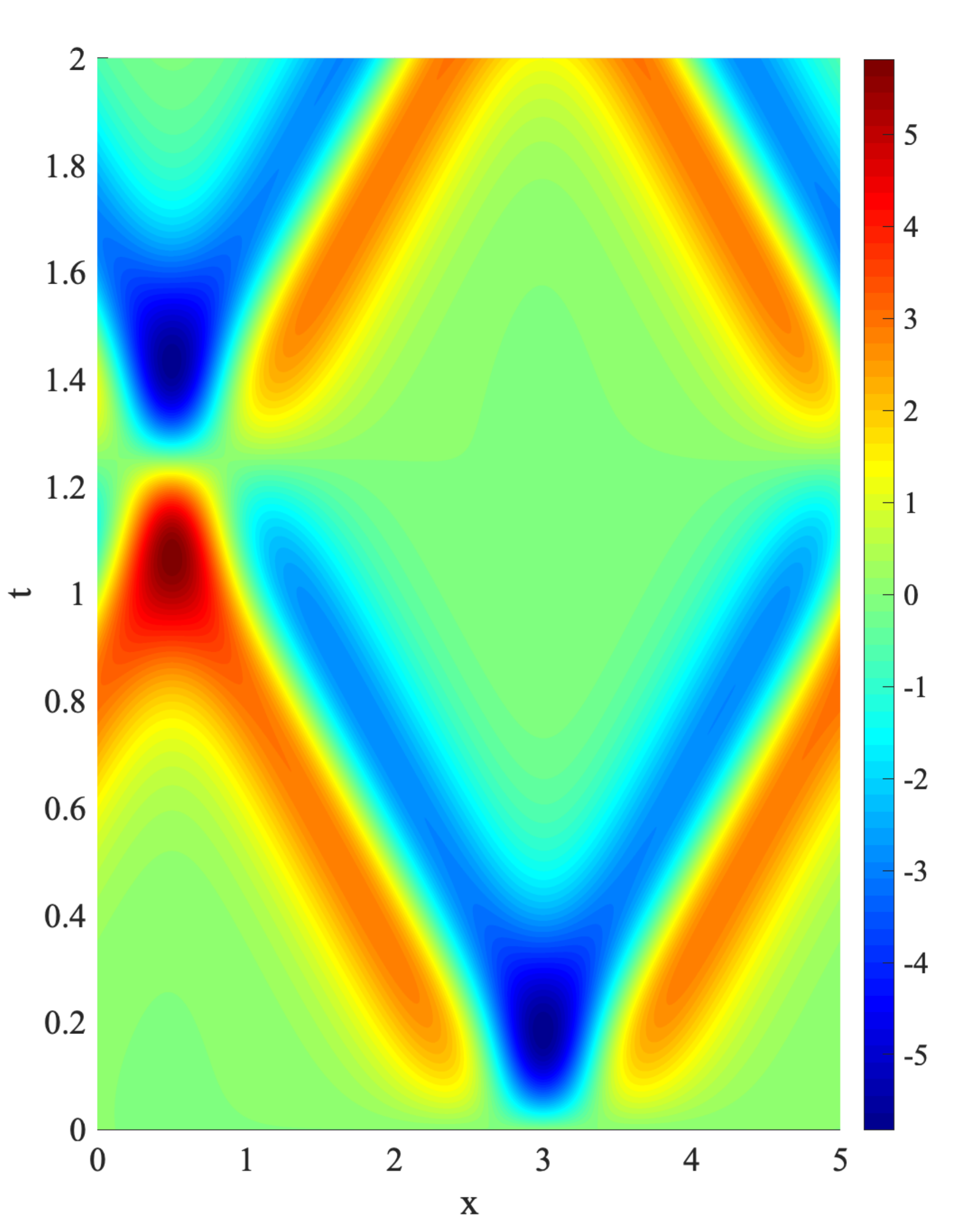}}\hspace{0.2em}
	\subfloat[{PINN} {absolute} error for $v$]{\includegraphics[width=0.25\linewidth]{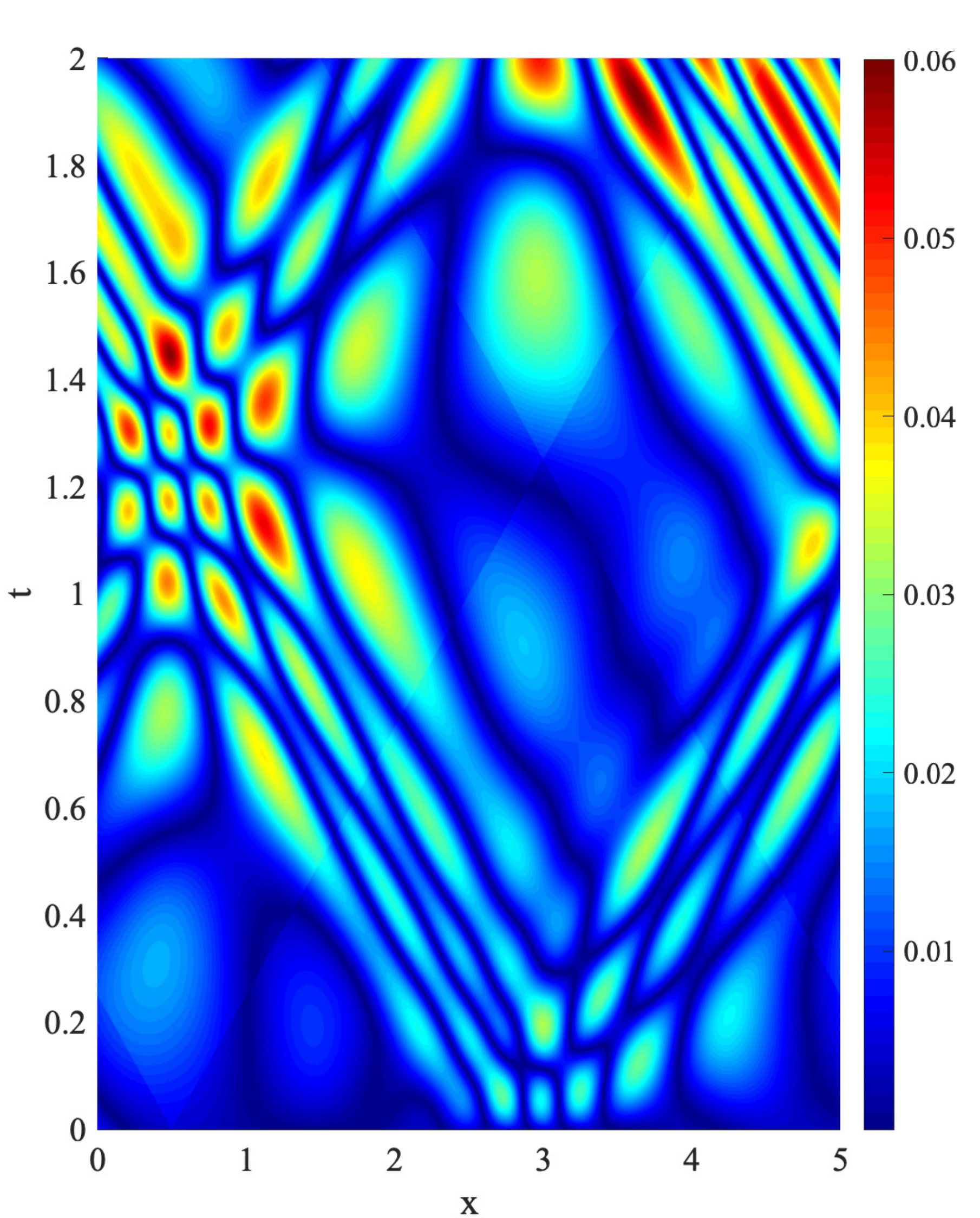}}
 }
	\caption{ Wave equation: Distributions of the True solutions, the PINN solutions and the PINN point-wise absolute errors for $u$ and $v$ in the spatial-temporal domain. $N=2000$ training points within the domain and on each of the domain boundaries.
 }
	\label{PINN_partpaper_Wave_fig1_1}
\end{figure}

We next test the PINN algorithm for solving the wave equation \eqref{wave} in one spatial dimension (plus time), under a configuration in accordance with that of~\cite{2021_JCP_Dong_modifiedbatch}. Consider the spatial-temporal domain, $ (x,t)\in D\times[0, T] = [0, 5] \times [0, 2] $, and the initial-boundary value problem with the wave equation on this domain,
\begin{subequations}\label{part_paper_wave_eq1}
	\begin{align}
		&\frac{\partial^2 u}{\partial t^2} - c^2 \frac{\partial^2 u}{\partial x^2} = 0, \\
		& u(0, t) = u(5, t), \qquad
		\frac{\partial u}{\partial x} (0, t) = \frac{\partial u}{\partial x} (5,t), \\
		& u(x,0) = 2\, {\rm sech}^3 \left(\frac{3}{\delta_0}(x-x_0)\right),\qquad
		\frac{\partial u}{\partial t}(x,0) = 0,
	\end{align}
\end{subequations}
where $u(x,t)$ is the wave field  to be solved for, $ c $ is the wave speed, $ x_0 $ is the initial peak location of the wave, $ \delta_0 $ is a constant that controls the width of the wave profile, and the periodic boundary conditions are imposed on $ x=0 $ and $ 5 $. In the simulations, we employ $ c=2 $, $ \delta_0=2 $, and $ x_0 = 3 $. Then the above problem has the  solution,
\begin{equation*}\left\{
\begin{split}
	&u(x, t) = {\rm sech}^3\left(\frac{3}{\delta_0}\left(-2.5 +\xi\right)\right) + {\rm sech}^3\left(\frac{3}{\delta_0}\left(-2.5 +\eta\right)\right),\\
	& \xi = {\rm mod}\left(x-x_0 + ct + 2.5, 5\right),\quad
	\eta = {\rm mod}\left(x-x_0 - ct + 2.5, 5\right),
 \end{split}
 \right.
\end{equation*}
where mod refers to the modulo operation. The two terms in $ u(x,t) $ represent the leftward- and rightward-traveling waves, respectively. 

We reformulate the problem~\eqref{part_paper_wave_eq1} into the following system,
\begin{subequations}\label{eq_62}
\begin{align}
&
u_t - v=0,\qquad
v_t -c^2u_{xx}=0, \label{eq_62a}
\\
& u(0,t) = u(5,t), \qquad
u_x(0,t) = u_x(5,t), \\
& u(x,0) = 2\, {\rm sech}^3 \left(\frac{3}{\delta_0}(x-x_0)\right),\qquad
		v(x, 0) = 0,
\end{align}
\end{subequations}
where $v(x,t)$ is an auxiliary field given by the first equation in~\eqref{eq_62a}.


To solve the system \eqref{eq_62} with PINN, we employ $ 90 $ and $ 60 $ neurons in the first and the second hidden layers of neural networks, respectively. 
We employ the following loss function in PINN in light of~\eqref{wave_T},
\begin{align}\label{eq_63}
	\text{Loss}
	=& \frac{W_1}{N} \sum_{n=1}^{N}\left[u_{\theta t}(x_{int}^n, t_{int}^n) - v_{\theta}(x_{int}^n, t_{int}^n)\right]^2 \notag \\
 &+ \frac{W_2}{N} \sum_{n=1}^{N} \left[v_{\theta t}(x_{int}^n, t_{int}^n) - u_{\theta xx}(x_{int}^n, t_{int}^n)\right]^2 \nonumber \\
	&+ \frac{W_3}{N} \sum_{n=1}^{N}\left[u_{\theta tx}(x_{int}^n, t_{int}^n) - v_{\theta x}(x_{int}^n, t_{int}^n)\right]^2 \notag \\
 &+ \frac{W_4}{N} \sum_{n=1}^{N}\left[ u_{\theta}(x_{tb}^n, 0) - 2\, {\rm sech}^3 \left(\frac{3}{\delta}_0(x_{tb}^n-x_0)\right)\right]^2 \nonumber \\
	&+ \frac{W_5}{N} \sum_{n=1}^{N} \left[v_{\theta}(x_{tb}^n, 0)\right]^2 +  \frac{W_6}{N} \sum_{n=1}^{N}\left[u_{\theta x}(x_{tb}^n, 0) + \frac{18 \sinh((3x_{tb}^n - 3x_0)/\delta_0)}{\delta_0\cosh^4((3x_{tb}^n- 3x_0)/\delta_0)}\right]^2 \nonumber \\ 
	& + \frac{W_7}{N} \sum_{n=1}^{N} \left[v_{\theta}(0, t_{sb}^n) - v_{\theta}(5, t_{sb}^n)\right]^2 
 + \frac{W_8}{N} \sum_{n=1}^{N} \left[u_{\theta x}(0, t_{sb}^n) - u_{\theta x}(5, t_{sb}^n)\right]^2.
\end{align}
Note that in the simulations we have employed the same number of collocation points ($N$) within the domain and on each of the domain boundaries. The above loss function differs slightly from the one in the error analysis~\eqref{wave_T}, in several aspects.
First, we have added a set of penalty coefficients $W_n>0$ ($1\leq n\leq 8$) for different loss terms in numerical simulations.
Second, the collocation points used in simulations (e.g.~$x_{int}^n$, $t_{int}^n$, $x_{sb}^n$, $t_{sb}^n$, $x_{tb}^n$) are generated randomly within the domain or on the domain boundaries from a uniform distribution. In addition, the averaging used here do not exactly correspond 
to the numerical quadrature rule (mid-point rule) used in the
theoretical analysis.

We have also considered another form (given below) for the loss function, as suggested by an alternate analysis as discussed in Remark~\ref{sec5_Remark1} (see equation~\eqref{wave_TT1}),
\begin{align}\label{eq_64}
	\text{Loss}
	=& \frac{W_1}{N} \sum_{n=1}^{N}\left[u_{\theta t}(x_{int}^n, t_{int}^n) - v_{\theta}(x_{int}^n, t_{int}^n)\right]^2 \notag \\
 &+ \frac{W_2}{N} \sum_{n=1}^{N} \left[v_{\theta t}(x_{int}^n, t_{int}^n) - u_{\theta xx}(x_{int}^n, t_{int}^n)\right]^2 \nonumber \\
	&+ \frac{W_3}{N} \sum_{n=1}^{N}\left[u_{\theta tx}(x_{int}^n, t_{int}^n) - v_{\theta x}(x_{int}^n, t_{int}^n)\right]^2 \notag \\
 &+ \frac{W_4}{N} \sum_{n=1}^{N}\left[u_{\theta}(x_{tb}^n, 0) - 2\, {\rm sech}^3 \left(\frac{3}{\delta}_0(x_{tb}^n-x_0)\right)\right]^2 \nonumber \\
	&+ \frac{W_5}{N} \sum_{n=1}^{N} \left[v_{\theta}(x_{tb}^n, 0)\right]^2 +  \frac{W_6}{N} \sum_{n=1}^{N}\left[u_{\theta x}(x_{tb}^n, 0) + \frac{18 \sinh((3x_{tb}^n - 3x_0)/\delta_0)}{\delta_0\cosh^4((3x_{tb}^n- 3x_0)/\delta_0)}\right]^2 \nonumber \\ 
	& + \frac{W_7}{N} \sum_{n=1}^{N} \left|v_{\theta}(0, t_{sb}^n) - v_{\theta}(5, t_{sb}^n)\right| 
 + \frac{W_8}{N} \sum_{n=1}^{N} \left|u_{\theta x}(0, t_{sb}^n) - u_{\theta x}(5, t_{sb}^n)\right|.
\end{align}
The difference between this form and the form~\eqref{eq_63} lies in the last two terms, with the terms here not squared.

The loss function~\eqref{eq_63} will be referred to as the loss form \#1 in subsequent discussions, and~\eqref{eq_64} will be referred to as the loss form \#2.
The PINN schemes that employ these two different loss forms will be referred to as
PINN-F1 and PINN-F2, respectively.

Figure \ref{PINN_partpaper_Wave_fig1_1} shows distributions of the exact solutions, the PINN solutions, and the PINN point-wise {absolute} errors for $u$ and $v$ in the spatial-temporal domain. 
Here the PINN solution is computed by PINN-F1, 
in which penalty coefficients are given by
$\bm{W}=(W_1,\dots,W_8)= (0.8,0.8,0.8,0.5,0.5,0.5,0.9,0.9) $.
One can observe that the  method has captured the wave fields for $u$ and $\frac{\partial u}{\partial t}$ reasonably well, with the error for $u$ notably smaller than that of $\frac{\partial u}{\partial t}$.

\begin{figure}[tb]
	\centering
	\subfloat[$ t=0.5 $]{
		\begin{minipage}[b]{0.25\textwidth}
			\includegraphics[scale=0.25]{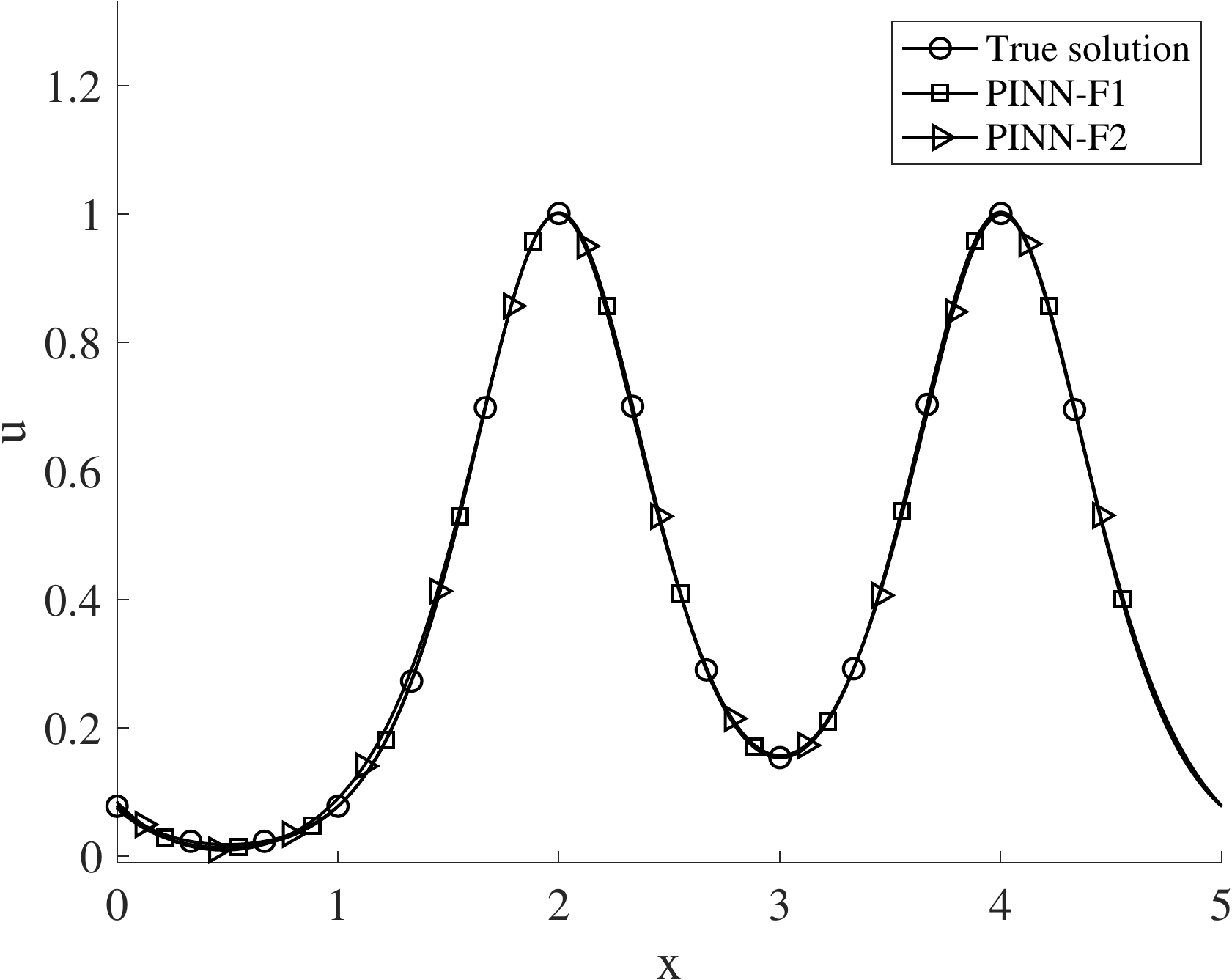}\\
			\includegraphics[scale=0.25]{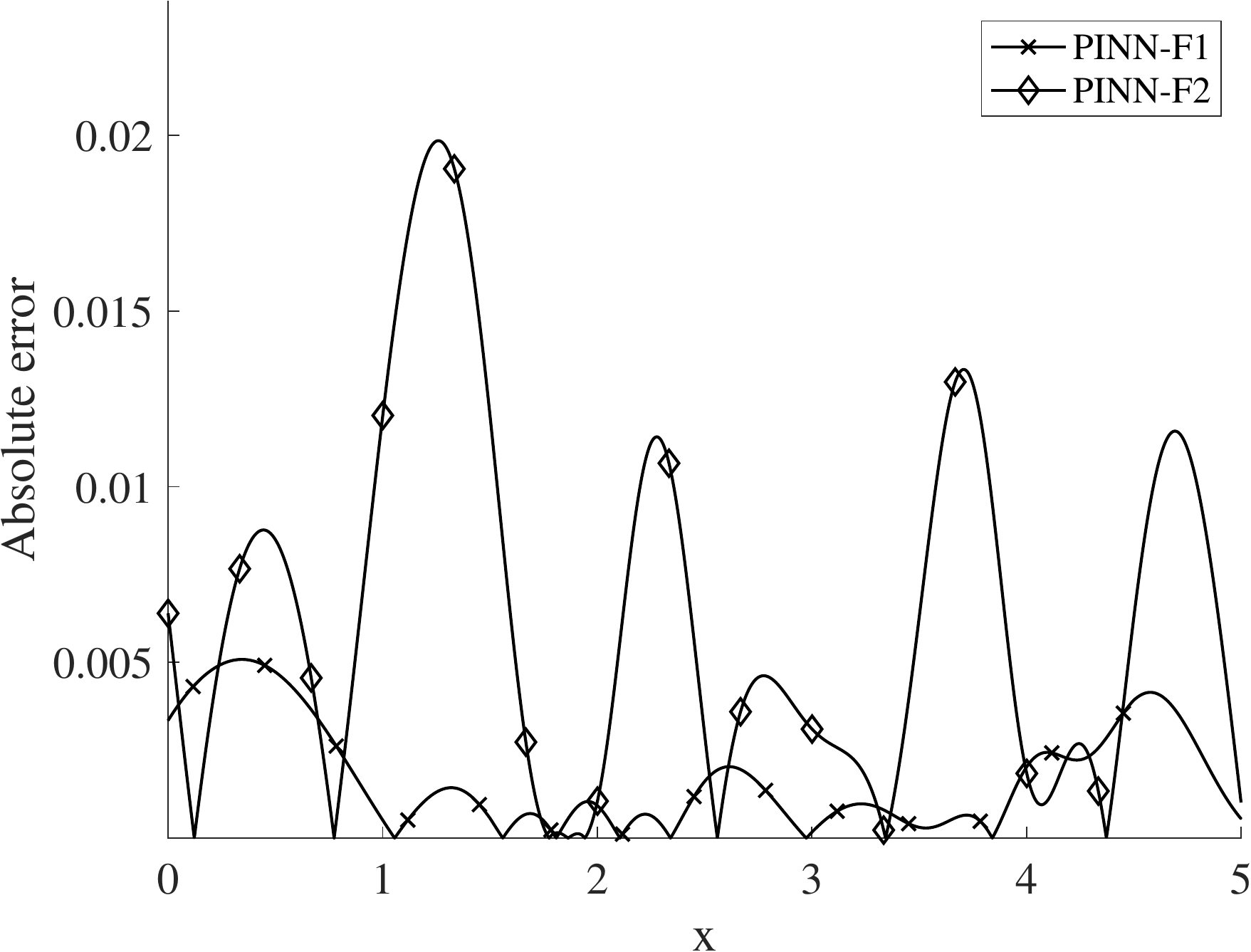}
		\end{minipage}
  	}\qquad
	\subfloat[$ t=1 $]{
		\begin{minipage}[b]{0.25\textwidth}
			\includegraphics[scale=0.25]{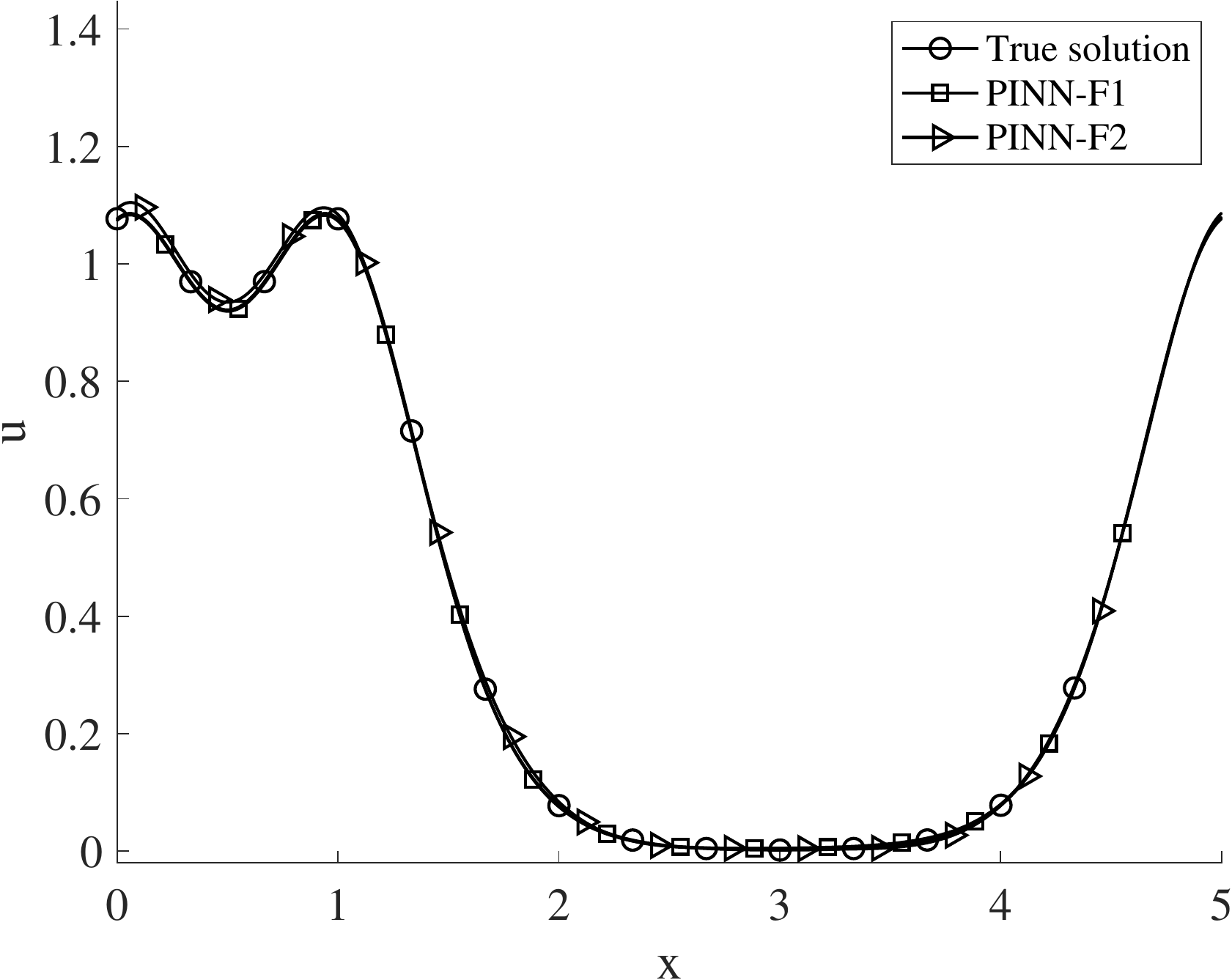}\\
			\includegraphics[scale=0.25]{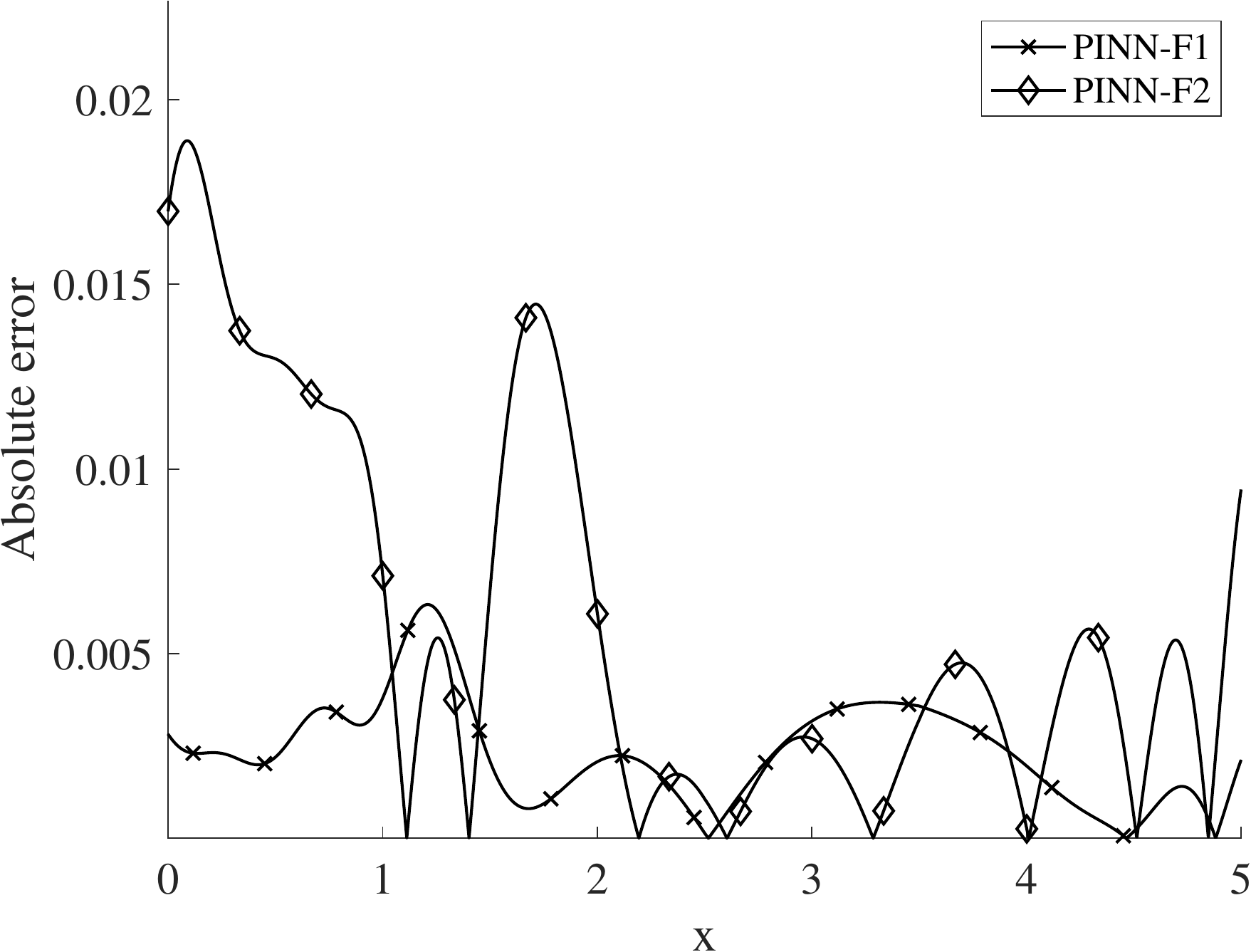}
		\end{minipage}
	}\qquad
	\subfloat[$ t=1.5 $]{
		\begin{minipage}[b]{0.25\textwidth}
			\includegraphics[scale=0.25]{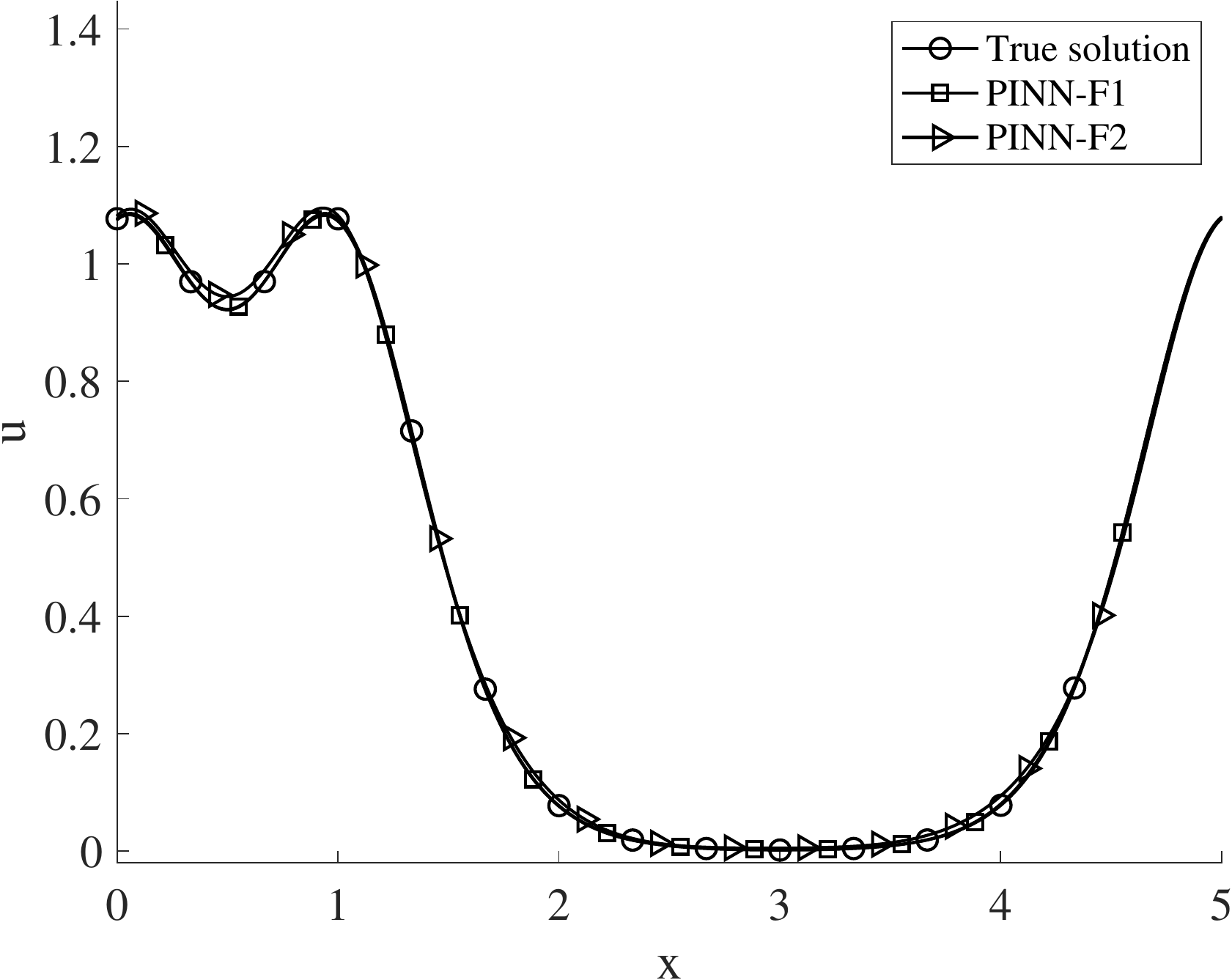}\\
			\includegraphics[scale=0.25]{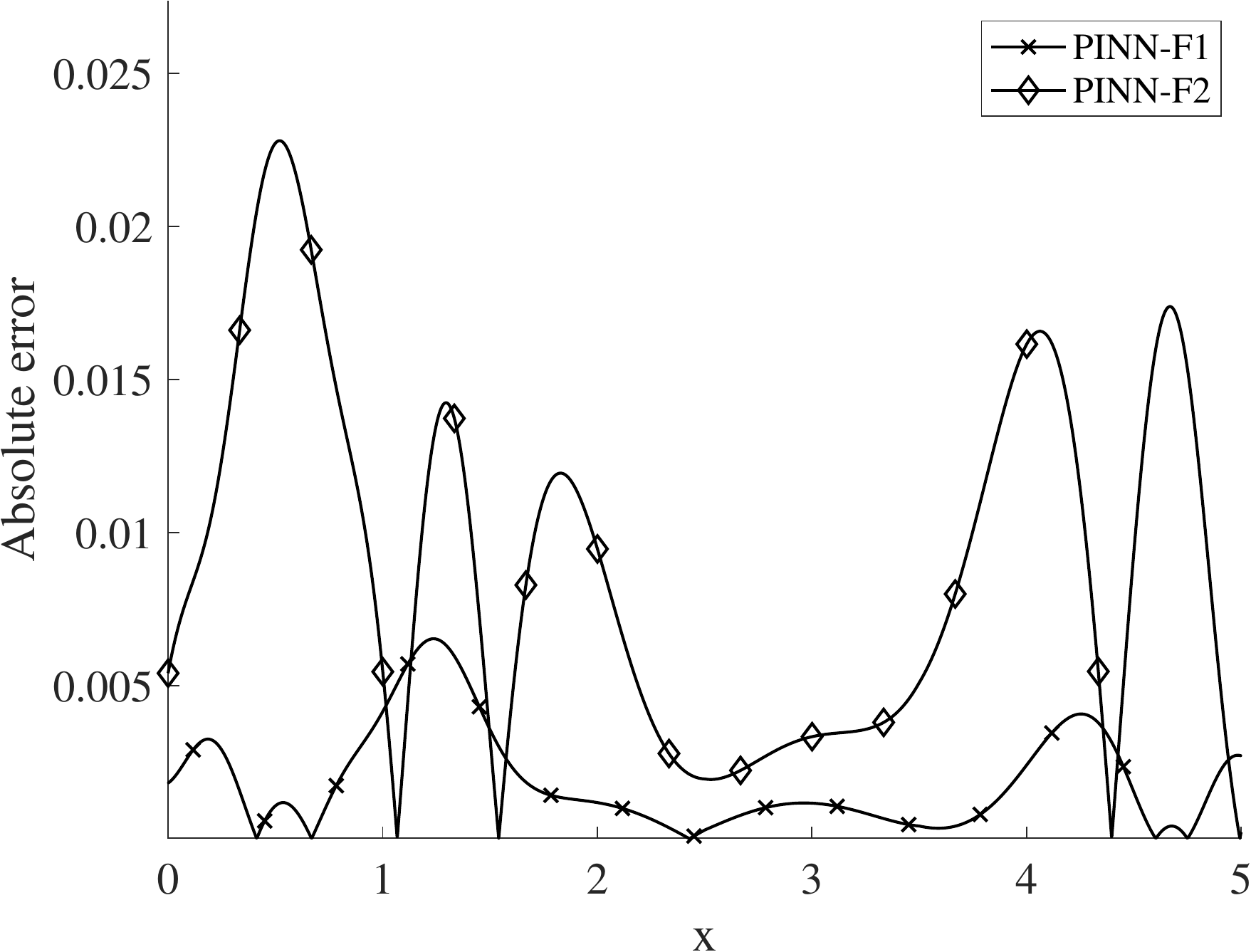}
		\end{minipage}
	}
	\caption{Wave equation: Comparison of profiles of $u$ (top row) and its absolute error (bottom row) between the PINN solutions (loss forms \#1 and \#2) and the exact solution at time instants (a) $t=0.5$, (b) $t=1.0$, and (c) $t=1.5$. 
 $N=2000$ training data points within the domain and on each of the domain boundaries ($x=0$ and $5$, and $t=0$).
 }
	\label{fg_2}
\end{figure}

\begin{figure}[tb]
	\centering
	\subfloat[$ t=0.5 $]{
		\begin{minipage}[b]{0.25\textwidth}
			\includegraphics[scale=0.25]{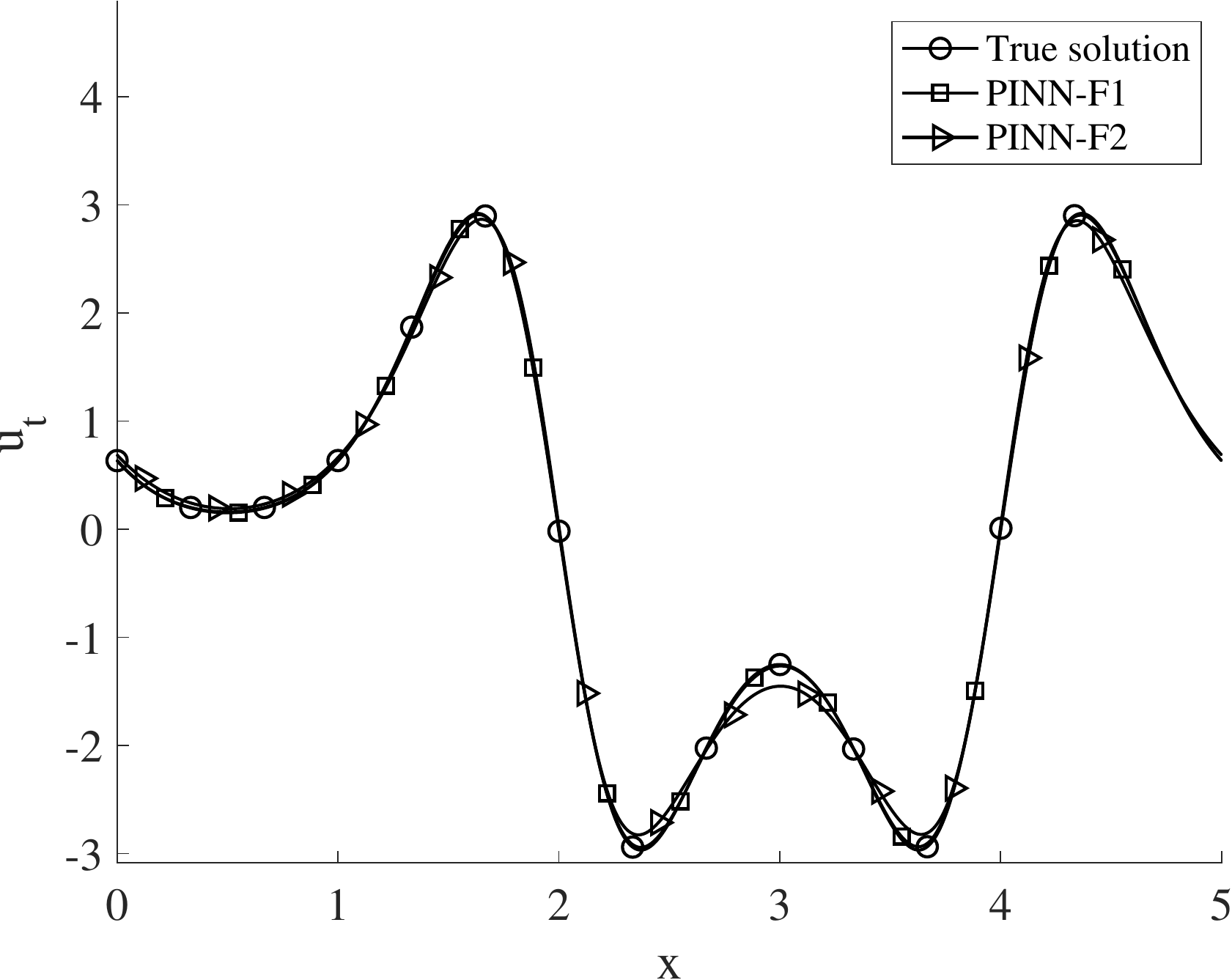}\\
			\includegraphics[scale=0.25]{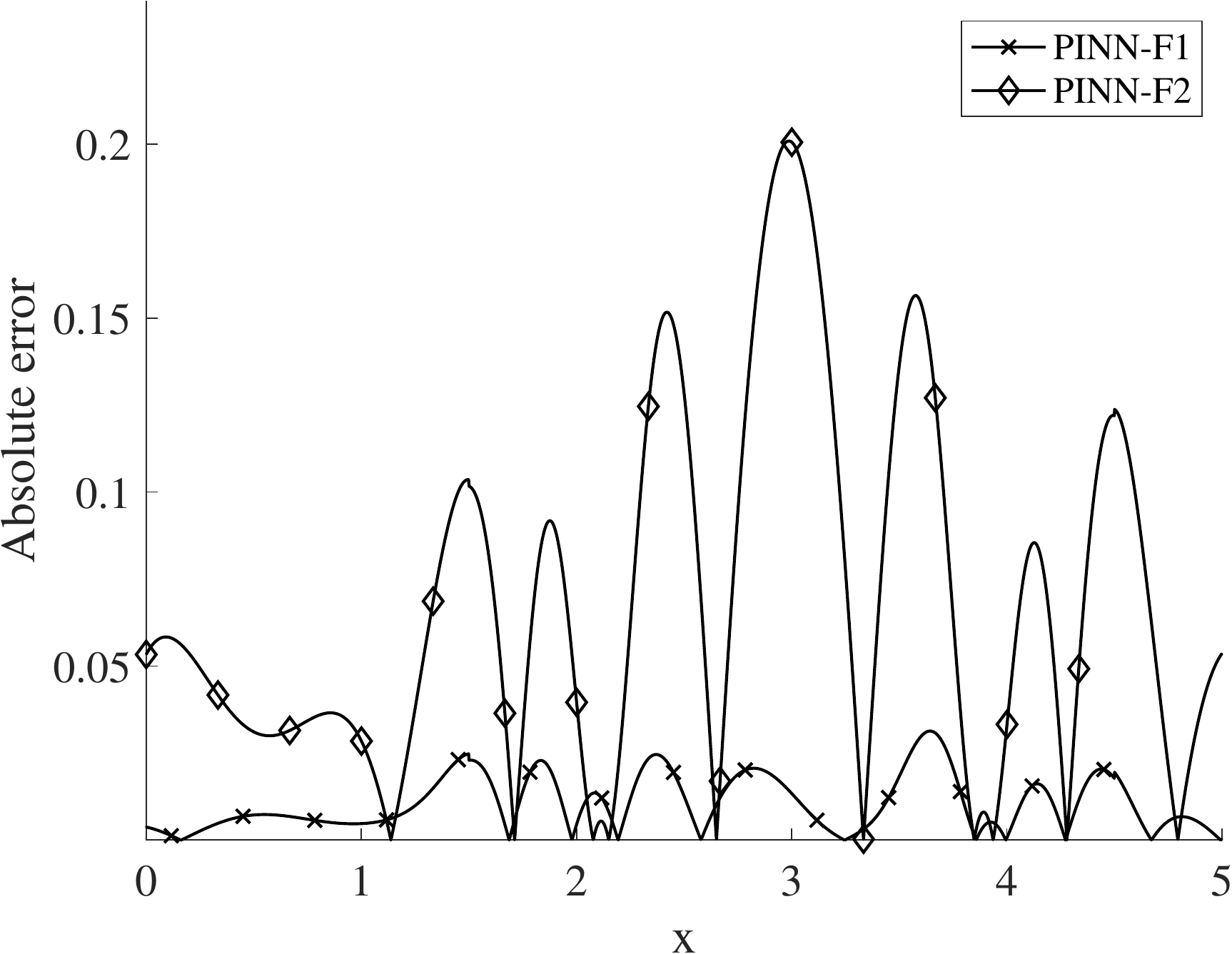}
		\end{minipage}
	}\qquad
	\subfloat[$ t=1 $]{
		\begin{minipage}[b]{0.25\textwidth}
			\includegraphics[scale=0.25]{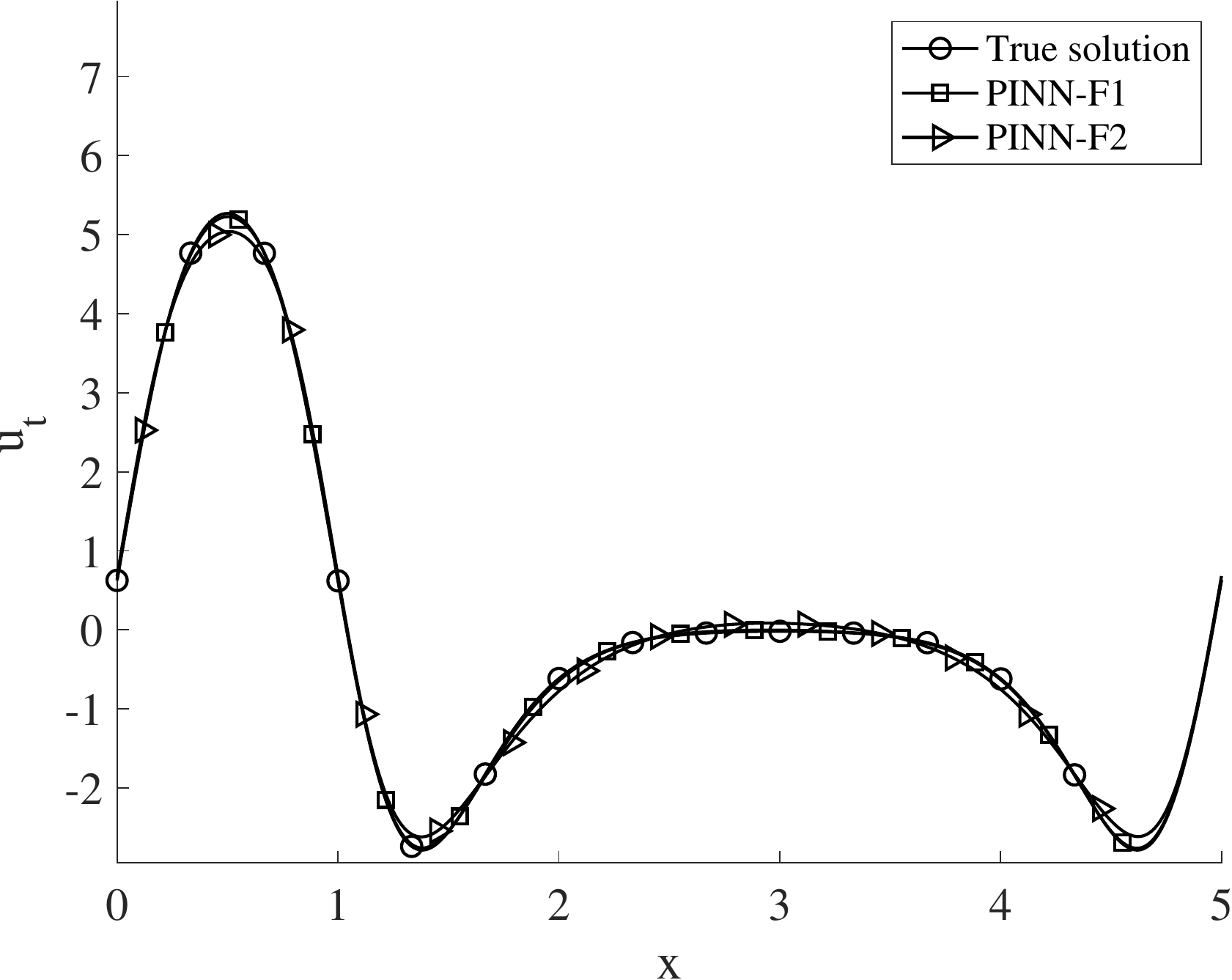}\\
			\includegraphics[scale=0.25]{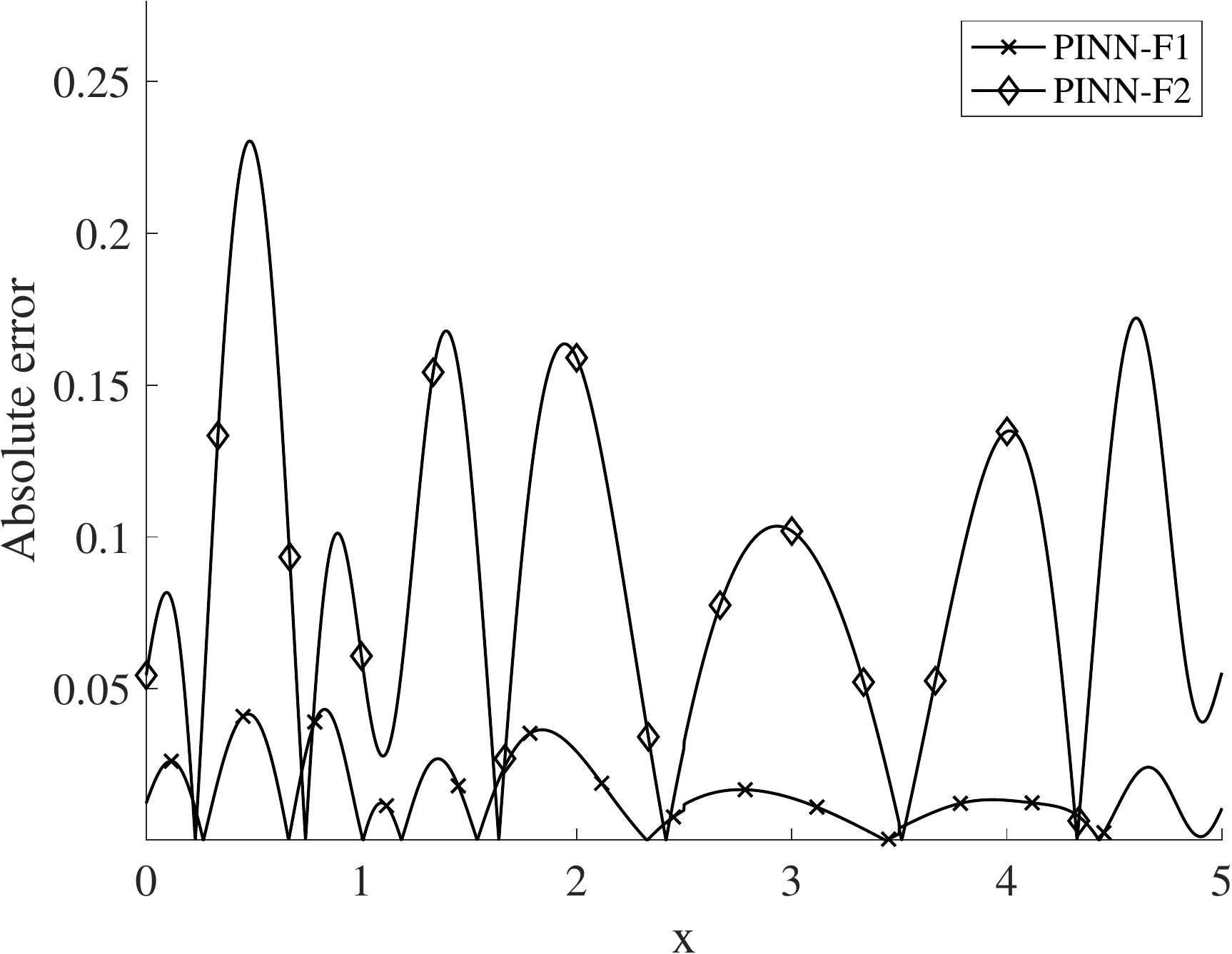}
		\end{minipage}
	}\qquad
	\subfloat[$ t=1.5 $]{
		\begin{minipage}[b]{0.25\textwidth}
			\includegraphics[scale=0.25]{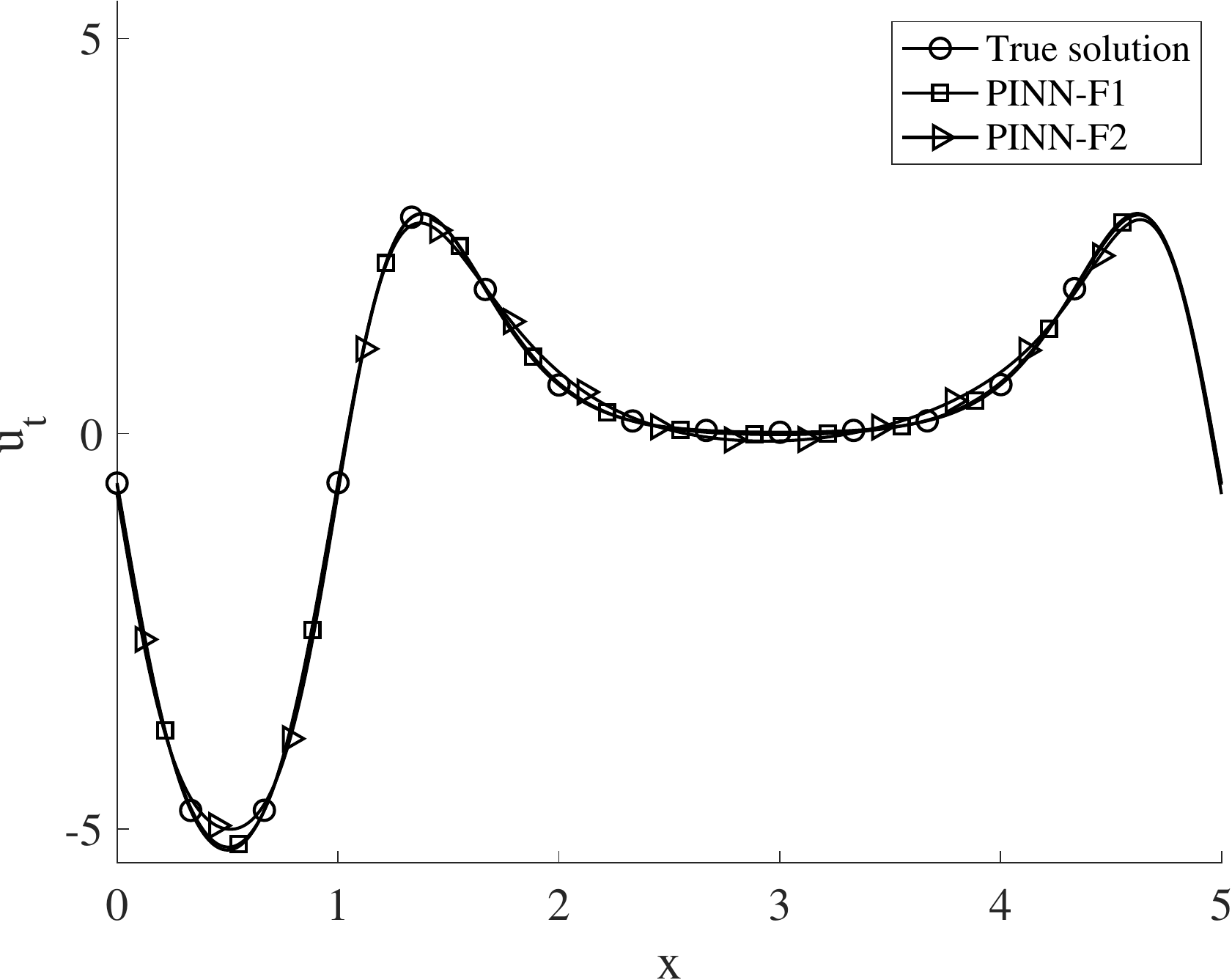}\\
			\includegraphics[scale=0.25]{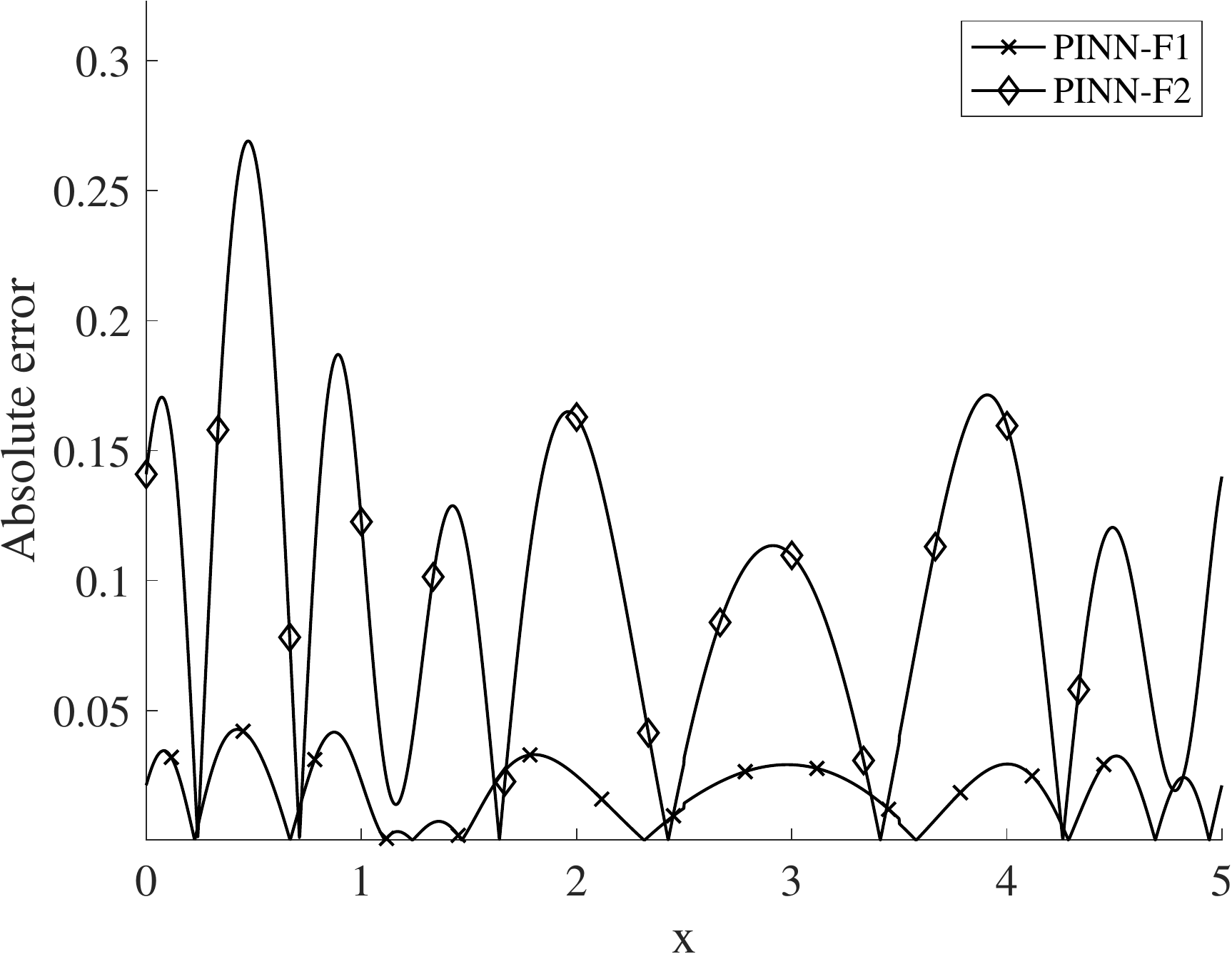}
		\end{minipage}
	}
	\caption{Wave equation: Comparison of the profiles of $v=\frac{\partial u}{\partial t}$ (top row) and its absolute error (bottom row) between the PINN solutions (loss forms \#1 and \#2) and the exact solution at time instants (a) $t=0.5$, (b) $t=1.0$, and (c) $t=1.5$. 
 $N=2000$ training data points within the domain and on each of the domain boundaries ($x=0$ and $5$, and $t=0$).
 }
	\label{fg_3}
\end{figure}

\begin{figure}[tb]
		\centering
  		\subfloat[PINN-F1]{\includegraphics[width=0.4\textwidth]{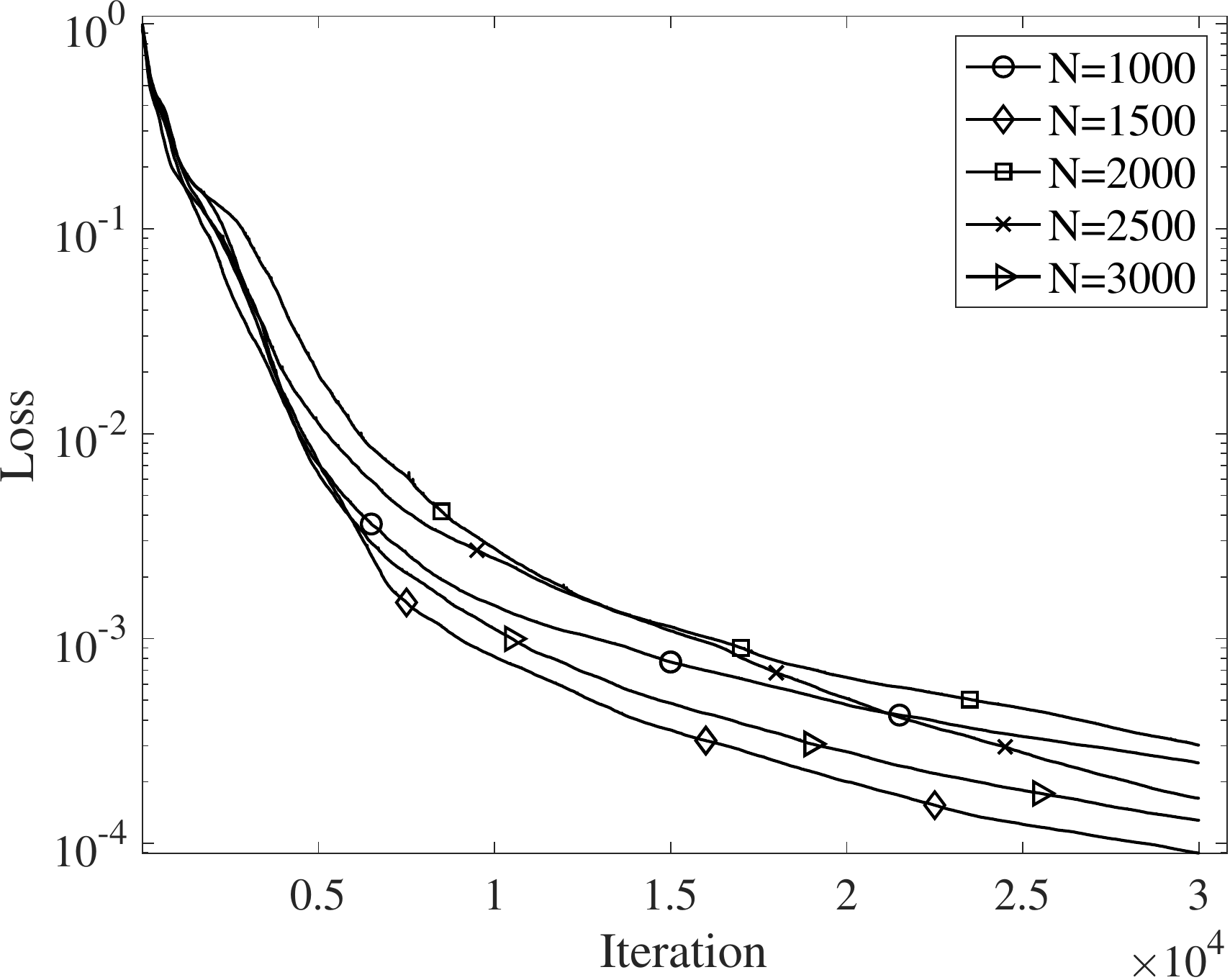}}
    \qquad
		\subfloat[PINN-F2]{\includegraphics[width=0.4\textwidth]{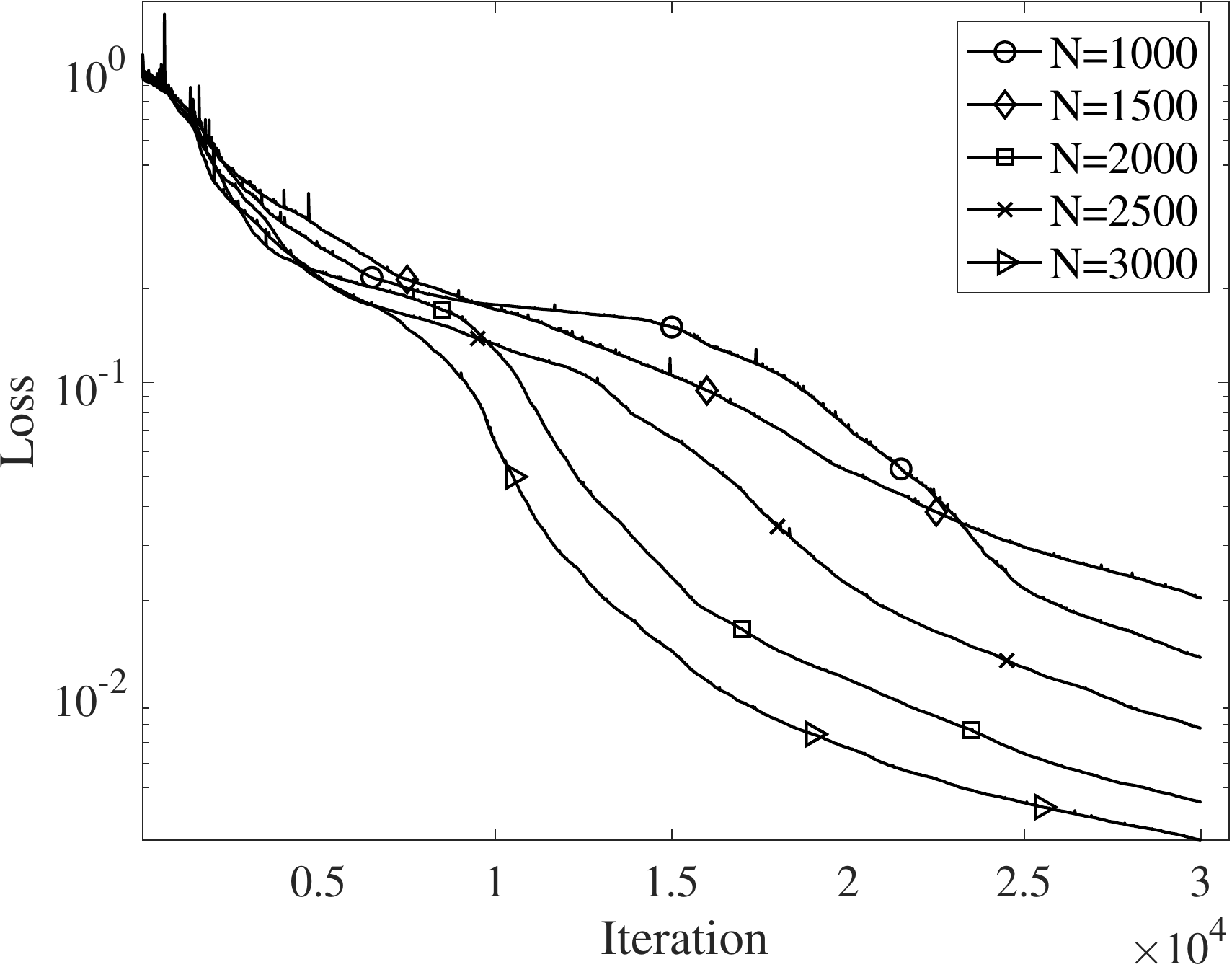}}
	\caption{Wave equation: Histories of the loss function
 versus the training iteration with PINN-F1 and PINN-F2, corresponding to different number of
 training data points ($N$). 
 }
	\label{fg_4}
\end{figure}

Figures~\ref{fg_2} and~\ref{fg_3} provide a comparison of the  solutions obtained 
using the two forms of loss functions.
Figure~\ref{fg_2} compares profiles of the PINN-F1 and PINN-F2 solutions, and the exact solution, for $u$ (top row) at three time instants ($t=0.5$, $1.0$, and $1.5$), as well as the error profiles (bottom row).
Figure~\ref{fg_3} shows the corresponding results for the field variable $v=\frac{\partial u}{\partial t}$.
These results are obtained by using $N=2000$ training data points in the domain and on each of the domain boundaries.
It is observed that both PINN schemes, with the loss functions given by~\eqref{eq_63} and~\eqref{eq_64} respectively, have captured the solution reasonably well. We further observe that the PINN-F1 scheme (with the loss form~\eqref{eq_63}) produces notably more accurate results than the PINN-F2 (with loss form~\eqref{eq_64}),
especially for the field $\frac{\partial u}{\partial t}$.

We have varied the number of training data points $N$ systematically and studied its effect on the PINN results.
Figure~\ref{fg_4} shows the loss histories of PINN-F1 and PINN-F2 corresponding to different number of training data points ($N$) in
the simulations, with a total of $30,000$ training iterations.
We can make two observations. First, the history curves with the loss function form \#1 is generally smoother, indicating that the loss function decreases almost monotonically as the training progresses.
On the other hand, significant fluctuations in
the loss history can be observed with the form \#2. Second, the eventual loss values produced by the loss form \#1 are significantly smaller,
by over an order of magnitude, than those produced by the loss form \#2.

Table \ref{tab_1} is a further comparison between the PINN-F1 and PINN-F2. Here the $l_2$ and $l_{\infty}$ errors of $u$ and $v$ computed by PINN-F1 and PINN-F2 corresponding to different training data points ($N$) have been listed.
There appears to be a general trend that 
the errors tend to decrease with increasing number of training points, but the decrease is not monotonic. 
It can be observed that the $u$ errors are notably smaller than those for $v=\frac{\partial u}{\partial t}$, as signified earlier in e.g.~Figure~\ref{PINN_partpaper_Wave_fig1_1}.
One can again observe that PINN-F1 results are notably more accurate than those of PINN-F2 for the wave equation.


\begin{table}[tb]
	\caption{Wave equation: The $u$ and $v$ errors versus the number of training data points $N$.
 }
	\label{tab_1}
	\centering
	\begin{tabular}[b]{ c | c|  c   c  | c   c   }
		\hline
  \multirow{2}{*}{method} &
		\multirow{2}{*}{$ N $}&\multicolumn{2}{c}{$ l_2 $-error} & \multicolumn{2}{c}{$ l_\infty $-error} \\ 
		\cline{3-6}
		& &$ u_\theta $  &$ v_\theta $ &$ u_\theta $  &$ v_\theta $ \\
		\hline
  & 1000&  5.7013e-03&  1.3531e-02&   1.8821e-02&  4.6631e-02\\
		\cline{3-6}
	&	1500&  2.1689e-03&   4.1035e-03&  6.7631e-03&   1.5109e-02\\
		\cline{3-6}
 PINN-F1	&	2000&  4.6896e-03&   9.6417e-03&  1.3828e-02&   3.3063e-02\\
		\cline{3-6}
	&	2500&  3.7879e-03&   9.8574e-03&  1.2868e-02&   3.3622e-02\\
		\cline{3-6}
	&	3000&  2.6588e-03&   6.0746e-03&  8.1457e-03&   1.9860e-02\\
  \hline
		& 1000&  4.7281e-02&  9.2431e-02&  1.4367e-01&  3.2764e-01\\
		\cline{3-6}
		& 1500&  4.9087e-02&  1.2438e-01&  2.1525e-01&  5.0601e-01\\
		\cline{3-6}
	PINN-F2	& 2000&  1.8554e-02&  4.9224e-02&  6.0780e-02&  1.6358e-01\\
		\cline{3-6}
		& 2500&  2.3526e-02&  5.4266e-02&  9.8690e-02&  1.9467e-01\\
		\cline{3-6}
		& 3000&  1.4164e-02&  3.7796e-02&  5.3045e-02&  1.4179e-01\\
		\hline
	\end{tabular}
\end{table}

\begin{figure}
	\centering
	\subfloat[PINN-F1]{\includegraphics[width=0.4\linewidth]{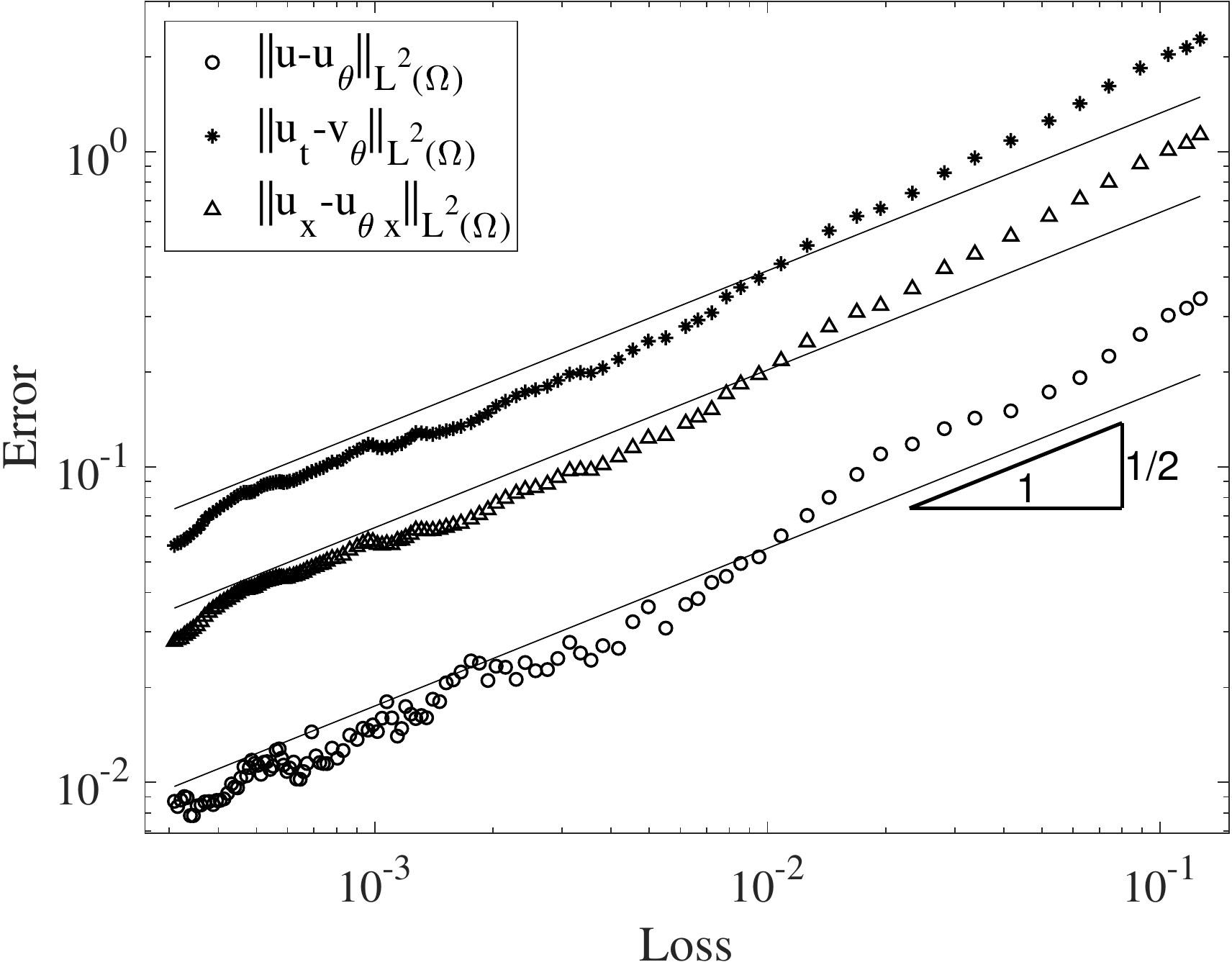}}
 \qquad
	\subfloat[PINN-F2]{\includegraphics[width=0.4\linewidth]{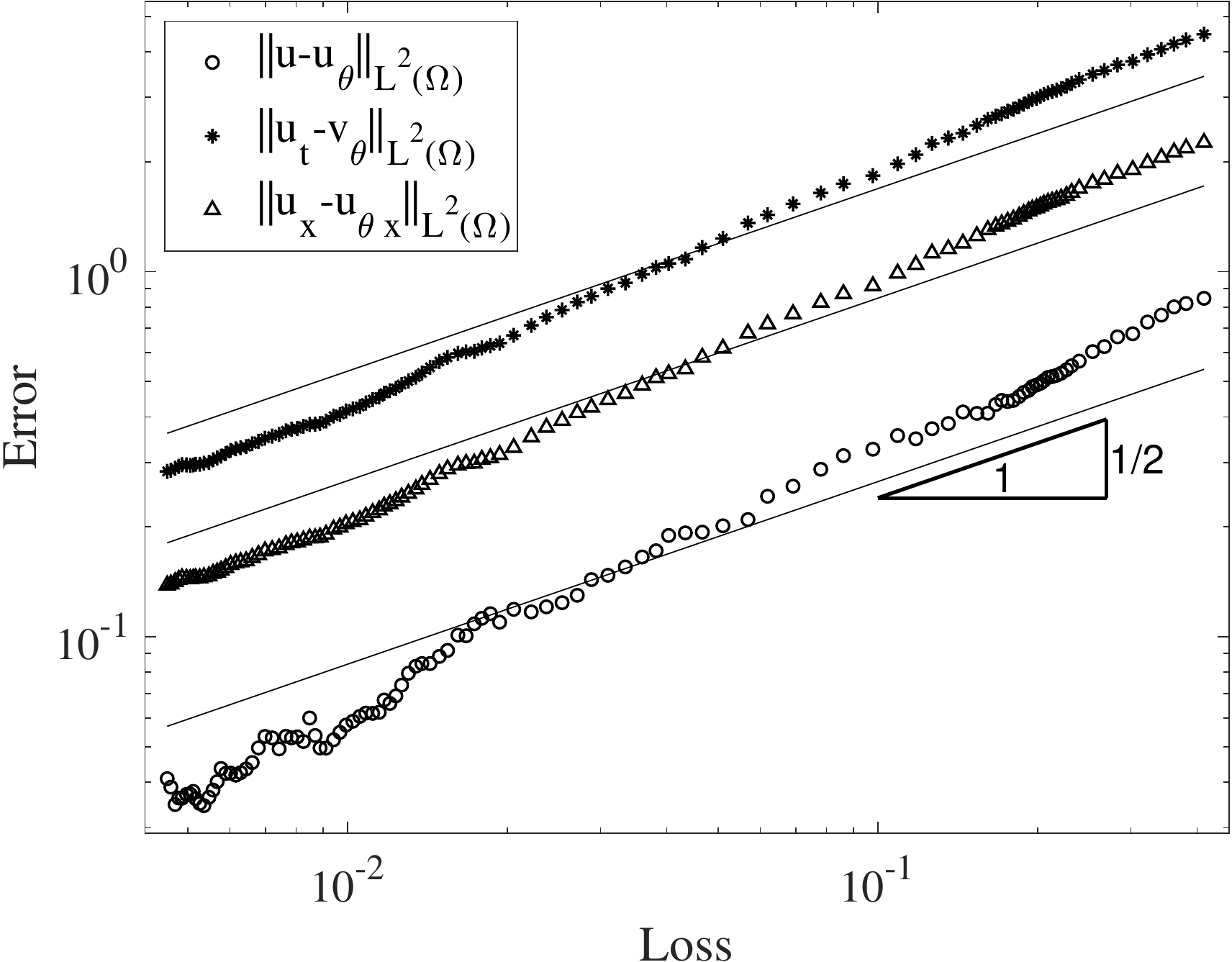}}
	\caption{Wave equation: The $l^2$ errors of $u$, $\frac{\partial u}{\partial t}$, and $\frac{\partial u}{\partial x}$ as a function of the training loss value.
$N=2000$ training data points.
 }
	\label{fg_5}
\end{figure}

Theorem \ref{sec5_Theorem3} suggests the solution errors for $u$, $v=\frac{\partial u}{\partial t}$, and $\nabla u$ approximately scale as the square root of
the training loss function. Figure~\ref{fg_5} provides some numerical evidence for this point. Here we plot the $l^2$ errors for $u$, $\frac{\partial u}{\partial t}$ and $\frac{\partial u}{\partial x}$ from our simulations as a function of the training loss value for PINN-F1 and PINN-F2 in logarithmic scales. It is evident that for PINN-F1 the scaling essentially follows the square root relation. For PINN-F2 the relation between the error and the training loss appears to scale with a power somewhat larger than $\frac12$.

\subsection{Sine-Gordon Equation}

\begin{figure}[tb]
	\centering
	\subfloat[True solution for $u$]{\includegraphics[width=0.2\linewidth]{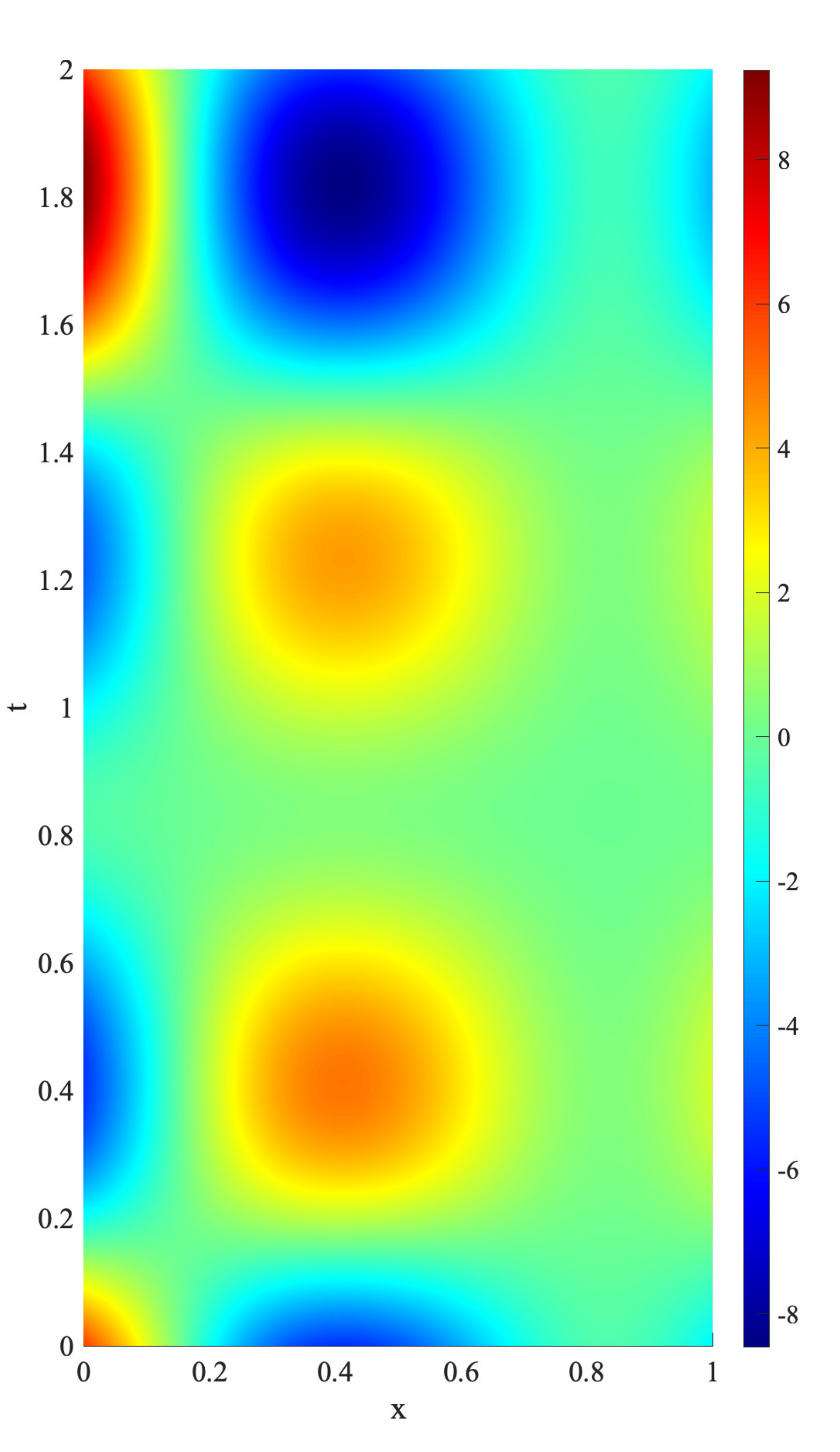}}\hspace{0.2em}
	\subfloat[PINN solution for $u$]{\includegraphics[width=0.2\linewidth]{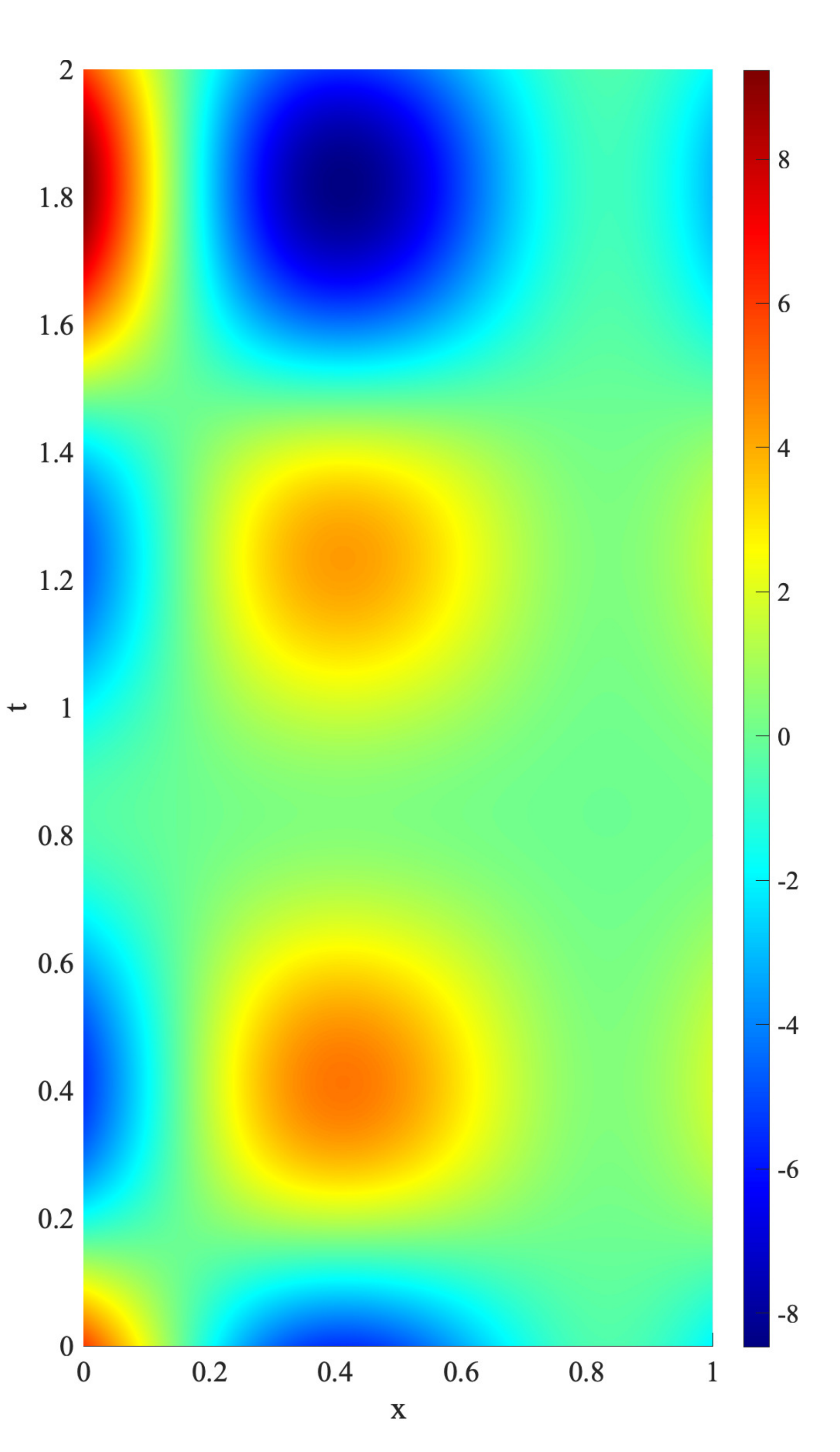}}\hspace{0.2em}
	\subfloat[Solution error for $u$]{\includegraphics[width=0.2\linewidth]{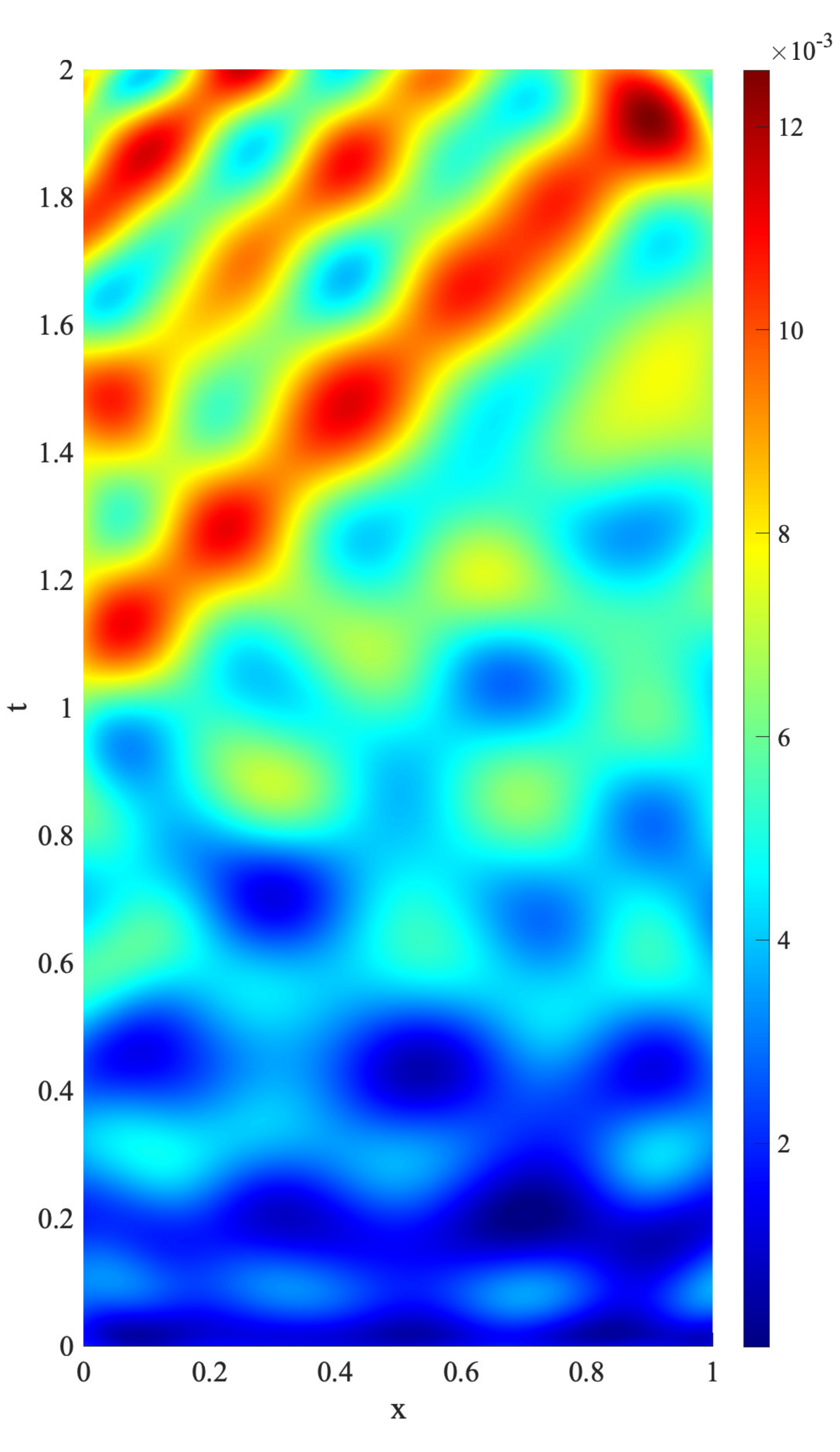}} \\
 	\subfloat[True solution for $v$]{\includegraphics[width=0.2\linewidth]{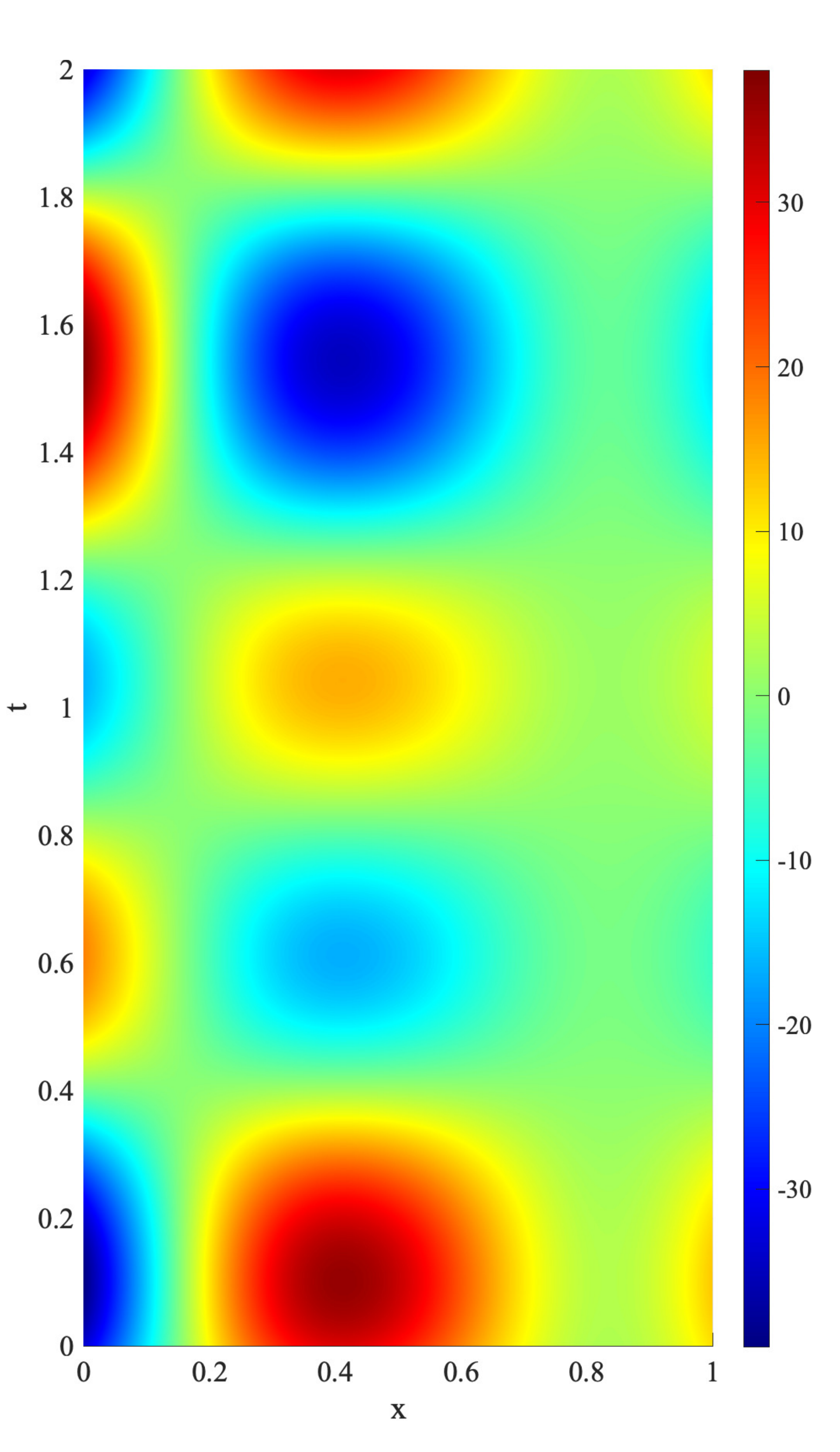}}\hspace{0.2em}
	\subfloat[PINN solution for $v$]{\includegraphics[width=0.2\linewidth]{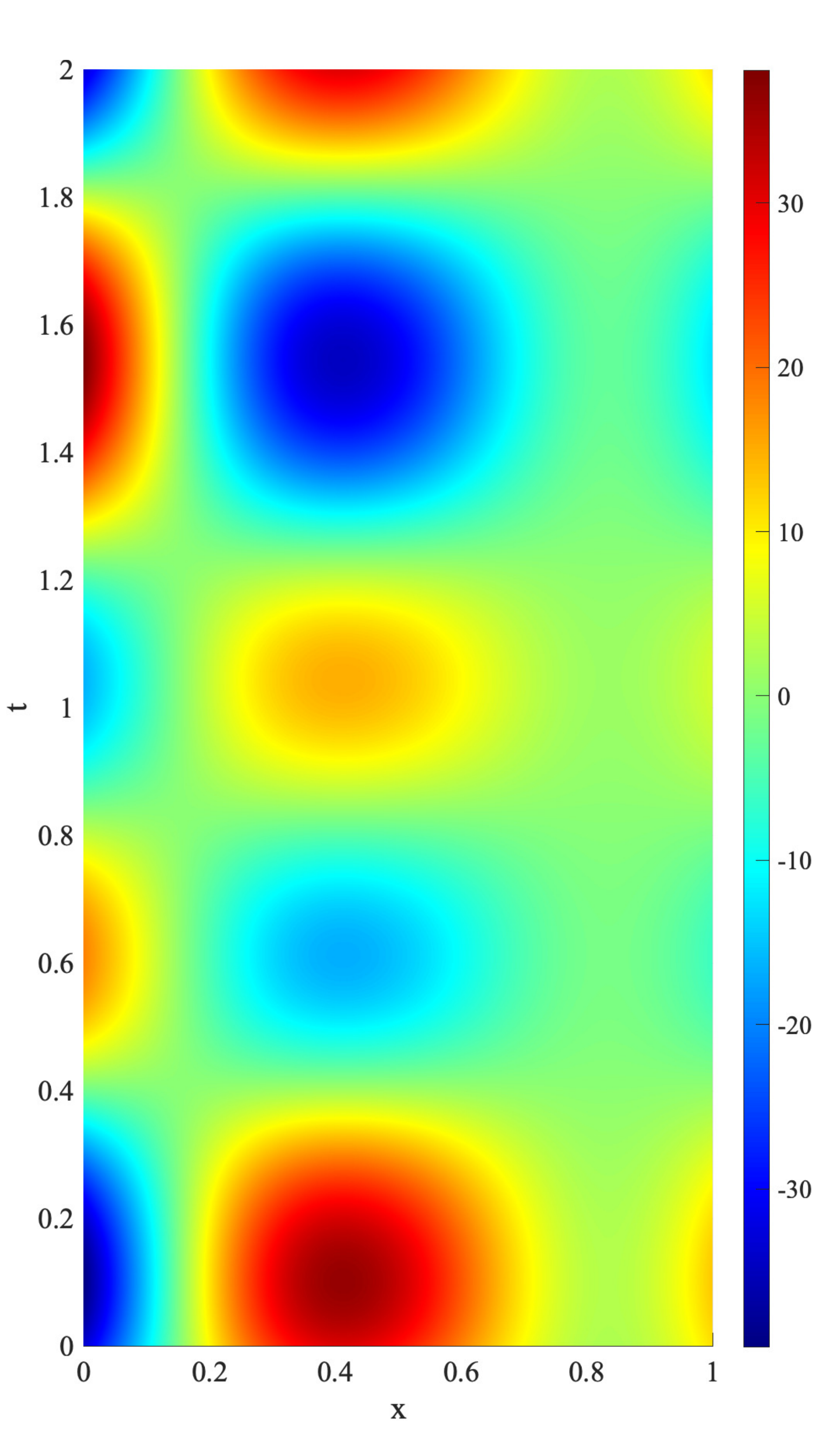}}\hspace{0.2em}
	\subfloat[Solution error for $v$]{\includegraphics[width=0.2\linewidth]{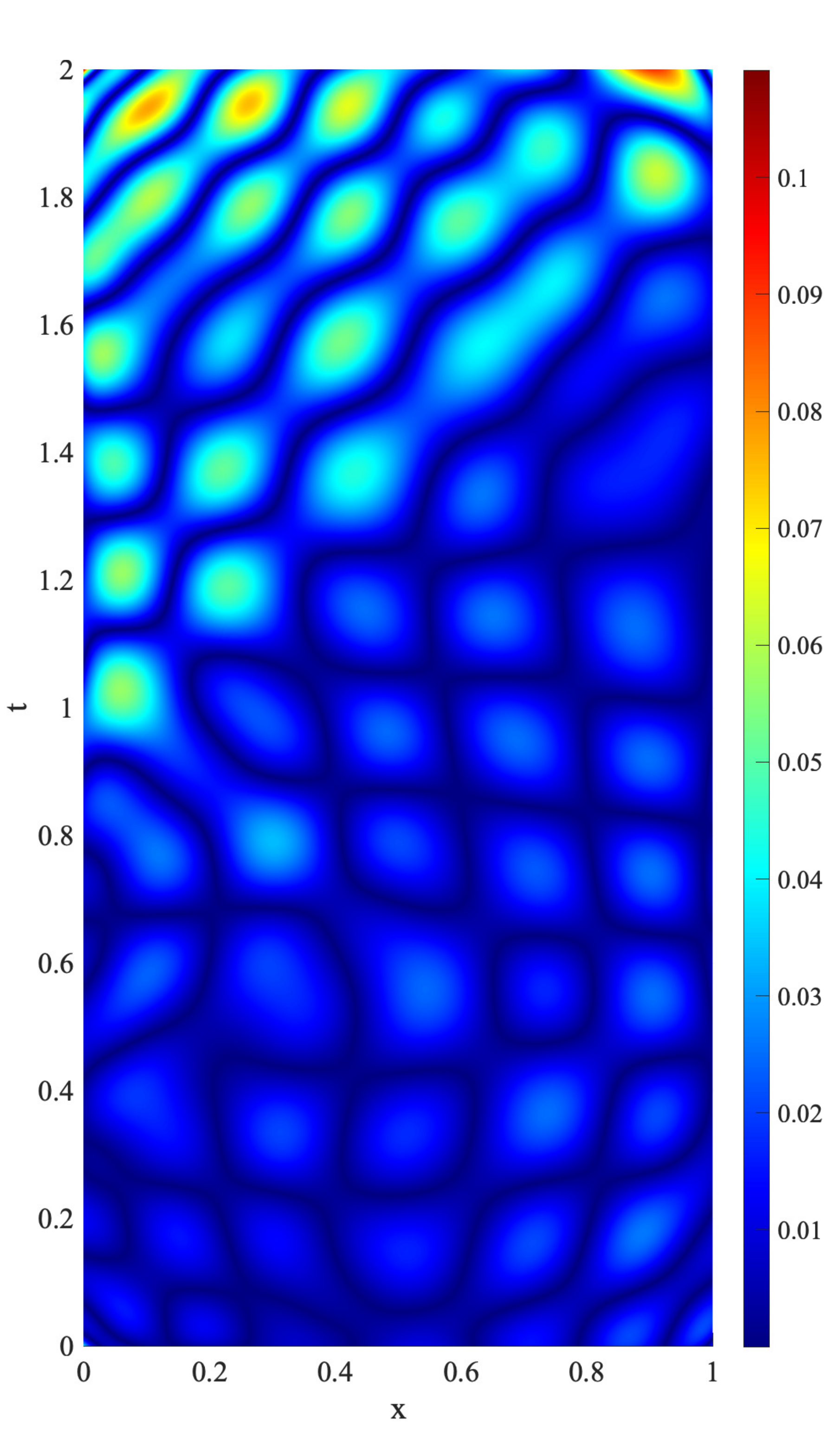}}
	\caption{ Sine-Gordon equation:
 Distributions of the exact solution (left column), the PINN solution (middle column) and the PINN absolute error (right column) for $u$ (top row) and for $v=\frac{\partial u}{\partial t}$ (bottom row).
 $N=2000$ collocation points within the domain and on the domain boundaries.
 }
 \label{fg_6}
 \end{figure}

We test the PINN algorithm suggested by the theoretical analysis for the Sine-Gordon equation~\eqref{SG} in this subsection. Consider the spatial-temporal domain
$(x,t)\in\Omega= D\times [0, T] = [0, 1] \times [0, 2]$, and the following initial/boundary value problem on this domain,
\begin{subequations}\label{eq_65}
\begin{align}
	& \frac{\partial^2u}{\partial t^2} - \frac{\partial^2 u}{\partial x^2} + u + \sin(u) = f(x,t), \\
 &u({0},t)=\phi_1(t),  \qquad u({1},t)=\phi_2(t), \\
	&u({x},0)=\psi_{1}({x}),\qquad \frac{\partial u}{\partial t}({x},0)=\psi_{2}({x}).
\end{align}
\end{subequations}
In these equations, $ u(x,t) $ is the field function to be solved for, $ f(x,t) $ is a source term, $ \psi_{1} $ and $ \psi_2 $ are the initial conditions, and $ \phi_1 $ and $ \phi_2 $ are the boundary conditions.
The source term, initial and boundary conditions appropriately are chosen by the following exact solution,  
\begin{align}\label{PINN_SG_eq1}
	u(x, t) = \left[2\cos\left(\pi x + \frac{\pi}{5}\right) + \frac{9}{5}\cos\left(2\pi x + \frac{7\pi}{20}\right)\right]\left[ 2\cos\left(\pi t + \frac{\pi}{5}\right) + \frac{9}{5}\cos\left(2\pi t + \frac{7\pi}{20}\right) \right].
\end{align}

To simulate this problem with PINN, we reformulate the problem as follows,
\begin{subequations}\label{eq_67}
\begin{align}
	&u_{t} - v = 0, \label{eq_67a} \\ 
 &  v_{t} - u_{xx} + u + \sin(u) = f(x,t), \\
 &u({0},t)=\phi_1(t),  \qquad u({1},t)=\phi_2(t),\\
	&u({x},0)=\psi_{1}({x}),\qquad v({x},0)=\psi_{2}({x}),
\end{align}
\end{subequations}
where $v$ is a variable defined by equation~\eqref{eq_67a}.

In light of~\eqref{SG_T}, we employ the following loss function in PINN,
\begin{align}\label{eq_68}
\text{Loss}=
     &\frac{W_1}{N}\sum_{n=1}^{N}\left[ u_{\theta t}(x_{int}^n, t_{int}^n) - v_{\theta}(x_{int}^n, t_{int}^n) \right]^2 \nonumber \\
 	& + \frac{W_2}{N}\sum_{n=1}^{N}\left[ v_{\theta t}(x_{int}^n, t_{int}^n) - u_{\theta xx}(x_{int}^n, t_{int}^n) + u_{\theta}(x_{int}^n, t_{int}^n) + \sin(u_{\theta}(x_{int}^n, t_{int}^n)) - f(x_{int}^n, t_{int}^n) \right]^2 \nonumber \\
 	& +\frac{W_3}{N}\sum_{n=1}^{N}\left[ u_{\theta t x}(x_{int}^n, t_{int}^n) - v_{\theta x}(x_{int}^n, t_{int}^n) \right]^2 + \frac{W_4}{N} \sum_{n=1}^{N}\left[ u_{\theta}(x_{tb}^n,0) - \psi_{1}(x_{tb}^n)\right]^2 \nonumber \\
 	& + \frac{W_5}{N} \sum_{n=1}^{N}\left[ v_{\theta}(x_{tb}^n,0) - \psi_{2}(x_{tb}^n)\right]^2+ \frac{W_6}{N}\sum_{n=1}^{N} \left[u_{\theta x}(x_{tb}^n,0) - \psi_{1 x}(x_{tb}^n)\right]^2 \nonumber \\
 	& + \frac{W_7}{N} \sum_{n=1}^{N}\left[ | v_{\theta}(0, t_{sb}^n) - \phi_{1t}({t_{sb}^n})| + | v_{\theta}(1, t_{sb}^n) - \phi_{2t}({t_{sb}^n}) |\right],
\end{align}
where $W_n>0$ ($1\leq n\leq 7$)
are the penalty coefficients for different loss terms added in the PINN implementation.
It should be noted that the loss terms with the coefficients $W_3$ and $W_6$ will be absent from the conventional PINN formulation (see~\cite{Raissi2019pinn}). These terms in the training loss are necessary based on the error analysis in Section~\ref{Sine-Gordon}. It should also be noted that the $W_7$ loss terms are not squared, as dictated by the theoretical analysis
of Section~\ref{Sine-Gordon}.

We have also implemented a PINN scheme with a variant form for the loss function,
\begin{align}\label{eq_69}
\text{Loss}=
     &\frac{W_1}{N}\sum_{n=1}^{N}\left[ u_{\theta t}(x_{int}^n, t_{int}^n) - v_{\theta}(x_{int}^n, t_{int}^n) \right]^2 \nonumber \\
 	& + \frac{W_2}{N}\sum_{n=1}^{N}\left[ v_{\theta t}(x_{int}^n, t_{int}^n) - u_{\theta xx}(x_{int}^n, t_{int}^n) + u_{\theta}(x_{int}^n, t_{int}^n) + \sin(u_{\theta}(x_{int}^n, t_{int}^n)) - f(x_{int}^n, t_{int}^n) \right]^2 \nonumber \\
 	& +\frac{W_3}{N}\sum_{n=1}^{N}\left[ u_{\theta t x}(x_{int}^n, t_{int}^n) - v_{\theta x}(x_{int}^n, t_{int}^n) \right]^2 + \frac{W_4}{N} \sum_{n=1}^{N}\left[ u_{\theta}(x_{tb}^n,0) - \psi_{1}(x_{tb}^n)\right]^2 \nonumber \\
 	& + \frac{W_5}{N} \sum_{n=1}^{N}\left[ v_{\theta}(x_{tb}^n,0) - \psi_{2}(x_{tb}^n)\right]^2+ \frac{W_6}{N}\sum_{n=1}^{N} \left[ u_{\theta x}(x_{tb}^n,0) - \psi_{1 x}(x_{tb}^n)\right]^2 \nonumber \\
 	& + \frac{W_7}{N} \sum_{n=1}^{N}\left[ ( v_{\theta}(0, t_{sb}^n) - \phi_{1t}({t_{sb}^n}))^2 + ( v_{\theta}(1, t_{sb}^n) - \phi_{2t}({t_{sb}^n}) )^2\right].
\end{align}
The difference between~\eqref{eq_69} and~\eqref{eq_68} lies in the $W_7$ terms. These $W_7$ terms in~\eqref{eq_69} are squared, and they are not in~\eqref{eq_68}.
We refer to the PINN scheme employing the loss function~\eqref{eq_68} as PINN-G1 and the scheme employing the loss function~\eqref{eq_69} as PINN-G2.

In the simulations we employ a feed-forward neural network with two input nodes (representing $x$ and $t$), two output nodes (representing $u$ and $v$), and two hidden layers, each having a width of $80$ nodes. The $\tanh$ activation function has been used for all the hidden nodes. We employ $N$ collocation points generated from a uniform random distribution within the domain, on 
 each of the domain boundary, and also on the initial boundary, where $N$ is varied systematically in the simulations. The penalty coefficients in the loss functions are taken to be
 $ \bm{W} = (W_1,\dots,W_7)=(0.5, 0.4, 0.5, 0.6, 0.6, 0.6, 0.8)$.

\begin{figure}[tb]
	\centering
	\subfloat[$ t=0.5 $]{
		\begin{minipage}[b]{0.25\textwidth}
			\includegraphics[scale=0.25]{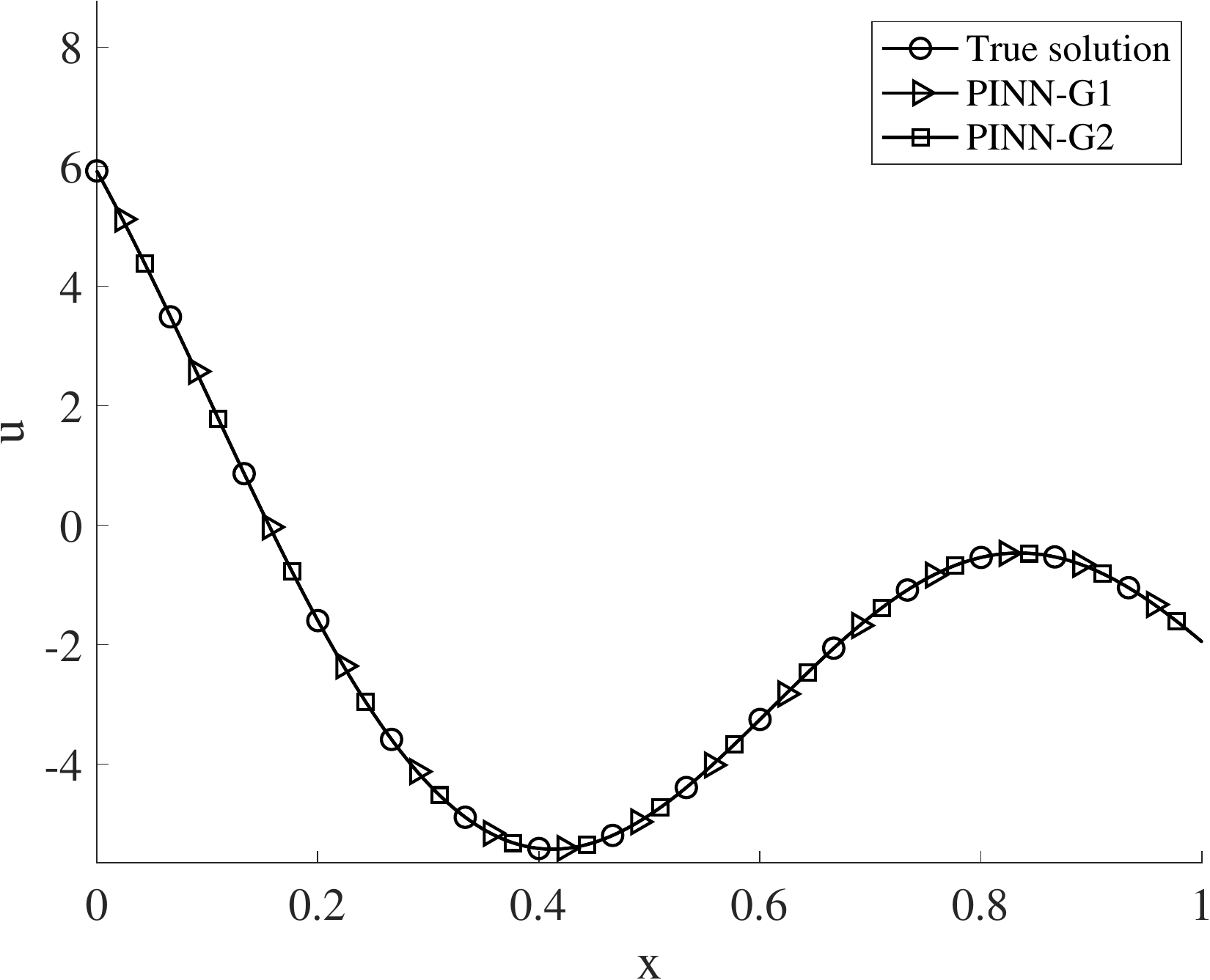}\\
			\includegraphics[scale=0.25]{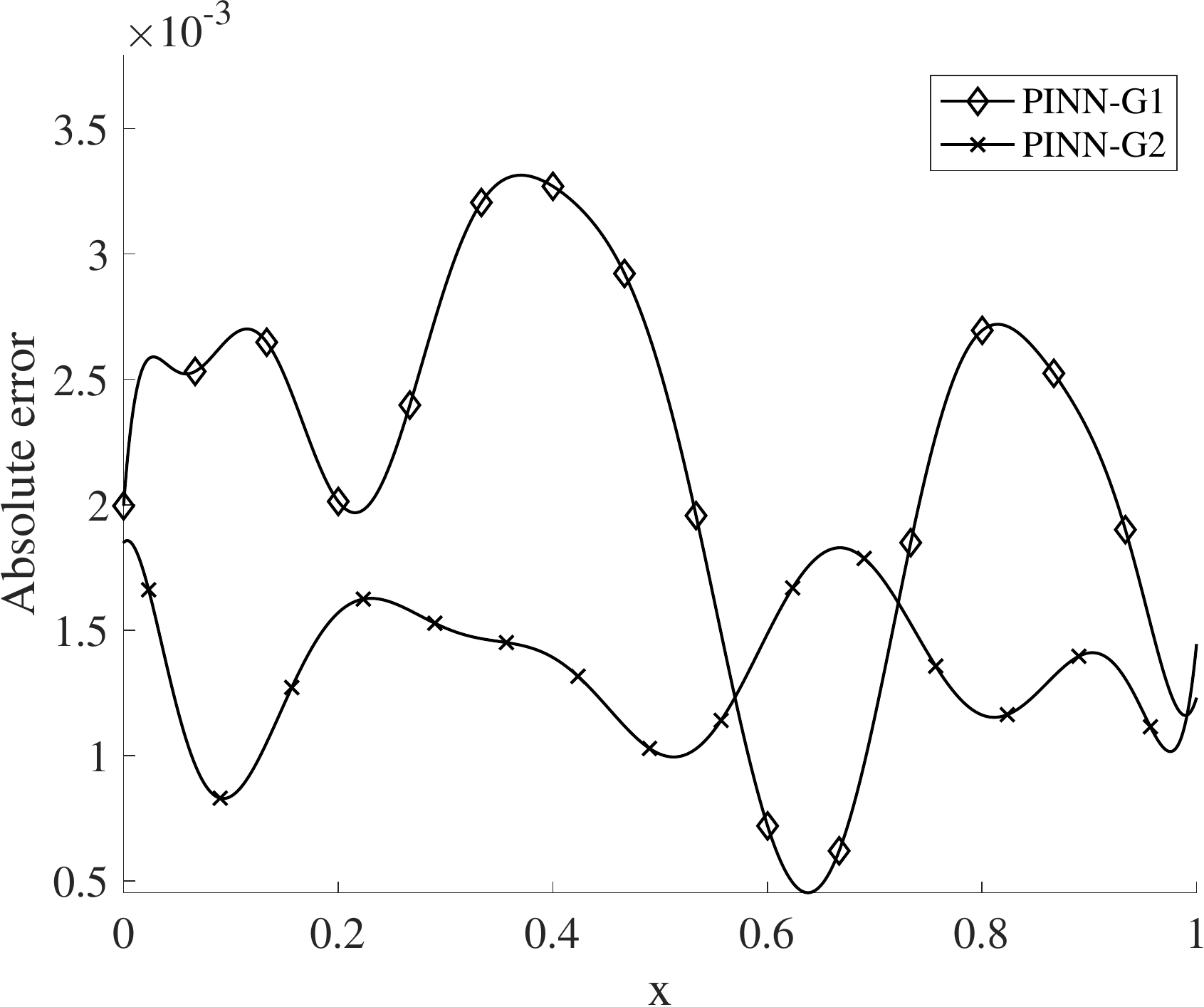}
		\end{minipage}
	}
	\subfloat[$ t=1 $]{
		\begin{minipage}[b]{0.25\textwidth}
			\includegraphics[scale=0.25]{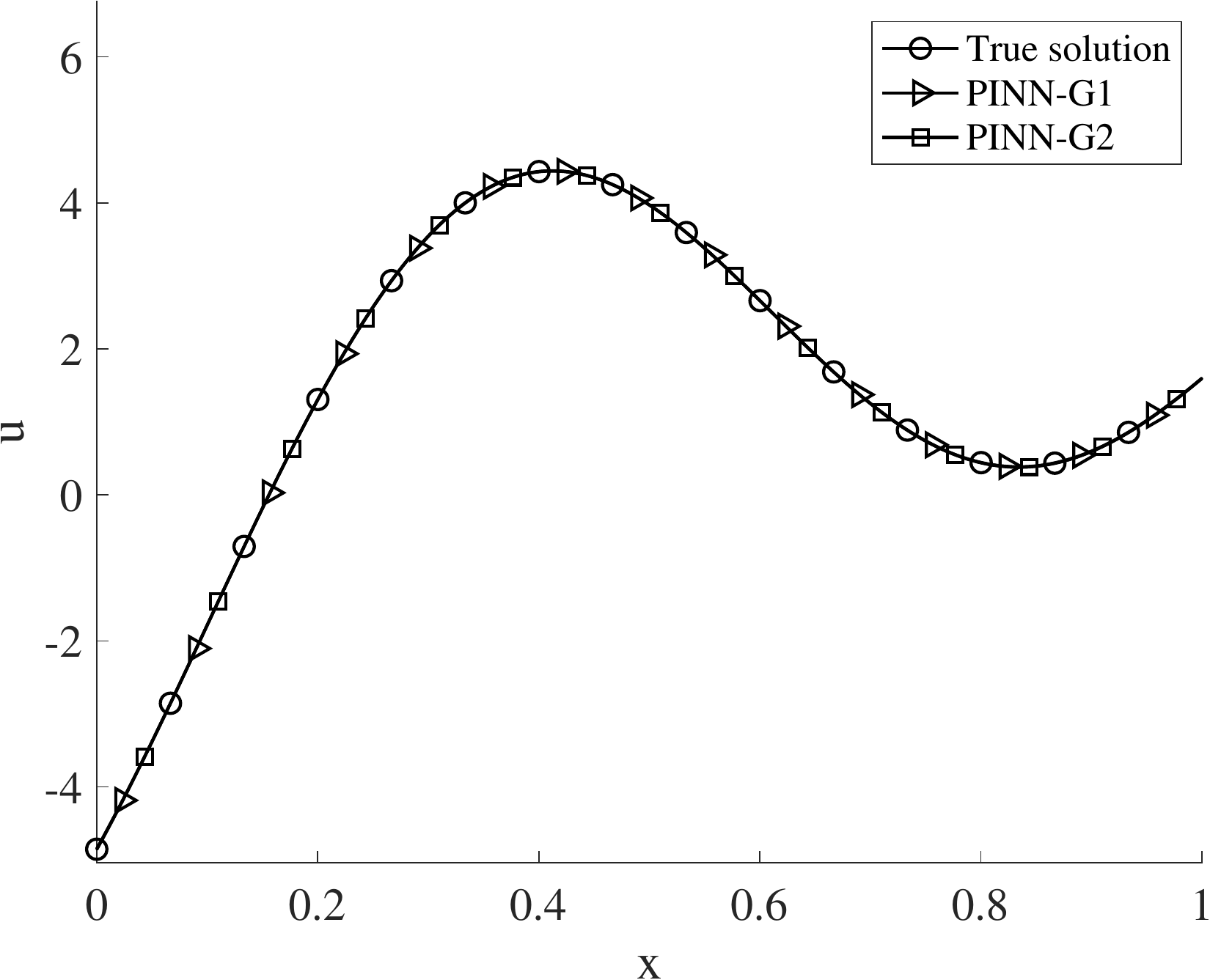}\\
			\includegraphics[scale=0.25]{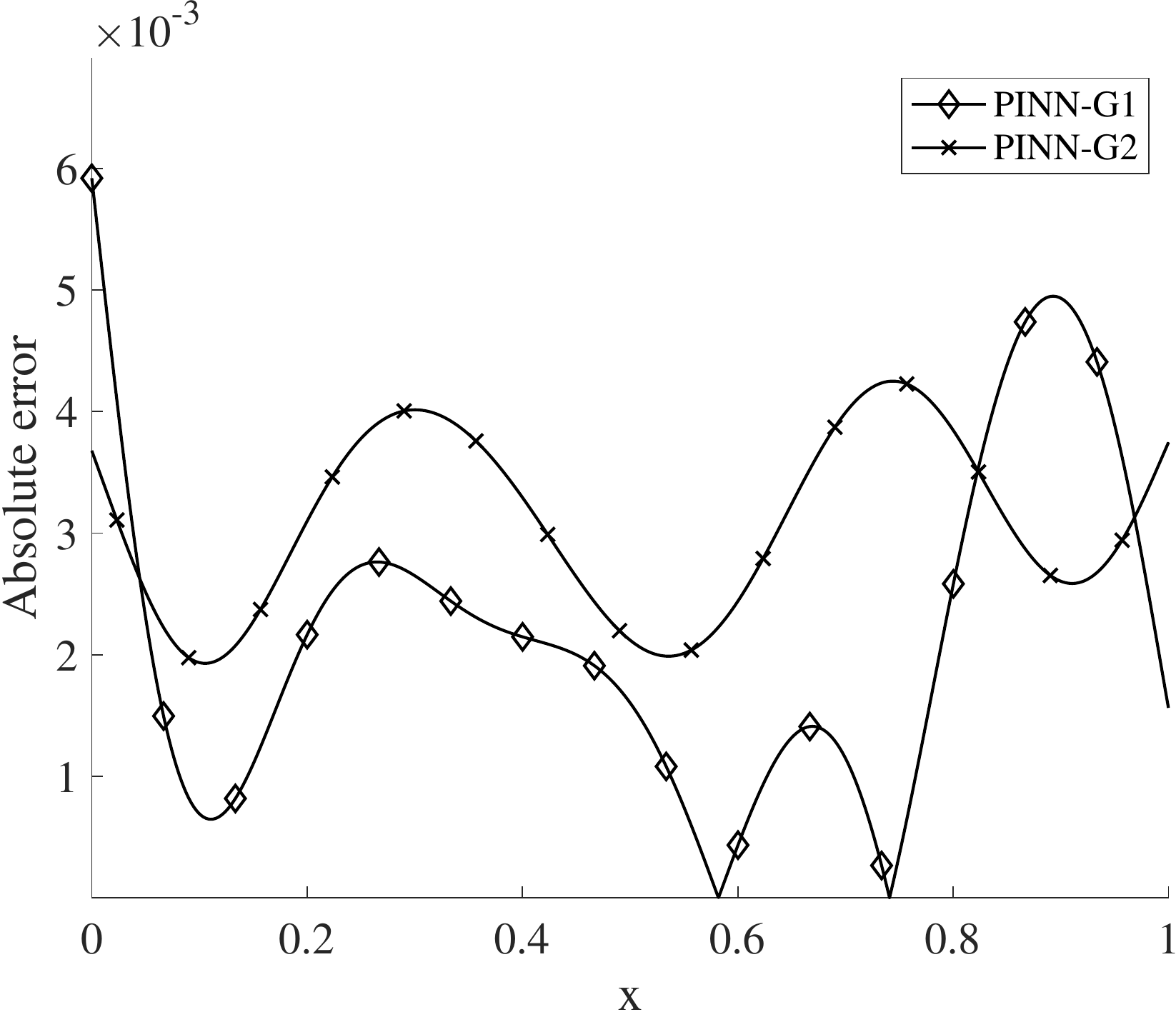}
		\end{minipage}
	}
	\subfloat[$ t=1.5 $]{
		\begin{minipage}[b]{0.25\textwidth}
			\includegraphics[scale=0.25]{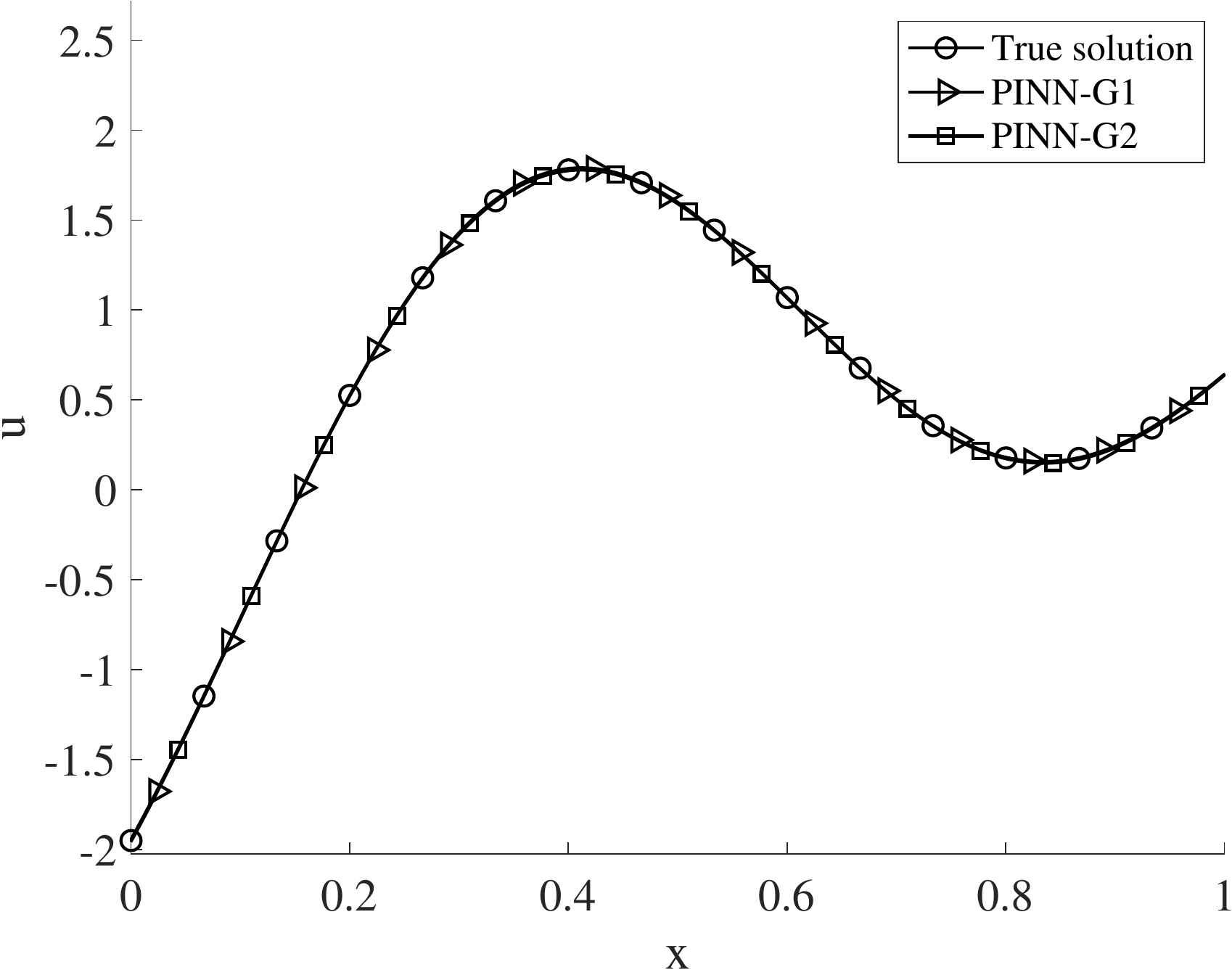}\\
			\includegraphics[scale=0.25]{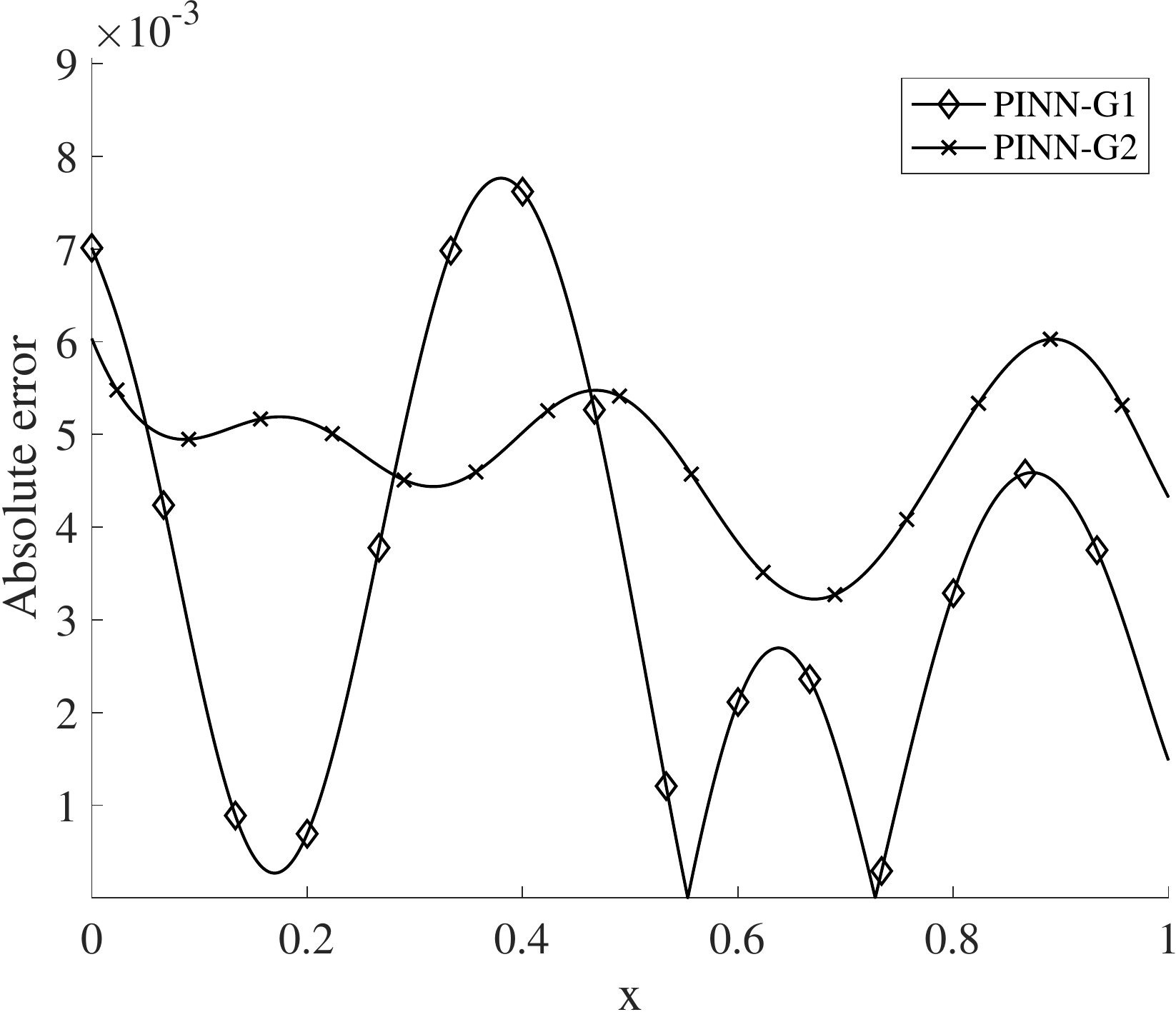}
		\end{minipage}
	}
	\caption{Sine-Gordon equation: Top row, comparison of profiles between the exact solution and PINN-G1/PINN-G2 solutions for $u$ at several time instants. Bottom row, profiles of the absolute error of the PINN-G1 and PINN-G2 solutions for $u$. $N=2000$ 
 training collocation points.
 }
	\label{fg_7}
\end{figure}

\begin{figure}[tb]
	\centering
	\subfloat[$ t=0.5 $]{
		\begin{minipage}[b]{0.25\textwidth}
			\includegraphics[scale=0.25]{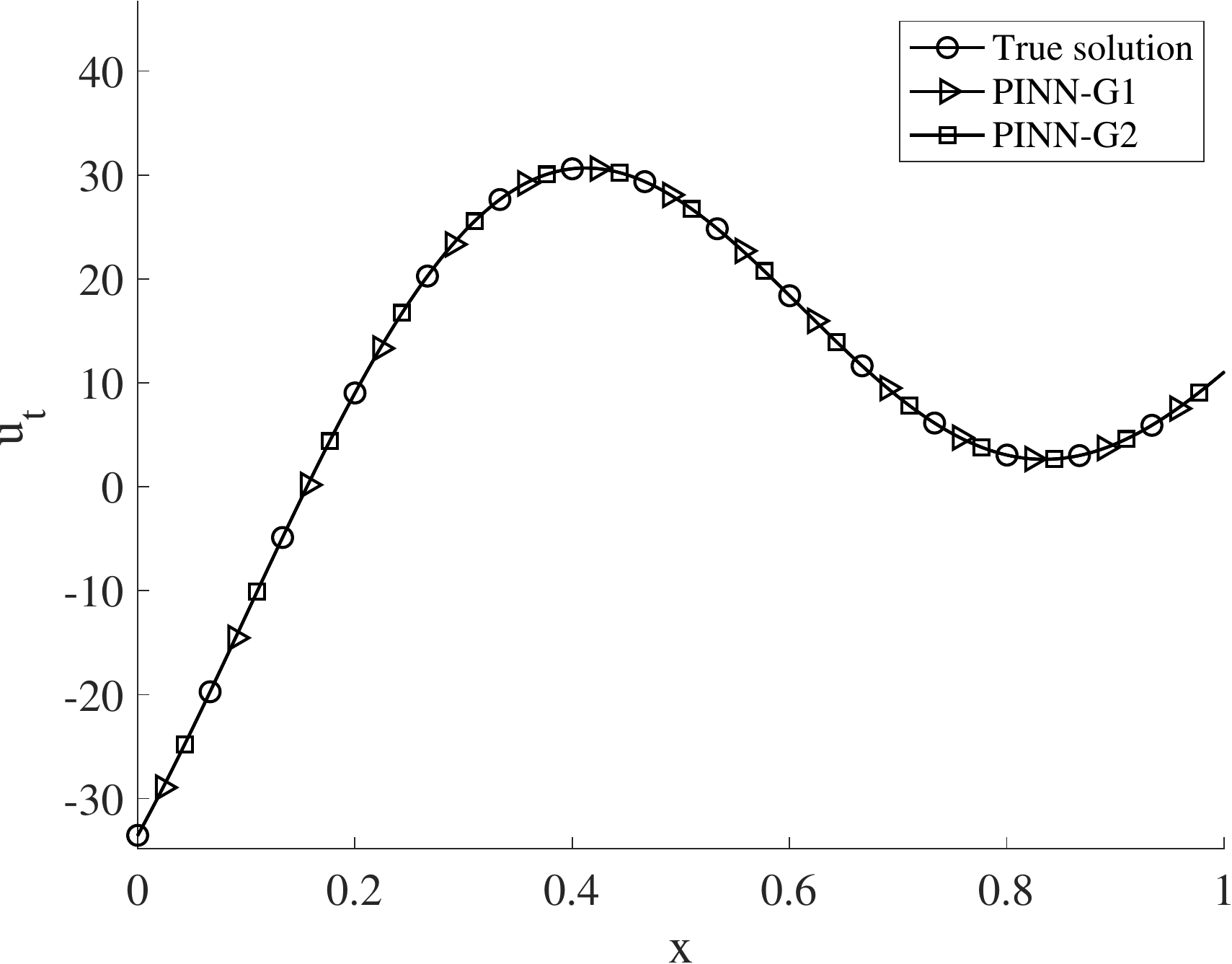}\\
			\includegraphics[scale=0.25]{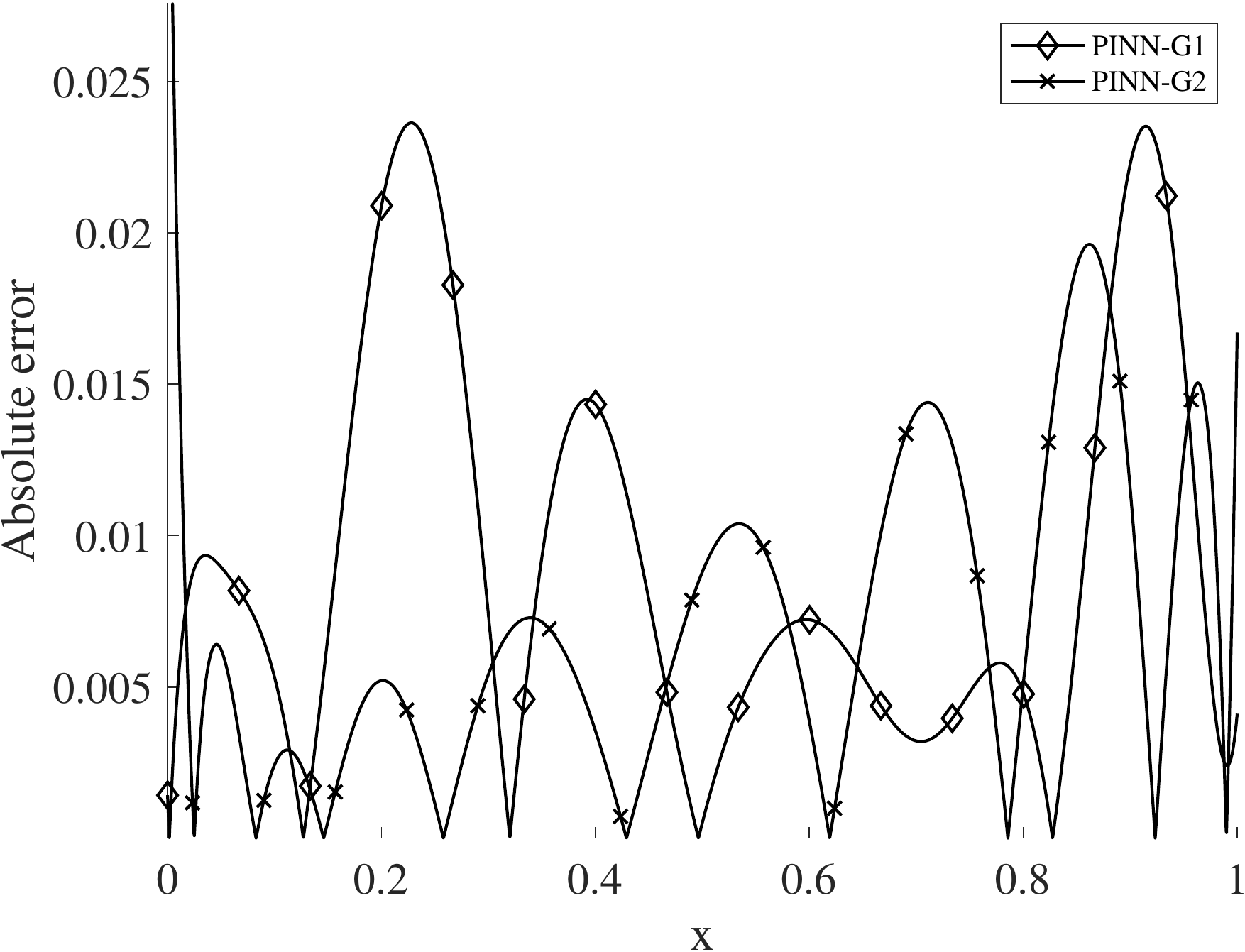}
		\end{minipage}
	}
	\subfloat[$ t=1 $]{
		\begin{minipage}[b]{0.25\textwidth}
			\includegraphics[scale=0.25]{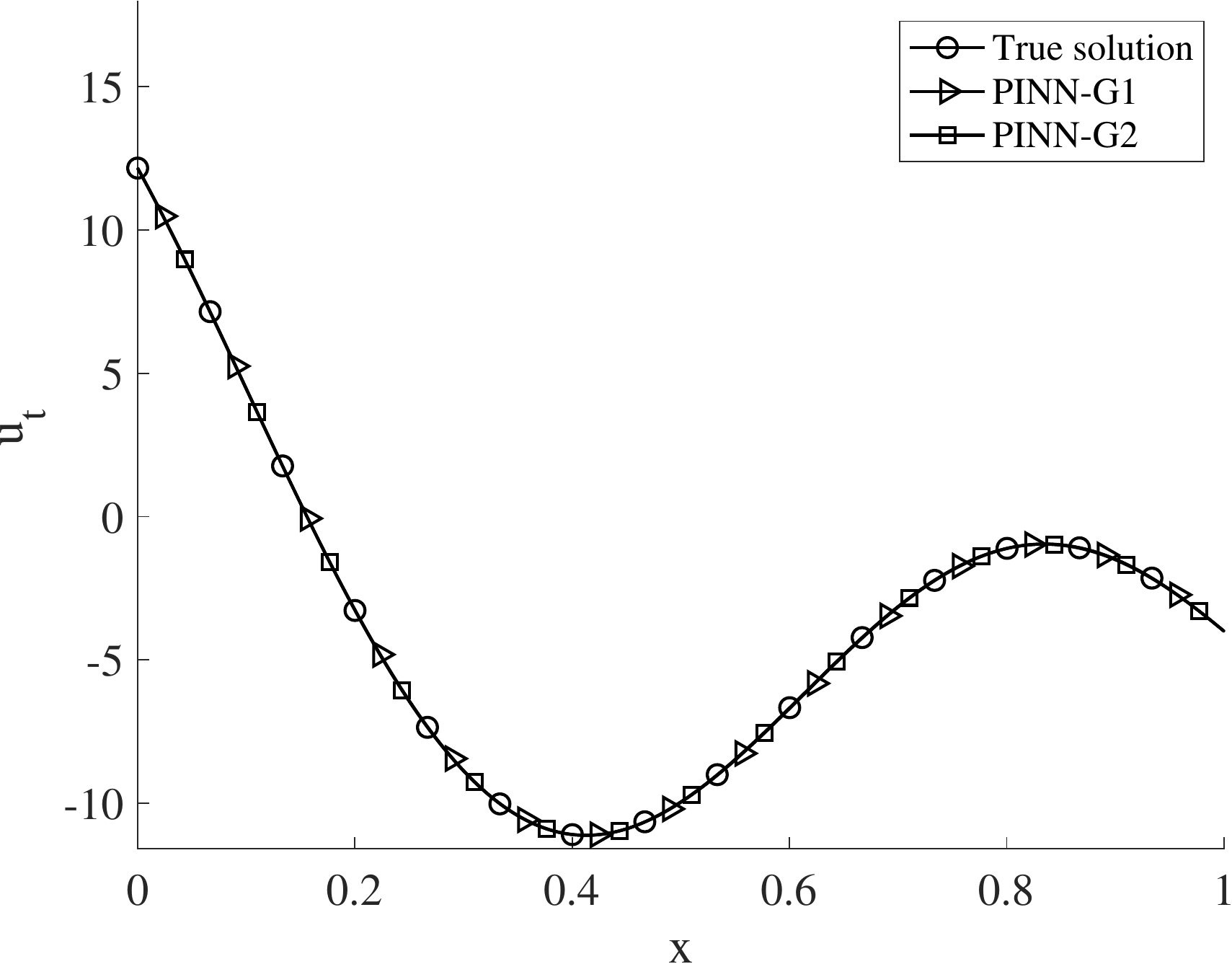}\\
			\includegraphics[scale=0.25]{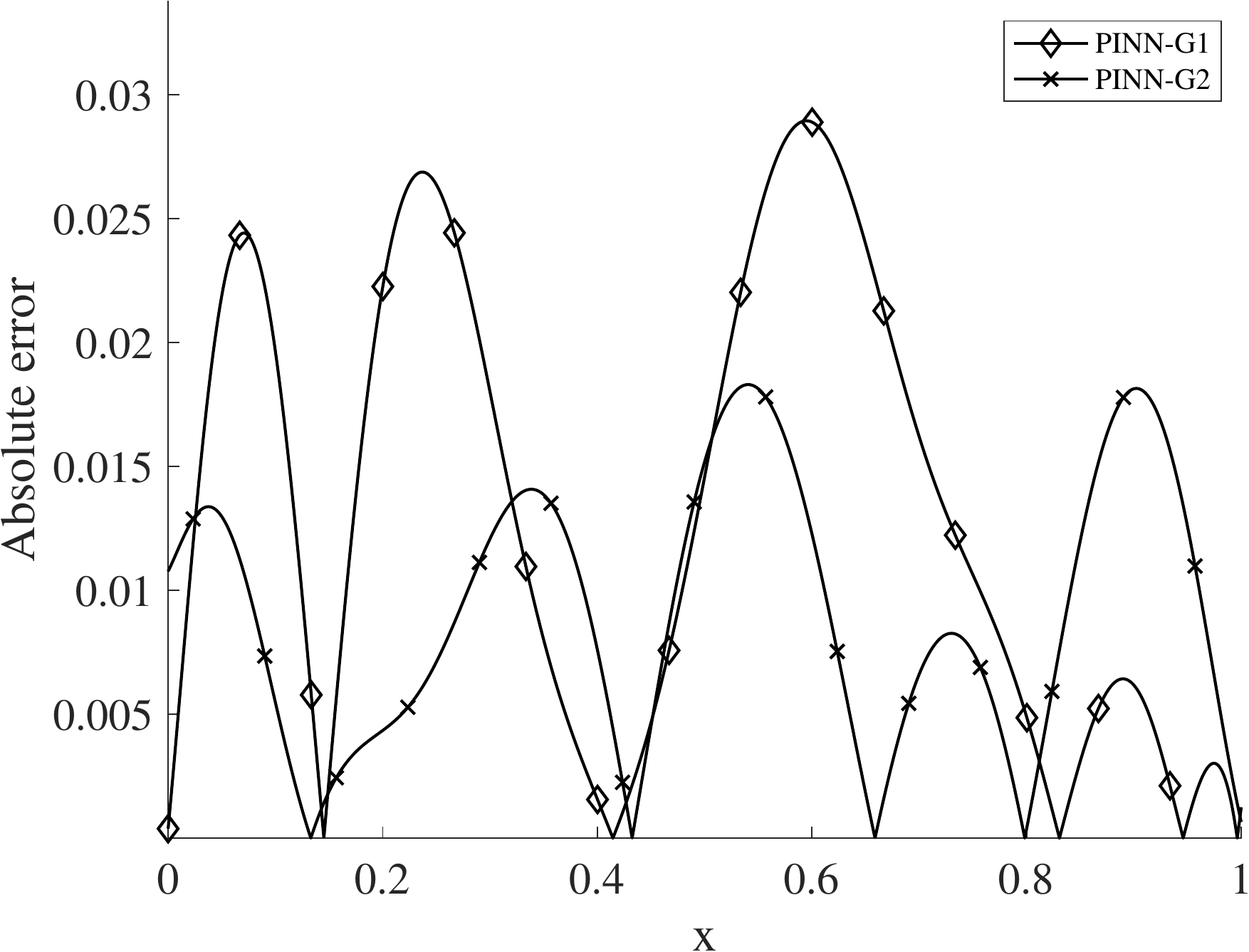}
		\end{minipage}
	}
	\subfloat[$ t=1.5 $]{
		\begin{minipage}[b]{0.25\textwidth}
			\includegraphics[scale=0.25]{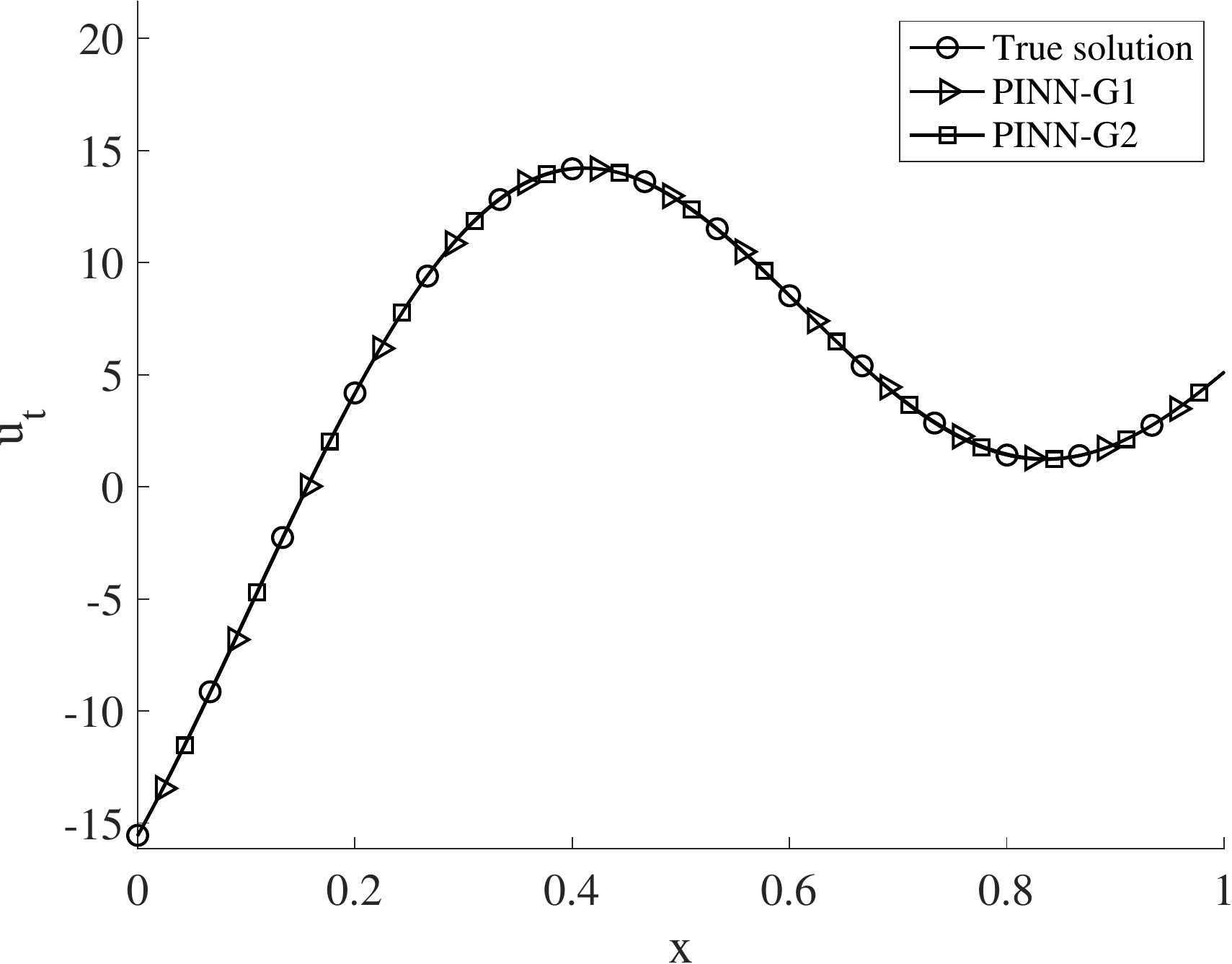}\\
			\includegraphics[scale=0.25]{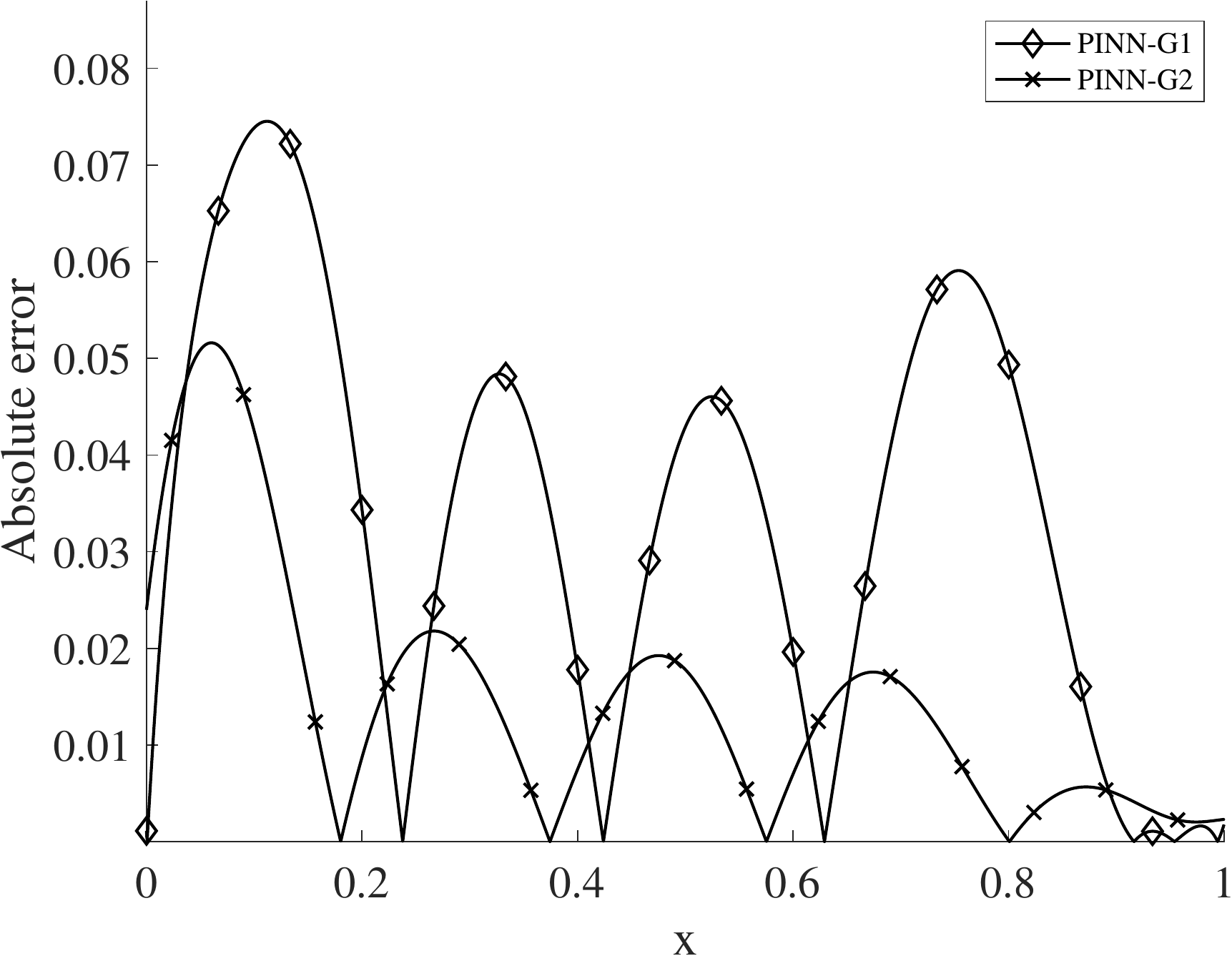}
		\end{minipage}
	}
	\caption{Sine-Gordon equation: Top row, comparison of profiles between the exact solution and PINN-G1/PINN-G2 solutions for $v=\frac{\partial u}{\partial t}$ at several time instants. Bottom row, profiles of the absolute error of the PINN-G1 and PINN-G2 solutions for $v$. $N=2000$ 
 training collocation points.
 }
	\label{fg_8}
\end{figure}

Figure~\ref{fg_6} shows distributions of of $u(x,t)$ and $v=\frac{\partial u}{\partial t}$ from the exact solution (left column) and 
the PINN solution (middle column),
as well as the point-wise absolute errors of the PINN solution for these fields (right column).
These results are obtained by PINN-G2 with $N=2000$ random collocation points within the domain and on each of the domain boundaries. The PINN solution is in good agreement with the true solution.

Figures~\ref{fg_7} and~\ref{fg_8} compare the profiles of $u$ and $v$ between the exact solution, and the solutions obtained by PINN-G1 and PINN-G2, at several time instants ($t=0.5$, $1$ and $1.5$). Profiles of the absolute errors of the PINN-G1/PINN-G2 solutions are also shown in these figures.
We observe that both PINN-G1 and PINN-G2 have captured the solution for $u$ quite accurately, and to a lesser extent, also for $v$.
Comparison of the error profiles between PINN-G1 and PINN-G2 suggests that the PINN-G2 error in general appears to be somewhat smaller than that of PINN-G1. But this seems not to be  true consistently in the entire domain.

\begin{figure}[tb]
	\centering
	\subfloat[PINN-G1]{\includegraphics[width=0.4\linewidth]{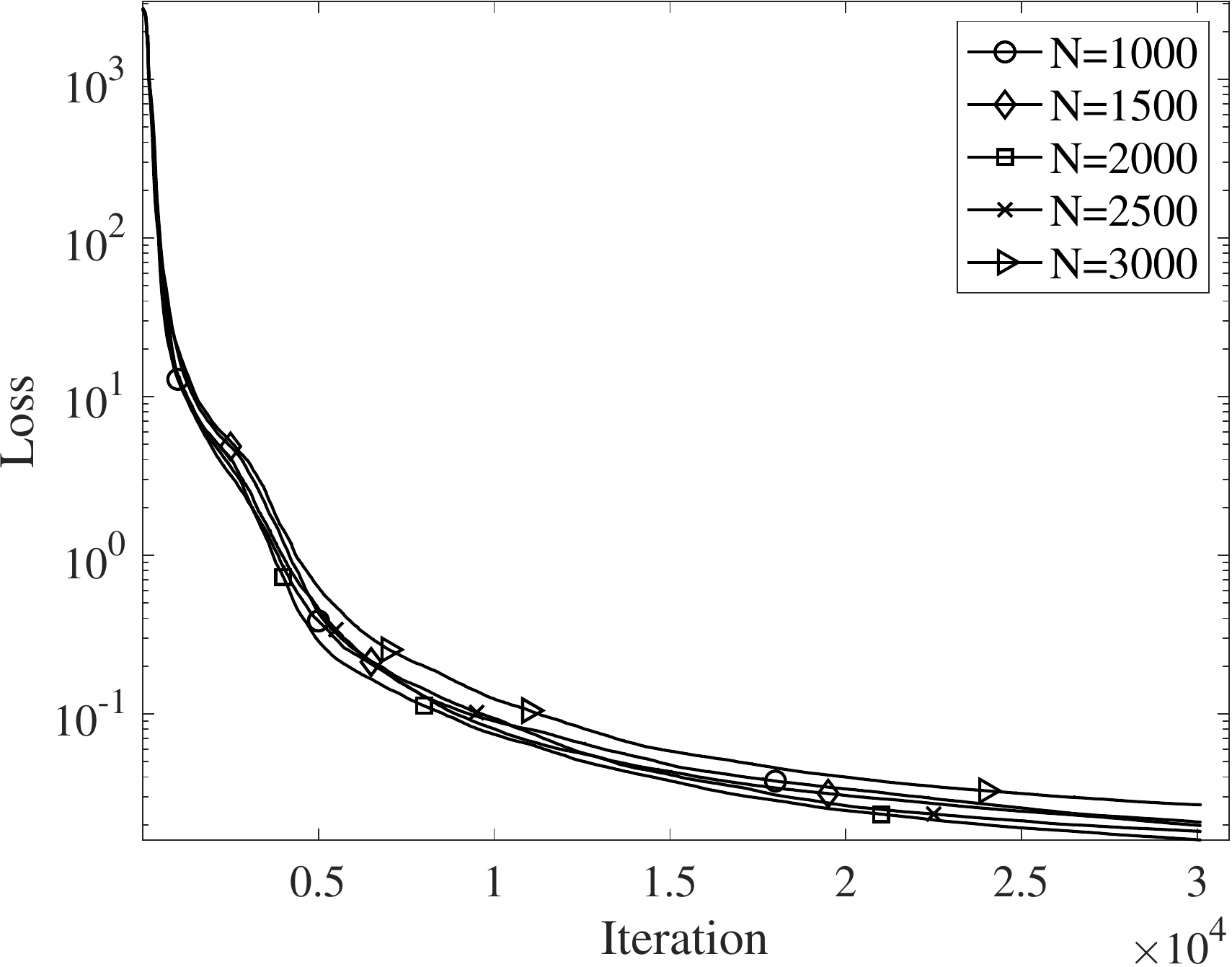}}\qquad
	\subfloat[PINN-G2]{\includegraphics[width=0.4\linewidth]{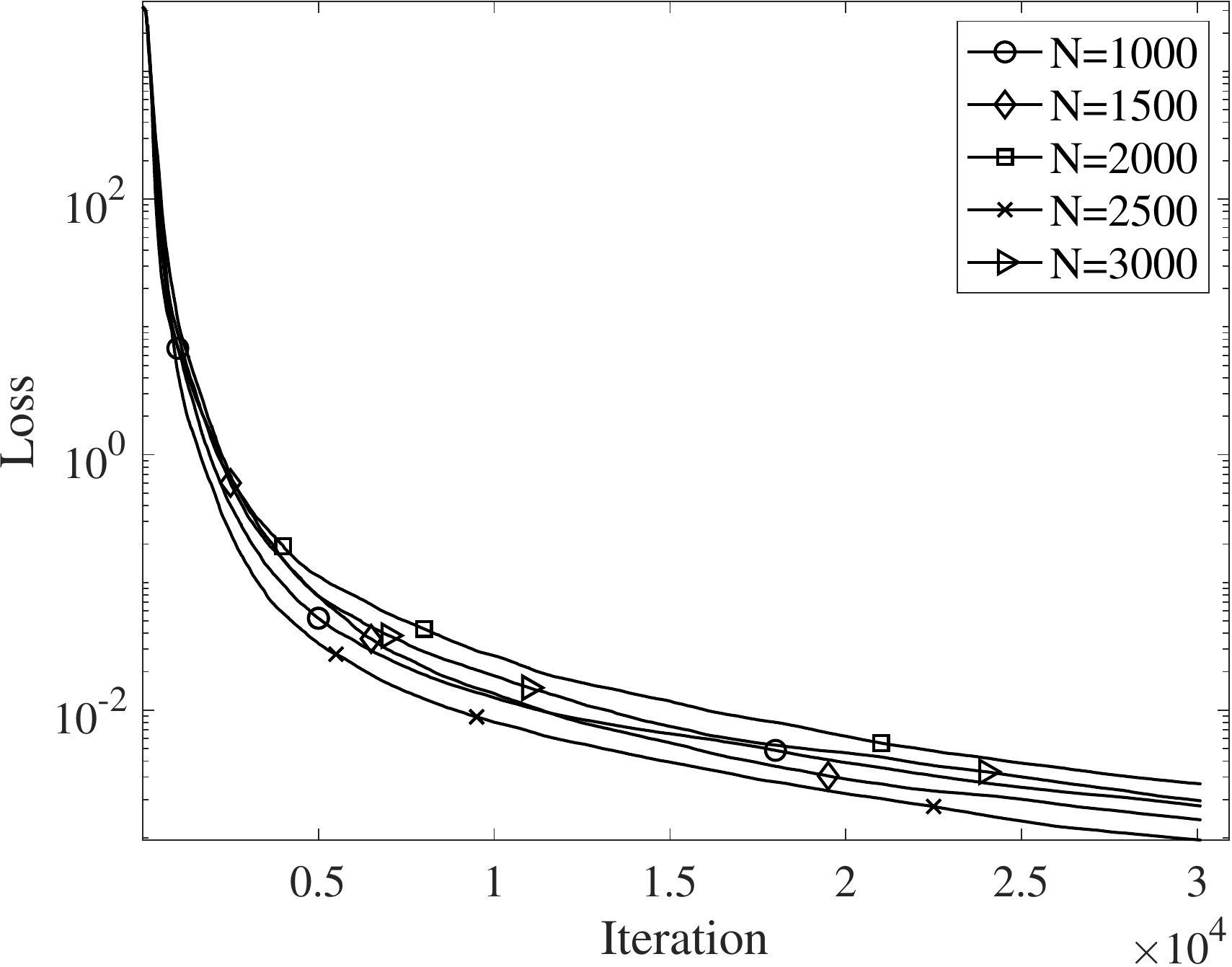}}\hspace{0.1em}
	\caption{Sine-Gordon equation: Loss histories of (a) PINN-G1 and (b) PINN-G2 corresponding to various numbers of training collocation points. 
 }
	\label{PINN_partpaper_SG_fig3}
\end{figure}

\begin{table}[tb]
    \caption{Sine-Gordon equation: The $l_2$ and $l_{\infty}$ errors for $u$ and $v$ versus the number of training collocation points $N$ corresponding to PINN-G1 and PINN-G2. }
    \label{tab_2}
	\centering
	\begin{tabular}[b]{ c | c@{\ \ }|  c  @{\ \ \ }  c  @{\ \ } | c  @{\ \ \ }  c  @{\ \ } }
		\hline
    \multirow{2}{*}{method}	&	
  \multirow{2}{*}{$ N $}&\multicolumn{2}{c}{$ l_2 $-error} & \multicolumn{2}{c}{$ l_\infty $-error} \\ 
		\cline{3-6}
		& &$ u_\theta $  &$ v_\theta $ &$ u_\theta $  &$ v_\theta $ \\
		\hline
		& 1000&  3.0818e-03&  4.3500e-03&  9.6044e-03&  1.8894e-02\\
		\cline{3-6}
		& 1500&  3.4335e-03&  4.8035e-03&  1.0566e-02&  1.7050e-02\\
		\cline{3-6}
	PINN-G1	& 2000&  2.1914e-03&  3.0055e-03&  7.5882e-03&  1.1099e-02\\
		\cline{3-6}
		& 2500&  3.0172e-03&  3.5698e-03&  9.2515e-03&  1.4645e-02\\
		\cline{3-6}
		& 3000&  2.5281e-03&  4.4858e-03&  7.2785e-03&  1.6213e-02\\
		\hline
		& 1000&  3.0674e-03&  2.0581e-03&  7.3413e-03&  1.1323e-02\\
		\cline{3-6}
		& 1500&  1.0605e-03&  1.4729e-03&  2.2914e-03&  6.2831e-03\\
		\cline{3-6}
	PINN-G2	& 2000&  2.2469e-03&  1.6072e-03&  4.8842e-03&  8.8320e-03\\
		\cline{3-6}
		& 2500&  6.6072e-04&  6.0509e-04&  1.4099e-03&  4.3423e-03\\
		\cline{3-6}
		& 3000&  6.6214e-04&  1.0830e-03&  1.9697e-03&  7.8866e-03\\
		\hline
	\end{tabular}
\end{table}


The effect of the collocation points on the PINN results has been studied by varying the number of training collocation points systematically between $N=1000$ and $N=3000$
within the domain and on each of the domain boundaries. The results are provided in Figure~\ref{PINN_partpaper_SG_fig3} and Table~\ref{tab_2}.
Figure~\ref{PINN_partpaper_SG_fig3} shows  histories of the loss function corresponding to different number of collocation points for PINN-G1 and PINN-G2.
Table~\ref{tab_2} provides the $l_2$ and $l_{\infty}$ errors of $u$ and $v$ versus the number of collocation points computed by PINN-G1 and PINN-G2.
The PINN errors in general tend to decrease with increasing number of collocation points, but this trend is not monotonic.
It can be observed that both PINN-G1 and PINN-G2 have captured the solutions quite accurately, with those errors from PINN-G2 in general slightly better.

\begin{figure}
	\centering
	\subfloat[PINN-G1]{\includegraphics[width=0.4\linewidth]{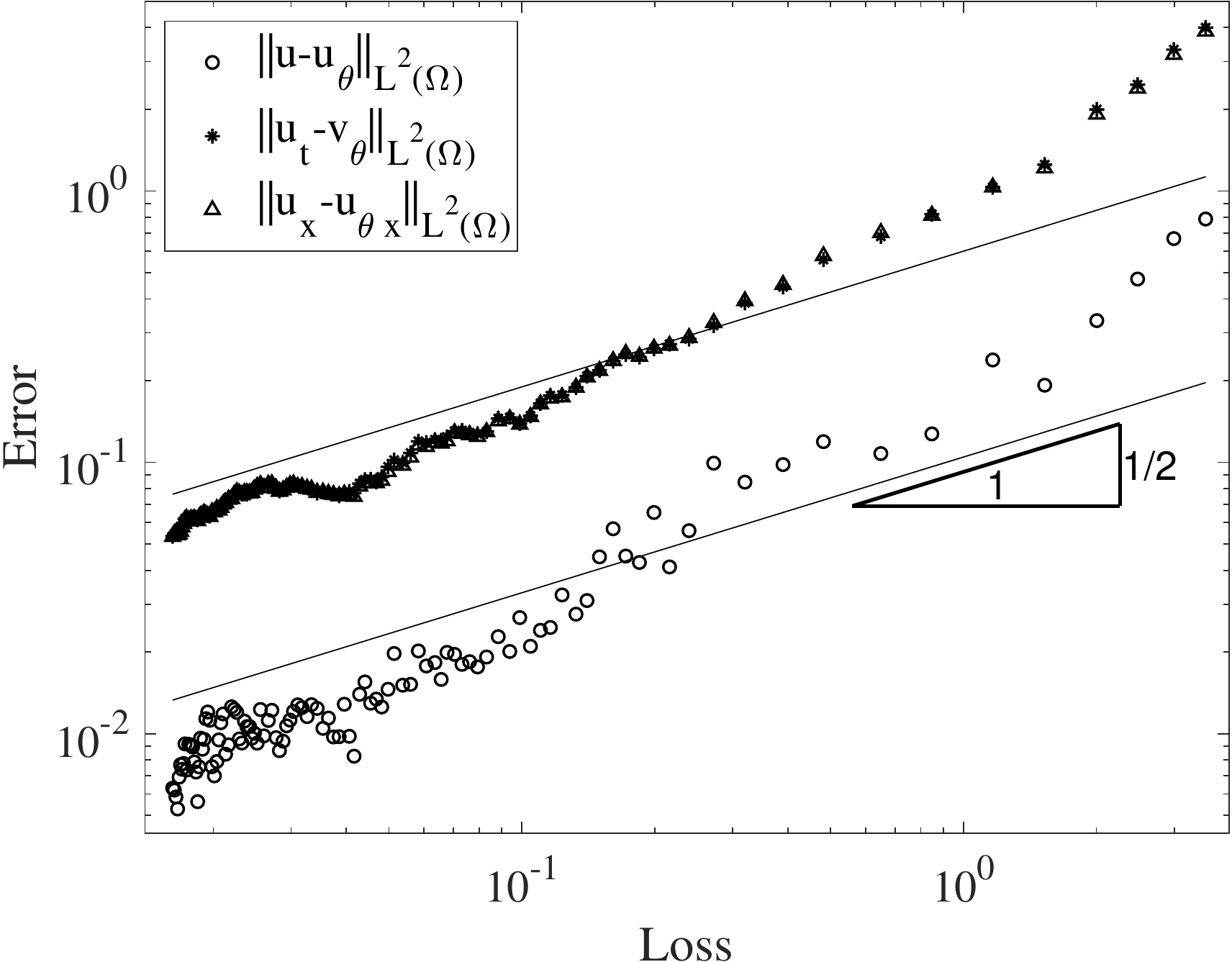}}\qquad
	\subfloat[PINN-G2]{\includegraphics[width=0.4\linewidth]{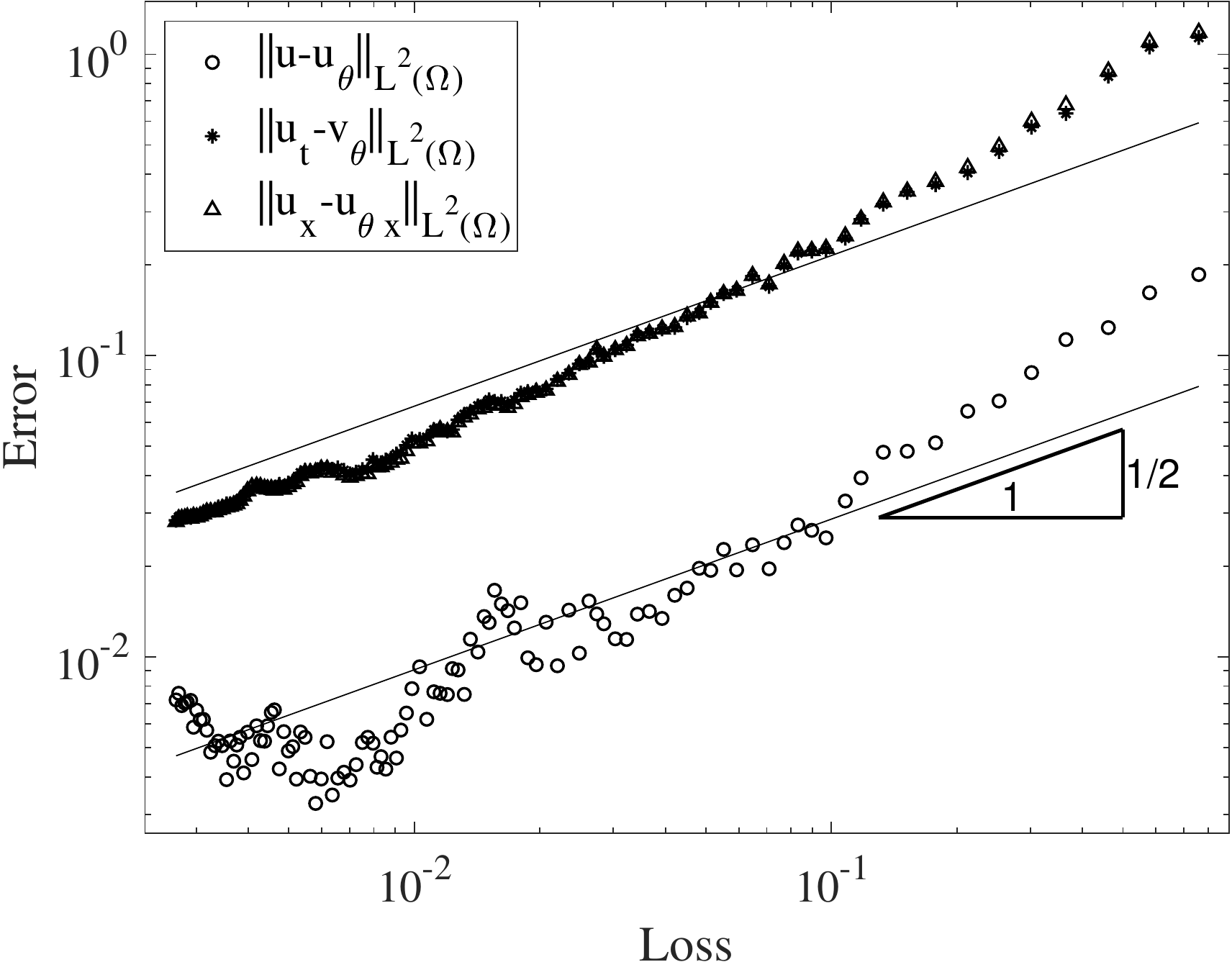}}\hspace{0.1em}
	\caption{Sine-Gordon equation: The $l^2$ errors of $u$, $\frac{\partial u}{\partial t}$, and $\frac{\partial u}{\partial x}$ as a function of the training loss value.
 }
	\label{PINN_partpaper_SG_errorrate}
\end{figure}


Figure \ref{PINN_partpaper_SG_errorrate} provides some numerical evidence for the relation between the total error and the training loss as suggested by Theorem \ref{sec6_Theorem3}.
Here we plot the $l_2$ errors for $u$, $v$ and $\frac{\partial u}{\partial x}$ as a function of the training loss value obtained by PINN-G1 and PINN-G2.  
The results indicate that the total error scales approximately as the square root of the training loss, which in some sense corroborates the error-loss relation as expressed in Theorem \ref{sec6_Theorem3}.

\subsection{Linear Elastodynamic Equation}

In this subsection we look into the linear elastodynamic equation (in two spatial dimensions plus time) and test the PINN algorithm as suggested by the theoretical analysis in Section~\ref{Elasto-dynamics} using this equation.
Consider the spatial-temporal domain $(x,y,t)\in \Omega = D \times [0, T]= [0, 1]\times[0,1]\times [0, 2] $, and the following initial/boundary value problem with the linear elastodynamics equation on $ \Omega $:
\begin{subequations}\label{num_ela_eq00}
\begin{align}\label{num_ela_eq1}
	&\rho\frac{\partial^2 \bm{u}}{\partial t^2} - 2\mu\nabla\cdot(\underline{\bm{\varepsilon}}(\bm{u})) -\lambda\nabla(\nabla\cdot\bm{u})= \bm{f}(\bm{x}, t),\\
	&\bm{u}|_{\Gamma_d}=\bm{\phi}_{d},\qquad \Big(2\mu\underline{\bm{\varepsilon}}(\bm{u}) +\lambda(\nabla\cdot\bm{u})\Big)|_{\Gamma_n}\bm{n}=\bm{\phi}_n,\\
	&\bm{u}(\bm{x}, 0)=\bm{\psi}_{1}, \qquad \frac{\partial\bm u}{\partial t}(\bm{x}, 0)=\bm{\psi}_{2},
\end{align}
\end{subequations}
where $ \bm u = (u_1(\bm x, t), u_2(\bm x, t))^T $ ($ \bm x=(x,y)\in D $, $ t\in[0, T] $) is the displacement field to be solved for,  $ \bm{f}(\bm{x}, t) $ is a  source term, and $\rho$, $\mu$ and $\lambda$ are material constants. $ \Gamma_d $ is the Dirichlet boundary and $ \Gamma_n $ is the Neumann boundary, with $ \partial D=\Gamma_d \cup \Gamma_n $ and $ \Gamma_d \cap \Gamma_n = \emptyset $, where $ \bm{n} $ is the outward-pointing unit normal vector. In our simulations we choose the left boundary ($x=0$) as the Dirichlet boundary, and the rest are Neumann boundaries. $ \bm{\phi}_d $ and $ \bm{\phi}_n $ are Dirichlet and Neumann boundary conditions, respectively. $ \bm{\psi}_1 $ and $ \bm{\psi}_2 $ are the initial conditions for the displacement and the velocity. We employ the material parameter values $ \mu = \lambda = \rho = 1 $, and the following manufactured solution (\cite{2018_CMAME_DGelastodynamics}) to this problem,
\begin{align}\label{num_ela_eq2}
	\bm u(\bm x, t) = \sin(\sqrt{2}\pi t) \begin{bmatrix}
		-\sin(\pi x)^2\sin(2\pi y)\\
		\sin(2\pi x)\sin(\pi y)^2
	\end{bmatrix}.
\end{align}
The source term $ \bm{f}(\bm{x}, t) $, the boundary/initial distributions $ \bm{\phi}_d $, $ \bm{\phi}_n $, $ \bm{\psi}_{1} $ and $ \bm{\psi}_{2} $ are chosen by the expression \eqref{num_ela_eq2}.

To simulate this problem using the PINN algorithm suggested by
the theoretical analysis from Section \ref{Elasto-dynamics}, we reformulate~\eqref{num_ela_eq00} into the following system
\begin{subequations}\label{num_ela_eq3}
 	\begin{align}
 		\label{num_ela_eq3_1}
 		&\bm{u}_{t} - \bm{v} = \bm{0},\qquad \bm{v}_{t} - 2\nabla\cdot(\underline{\bm{\varepsilon}}(\bm{u})) -\nabla(\nabla\cdot\bm{u})= \bm{f}(\bm{x}, t),\\
 		&\bm{u}|_{\Gamma_d}=\bm{\phi}_{d},\qquad \Big(2\underline{\bm{\varepsilon}}(\bm{u}) +(\nabla\cdot\bm{u})\Big)|_{\Gamma_n}\bm{n}=\bm{\phi}_n,\\
 		&\bm{u}(\bm{x}, 0)=\bm{\psi}_{1}, \qquad \bm{v}(\bm{x}, 0)=\bm{\psi}_{2},
 	\end{align}
\end{subequations}
where $ \bm{v}(\bm{x}, t) $ is an intermediate variable (representing the velocity) as given by \eqref{num_ela_eq3_1}.

In light of~\eqref{elast_T},
we employ the following loss function for PINN, 
\begin{align}\label{num_ela_eq4}
 	\text{Loss}
    &= \frac{W_1}{N} \sum_{n=1}^{N} \left[ \bm{u}_{\theta t}(\bm{x}_{int}^n, t_{int}^n) - \bm{v}_{\theta}(\bm{x}_{int}^n, t_{int}^n)\right]^2 \nonumber\\
 	& +  \frac{W_2}{N} \sum_{n=1}^{N} \left[ \bm{v}_{\theta t}(\bm{x}_{int}^n, t_{int}^n) - 2\nabla\cdot (\underline{\bm{\varepsilon}}(\bm{u}_{\theta}(\bm{x}_{int}^n, t_{int}^n))) - \nabla(\nabla\cdot \bm{u}_{\theta}(\bm{x}_{int}^n, t_{int}^n)) -\bm{f}(\bm{x}_{int}^n, t_{int}^n)) \right]^2 \nonumber\\
 	& + \frac{W_3}{N} \sum_{n=1}^{N}\left[ \underline{\bm{\varepsilon}}(\bm{u}_{\theta t}(\bm{x}_{int}^n, t_{int}^n) - \bm{v}_{\theta}(\bm{x}_{int}^n, t_{int}^n))\right]^2 + \frac{W_4}{N} \sum_{n=1}^{N} \left[ \nabla\cdot(\bm{u}_{\theta t}(\bm{x}_{int}^n, t_{int}^n) - \bm{v}_{\theta}(\bm{x}_{int}^n, t_{int}^n)) \right]^2 \nonumber\\
 	&+ \frac{W_5}{N} \sum_{n=1}^{N} \left[ \bm{u}_{\theta}(\bm{x}_{tb}^n, 0) -\bm{\psi}_1(\bm{x}_{tb}^n) \right]^2 + \frac{W_6}{N} \sum_{n=1}^{N} \left[ \bm{v}_{\theta}(\bm{x}_{tb}^n, 0) -\bm{\psi}_2(\bm{x}_{tb}^n) \right]^2 \nonumber \\
 	&+ \frac{W_7}{N} \sum_{n=1}^{N} \left[ \underline{\bm{\varepsilon}}(\bm{u}_{\theta}(\bm{x}_{tb}^n, 0) -\bm{\psi}_1(\bm{x}_{tb}^n))  \right]^2 + \frac{W_8}{N} \sum_{n=1}^{N} \left[\nabla\cdot(\bm{u}_{\theta}(\bm{x}_{tb}^n, 0) -\bm{\psi}_1(\bm{x}_{tb}^n))  \right]^2  \nonumber \\
 	&+\frac{W_9}{N} \sum_{n=1}^{N} |\bm{v}_{\theta}(\bm{x}_{sb1}^n, t_{sb1}^n) - \bm{\phi}_{dt}(\bm{x}_{sb1}^n, t_{sb1}^n)| \nonumber\\
 	&+ \frac{W_{10}}{N} \sum_{n=1}^{N} | 2\underline{\bm{\varepsilon}}(\bm{u}_{\theta}(\bm{x}_{sb2}^n, t_{sb2}^n))\bm{n} +(\nabla\cdot\bm{u}_{\theta}(\bm{x}_{sb2}^n, t_{sb2}^n))\bm{n} - \bm{\phi}_n(\bm{x}_{sb2}^n, t_{sb2}^n) |,
\end{align}
where we have added the penalty coefficients, $W_n>0$ ($ 1\leq n\leq 10 $), for different loss terms in the implementation, and $N$ denotes the number of collocation points within the domain and on the domain boundaries. 
In the numerical tests we have also implemented another form for the loss function as follows,
\begin{align}\label{num_ela_eq5}
	\text{Loss}&
	=\frac{W_1}{N} \sum_{n=1}^{N} \left[ \bm{u}_{\theta t}(\bm{x}_{int}^n, t_{int}^n) - \bm{v}_{\theta}(\bm{x}_{int}^n, t_{int}^n)\right]^2 \nonumber\\
	& +  \frac{W_2}{N} \sum_{n=1}^{N} \left[ \bm{v}_{\theta t}(\bm{x}_{int}^n, t_{int}^n) - 2\nabla\cdot (\underline{\bm{\varepsilon}}(\bm{u}_{\theta}(\bm{x}_{int}^n, t_{int}^n))) - \nabla(\nabla\cdot \bm{u}_{\theta}(\bm{x}_{int}^n, t_{int}^n)) -\bm{f}(\bm{x}_{int}^n, t_{int}^n)) \right]^2 \nonumber\\
	& + \frac{W_3}{N} \sum_{n=1}^{N}\left[ \underline{\bm{\varepsilon}}(\bm{u}_{\theta t}(\bm{x}_{int}^n, t_{int}^n) - \bm{v}_{\theta}(\bm{x}_{int}^n, t_{int}^n))\right]^2 + \frac{W_4}{N} \sum_{n=1}^{N} \left[ \nabla\cdot(\bm{u}_{\theta t}(\bm{x}_{int}^n, t_{int}^n) - \bm{v}_{\theta}(\bm{x}_{int}^n, t_{int}^n)) \right]^2 \nonumber\\
	&+ \frac{W_5}{N} \sum_{n=1}^{N} \left[ \bm{u}_{\theta}(\bm{x}_{tb}^n, 0) -\bm{\psi}_1(\bm{x}_{tb}^n) \right]^2 + \frac{W_6}{N} \sum_{n=1}^{N} \left[ \bm{v}_{\theta}(\bm{x}_{tb}^n, 0) -\bm{\psi}_2(\bm{x}_{tb}^n) \right]^2 \nonumber \\
	&+ \frac{W_7}{N} \sum_{n=1}^{N} \left[ \underline{\bm{\varepsilon}}(\bm{u}_{\theta}(\bm{x}_{tb}^n, 0) -\bm{\psi}_1(\bm{x}_{tb}^n))  \right]^2 + \frac{W_8}{N} \sum_{n=1}^{N} \left[ \nabla\cdot(\bm{u}_{\theta}(\bm{x}_{tb}^n, 0) -\bm{\psi}_1(\bm{x}_{tb}^n))  \right]^2  \nonumber \\
	&+\frac{W_9}{N} \sum_{n=1}^{N} \left[\bm{v}_{\theta}(\bm{x}_{sb1}^n, t_{sb1}^n) - \bm{\phi}_{dt}(\bm{x}_{sb1}^n, t_{sb1}^n)\right]^2 \nonumber\\
	&+ \frac{W_{10}}{N} \sum_{n=1}^{N} \left[2\underline{\bm{\varepsilon}}(\bm{u}_{\theta}(\bm{x}_{sb2}^n, t_{sb2}^n))\bm{n} +(\nabla\cdot\bm{u}_{\theta}(\bm{x}_{sb2}^n, t_{sb2}^n))\bm{n} - \bm{\phi}_n(\bm{x}_{sb2}^n, t_{sb2}^n) \right]^2.
\end{align}
The difference between these two forms for the loss function lies in the $W_9$ and $W_{10}$ terms.
It should be noted that 
the $W_9$ and $ W_{10} $ terms in \eqref{num_ela_eq4} are not squared, in light of the error terms \eqref{elast_T1}$ - $\eqref{elast_T8} from the theoretical analysis. In contrast, these terms are squared in \eqref{num_ela_eq5}. The PINN scheme utilizing the loss function \eqref{num_ela_eq4} is henceforth referred to as PINN-H1, and the scheme that employs the loss function \eqref{num_ela_eq5} shall be referred to as PINN-H2.

In the  simulations, we employ a feed-forward neural network with three input nodes, which represent $ \bm{x}=(x,y) $ and the time variable t, and four output nodes, which represent $ \bm{u}=(u_1, u_2) $ and $ \bm{v}=(v_1, v_2) $. The neural network has two hidden layers, with widths of 90 and 60 nodes, respectively, and the $\tanh$ activation function for all the hidden nodes. 
For the network training, $ N $ collocation points are generated from a uniform random distribution within the domain, on each of the domain boundary, as well as on the initial boundary. $ N $ is systematically varied in the simulations. We employ the penalty coefficients $ \bm{W} = (W_1, ..., W_{10}) = (0.9, 0.9, 0.9, 0.9, 0.5, 0.5, 0.5, 0.5, 0.9, 0.9) $ in the simulations.

\begin{figure}[htpb]
	\centering
	\subfloat[$ t=0.5 $]{
		\begin{minipage}[b]{0.25\textwidth}
			\includegraphics[scale=0.24]{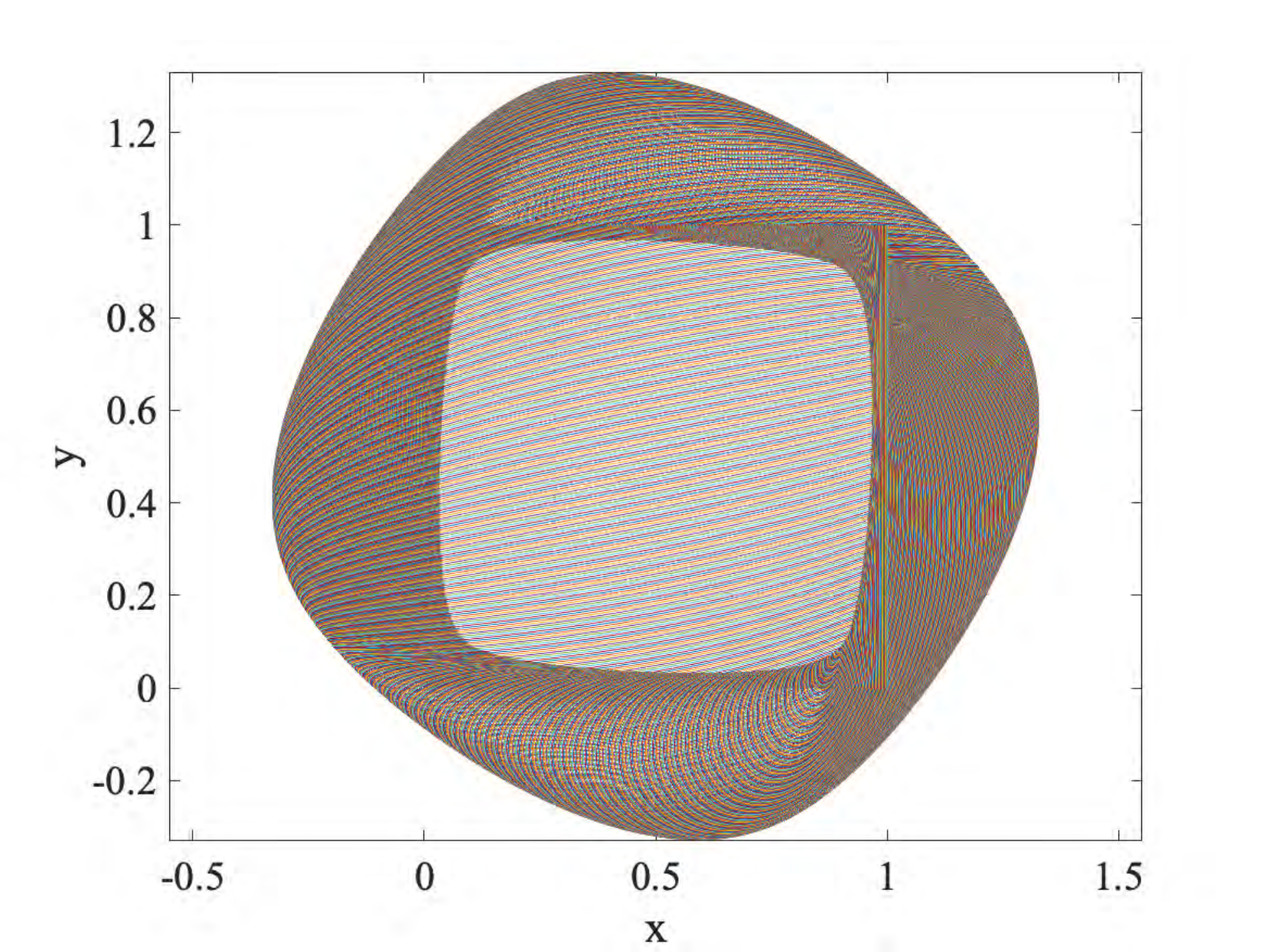}\\
			\includegraphics[scale=0.24]{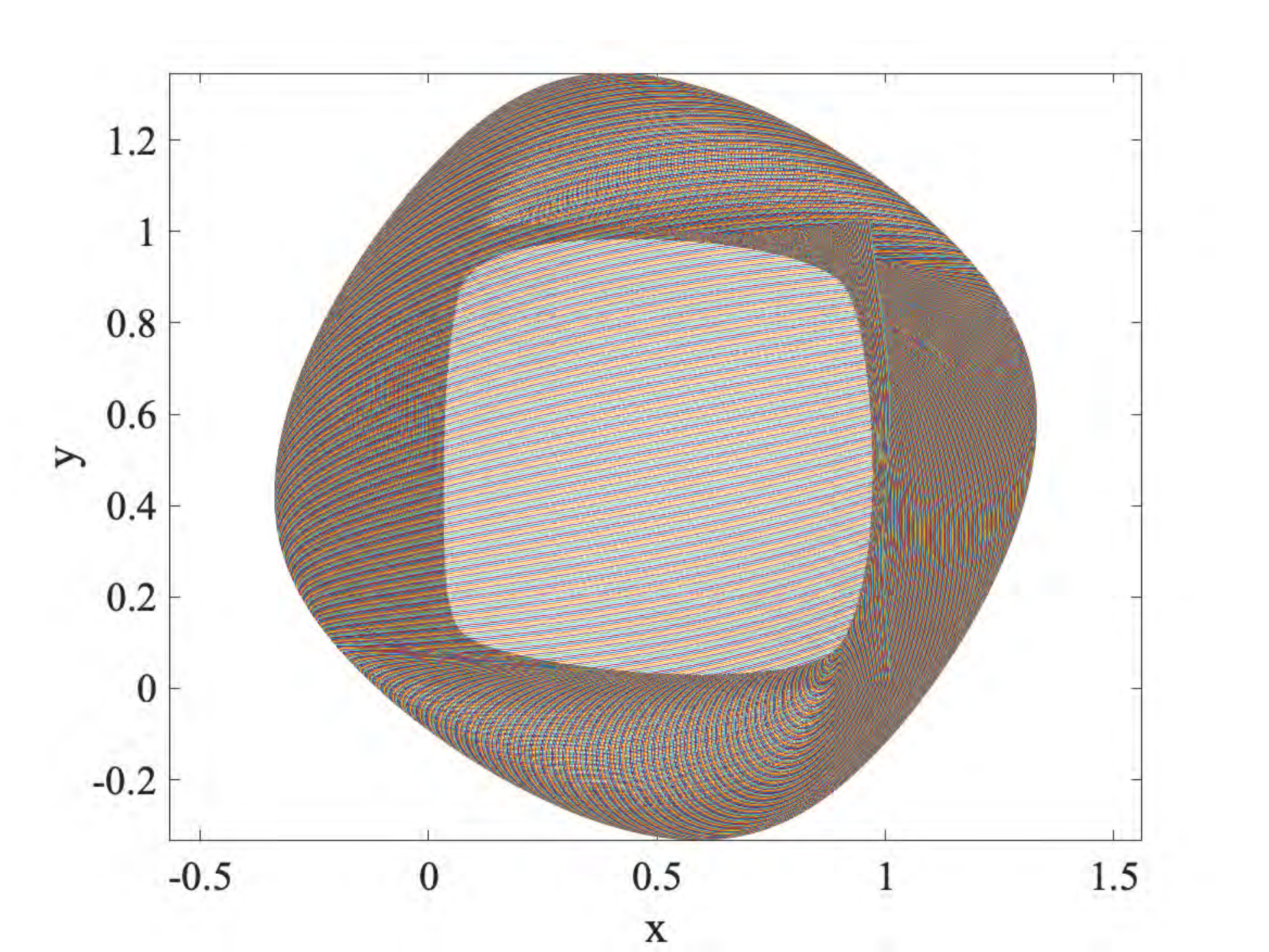}\\
			\includegraphics[scale=0.24]{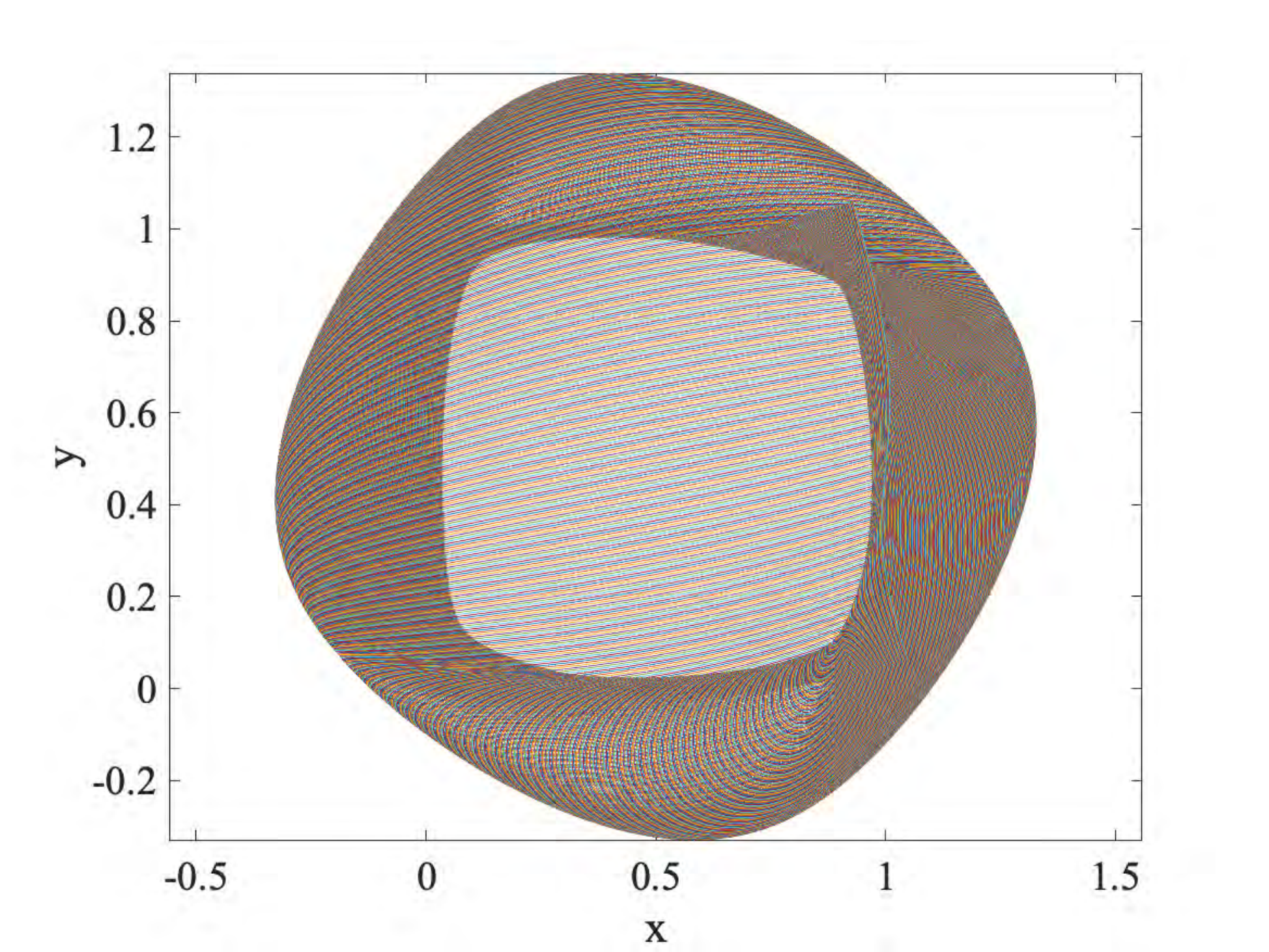}
		\end{minipage}
	}
	\subfloat[$ t=1 $]{
		\begin{minipage}[b]{0.25\textwidth}
			\includegraphics[scale=0.24]{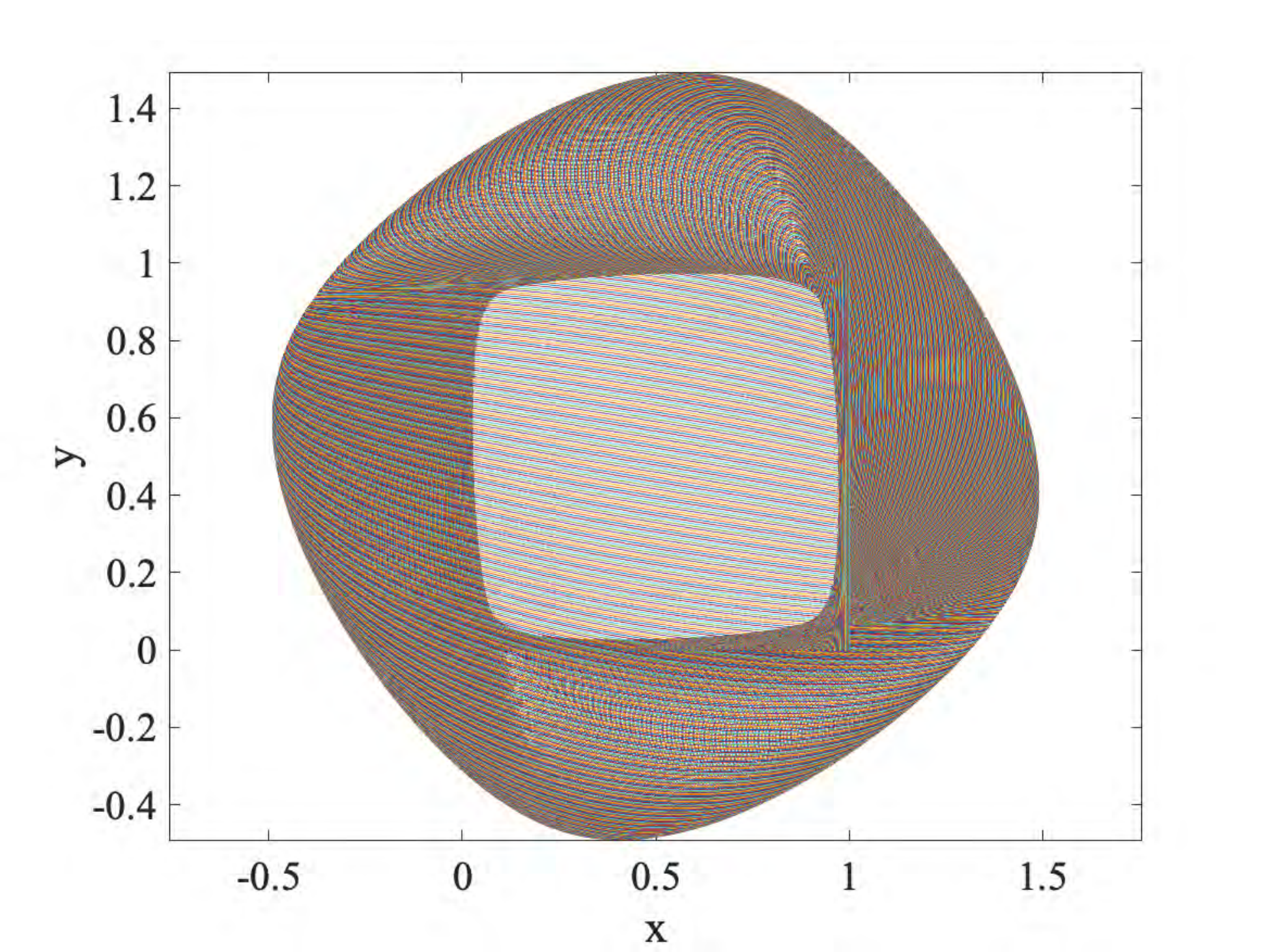}\\
			\includegraphics[scale=0.24]{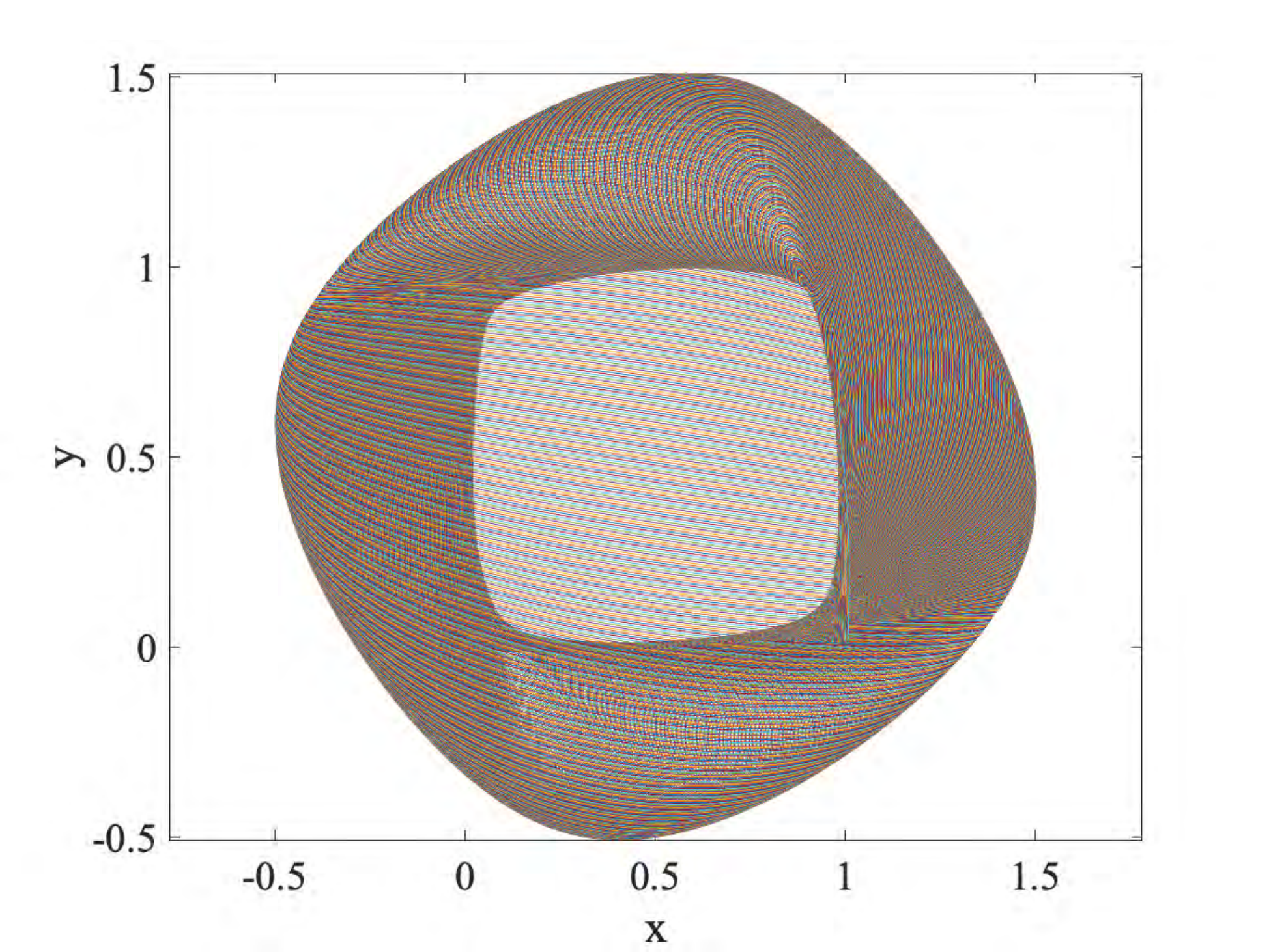}\\
			\includegraphics[scale=0.24]{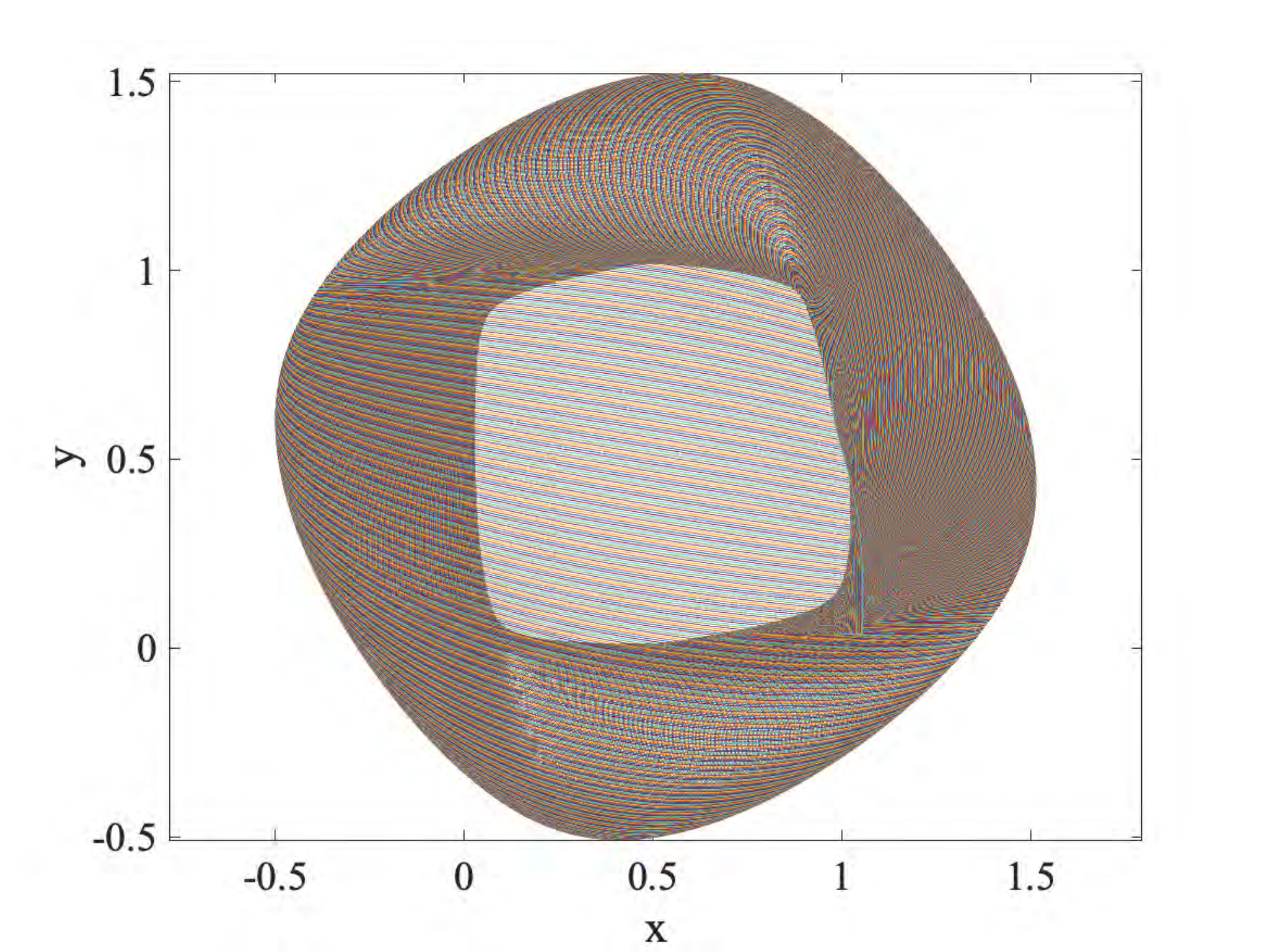}
		\end{minipage}
	}
	\subfloat[$ t=1.5 $]{
		\begin{minipage}[b]{0.25\textwidth}
			\includegraphics[scale=0.24]{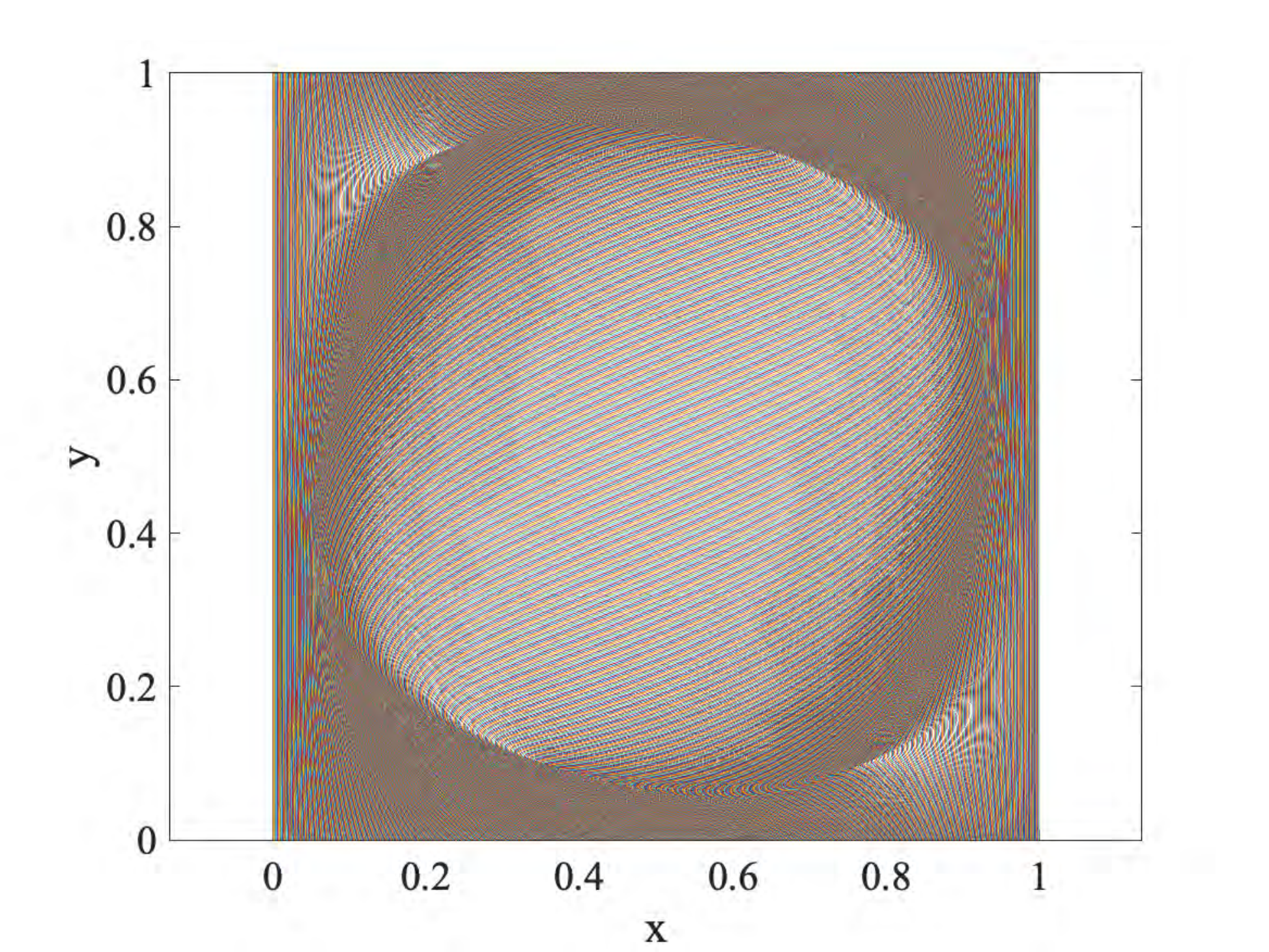}\\
			\includegraphics[scale=0.24]{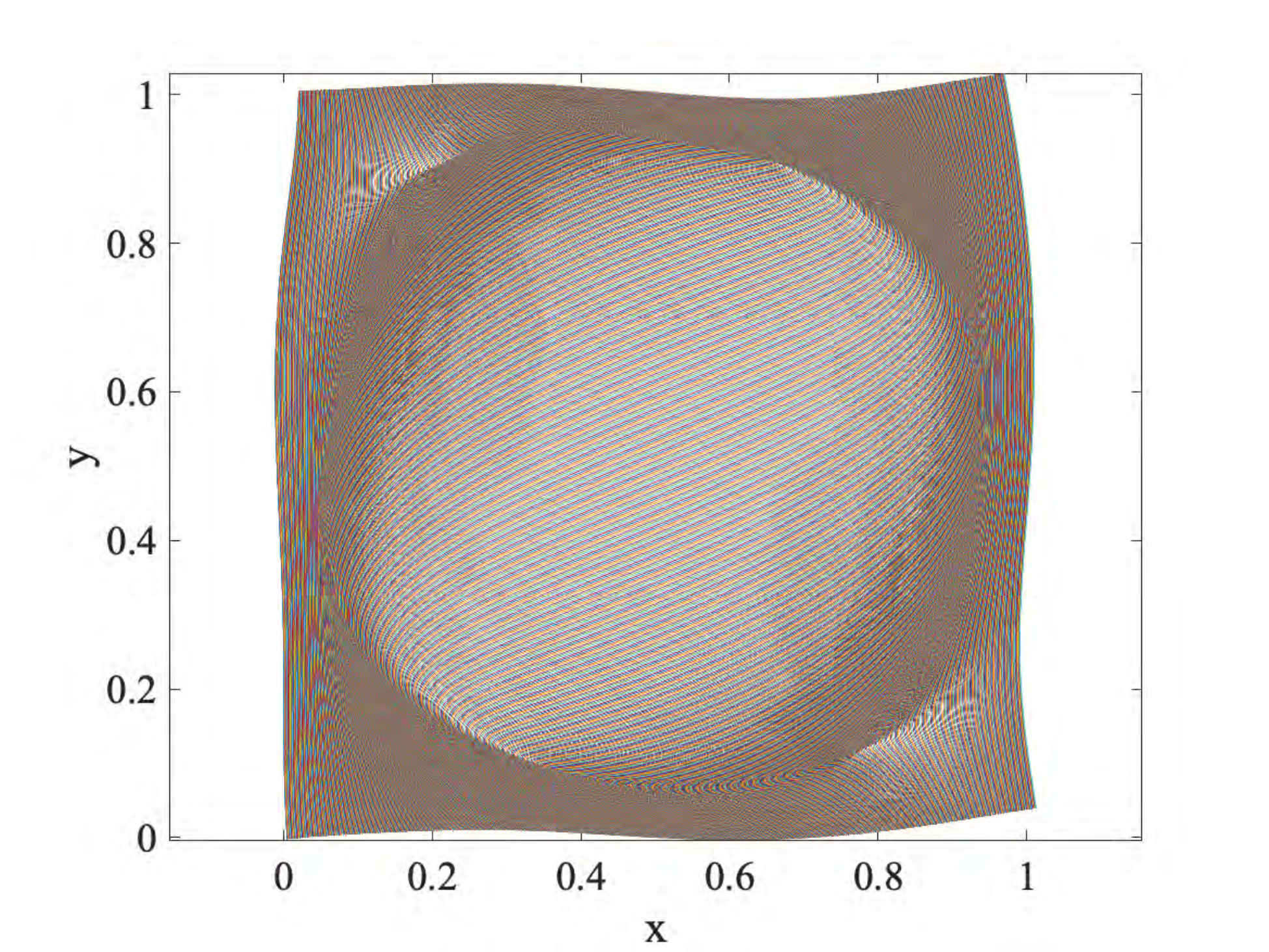}\\
			\includegraphics[scale=0.24]{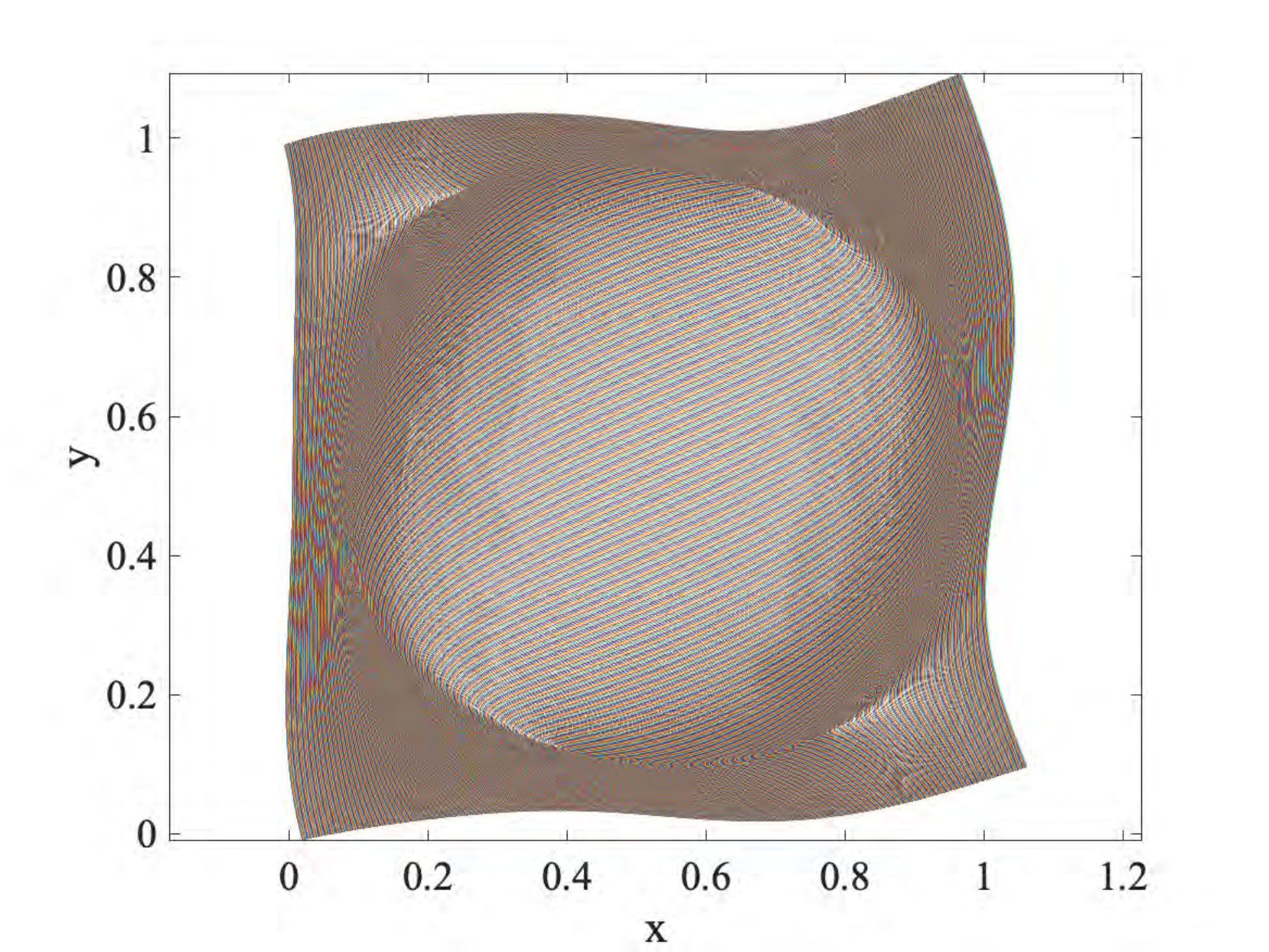}
		\end{minipage}
	}
\caption{Linear elastodynamic equation: 
Visualization of the deformed configuration at time instants (a) $t=0.5$, (b) $t=1.0$, and (c) $t=1.5$ from the exact solution (top row), the PINN-H1 solution (middle row) and the PINN-H2 solution (bottom row). 
Plotted here are the deformed field, $\bm x+\bm u(\bm x,t)$, for a set of grid points $\bm x\in D=[0,1]\times [0,1]$.
$N=2000$ training collocation points within domain and on the domain boundaries.
}
	\label{fg_11}
\end{figure}

\begin{figure}[htpb]
	\centering
	\subfloat[$ t=0.5 $]{
		\begin{minipage}[b]{0.25\textwidth}
			\includegraphics[scale=0.21]{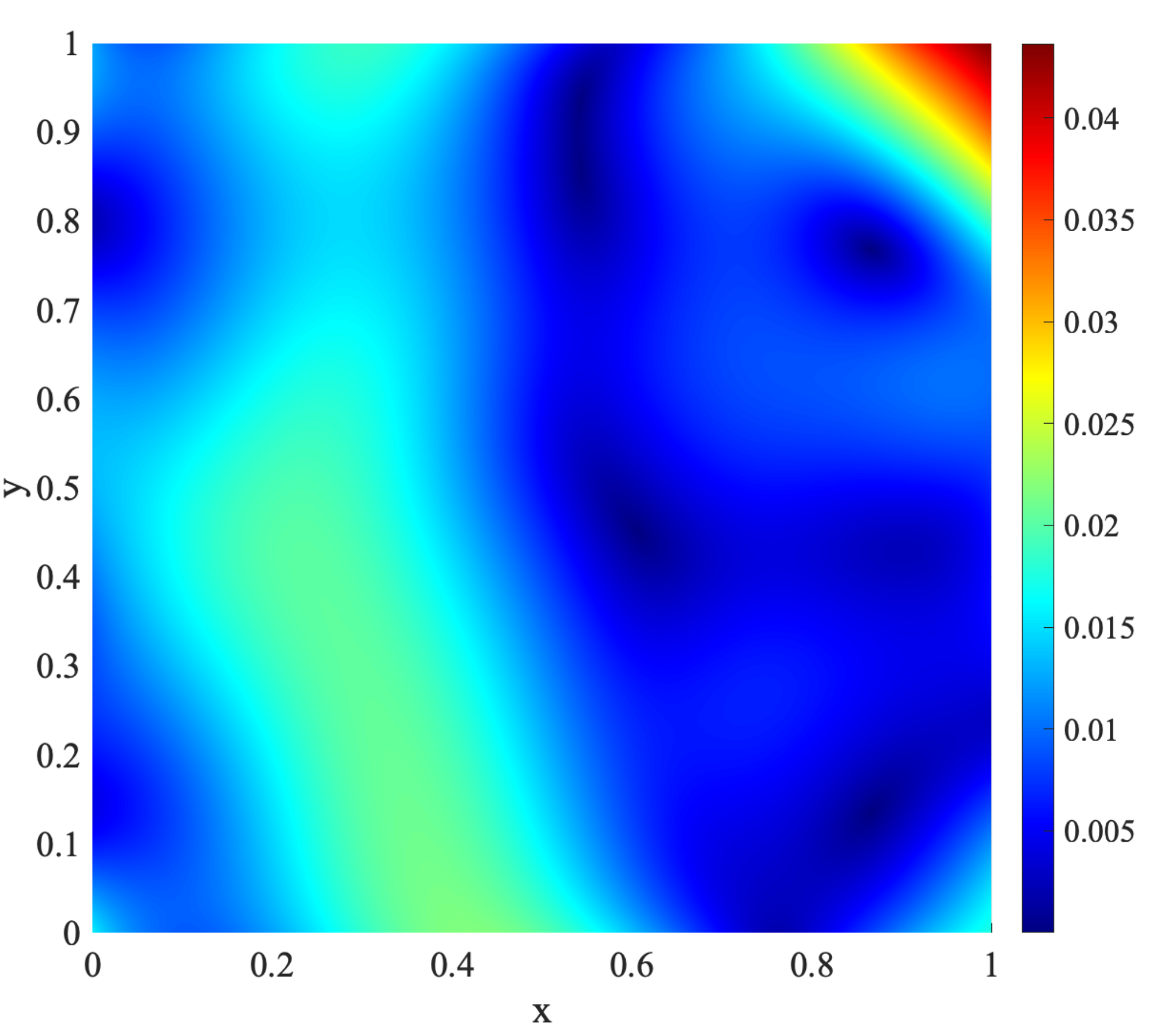}\\
			\includegraphics[scale=0.21]{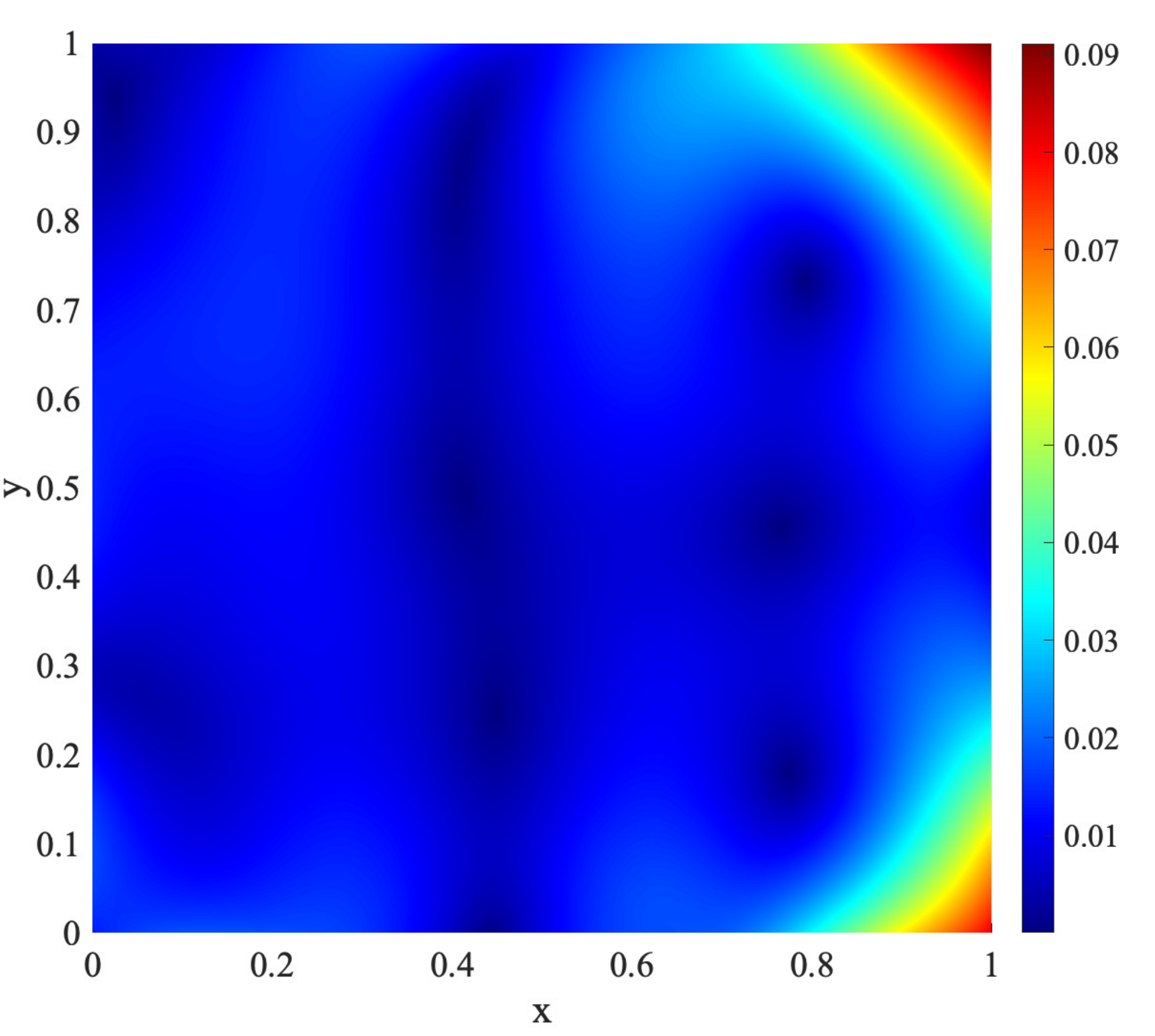}
		\end{minipage}
	}
	\subfloat[$ t=1 $]{
		\begin{minipage}[b]{0.25\textwidth}
			\includegraphics[scale=0.21]{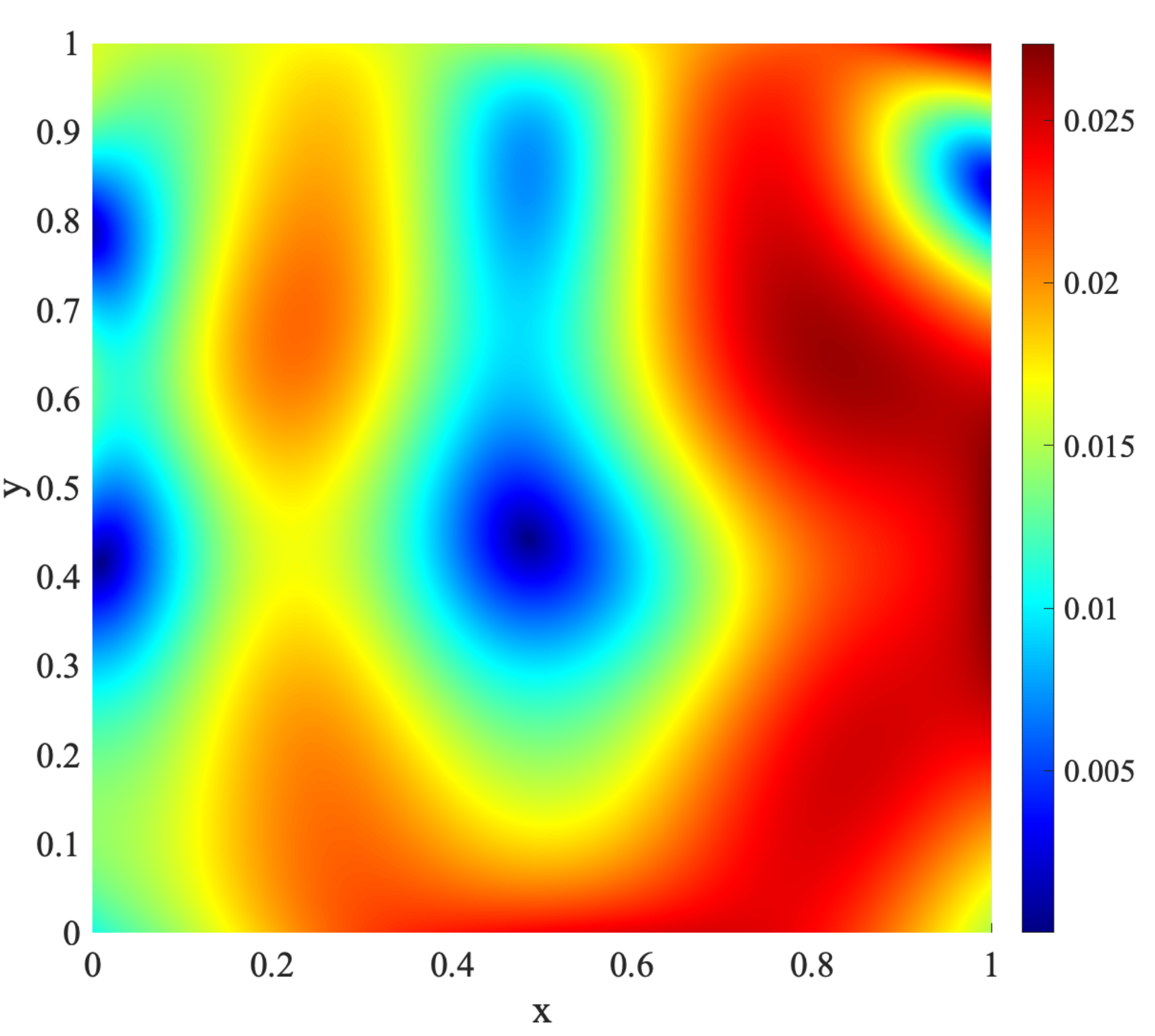}\\
			\includegraphics[scale=0.21]{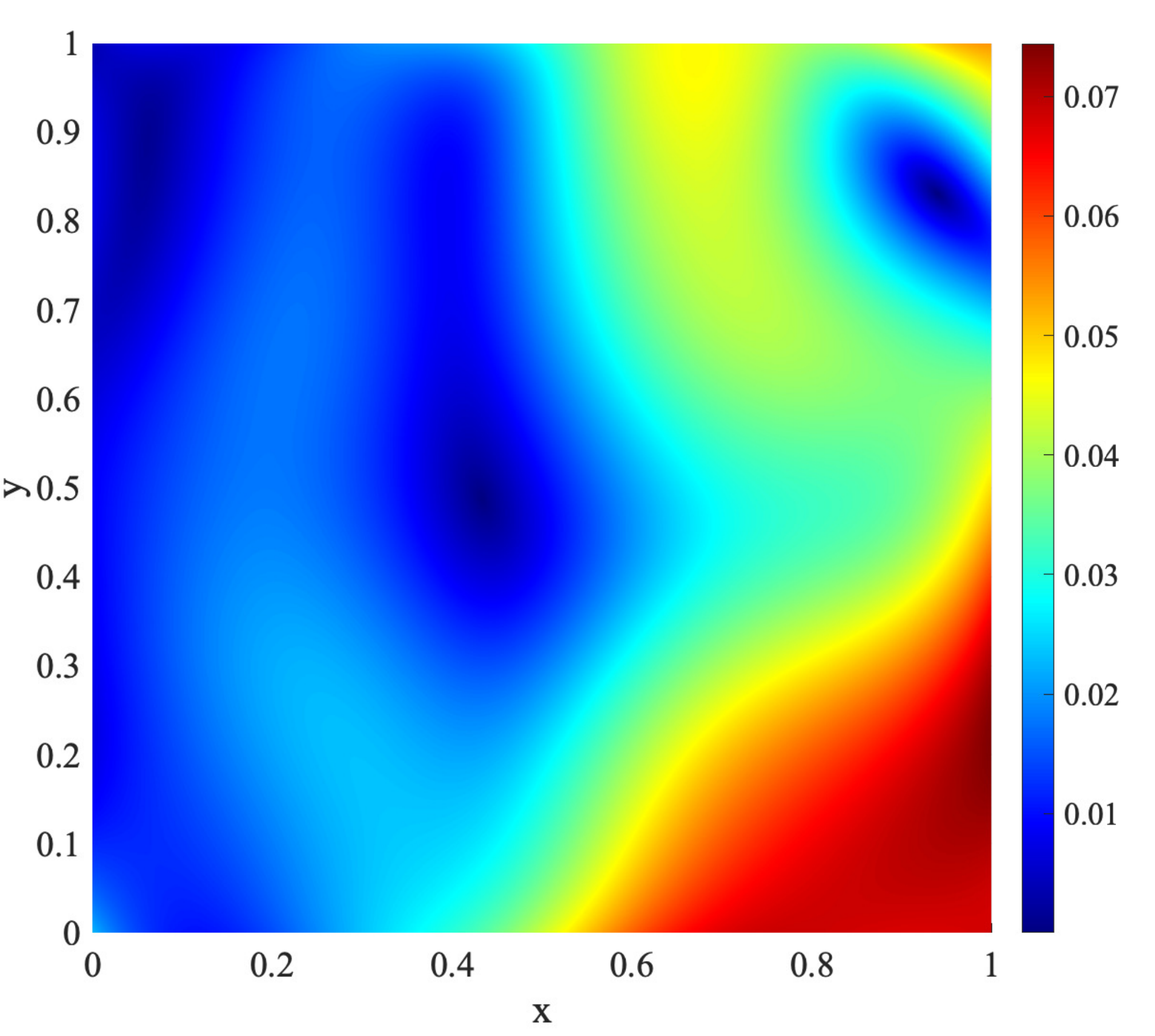}
		\end{minipage}
	}
	\subfloat[$ t=1.5 $]{
		\begin{minipage}[b]{0.25\textwidth}
			\includegraphics[scale=0.21]{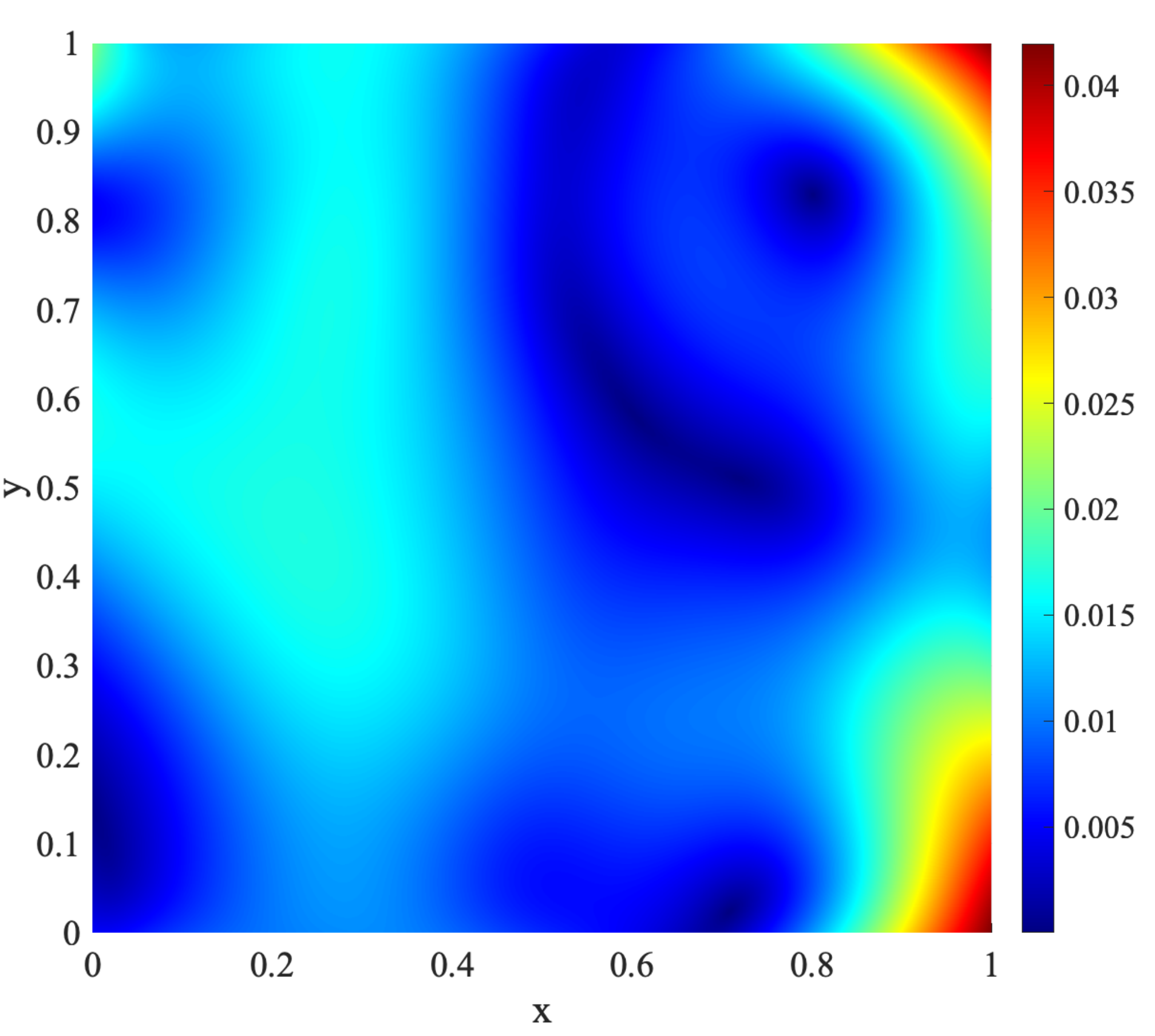}\\
			\includegraphics[scale=0.21]{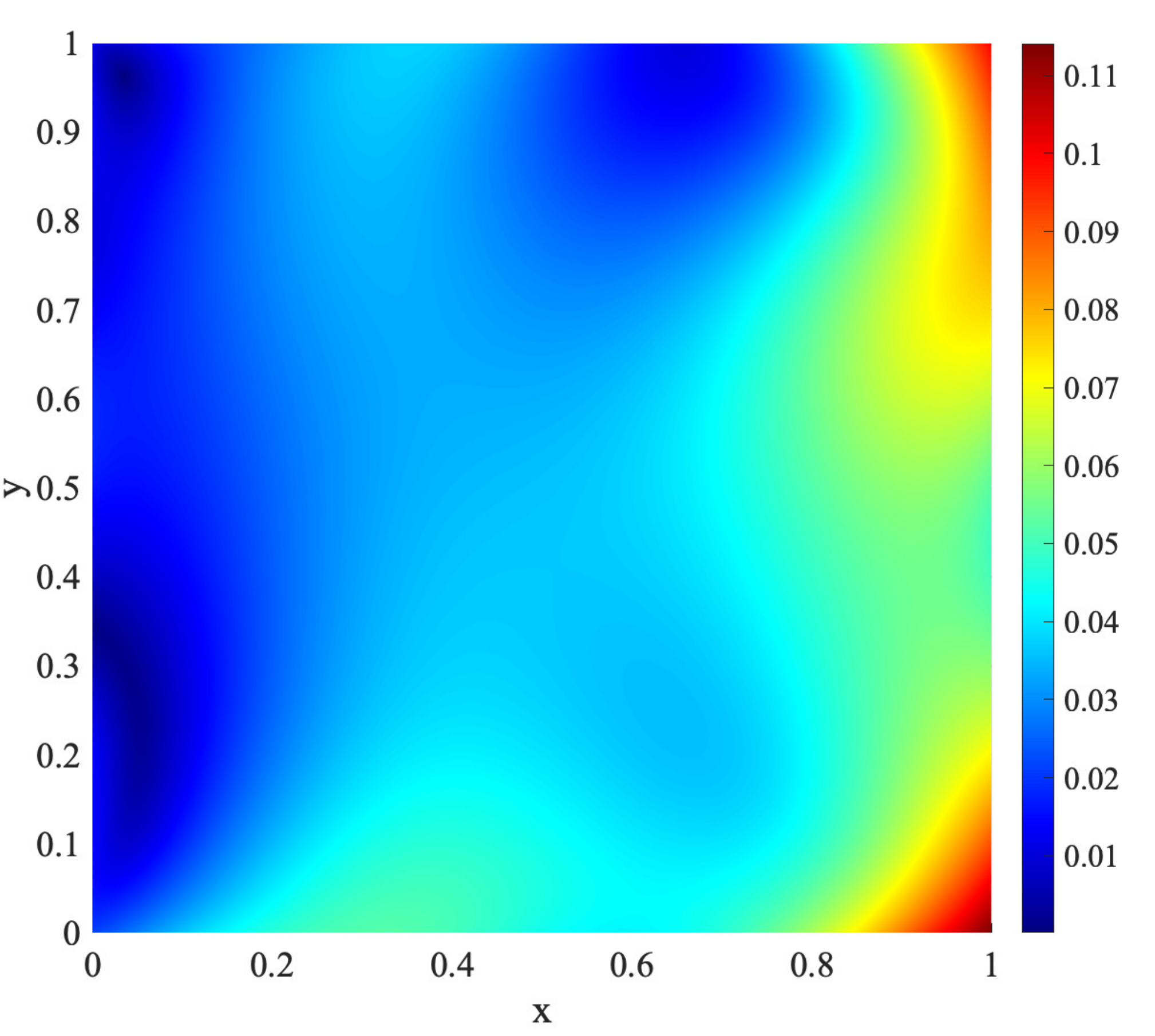}
		\end{minipage}
	}
\caption{Linear elastodynamic equation: 
Distributions of the point-wise absolute error,
$\|\bm u_{\theta}-\bm u \|$, of the PINN-H1 solution (top row) and the PINN-H2 solution (bottom row) at three time instants (a) $t=0.5$,
(b) $t=1.0$, and (c) $t=1.5$.
$N=2000$ training collocation points within domain and on the domain boundaries.
}
	\label{fg_a11}
\end{figure}

In Figures~\ref{fg_11} and~\ref{fg_a11} we compare the PINN-H1/PINN-H2 solutions with the exact solution and provide an overview of their errors.
Figure~\ref{fg_11} is a visualization of the deformed configuration of the domain.
Here we have plotted the deformed field, $\bm x+\bm u(\bm x,t)$, for a set of grid points
$\bm x\in D$ at three time instants from the exact solution, the PINN-H1 and PINN-H2 solutions.
Figure~\ref{fg_a11} shows distributions of 
the point-wise absolute error of the PINN-H1/PINN-H2 solutions,
$\|\bm u_{\theta}-\bm u \|=\sqrt{(u_{\theta 1}(\bm x,t) - u_1(\bm x,t))^2+(u_{\theta 2}(\bm x,t) - u_2(\bm x,t))^2}$, at the same three time instants.
Here $ \bm u_\theta = (u_{\theta 1}, u_{\theta 2}) $ denotes the PINN solution. 
While both  PINN schemes capture the solution fairly well at $ t=0.5$ and $1$, at $ t=1.5 $ both schemes show larger deviations from the true solution. In general, the PINN-H1 scheme appears to produce a better approximation to the solution than PINN-H2.

\begin{table}[tb]\small
    \caption{Linear elastodynamic equation: The $l_2$ and $l_\infty$ errors for $\bm{u}=(u_1, u_2)$ and $\bm{v}=(v_1,v_2)$ versus the number of training data points $N$ from the PINN-H1 and PINN-H2 solutions.}
	\label{PINN_partpaper_Ela_tab_1}
    \centering
	\begin{tabular}[b]{ c@{\ \ }| c  @{\ \ \ }  c  @{\ \ }  c  @{\ \ \ }  c  | c  @{\ \ \ }  c  @{\ \ } c  @{\ \ \ }  c }
		\hline
		\multirow{2}{*}{$ N $}&\multicolumn{4}{c}{$ l_2 $-error} & \multicolumn{4}{c}{$ l_\infty $-error} \\ 
		\cline{2-9}
		{} &$ u_{\theta 1} $  &$ u_{\theta 2} $  &$ v_{\theta 1} $  &$ v_{\theta 2} $ &$ u_{\theta 1} $  &$ u_{\theta 2} $  &$ v_{\theta 1} $  &$ v_{\theta 2} $ \\
		\hline
		{}&\multicolumn{8}{c}{PINN-H1}\\
		\cline{2-9}
		1000&  4.8837e-02&  6.0673e-02&   4.7460e-02&  5.1640e-02&  1.7189e-01&  2.1201e-01&  6.9024e-01&  6.1540e-01\\
		\cline{2-9}
		1500&  2.8131e-02&  3.1485e-02&   4.1104e-02&  4.1613e-02&  1.9848e-01&  2.4670e-01&  3.4716e-01&  4.0582e-01\\
		\cline{2-9}
		2000&  2.7796e-02&  4.0410e-02&   3.5891e-02&  4.6334e-02&  1.4704e-01&  1.7687e-01&  4.0678e-01&  5.0022e-01\\
		\cline{2-9}
		2500&  3.0909e-02&  4.0215e-02&   3.3966e-02&  4.4024e-02&  1.7589e-01&  2.4211e-01&  4.1403e-01&  3.9570e-01\\
		\cline{2-9}
		3000&  2.6411e-02&  3.5600e-02&   4.3209e-02&  5.2802e-02&  1.4289e-01&  1.3625e-01&  5.1167e-01&  5.3298e-01\\
		\hline 
		{}&\multicolumn{8}{c}{PINN-H2}\\
		\cline{2-9}
		1000&  4.9869e-02&  1.3451e-01&   5.6327e-02&  5.4796e-02&  3.2314e-01&  3.4978e-01&  6.7624e-01&  5.7277e-01\\
		\cline{2-9}
		1500&  5.4708e-02&  1.3987e-01&   4.5871e-02&  5.1622e-02&  2.8609e-01&  5.2598e-01&  4.9343e-01&  2.3518e-01\\
		\cline{2-9}
		2000&  6.2114e-02&  1.0190e-01&   6.4477e-02&  5.0011e-02&  2.5745e-01&  3.1642e-01&  5.9057e-01&  5.8411e-01\\
		\cline{2-9}
		2500&  3.7887e-02&  6.0630e-02&   5.4363e-02&  5.0659e-02&  2.2212e-01&  2.4774e-01&  5.3681e-01&  3.5427e-01\\
		\cline{2-9}
		3000&  5.4862e-02&  6.3407e-02&   5.5208e-02&  6.0082e-02&  3.4102e-01&  2.1308e-01&  5.1894e-01&  4.4995e-01\\
		\hline
	\end{tabular}
\end{table}

\begin{figure}[tb]
	\centering
	\subfloat[PINN-H1]{\includegraphics[width=0.4\linewidth]{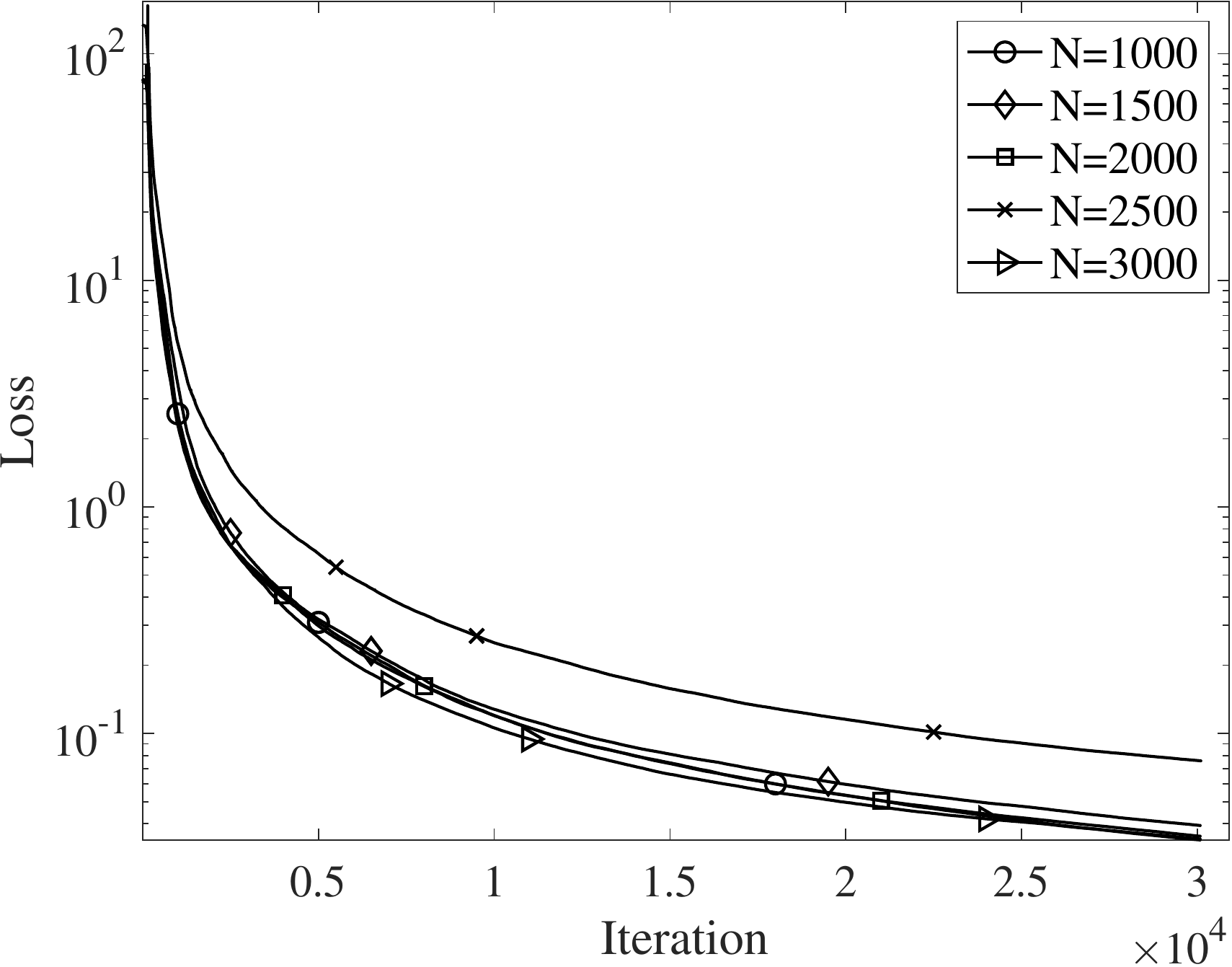}}\qquad
	\subfloat[PINN-H2]{\includegraphics[width=0.4\linewidth]{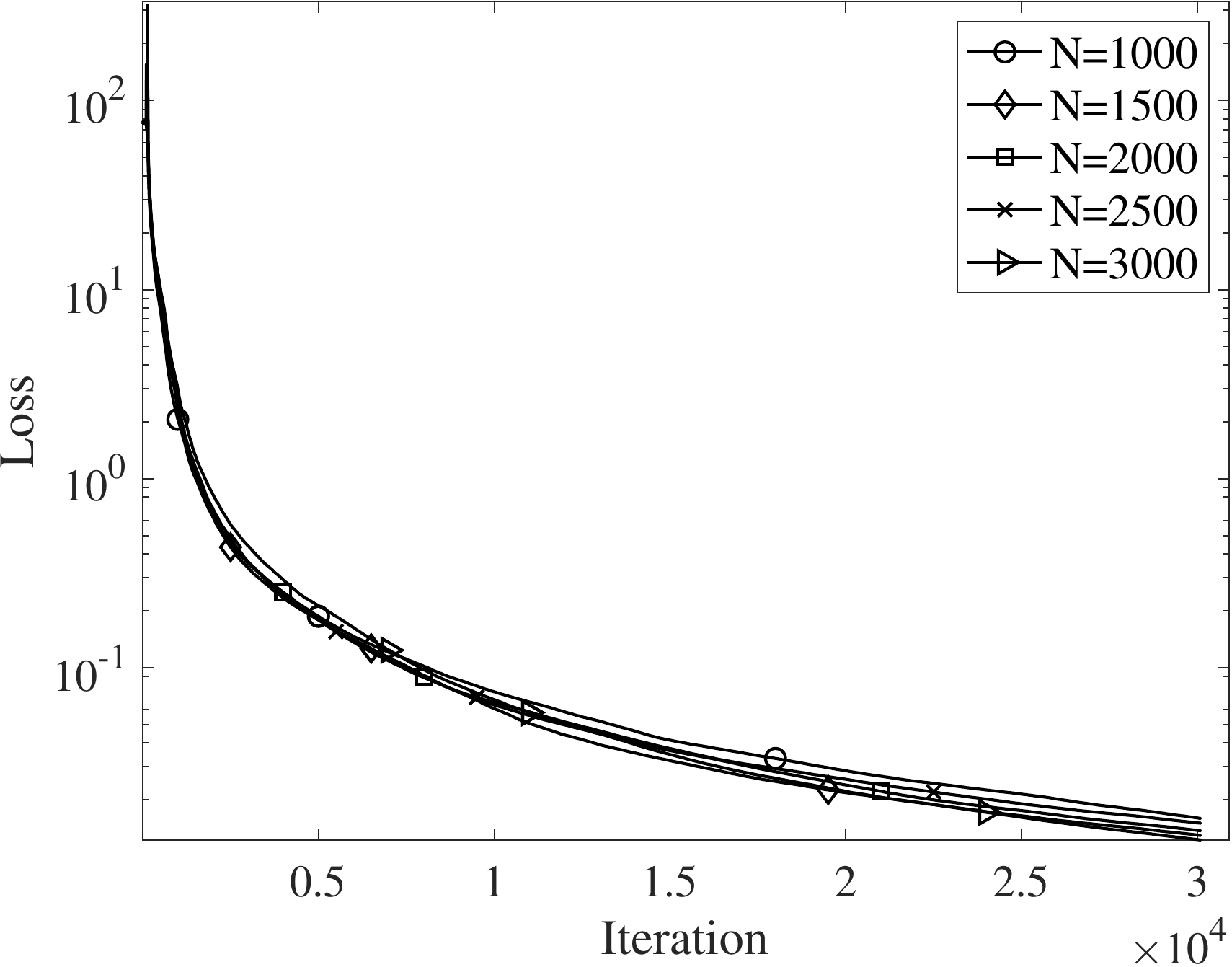}}
	\caption{Linear elastodynamic equation: Training loss histories of PINN-H1 and PINN-H2 corresponding to different numbers of collocation points ($N$) in the simulation.
 }
	\label{PINN_partpaper_Ela_fig2}
\end{figure}

The effect of the number of collocation points ($N$) on the PINN results has been studied in Figure~\ref{PINN_partpaper_Ela_fig2} and Table~\ref{PINN_partpaper_Ela_tab_1}, where $N$ is systematically varied in the range $N=1000$ to $N=3000$. Figure~\ref{PINN_partpaper_Ela_fig2} shows the histories of the loss function for training PINN-H1 and PINN-H2 under different collocation points. Table~\ref{PINN_partpaper_Ela_tab_1} lists the corresponding $l_2$ and $l_{\infty}$ errors of $\bm u$ and $\bm v$ obtained from PINN-H1 and PINN-H2. 
One can observe that the PINN errors in general tend to improve with increasing number of collocation points. It can also be observed that the PINN-H1 errors in general appear better than those of PINN-H2 for this problem.

\begin{figure}[tb]
	\centering
	\subfloat[PINN-H1]{\includegraphics[width=0.4\linewidth]{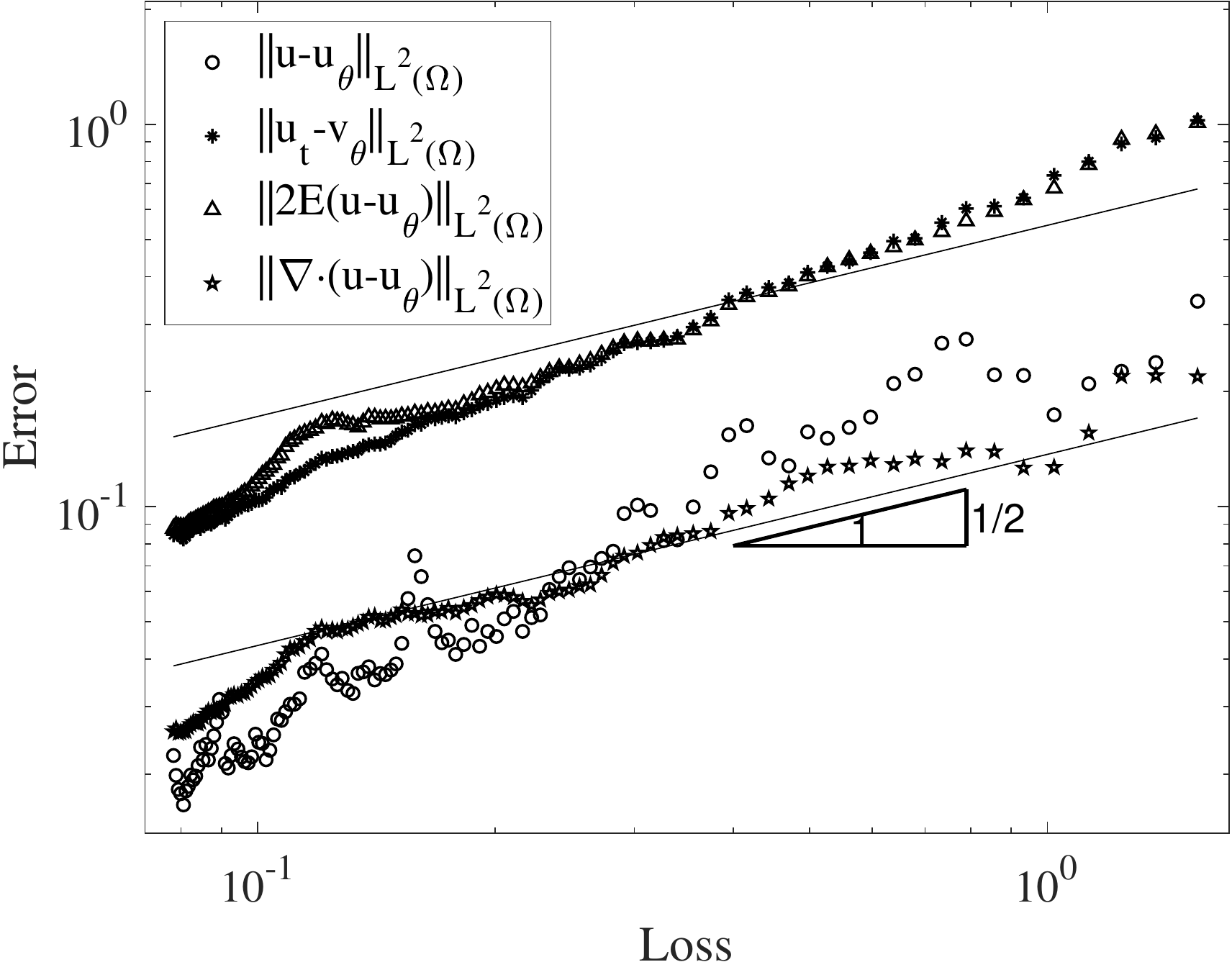}}\qquad
	\subfloat[PINN-H2]{\includegraphics[width=0.4\linewidth]{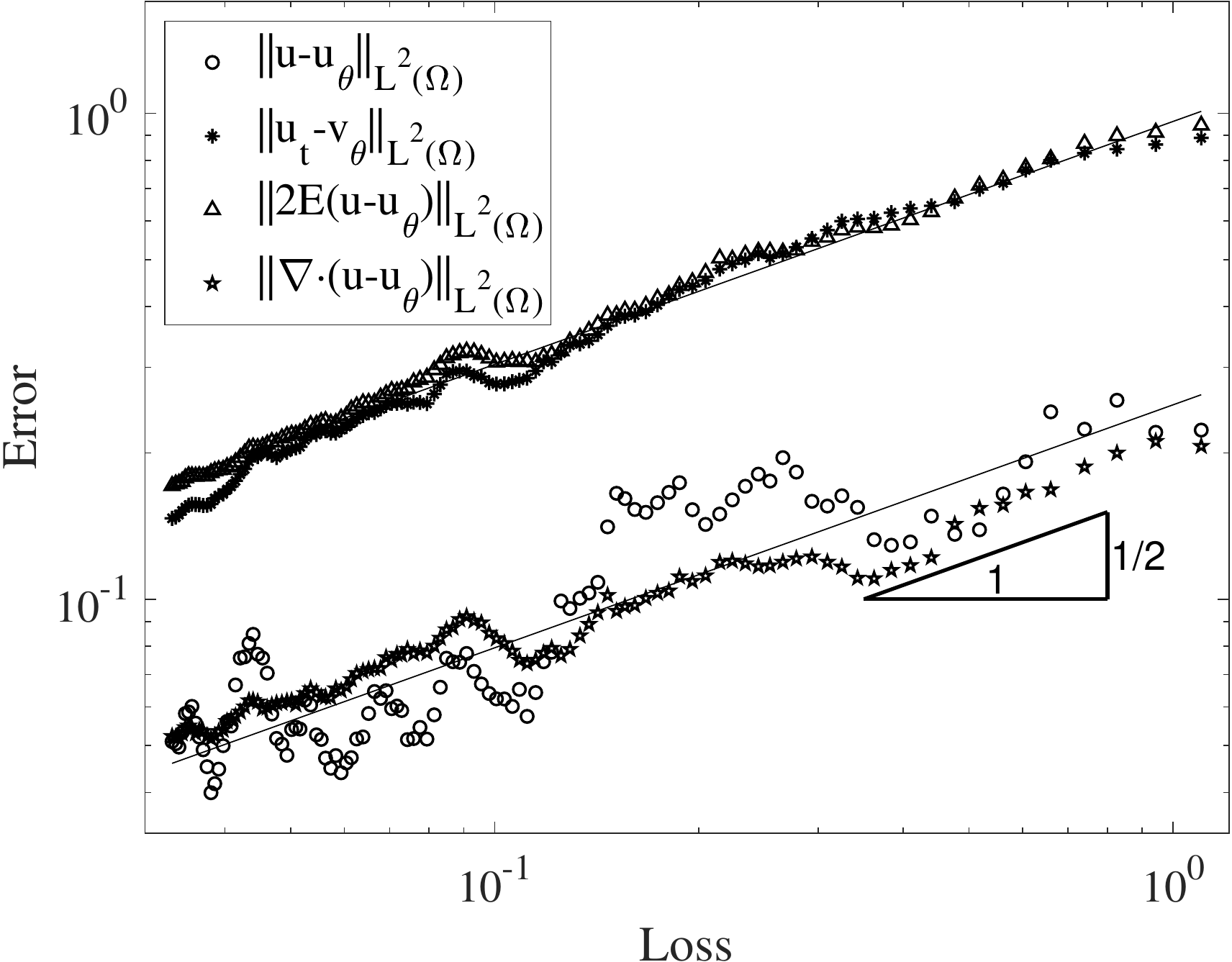}}\hspace{0.1em}
	\caption{Linear elastodynamic equation: The errors for $\bm u$, $\bm u_t$, $\underline{\bm\varepsilon}(\bm u)$ and $\nabla\cdot\bm u$ versus the training loss value obtained by PINN-H1 and PINN-H2.
 }
\label{PINN_partpaper_LinearElasticity_errorrate}
\end{figure}

Figure~\ref{PINN_partpaper_LinearElasticity_errorrate} shows the errors of $\bm u$, $\bm u_t$, $\underline{\bm\varepsilon}(\bm u)$ and
$\nabla\cdot\bm u$ as a function of the loss function value in the network training of PINN-H1 and PINN-H2. The data indicates that these errors approximately scale as the square root of the training loss, which is consistent with the relation as given by Theorem~\ref{sec9_Theorem3}. This in a sense provides numerical evidence for the theoretical analysis in Section~\ref{Elasto-dynamics}.

	\section{Concluding Remarks}
\label{Conclusion}

In the present paper we have considered the approximation of a class of dynamic PDEs of second order in time by physics-informed neural networks (PINN). We provide an analysis of the convergence and the error of PINN for approximating the wave equation, the Sine-Gordon equation, and the linear elastodynamic equation. Our analyses show that, with feed-forward neural networks having two hidden layers and the $\tanh$ activation function for all the hidden nodes, the PINN approximation errors for the solution field, its time derivative and its gradient can be bounded by the PINN training loss and the number of training data points (quadrature points). 

Our theoretical analyses further suggest new forms for the PINN training loss function, which contain certain residuals that are crucial to the error estimate but would be absent from the canonical PINN formulation of the loss function. These typically include the gradient of the equation residual, the gradient of the initial-condition residual, and the time derivative of the boundary-condition residual. In addition, depending on the type of boundary conditions involved in the problem, our analyses suggest that a norm other than the commonly-used $L^2$ norm may be more appropriate for the boundary residuals in the loss function. 
Adopting these new forms of the loss function suggested by the theoretical analyses leads to a variant PINN algorithm. We have implemented the new algorithm and presented a number of numerical experiments on the wave equation, the Sine-Gordon equation and the linear elastodynamic equation. The simulation results demonstrate that the method can capture the solution field well for these PDEs. The numerical data corroborate the theoretical analyses.  


\section*{Declarations}
The authors declare that they have no known competing financial interests or personal relationships that could have appeared to influence the work reported in this paper. 

\section*{Availability of data/code and material }
Data will be made available on reasonable request.

\section*{Acknowledgements}
The work was partially supported by the China Postdoctoral Science Foundation (No.2021M702747), Natural Science Foundation of Hunan Province (No.2022JJ40422), NSF of China (No.12101495), General Special Project of Education Department of Shaanxi Provincial Government (No.21JK0943), and the US National Science Foundation (DMS-2012415).

        \section{Appendix: Auxiliary Results and Proofs of Main Theorems from Sections~\ref{Sine-Gordon} and~\ref{Elasto-dynamics}}\label{Appendix}

\subsection{Notation}\label{Notation}

Let a $d$-tuple of non-negative integers $\alpha\in \mathbb{N}_0^d$ be multi-index with $d \in\mathbb{N}$. 
For given two multi-indices $\alpha,\beta \in \mathbb{N}_0^d$, we say that $\alpha,\beta$, if, and only if, $\alpha_i\leq \beta_i$ for all $i=1,\cdots,d$. And then, denote 
\[
|\alpha| = \sum_{i=1}^d\alpha_i,\qquad \alpha!=\prod_{i=1}^d\alpha_i!, 
\qquad \begin{pmatrix}
	\alpha\\	
	\beta
\end{pmatrix}=\frac{\alpha!}{\beta!(\alpha-\beta)!}.
\]
Let $P_{m,n}=\{\alpha\in \mathbb{N}_0^n,  |\alpha|=m\}$, for which it holds
\[
|P_{m,n}| =\begin{pmatrix}
	m+n-1\\	
	m
\end{pmatrix}.
\]

\subsection{Some Auxiliary Results}\label{Auxiliary lemmas}

\begin{Lemma}\label{Ar_1} Let $d\in \mathbb{N}, k,l\in \mathbb{N}_0$ with $k>l+\frac{d}{2}$ and $\Omega \subset \mathbb{R}^d$ be an open set. Every function $f\in H^k(\Omega)$ has a continuous representative belonging to $C^l(\Omega)$.	
\end{Lemma}

\begin{Lemma}\label{Ar_01} Let $d\in \mathbb{N}, k\in \mathbb{N}_0$, $f\in H^{k}(\Omega)$ and $g\in W^{k,\infty}(\Omega)$ with $\Omega \subset \mathbb{R}^d$, then	
	\[
	\|fg\|_{H^k(\Omega)}\leq 2^k\|f\|_{H^k(\Omega)}\|g\|_{W^{k,\infty}(\Omega)}.
	\]
\end{Lemma}

\begin{Lemma}[Multiplicative trace inequality, e.g. \cite{DeRyck2021On}]\label{Ar_2} 
Let $d\geq 2$, $\Omega \subset \mathbb{R}^d$ be a Lipschitz domain and let $\gamma_0: H^1(\Omega)\rightarrow L^2(\partial\Omega): u \mapsto u|_{\partial\Omega}$ be the trace operator. Denote by $h_{\Omega}$ the diameter of $\Omega$ and by $\rho_{\Omega}$ the radius of the largest $d$-dimensional ball that can be inscribed into $\Omega$. Then it holds that
	\begin{equation}\label{lem_aux1}
		\|\gamma_0 u\|_{L^2(\partial \Omega)}\leq C_{h_{\Omega},d,\rho_{\Omega}}\|u\|_{H^1(\Omega)},
	\end{equation}
where $C_{h_{\Omega},d,\rho_{\Omega}}=\sqrt{\frac{2\max\{2h_{\Omega},d\}}{\rho_{\Omega}}}$.
\end{Lemma}

\begin{Lemma}[\cite{2023_IMA_Mishra_NS}]\label{Ar_3}
Let $d, n, L, W\in \mathbb{N}$ and let $u_{\theta}:\mathbb{R}^{d}\rightarrow \mathbb{R}^{d}$ be a neural network with $\theta\in \Theta$ for $L\geq 2, R, W\geq 1$, c.f. Definition \ref{pre_lem1}. Assume that $\|\sigma\|_{C^n}\geq 1$. Then it holds for $1\leq j\leq d$ that
	\begin{equation}\label{lem_aux2}
		\|(u_{\theta})_j\|_{C^n(\Omega)}\leq 16^Ld^{2n}(e^2n^4W^3R^n\|\sigma\|_{C^n(\Omega)})^{nL}.
	\end{equation}
\end{Lemma}

\begin{Lemma}[\cite{2023_IMA_Mishra_NS}]\label{Ar_4}
Let $d\geq2, m\geq 3, \sigma>0, a_i, b_i \in \mathbb{Z}$ with $a_i <b_i$ for $1\leq i\leq d$, $\Omega=\prod_{i=1}^d[a_i,b_i]$ and $f\in H^m(\Omega)$. Then for every $N\in \mathbb{N}$ with $N>5$ there exists a tanh neural network $\hat{f}^N$ with two hidden layers, one of width at most $3[\frac{m}{2}]|P_{m-1,d+1}|+\sum_{i=1}^d(b_i-a_i)(N-1)$ and another of width at most $3[\frac{d+2}{2}]|P_{d+1,d+1}|N^d\prod_{i=1}^d(b_i-a_i)$, such that for $k=0,1,2$ it holds that
	\begin{equation}\label{lem_aux3}
		\|f-\hat{f}^N\|_{H^k(\Omega)}\leq 2^k3^dC_{k,m,d,f}(1+\sigma){\rm ln}^k\left(\beta_{k,\sigma,d,f}N^{d+m+2}\right)N^{-m+k},
	\end{equation}
	and where
	\begin{align*}
		&C_{k,m,d,f}=\max_{0\leq l \leq k}\left(
		\begin{array}{c}
			d+l-1 \\
			l \\
		\end{array}
		\right)^{1/2}\frac{((m-l)!)^{1/2}}{([\frac{m-l}{d}]!)^{d/2}}\left(\frac{3\sqrt{d}}{\pi}\right)^{m-l}|f|_{H^m(\Omega)},\\
		&\beta_{k,\sigma,d,f}=\frac{5\cdot2^{kd}\max\{\prod_{i=1}^d(b_i-a_i),d\}\max\{\|f\|_{W^{k,\infty}(\Omega)},1\}}{3^d\sigma\min\{1,C_{k,m,d,f}\}}.	
	\end{align*}
	Moreover, the weights of $\hat{f}^N$ scale as $O(N^{\gamma})$ with $\gamma=\max\{m^2/2,d(1+m/2+d/2)\}$.
\end{Lemma}

\subsection{Proof of Main Theorems from Section~\ref{Sine-Gordon}: Sine-Gordon Equation}\label{Proof_sinegordon}

\vspace{0.1in}
\noindent\underline{\bf Theorem \ref{sec6_Theorem1}:}
	Let $d$, $r$, $k \in \mathbb{N}$ with $k\geq 3$. Assume that $g(u)$ is Lipschitz continuous,  $u \in C^k(D\times[0,T])$ and $v \in C^{k-1}(D\times[0,T])$. Then for every integer $N>5$, there exist $\tanh$ neural networks $u_{\theta}$ and $v_{\theta}$, each having two hidden layers, of widths at most $3\lceil\frac{k}{2}\rceil|P_{k-1,d+2}| + \lceil NT\rceil+ d(N-1)$ and $3\lceil\frac{d+3}{2}\rceil|P_{d+2,d+2}| \lceil NT\rceil N^d$, such that
		\begin{align*}
			&\|R_{int1}\|_{L^2(\Omega)},\|R_{tb1}\|_{L^2(D)}\lesssim {\rm ln}NN^{-k+1},\\
			&\|R_{int2}\|_{L^2(\Omega)},\|\nabla R_{int1}\|_{L^2(\Omega)}, \|\nabla R_{tb1}\|_{L^2(D)}\lesssim {\rm ln}^2NN^{-k+2},\\
			&\|R_{tb2}\|_{L^2(D)},\|R_{sb}\|_{L^2(\partial D\times [0,t])}\lesssim {\rm ln}NN^{-k+2}.
		\end{align*}
\begin{proof} Based on $u \in C^k(D\times[0,T])$, $v \in C^{k-1}(D\times[0,T])$ and Lemma \ref{Ar_4}, there exist neural networks $u_{\theta}$ and $v_{\theta}$, with the same two hidden layers and widths $3\lceil\frac{k}{2}\rceil|P_{k-1,d+2}| + \lceil NT\rceil+ d(N-1)$ and $3\lceil\frac{d+3}{2}\rceil|P_{d+2,d+2}| \lceil NT\rceil N^d$, such that for every $0 \leq l \leq 2$ and $0\leq s\leq 2$,
	\begin{align*}
		&\|u_{\theta}-u\|_{H^l(\Omega)}\leq C_{l,k,d+1,u}\lambda_{l,u}(N)N^{-k+l},\\
		&\|v_{\theta}-v\|_{H^{s}(\Omega)}\leq C_{s,k-1,d+1,v}\lambda_{s,v}(N)N^{-k+1+s}.
	\end{align*}
	It is now straightforward to bound the PINN residual.
	\begin{align*}
		&\|\hat{u}_t\|_{L^2(\Omega)}\leq\|\hat{u}\|_{H^1(\Omega)},\qquad \|\hat{v}_t\|_{L^2(\Omega)}\leq\|\hat{v}\|_{H^1(\Omega)},\\
		&\|\Delta\hat{u}\|_{L^2(\Omega)}\leq\|\hat{u}\|_{H^2(\Omega)}\qquad \|\nabla\hat{u}_t\|_{L^2(\Omega)}\leq\|\hat{u}\|_{H^2(\Omega)},\\
		&\|\nabla\hat{v}\|_{L^2(\Omega)}\leq\|\hat{v}\|_{H^1(\Omega)},\\
		&\|\hat{u}\|_{L^2(D)}\leq \|\hat{u}\|_{L^2(\partial\Omega)}\leq C_{h_{\Omega},d+1,\rho_{\Omega}}\|\hat{u}\|_{H^1(\Omega)},\\
		&\|\hat{v}\|_{L^2(D)}\leq \|\hat{v}\|_{L^2(\partial\Omega)}\leq C_{h_{\Omega},d+1,\rho_{\Omega}}\|\hat{v}\|_{H^1(\Omega)},\\
		&\|\nabla\hat{u}\|_{L^2(D)}\leq \|\nabla\hat{u}\|_{L^2(\partial\Omega)}\leq C_{h_{\Omega},d+1,\rho_{\Omega}}\|\hat{u}\|_{H^2(\Omega)},\\
		&\|\hat{v}\|_{L^2(\partial D\times [0,t])}\leq \|\hat{v}\|_{L^2(\partial\Omega)}\leq C_{h_{\Omega},d+1,\rho_{\Omega}}\|\hat{v}\|_{H^1(\Omega)}.
	\end{align*}
	 Similar to Theorem \ref{sec5_Theorem1}, we can obtain                         
	\begin{align*}
		&\|R_{int1}\|_{L^2(\Omega)}=\|\hat{u}_t-\hat{v}\|_{L^2(\Omega)}\leq
		\|\hat{u}\|_{H^1(\Omega)}+\|\hat{v}\|_{L^2(\Omega)}\lesssim {\rm ln}NN^{-k+1},\\
		&\|R_{int2}\|_{L^2(\Omega)}=\|\varepsilon^2\hat{v}_{t}-a^2\Delta \hat{u} +\varepsilon_1^2\hat{u}+g(u_{\theta})-g(u)\|_{L^2(\Omega)}\\
		&\qquad\leq
		\varepsilon^2\|\hat{v}\|_{H^1(\Omega)}+a^2\|\hat{u}\|_{H^2(\Omega)}+\varepsilon_1^2\|\hat{u}\|_{L^2(\Omega)}+L\|\hat{u}\|_{L^2(\Omega)}\lesssim {\rm ln}^2NN^{-k+2},\\
		&\|\nabla R_{int1}\|_{L^2(\Omega)}=\|\nabla(\hat{u}_t-\hat{v})\|_{L^2(\Omega)}\leq
		\|\hat{u}\|_{H^2(\Omega)}+\|\hat{v}\|_{H^1(\Omega)}\lesssim {\rm ln}^2NN^{-k+2},\\
		&\|R_{tb1}\|_{L^2(D)}\leq C_{h_{\Omega},d+1,\rho_{\Omega}}\|\hat{u}\|_{H^1(\Omega)}\lesssim {\rm ln}NN^{-k+1},\\
		&\|R_{tb2}\|_{L^2(D)}, \|R_{sb}\|_{L^2(\partial D\times [0,t])}\leq C_{h_{\Omega},d+1,\rho_{\Omega}}\|\hat{v}\|_{H^1(\Omega)}\lesssim {\rm ln}NN^{-k+2},\\
		&\|\nabla R_{tb1}\|_{L^2(D)}\leq C_{h_{\Omega},d+1,\rho_{\Omega}}\|\hat{u}\|_{H^2(\Omega)}\lesssim {\rm ln}^2NN^{-k+2}.
	\end{align*}
\end{proof}

\vspace{0.1in}
\noindent\underline{\bf Theorem \ref{sec6_Theorem2}:}
Let $d\in \mathbb{N}$, $u\in C^1(\Omega)$ and $v\in C^0(\Omega)$ be the classical solution to
	the Sine-Gordon equation \eqref{SG}. Let $(u_{\theta},v_{\theta})$ denote the PINN approximation with the parameter $\theta$. Then the following relation holds,
	\[
		\int_0^{T}\int_{D}(|\hat{u}(\bm{x},t)|^2+a^2|\nabla \hat{u}(\bm{x},t)|^2+\varepsilon^2|\hat{v}(\bm{x},t)|^2)\dx\dt
		\leq C_GT\exp\left((2+\varepsilon_1^2+L+a^2)T\right),
	\]
	where $C_G$ is defined in the proof.
\begin{proof} By taking the inner product of \eqref{SG_error_eq1} and \eqref{SG_error_eq2} with $\hat{u}$ and $\hat{v}$ over $D$, respectively, we have
	\begin{align}
		\label{sec6_eq0}
		\frac{d}{2dt}\int_{D} |\hat{u}|^2\dx &= \int_{D}\hat{u}\hat{v}\dx+\int_{D} R_{int1}\hat{u}\dx\leq \int_{D} |\hat{u}|^2\dx+\frac{1}{2}\int_{D} |R_{int1}|^2\dx+\frac{1}{2}\int_{D} |\hat{v}|^2\dx,\\
		\label{sec6_eq1}
		\varepsilon^2\frac{d}{2dt}\int_{D} |\hat{v}|^2\dx &=- a^2\int_{D}\nabla\hat{u}\cdot\nabla\hat{v}\dx+a^2\int_{\partial D} R_{sb}\nabla\hat{u}\cdot\bm{n}\ds-\varepsilon_1^2\int_{D}\hat{u}\hat{v}\dx
		\nonumber\\
		&\qquad-\int_{D} (g(u_{\theta})-g(u))\hat{v}\dx
		+\int_{D} R_{int2}\hat{v}\dx
		\nonumber\\
		&=- a^2\int_{D}\nabla\hat{u}\cdot\nabla\hat{u}_t\dx
		+a^2\int_{D}\nabla\hat{u}\cdot\nabla R_{int1}\dx
		+a^2\int_{\partial D} R_{sb}\nabla\hat{u}\cdot\bm{n}\ds-\varepsilon_1^2\int_{D}\hat{u}\hat{v}\dx
		\nonumber\\
		&\qquad-\int_{D} (g(u_{\theta})-g(u))\hat{v}\dx
		+\int_{D} R_{int2}\hat{v}\dx
		\nonumber\\
		&=-a^2\frac{d}{2dt}\int_{D} |\nabla\hat{u}|^2\dx 
		+a^2\int_{D}\nabla\hat{u}\cdot\nabla R_{int1}\dx
		+a^2\int_{\partial D} R_{sb}\nabla\hat{u}\cdot\bm{n}\ds-\varepsilon_1^2\int_{D}\hat{u}\hat{v}\dx
		\nonumber\\
		&\qquad-\int_{D} (g(u_{\theta})-g(u))\hat{v}\dx
		+\int_{D} R_{int2}\hat{v}\dx
		\nonumber\\
		&\leq-a^2\frac{d}{2dt}\int_{D} |\nabla\hat{u}|^2\dx +\frac{a^2}{2}\int_{D} |\nabla\hat{u}|^2\dx+\frac{a^2}{2}\int_{D} |\nabla R_{int1}|^2\dx+C_{\partial D}\left(\int_{\partial D}|R_{sb}|^2\ds\right)^{\frac{1}{2}}
		\nonumber\\
		&\qquad+\frac{1}{2}(\varepsilon_1^2+L)\int_{D} |\hat{u}|^2\dx+\frac{1}{2}(\varepsilon_1^2+L+1)\int_{D} |\hat{v}|^2\dx+\frac{1}{2}\int_{D} |R_{int2}|^2\dx,
	\end{align}
	where $C_{\partial D}=a^2|\partial D|^{\frac{1}{2}}(\|u\|_{C^1(\partial D\times[0,t])}+||u_{\theta}||_{C^1(\partial D\times[0,t])})$ and $\hat{v}=\hat{u}_t-R_{int1}$ have been used. 
	
 Add \eqref{sec6_eq0} to \eqref{sec6_eq1}, and we get
	\begin{align}\label{sec6_eq2}
		&\frac{d}{2dt}\int_{D} |\hat{u}|^2\dx+a^2\frac{d}{2dt}\int_{D} |\nabla\hat{u}|^2\dx +\varepsilon^2\frac{d}{2dt}\int_{D} |\hat{v}|^2\dx
		\nonumber\\
		&\qquad \leq \frac{1}{2}(\varepsilon_1^2+L+2)\int_{D} |\hat{u}|^2\dx+\frac{a^2}{2}\int_{D} |\nabla\hat{u}|^2\dx+\frac{1}{2}(\varepsilon_1^2+L+2)\int_{D} |\hat{v}|^2\dx+\frac{1}{2}\int_{D} |R_{int1}|^2\dx+\frac{1}{2}\int_{D} |R_{int2}|^2\dx
		\nonumber\\
		&\qquad+\frac{a^2}{2}\int_{D} |\nabla R_{int1}|^2\dx+C_{\partial D}\left(\int_{\partial D}|R_{sb}|^2\ds\right)^{\frac{1}{2}}.
	\end{align}
	Integrating \eqref{sec6_eq2} over $[0,
	\tau]$ for any $\tau \leq T$ and applying the Cauchy–Schwarz inequality, we obtain
	\begin{align*}
		&\int_{D} |\hat{u}(\bm{x},\tau)|^2\dx 
		+a^2\int_{D} |\nabla\hat{u}(\bm{x},\tau)|^2\dx
		+\varepsilon^2\int_{D} |\hat{v}(\bm{x},\tau)|^2\dx\\
		&\qquad\leq\int_{D}|R_{tb1}|^2\dx +a^2\int_{D}|\nabla R_{tb1}|^2\dx+\varepsilon^2\int_{D}|R_{tb2}|^2\dx+
		(2+\varepsilon_1^2+L+a^2)\int_{0}^{\tau}\int_{D} \left(|\hat{u}|^2+ |\nabla\hat{u}|^2+ |\hat{v}|^2 \right)\dx\dt
		\\
		&\qquad+ \int_{0}^{T}\int_{D}\left(|R_{int1}|^2+a^2|\nabla R_{int1}|^2+|R_{int2}|^2\right)\dx\dt+2C_{\partial D}|T|^{\frac{1}{2}}\left(\int_{0}^{T}\int_{\partial D}|R_{sb}|^2\ds\dt\right)^{\frac{1}{2}}.
	\end{align*}
	Applying the integral form of the Gr${\rm\ddot{o}}$nwall inequality to the above inequality leads to,
	\begin{equation}\label{sec6_eq3}
		\int_{D} |\hat{u}(\bm{x},\tau)|^2\dx 
		+a^2\int_{D} |\nabla\hat{u}(\bm{x},\tau)|^2\dx
		+\varepsilon^2\int_{D} |\hat{v}(\bm{x},\tau)|^2\dx
		\leq C_G\exp\left((2+\varepsilon_1^2+L+a^2)T\right),
	\end{equation}
	where
	\begin{align*}
		&C_G=\int_{D}(|R_{tb1}|^2+a^2|\nabla R_{tb1}|^2+\varepsilon^2|R_{tb2}|^2)\dx + \int_{0}^{T}\int_{D}(|R_{int1}|^2+|R_{int2}|^2+a^2|\nabla R_{int1}|^2)\dx\dt
		\\
		&\qquad +  2C_{\partial D}|T|^{\frac{1}{2}}\left(\int_{0}^{T}\int_{\partial D}|R_{sb}|^2\ds\dt\right)^{\frac{1}{2}}.
	\end{align*}
	Then, we integrate \eqref{sec6_eq3} over $[0,T]$ to end the proof.
\end{proof} 

\vspace{0.1in}
\noindent\underline{\bf Theorem \ref{sec6_Theorem3}:}
Let $d\in \mathbb{N}$ and $T>0$. Let $u\in C^4(\Omega)$ and $v\in C^3(\Omega)$ be the classical solution to	the Sine-Gordon equation \eqref{SG}. Let $(u_{\theta},v_{\theta})$ denote the PINN approximation with the parameter $\theta \in \Theta$. Then the following relation holds,
	\begin{align*}
		&\int_0^{T}\int_{D}(|\hat{u}(\bm{x},t)|^2+a^2|\nabla \hat{u}(\bm{x},t)|^2+\varepsilon^2|\hat{v}(\bm{x},t)|^2)\dx\dt\leq C_TT\exp\left((2+\varepsilon_1^2+L+a^2)T\right) 
		\nonumber\\
		&\qquad=\mathcal{O}(\mathcal{E}_T(\theta)^2 + M_{int}^{-\frac{2}{d+1}} +M_{tb}^{-\frac{2}{d}}+M_{sb}^{-\frac{1}{d}}),   
    \end{align*}
	where the constant $C_T$ is given in the proof.
\begin{proof} We can combine Theorem \ref{sec6_Theorem2} with the quadrature error formula \eqref{int1} to obtain the  error estimate,
	\begin{align*}
		\int_{D}|R_{tb1}|^2\dx&=\int_{D}|R_{tb1}|^2\dx-\mathcal{Q}_{M_{tb}}^{D}(R_{tb1}^2)+\mathcal{Q}_{M_{tb}}^{D}(R_{tb1}^2)\\
		&\leq C_{({R_{tb1}^2})}M_{tb}^{-\frac{2}{d}}+\mathcal{Q}_{M_{tb}}^{D}(R_{tb1}^2),\\
		\int_{D}|R_{tb2}|^2\dx&=\int_{D}|R_{tb2}|^2\dx-\mathcal{Q}_{M_{tb}}^{D}(R_{tb2}^2)+\mathcal{Q}_{M_{tb}}^{D}(R_{tb2}^2)\\
		&\leq C_{({R_{tb2}^2})}M_{tb}^{-\frac{2}{d}}+\mathcal{Q}_{M_{tb}}^{D}(R_{tb2}^2),\\
		\int_{D}|\nabla R_{tb1}|^2\dx&=\int_{D}|\nabla R_{tb1}|^2\dx-\mathcal{Q}_{M_{tb}}^{D}(|\nabla R_{tb1}|^2)+\mathcal{Q}_{M_{tb}}^{D}(|\nabla R_{tb1}|^2)\\
		&\leq C_{(|\nabla R_{tb1}|^2)}M_{tb}^{-\frac{2}{d}}+\mathcal{Q}_{M_{tb}}^{D}(|\nabla R_{tb1}|^2),\\
		\int_{\Omega}|R_{int1}|^2\dx\dt&=\int_{\Omega}|R_{int1}|^2\dx\dt-\mathcal{Q}_{M_{int}}^{\Omega}(R_{int1}^2)+\mathcal{Q}_{M_{int}}^{\Omega}(R_{int1}^2)\\
		&\leq C_{({R_{int1}^2})}M_{int}^{-\frac{2}{d+1}}+\mathcal{Q}_{M_{int}}^{\Omega}(R_{int1}^2),\\
		\int_{\Omega}|R_{int2}|^2\dx\dt&=\int_{\Omega}|R_{int2}|^2\dx\dt-\mathcal{Q}_{M_{int}}^{\Omega}(R_{int2}^2)+\mathcal{Q}_{M_{int}}^{\Omega}(R_{int2}^2)\\
		&\leq C_{({R_{int2}^2})}M_{int}^{-\frac{2}{d+1}}+\mathcal{Q}_{M_{int}}^{\Omega}(R_{int2}^2),\\
		\int_{\Omega}|\nabla R_{int1}|^2\dx\dt&=\int_{\Omega}|\nabla R_{int1}|^2\dx\dt-\mathcal{Q}_{M_{int}}^{\Omega}(|\nabla R_{int1}|^2)+\mathcal{Q}_{M_{int}}^{\Omega}(|\nabla R_{int1}|^2)\\
		&\leq C_{(|\nabla R_{int1}|^2)}M_{int}^{-\frac{2}{d+1}}+\mathcal{Q}_{M_{int}}^{\Omega}(|\nabla R_{int1}|^2),\\
		\int_{\Omega_*}|R_{sb}|^2\ds\dt&=\int_{\Omega_*}|R_{sb}|^2\ds\dt-\mathcal{Q}_{M_{sb}}^{\Omega_*}(R_{sb}^2)+\mathcal{Q}_{M_{sb}}^{\Omega_*}(R_{sb}^2)\\
		&\leq C_{({R_{sb}^2})}M_{sb}^{-\frac{2}{d}}+\mathcal{Q}_{M_{sb}}^{\Omega_*}(R_{sb}^2).
	\end{align*}
	In light of \eqref{sec6_eq3} and the above inequalities, we have
	\begin{equation*}
		\int_0^{T}\int_{D}(|\hat{u}(\bm{x},t)|^2+a^2|\nabla \hat{u}(\bm{x},t)|^2+\varepsilon^2|\hat{v}(\bm{x},t)|^2)\dx\dt
		\leq TC_T\exp\left((2+\varepsilon_1^2+L+a^2)T\right),
	\end{equation*}
	where 
	\begin{align*}
		C_T=&C_{({R_{tb1}^2})}M_{tb}^{-\frac{2}{d}}+\mathcal{Q}_{M_{tb}}^{D}(R_{tb1}^2)+\varepsilon^2\left(C_{({R_{tb2}^2})}M_{tb}^{-\frac{2}{d}}+\mathcal{Q}_{M_{tb}}^{D}(R_{tb2}^2) \right)\\
		&+a^2\left( C_{(|\nabla R_{tb1}|^2)}M_{tb}^{-\frac{2}{d}}+\mathcal{Q}_{M_{tb}}^{D}(|\nabla R_{tb1}|^2) \right)+C_{({R_{int1}^2})}M_{int}^{-\frac{2}{d+1}}+\mathcal{Q}_{M_{int}}^{\Omega}(R_{int1}^2)\\
		&+C_{({R_{int2}^2})}M_{int}^{-\frac{2}{d+1}}+\mathcal{Q}_{M_{int}}^{\Omega}(R_{int2}^2)
		+a^2\left(C_{(|\nabla R_{int1}|^2)}M_{int}^{-\frac{2}{d+1}}+\mathcal{Q}_{M_{int}}^{\Omega}(|\nabla R_{int1}|^2)\right),\\
		&+2C_{\partial D}|T|^{\frac{1}{2}}\left(C_{({R_{sb}^2})}M_{sb}^{-\frac{2}{d}}+\mathcal{Q}_{M_{sb}}^{\Omega_*}(R_{sb}^2)\right)^{\frac{1}{2}},
	\end{align*}
	and
	\begin{align*}
		&C_{({R_{tb1}^2})}\lesssim\|\hat{u}\|_{C^2}^2, \quad C_{({R_{tb2}^2})}\lesssim \|\hat{v}\|_{C^2}^2, \quad C_{(|\nabla R_{tb1}|^2)}\lesssim \|\hat{u}\|_{C^3}^2, \quad C_{({R_{int1}^2})}\lesssim \|\hat{u}\|_{C^3}^2+\|\hat{v}\|_{C^2}^2,\\
        &\qquad \qquad C_{({R_{int2}^2})}, C_{(|\nabla R_{int1}|^2)}\lesssim \|\hat{u}\|_{C^4}^2+\|\hat{v}\|_{C^3}^2,\quad C_{({R_{sb}^2})}\lesssim \|\hat{v}\|_{C^3}^2.
	\end{align*}
    Here, the boundedness $\|u_{\theta}\|_{C^n}$ and $\|v_{\theta}\|_{C^n}$ ($n \in \mathbb{N}$) of the above constants can be obtained by Lemma \ref{Ar_3} and $\|R_q^2\|_{C^n}\leq 2^n\|R_q\|_{C^n}^2$ for $R_q=R_{tb1}$, $R_{tb2}$, $\nabla R_{tb1}$, $R_{int1}$, $R_{int2}$, $\nabla R_{int1}$ and $R_{sb}$.  
\end{proof} 

\subsection{Proof of Main Theorems from Section~\ref{Elasto-dynamics}: Linear Elastodynamic Equation}\label{Proof_elasto}

\vspace{0.1in}
\noindent\underline{\bf Theorem \ref{sec9_Theorem1}:}
	Let $d$, $r$, $k \in \mathbb{N}$ with $k\geq 3$. Let $\bm{\psi}_{1}\in H^{r}(D)$, $\bm{\psi}_{2}\in H^{r-1}(D)$ and $\bm{f}\in H^{r-1}(D\times [0,T])$ with $r>\frac{d}{2}+k$. 
	For every integer $N>5$, there exist $\tanh$ neural networks $(\bm{u}_j)_{\theta}$ and $(\bm{v}_j)_{\theta}$, with $j=1,2,\cdots,d$, each with two hidden layers, of widths at most $3\lceil\frac{k}{2}\rceil|P_{k-1,d+2}| + \lceil NT\rceil+ d(N-1)$ and $3\lceil\frac{d+3}{2}\rceil|P_{d+2,d+2}| \lceil NT\rceil N^d$, such that
\begin{align*}
			&\|\bm{R}_{int1}\|_{L^2(\Omega)},\|\bm{R}_{tb1}\|_{L^2(\Omega)}\lesssim {\rm ln}NN^{-k+1},\\
			&\|\bm{R}_{int2}\|_{L^2(\Omega)},\|\underline{\bm{\varepsilon}}(\bm{R}_{int1})\|_{L^2(\Omega)},\|\nabla\cdot\bm{R}_{int1}\|_{L^2(\Omega)}\lesssim {\rm ln}^2NN^{-k+2},\\
			&\|\underline{\bm{\varepsilon}}(\bm{R}_{tb1})\|_{L^2(D)},\|\nabla\cdot\bm{R}_{tb1}\|_{L^2(D)},\|\bm{R}_{sb2}\|_{L^2(\Gamma_N\times [0,t])}\lesssim {\rm ln}^2NN^{-k+2},\\
			&\|\bm{R}_{tb2}\|_{L^2(D)}, \|\bm{R}_{sb1}\|_{L^2(\Gamma_D\times [0,t])}\lesssim {\rm ln}NN^{-k+2}.
\end{align*}
\begin{proof} Lemma \ref{sec9_Lemma2} implies that, 
	\[	
	\bm{u} \in C^k(D\times [0,T]),\qquad \bm{v} \in C^{k-1}(D\times [0,T]).
	\]
	Let $\bm{u}_{\theta}=((u_1)_{\theta},(u_2)_{\theta},\cdots,(u_d)_{\theta})$ and $\bm{v}_{\theta}=((v_1)_{\theta},(v_2)_{\theta},\cdots,(v_d)_{\theta})$. Based on Lemma \ref{Ar_4}, there exists $\tanh$ neural networks $(u_i)_{\theta}$ and $(v_i)_{\theta}$, with $i=1,2,\cdots,d$, each having two hidden layers, of widths at most $3\lceil\frac{k}{2}\rceil|P_{k-1,d+2}| + \lceil NT\rceil+ d(N-1)$ and $3\lceil\frac{d+3}{2}\rceil|P_{d+2,d+2}| \lceil NT\rceil N^d$, such that for every $0 \leq l \leq 2$ and $0\leq s\leq 2$, 
	\begin{align}
		\label{sec9_pinn_error1}
		&\|u_i-(u_i)_{\theta}\|_{H^l(\Omega)}\leq C_{l,k,d+1,u_i}\lambda_{l,u_i}(N)N^{-k+l},\\
		\label{sec9_pinn_error2}
		&\|v_i-(v_i)_{\theta}\|_{H^{s}(\Omega)}\leq C_{s,k-1,d+1,v_i}\lambda_{s,v_i}(N)N^{-k+1+s}.
	\end{align}
	Let $\partial_i$ represent the derivative with respect to the $i$-th dimension. For $1\leq i,\ j \leq d$, we have
	\begin{align*}
		&\|(\hat{u}_t)_i\|_{L^2(\Omega)}\leq\|\hat{u}_i\|_{H^1(\Omega)},\qquad \|(\hat{v}_t)_i\|_{L^2(\Omega)}\leq\|\hat{v}_i\|_{H^1(\Omega)},\\
        &\|\partial_i\partial_j\hat{u}_i\|_{L^2(\Omega)}, \|\partial_i\partial_i\hat{u}_i\|_{L^2(\Omega)}, \|\partial_j\partial_j\hat{u}_i\|_{L^2(\Omega)}\leq\|\hat{u}_i\|_{H^2(\Omega)},\\
		&\|\partial_j(\hat{u}_t)_i\|_{L^2(\Omega)}\leq\|(\hat{u}_t)_i\|_{H^1(\Omega)}\leq\|\hat{u}_i\|_{H^2(\Omega)},\qquad \|\partial_j\hat{v}_i\|_{L^2(\Omega)}\leq\|\hat{v}_i\|_{H^1(\Omega)},\\
		&\|\hat{u}_i\|_{L^2(D)}\leq \|\hat{u}_i\|_{L^2(\partial\Omega)}\leq C_{h_{\Omega},d+1,\rho_{\Omega}}
		\|\hat{u}_i\|_{H^1(\Omega)},\\		
		&\|\hat{v}_i\|_{L^2(D)}\leq \|\hat{v}_i\|_{L^2(\partial\Omega)}\leq C_{h_{\Omega},d+1,\rho_{\Omega}}\|\hat{v}_i\|_{H^1(\Omega)},\\
        &\|\partial_j(\hat{u})_i\|_{L^2(D)}\leq\|\partial_j(\hat{u})_i\|_{L^2(\partial\Omega)}\leq C_{h_{\Omega},d+1,\rho_{\Omega}}
		\|\hat{u}_i\|_{H^2(\Omega)},\\
		&\|\hat{v}_i\|_{L^2(\Gamma_D\times [0,t])}\leq \|\hat{v}_i\|_{L^2(\partial\Omega)}\leq C_{h_{\Omega},d+1,\rho_{\Omega}}\|\hat{v}_i\|_{H^1(\Omega)},\\	
		&\|\partial_i\hat{u}_{i}n_{i}\|_{L^2(\Gamma_N\times [0,t])},\|\partial_j\hat{u}_{i}n_{i}\|_{L^2(\Gamma_N\times [0,t])},\|\partial_j\hat{u}_{i}n_{j}\|_{L^2(\Gamma_N\times [0,t])}
        \leq C_{h_{\Omega},d+1,\rho_{\Omega}}\|\hat{u}_{i}\|_{H^2(\Omega)}.
	\end{align*}
 Using \eqref{sec9_pinn_error1} and \eqref{sec9_pinn_error2} and the above relations, we can now bound the PINN residuals,
	\begin{align*}
		&\|\bm{R}_{int1}\|_{L^2(\Omega)}\leq\|\hat{\bm{u}}_t-\hat{\bm{v}}\|_{L^2(\Omega)}\leq
		\|\hat{\bm{u}}\|_{H^1(\Omega)}+\|\hat{\bm{v}}\|_{L^2(\Omega)}\lesssim {\rm ln}NN^{-k+1},\\
		&\|\bm{R}_{int2}\|_{L^2(\Omega)}\leq\|\rho\hat{\bm{v}}_t-2\mu\nabla\cdot(\underline{\bm{\varepsilon}}(\hat{\bm{u}})) -\lambda\nabla(\nabla\cdot\hat{\bm{u}})\|_{L^2(\Omega)}\\
		&\qquad\lesssim\|\hat{\bm{v}}\|_{H^1(\Omega)}+\|\hat{\bm{u}} \|_{H^2(\Omega)}\lesssim{\rm ln}^2NN^{-k+2},\\
		&\|\underline{\bm{\varepsilon}}(\bm{R}_{int1})\|_{L^2(\Omega)},\|\nabla\cdot\bm{R}_{int1}\|_{L^2(\Omega)}\lesssim
		\|\hat{\bm{u}}\|_{H^2(\Omega)}+\|\hat{\bm{v}}\|_{H^1(\Omega)}\lesssim {\rm ln}^2NN^{-k+2},\\
		&\|\bm{R}_{tb1}\|_{L^2(D)}\leq\|\hat{\bm{u}}\|_{L^2(\partial\Omega)}\lesssim\|\hat{\bm{u}}\|_{H^1(\Omega)}\lesssim{\rm ln}NN^{-k+1},\\
		&\|\bm{R}_{tb2}\|_{L^2(D)}\leq\|\hat{\bm{v}}\|_{L^2(\partial\Omega)}\lesssim\|\hat{\bm{v}}\|_{H^1(\Omega)}\lesssim{\rm ln}NN^{-k+2},\\
		&\|\underline{\bm{\varepsilon}}(\bm{R}_{tb1})\|_{L^2(D)},\|\nabla\cdot\bm{R}_{tb1}\|_{L^2(D)}\lesssim\|\hat{\bm{u}}\|_{H^2(D)}\lesssim{\rm ln}^2NN^{-k+2},\\
		&\|\bm{R}_{sb1}\|_{L^2(\Gamma_D\times [0,t])}
		\leq\|\hat{\bm{v}}\|_{L^2(\partial\Omega)}
		\lesssim\|\hat{\bm{v}}\|_{H^1(\Omega)}\lesssim{\rm ln}NN^{-k+2},\\
		&\|\bm{R}_{sb2}\|_{L^2(\Gamma_N\times [0,t])}
		\leq\|2\mu\underline{\bm{\varepsilon}}(\hat{\bm{u}})\bm{n} +\lambda(\nabla\cdot\hat{\bm{u}})\bm{n}\|_{\partial\Omega}\lesssim\|\hat{\bm{u}}\|_{H^2(\Omega)}\lesssim{\rm ln}^2NN^{-k+2}.
	\end{align*}
\end{proof}

\vspace{0.1in}
\noindent\underline{\bf Theorem \ref{sec9_Theorem2}:}
Let $d\in \mathbb{N}$, $\bm{u} \in C^1(\Omega)$ and $\bm{v}\in C(\Omega)$ be the classical solution to the linear elastodynamic equation \eqref{elast}. Let $(\bm{u}_{\theta},\bm{v}_{\theta})$ denote the PINN approximation with the parameter $\theta$. then the following relation holds,
	\begin{equation*}
		\int_0^{T}\int_{D}( |\hat{\bm{u}}(\bm{x},t)|^2+2\mu|\underline{\bm{\varepsilon}}(\hat{\bm{u}}(\bm{x},t))|^2
		+\lambda|\nabla\cdot\hat{\bm{u}}(\bm{x},t)|^2+\rho|\hat{\bm{v}}(\bm{x},t)|^2)\dx\dt
		\leq C_GT\exp\left((2+2\mu+\lambda)T\right),
	\end{equation*}
	where $C_G$ is given in the proof.
\begin{proof} By taking the inner product of \eqref{elast_error_eq1} and \eqref{elast_error_eq2} with $\hat{\bm{u}}$ and $\hat{\bm{v}}$ and integrating over $D$, respectively, we have
	\begin{align}
		\label{sec9_eq0}
		&\frac{d}{2dt}\int_{D} |\hat{\bm{u}}|^2\dx = \int_{D}\hat{\bm{u}}\hat{\bm{v}}\dx+\int_{D} \bm{R}_{int1}\hat{\bm{u}}\dx\leq \int_{D} |\hat{\bm{u}}|^2\dx+\frac{1}{2}\int_{D} |\bm{R}_{int1}|^2\dx+\frac{1}{2}\int_{D} |\hat{\bm{v}}|^2\dx,\\
		\label{sec9_eq1}
		&\rho\frac{d}{2dt}\int_{D} |\hat{\bm{v}}|^2\dx 
		=- 2\mu\int_{D}\underline{\bm{\varepsilon}}(\hat{\bm{u}}):\nabla\hat{\bm{v}}\dx
		-\lambda\int_{D}(\nabla\cdot\hat{\bm{u}})(\nabla\cdot\hat{\bm{v}})\dx
		+\int_{\partial D}(2\mu\underline{\bm{\varepsilon}}(\hat{\bm{u}})\bm{n} +\lambda(\nabla\cdot\hat{\bm{u}})\bm{n})\cdot\hat{\bm{v}}\ds
		\nonumber\\
		&\qquad+\int_{D}\bm{R}_{int2}\hat{\bm{v}}\dx
		\nonumber\\
		&=- 2\mu\int_{D}\underline{\bm{\varepsilon}}(\hat{\bm{u}}):\nabla\hat{\bm{u}}_t\dx
		+2\mu\int_{D}\underline{\bm{\varepsilon}}(\hat{\bm{u}}):\nabla \bm{R}_{int1}\dx
		-\lambda\int_{D}(\nabla\cdot\hat{\bm{u}})(\nabla\cdot\hat{\bm{u}}_t)\dx+\lambda\int_{D}(\nabla\cdot \bm{R}_{int1})(\nabla\cdot\hat{\bm{v}})\dx
		\nonumber\\
		&\qquad
		+\int_{\Gamma_D}(2\mu\underline{\bm{\varepsilon}}(\hat{\bm{u}})\bm{n} +\lambda(\nabla\cdot\hat{\bm{u}})\bm{n})\cdot \bm{R}_{sb1}\ds
		+\int_{\Gamma_N}\bm{R}_{sb2}\cdot\hat{\bm{v}}\ds
		+\int_{D}\bm{R}_{int2}\hat{\bm{v}}\dx
		\nonumber\\
		&=-\frac{d}{dt}\int_{D}\mu|\underline{\bm{\varepsilon}}(\hat{\bm{u}})|^2\dx
		-\frac{d}{dt}\int_{D}\frac{\lambda}{2}|\nabla\cdot\hat{\bm{u}}|^2\dx+2\mu\int_{D}\underline{\bm{\varepsilon}}(\hat{\bm{u}}):\nabla \bm{R}_{int1}\dx+\lambda\int_{D}(\nabla\cdot \bm{R}_{int1})(\nabla\cdot\hat{\bm{v}})\dx
		\nonumber\\
		&\qquad
		+\int_{\Gamma_D}(2\mu\underline{\bm{\varepsilon}}(\hat{\bm{u}})\bm{n} +\lambda(\nabla\cdot\hat{\bm{u}})\bm{n})\cdot \bm{R}_{sb1}\ds
		+\int_{\Gamma_N}\bm{R}_{sb2}\cdot\hat{\bm{v}}\ds
		+\int_{D}\bm{R}_{int2}\hat{\bm{v}}\dx
		\nonumber\\
		&\leq-\frac{d}{dt}\int_{D}\mu|\underline{\bm{\varepsilon}}(\hat{\bm{u}})|^2\dx
		-\frac{d}{dt}\int_{D}\frac{\lambda}{2}|\nabla\cdot\hat{\bm{u}}|^2\dx +\mu\int_{D}|\underline{\bm{\varepsilon}}(\hat{\bm{u}})|^2\dx+\mu\int_{D} |\underline{\bm{\varepsilon}}(\bm{R}_{int1})|^2\dx
		\nonumber\\
		&\qquad
		+\frac{\lambda}{2}\int_{D}|\nabla\cdot \bm{R}_{int1}|\dx
		+\frac{\lambda}{2}\int_{D}|\nabla\cdot \hat{\bm{v}}|\dx
		+\frac{1}{2}\int_{D} |\hat{\bm{v}}|^2\dx+\frac{1}{2}\int_{D} |\bm{R}_{int2}|^2\dx
		\nonumber\\
		&\qquad +C_{\Gamma_D}\left(\int_{\Gamma_D}|\bm{R}_{sb1}|^2\ds\right)^{\frac{1}{2}}+C_{\Gamma_N}\left(\int_{\Gamma_N}|\bm{R}_{sb2}|^2\ds\right)^{\frac{1}{2}}.
	\end{align}
	Here we have used $\hat{\bm{v}}=\hat{\bm{u}}_t-\bm{R}_{int1}$, and the constants are given by $C_{\Gamma_D}=(2\mu+\lambda)|\Gamma_D|^{\frac{1}{2}}\|\bm{u}\|_{C^1(\Gamma_D\times [0,T])}+(2\mu+\lambda)|\Gamma_D|^{\frac{1}{2}}||\bm{u}_{\theta}||_{C^1(\Gamma_D\times [0,T])}$ and $C_{\Gamma_N}=|\Gamma_N|^{\frac{1}{2}}(\|\bm{v}\|_{C(\Gamma_N\times [0,T])}+||\bm{v}_{\theta}||_{C(\Gamma_N\times [0,T])})$. \\

 Add \eqref{sec9_eq0} to \eqref{sec9_eq1}, and we get, 
	\begin{align}\label{sec9_eq2}
		&\frac{d}{2dt}\int_{D} |\hat{\bm{u}}|^2\dx+\frac{d}{dt}\int_{D}\mu|\underline{\bm{\varepsilon}}(\hat{\bm{u}})|^2\dx
		+\frac{d}{2dt}\int_{D}\lambda|\nabla\cdot\hat{\bm{u}}|^2\dx+\rho\frac{d}{2dt}\int_{D} |\hat{\bm{v}}|^2\dx
		\nonumber\\
		&\qquad \leq\int_{D}|\hat{\bm{u}}|^2\dx+\mu\int_{D}|\underline{\bm{\varepsilon}}(\hat{\bm{u}})|^2\dx
		+\frac{\lambda}{2}\int_{D}|\nabla\cdot \hat{\bm{v}}|\dx+\int_{D}|\hat{\bm{v}}|^2\dx+\frac{1}{2}\int_{D}(|\bm{R}_{int1}|^2+|\bm{R}_{int2}|^2)\dx
		\nonumber\\
		&\qquad\qquad+\mu\int_{D} |\underline{\bm{\varepsilon}}(\bm{R}_{int1})|^2\dx+\frac{\lambda}{2}\int_{D}|\nabla\cdot \bm{R}_{int1}|\dx+C_{\Gamma_D}\left(\int_{\Gamma_D}|\bm{R}_{sb1}|^2\ds\right)^{\frac{1}{2}}+C_{\Gamma_N}\left(\int_{\Gamma_N}|\bm{R}_{sb2}|^2\ds\right)^{\frac{1}{2}}.
	\end{align}
	
 Integrating \eqref{sec9_eq2} over $[0,\tau]$ for any $\tau \leq T$ and applying Cauchy–Schwarz inequality, we obtain,
	\begin{align*}
		&\int_{D} |\hat{\bm{u}}(\bm{x},\tau)|^2\dx 
		+\int_{D}2\mu|\underline{\bm{\varepsilon}}(\hat{\bm{u}}(\bm{x},\tau))|^2\dx
		+\int_{D}\lambda|\nabla\cdot\hat{\bm{u}}(\bm{x},\tau)|^2\dx
		+\rho\int_{D} |\hat{\bm{v}}(\bm{x},\tau)|^2\dx\\
		&\qquad\leq\int_{D}|\bm{R}_{tb1}|^2\dx +\int_{D}2\mu|\underline{\bm{\varepsilon}}(\bm{R}_{tb1})|^2\dx
		+\int_{D}\lambda|\nabla\cdot \bm{R}_{tb1}|^2\dx+\rho\int_{D}|\bm{R}_{tb2}|^2\dx\\
		&\qquad\qquad
		+(2+2\mu+\lambda)\int_{0}^{\tau}\int_{D}\left(|\hat{\bm{u}}|^2+|\underline{\bm{\varepsilon}}(\hat{\bm{u}})|^2+|\nabla\cdot\hat{\bm{u}}|^2+ |\hat{\bm{v}}|^2 \right)\dx\dt
		\\
		&\qquad\qquad+ \int_{0}^{T}\int_{D}\left(|\bm{R}_{int1}|^2+2\mu|\underline{\bm{\varepsilon}}(\bm{R}_{int1})|^2+\lambda|\nabla\cdot \bm{R}_{int1}|^2+|\bm{R}_{int2}|^2\right)\dx\dt
		\\
		&\qquad\qquad+2|T|^{\frac{1}{2}}C_{\Gamma_D}\left(\int_{0}^{T}\int_{\Gamma_D}|\bm{R}_{sb1}|^2\ds\dt\right)^{\frac{1}{2}}+2|T|^{\frac{1}{2}}C_{\Gamma_N}\left(\int_{0}^{T}\int_{\Gamma_N}|\bm{R}_{sb2}|^2\ds\dt\right)^{\frac{1}{2}}.
	\end{align*}
 By applying the integral form of the Gr${\rm\ddot{o}}$nwall inequality to the above inequality, we have
	\begin{equation}\label{sec9_eq3}
		\int_{D} (|\hat{\bm{u}}(\bm{x},\tau)|^2+2\mu|\underline{\bm{\varepsilon}}(\hat{\bm{u}}(\bm{x},\tau))|^2+\lambda|\nabla\cdot\hat{\bm{u}}(\bm{x},\tau)|^2
		+\rho\int_{D} |\hat{\bm{v}}(\bm{x},\tau)|^2)\dx
		\leq C_G\exp\left((2+2\mu+\lambda)T\right),
	\end{equation}
	where
	\begin{align*}
		&C_G=\int_{D}|\bm{R}_{tb1}|^2\dx+\int_{D}2\mu|\underline{\bm{\varepsilon}}(\bm{R}_{tb1})|^2\dx
		+\int_{D}\lambda|\nabla\cdot \bm{R}_{tb1}|^2\dx+\rho\int_{D}|\bm{R}_{tb2}|^2\dx
		\\
		&\qquad
		+\int_{0}^{T}\int_{D}\left(|\bm{R}_{int1}|^2+2\mu|\underline{\bm{\varepsilon}}(\bm{R}_{int1})|^2+\lambda|\nabla\cdot \bm{R}_{int1}|^2+|\bm{R}_{int2}|^2\right)\dx\dt
		\\
		&\qquad +2|T|^{\frac{1}{2}}C_{\Gamma_D}\left(\int_{0}^{T}\int_{\Gamma_D}|\bm{R}_{sb1}|^2\ds\dt\right)^{\frac{1}{2}}+2|T|^{\frac{1}{2}}C_{\Gamma_N}\left(\int_{0}^{T}\int_{\Gamma_N}|\bm{R}_{sb2}|^2\ds\dt\right)^{\frac{1}{2}}.
	\end{align*}
	Then, we finish the proof by integrating \eqref{sec9_eq3} over $[0,T]$.
\end{proof} 

\vspace{0.1in}
\noindent\underline{\bf Theorem \ref{sec9_Theorem3}:} Let $d\in \mathbb{N}$, $\bm{u}\in C^4(\Omega)$ and $\bm{v}\in C^3(\Omega)$ be the classical solution to
	the linear elastodynamic equation \eqref{elast}. Let $(\bm{u}_{\theta},\bm{v}_{\theta})$ denote the PINN approximation with the parameter $\theta$. Then the following relation holds,
	\begin{align*}
	&\int_0^{T}\int_{D}( |\hat{\bm{u}}(\bm{x},t)|^2+2\mu|\underline{\bm{\varepsilon}}(\hat{\bm{u}}(\bm{x},t))|^2
	+\lambda|\nabla\cdot\hat{\bm{u}}(\bm{x},t)|^2+\rho|\hat{\bm{v}}(\bm{x},t)|^2)\dx\dt
	\leq C_TT\exp\left((2+2\mu+\lambda)T\right)
		\nonumber\\
	&\qquad=\mathcal{O}(\mathcal{E}_T(\theta)^2 + M_{int}^{-\frac{2}{d+1}} +M_{tb}^{-\frac{2}{d}}+M_{sb}^{-\frac{1}{d}}),       
    \end{align*}
	where $C_T$ is defined in the following proof.
\begin{proof} By the definitions of different components of the training error \eqref{elast_TT} and applying the estimate \eqref{int1} on the quadrature error, we have
	\begin{align*}
		\int_{D}|\bm{R}_{tb1}|^2\dx&=\int_{D}|\bm{R}_{tb1}|^2\dx-\mathcal{Q}_{M_{tb}}^{D}(\bm{R}_{tb1}^2)+\mathcal{Q}_{M_{tb}}^{D}(\bm{R}_{tb1}^2)\\
		&\leq C_{({\bm{R}_{tb1}^2})}M_{tb}^{-\frac{2}{d}}+\mathcal{Q}_{M_{tb}}^{D}(\bm{R}_{tb1}^2),\\
		\int_{D}|\bm{R}_{tb2}|^2\dx&=\int_{D}|\bm{R}_{tb2}|^2\dx-\mathcal{Q}_{M_{tb}}^{D}(\bm{R}_{tb2}^2)+\mathcal{Q}_{M_{tb}}^{D}(\bm{R}_{tb2}^2)\\
		&\leq C_{({\bm{R}_{tb2}^2})}M_{tb}^{-\frac{2}{d}}+\mathcal{Q}_{M_{tb}}^{D}(\bm{R}_{tb2}^2),\\
		\int_{D}|\underline{\bm{\varepsilon}}(\bm{R}_{tb1})|^2\dx
		&=\int_{D}|\underline{\bm{\varepsilon}}(\bm{R}_{tb1})|^2\dx-\mathcal{Q}_{M_{tb}}^{D}(|\underline{\bm{\varepsilon}}(\bm{R}_{tb1})|^2)+\mathcal{Q}_{M_{tb}}^{D}(|\underline{\bm{\varepsilon}}(\bm{R}_{tb1})|^2)\\
		&\leq C_{(|\underline{\bm{\varepsilon}}(\bm{R}_{tb1})|^2)}M_{tb}^{-\frac{2}{d}}+\mathcal{Q}_{M_{tb}}^{D}(|\underline{\bm{\varepsilon}}(\bm{R}_{tb1})|^2),\\
		\int_{D}|\nabla\cdot \bm{R}_{tb1}|^2\dx
		&=\int_{D}|\nabla\cdot \bm{R}_{tb1}|^2\dx-\mathcal{Q}_{M_{tb}}^{D}(|\nabla\cdot \bm{R}_{tb1}|^2)+\mathcal{Q}_{M_{tb}}^{D}(|\nabla\cdot \bm{R}_{tb1}|^2)\\
		&\leq C_{(|\nabla\cdot \bm{R}_{tb1}|^2)}M_{tb}^{-\frac{2}{d}}+\mathcal{Q}_{M_{tb}}^{D}(|\nabla\cdot \bm{R}_{tb1}|^2),\\
		\int_{\Omega}|\bm{R}_{int1}|^2\dx\dt&=\int_{\Omega}|\bm{R}_{int1}|^2\dx\dt-\mathcal{Q}_{M_{int}}^{\Omega}(\bm{R}_{int1}^2)+\mathcal{Q}_{M_{int}}^{\Omega}(\bm{R}_{int1}^2)\\
		&\leq C_{({\bm{R}_{int1}^2})}M_{int}^{-\frac{2}{d+1}}+\mathcal{Q}_{M_{int}}^{\Omega}(\bm{R}_{int1}^2),\\
		\int_{\Omega}|\bm{R}_{int2}|^2\dx\dt&=\int_{\Omega}|\bm{R}_{int2}|^2\dx\dt-\mathcal{Q}_{M_{int}}^{\Omega}(\bm{R}_{int2}^2)+\mathcal{Q}_{M_{int}}^{\Omega}(\bm{R}_{int2}^2)\\
		&\leq C_{({\bm{R}_{int2}^2})}M_{int}^{-\frac{2}{d+1}}+\mathcal{Q}_{M_{int}}^{\Omega}(\bm{R}_{int2}^2),\\
		\int_{\Omega}|\underline{\bm{\varepsilon}}(\bm{R}_{int1})|^2\dx\dt&=\int_{\Omega}|\underline{\bm{\varepsilon}}(\bm{R}_{int1})|^2\dx\dt-\mathcal{Q}_{M_{int}}^{\Omega}(|\underline{\bm{\varepsilon}}(\bm{R}_{int1})|^2)+\mathcal{Q}_{M_{int}}^{\Omega}(|\underline{\bm{\varepsilon}}(\bm{R}_{int1})|^2)\\
		&\leq C_{(|\underline{\bm{\varepsilon}}(\bm{R}_{int1})|^2)}M_{int}^{-\frac{2}{d+1}}+\mathcal{Q}_{M_{int}}^{\Omega}(|\underline{\bm{\varepsilon}}(\bm{R}_{int1})|^2),\\
		\int_{\Omega}|\nabla\cdot\bm{R}_{int1}|^2\dx\dt&=\int_{\Omega}|\nabla\cdot\bm{R}_{int1}|^2\dx\dt-\mathcal{Q}_{M_{int}}^{\Omega}(|\nabla\cdot\bm{R}_{int1}|^2)+\mathcal{Q}_{M_{int}}^{\Omega}(|\nabla\cdot\bm{R}_{int1}|^2)\\
		&\leq C_{(|\nabla\cdot\bm{R}_{int1}|^2)}M_{int}^{-\frac{2}{d+1}}+\mathcal{Q}_{M_{int}}^{\Omega}(|\nabla\cdot\bm{R}_{int1}|^2),\\
		\int_{\Omega_D}|\bm{R}_{sb1}|^2\ds\dt&=\int_{\Omega_D}|\bm{R}_{sb1}|^2\ds\dt-\mathcal{Q}_{M_{sb1}}^{\Omega_D}(\bm{R}_{sb1}^2)+\mathcal{Q}_{M_{sb1}}^{\Omega_D}(\bm{R}_{sb1}^2)\\
		&\leq C_{({\bm{R}_{sb1}^2})}M_{sb1}^{-\frac{2}{d}}+\mathcal{Q}_{M_{sb1}}^{\Omega_D}(\bm{R}_{sb1}^2),\\
		\int_{\Omega_N}|\bm{R}_{sb2}|^2\ds\dt&=\int_{\Omega_N}|\bm{R}_{sb2}|^2\ds\dt-\mathcal{Q}_{M_{sb2}}^{\Omega_N}(\bm{R}_{sb2}^2)+\mathcal{Q}_{M_{sb2}}^{\Omega_N}(\bm{R}_{sb2}^2)\\
		&\leq C_{({\bm{R}_{sb2}^2})}M_{sb2}^{-\frac{2}{d}}+\mathcal{Q}_{M_{sb2}}^{\Omega_N}(\bm{R}_{sb2}^2).
	\end{align*}
	
 In light of the above inequalities and \eqref{sec9_eq3}, we obtain
	\begin{equation*}
		\int_0^{T}\int_{D}( |\hat{\bm{u}}(\bm{x},t)|^2+2\mu|\underline{\bm{\varepsilon}}(\hat{\bm{u}}(\bm{x},t))|^2
		+\lambda|\nabla\cdot\hat{\bm{u}}(\bm{x},t)|^2+\rho|\hat{\bm{v}}(\bm{x},t)|^2)\dx\dt
		\leq TC_T\exp\left((2+2\mu+\lambda)T\right),
	\end{equation*}
	where 
	\begin{align*}
		C_T=&C_{({\bm{R}_{tb1}^2})}M_{tb}^{-\frac{2}{d}}+\mathcal{Q}_{M_{tb}}^{D}(\bm{R}_{tb1}^2)
		+\rho\left(C_{({\bm{R}_{tb2}^2})}M_{tb}^{-\frac{2}{d}}+\mathcal{Q}_{M_{tb}}^{D}(\bm{R}_{tb2}^2)\right)
		+ 2\mu\left(C_{(|\underline{\bm{\varepsilon}}(\bm{R}_{tb1})|^2)}M_{tb}^{-\frac{2}{d}}+\mathcal{Q}_{M_{tb}}^{D}(|\underline{\bm{\varepsilon}}(\bm{R}_{tb1})|^2)\right)\\
		&+\lambda\left(C_{(|\nabla\cdot \bm{R}_{tb1}|^2)}M_{tb}^{-\frac{2}{d}}+\mathcal{Q}_{M_{tb}}^{D}(|\nabla\cdot \bm{R}_{tb1}|^2)\right)+C_{({\bm{R}_{int1}^2})}M_{int}^{-\frac{2}{d+1}}+\mathcal{Q}_{M_{int}}^{\Omega}(\bm{R}_{int1}^2)\\
		&+ C_{({\bm{R}_{int2}^2})}M_{int}^{-\frac{2}{d+1}}+\mathcal{Q}_{M_{int}}^{\Omega}(\bm{R}_{int2}^2)
		+2\mu\left(C_{(|\underline{\bm{\varepsilon}}(\bm{R}_{int1})|^2)}M_{int}^{-\frac{2}{d+1}}+\mathcal{Q}_{M_{int}}^{\Omega}(|\underline{\bm{\varepsilon}}(\bm{R}_{int1})|^2)\right)\\
		&+\lambda\left(C_{(|\nabla\cdot\bm{R}_{int1}|^2)}M_{int}^{-\frac{2}{d+1}}+\mathcal{Q}_{M_{int}}^{\Omega}(|\nabla\cdot\bm{R}_{int1}|^2)\right)
		+ 2|T|^{\frac{1}{2}}C_{\Gamma_D}\left(C_{({\bm{R}_{sb1}^2})}M_{sb1}^{-\frac{2}{d}}+\mathcal{Q}_{M_{sb1}}^{\Omega_D}(\bm{R}_{sb1}^2)\right)^{\frac{1}{2}}\\
		&+2|T|^{\frac{1}{2}}C_{\Gamma_N}\left(C_{({\bm{R}_{sb2}^2})}M_{sb2}^{-\frac{2}{d}}+\mathcal{Q}_{M_{sb2}}^{\Omega_N}(\bm{R}_{sb2}^2)\right)^{\frac{1}{2}}.
	\end{align*}
	 The boundedness of the constants $C(\bm{R}_q^2)$ can be obtained from Lemma \ref{Ar_3} and $\|\bm{R}_q^2\|_{C^n}\leq 2^n\|\bm{R}_q\|_{C^n}^2$, with $\bm{R}_q=\bm{R}_{tb1}$, $\bm{R}_{tb2}$, $\underline{\bm{\varepsilon}}(\bm{R}_{tb1})$, $\nabla\cdot \bm{R}_{tb1}$, $\bm{R}_{int1}$, $\bm{R}_{int2}$, $\underline{\bm{\varepsilon}}(\bm{R}_{int1})$, $\nabla\cdot\bm{R}_{int1}$, $\bm{R}_{sb1}$ and $\bm{R}_{sb2}$. 
\end{proof}  

\bibliographystyle{plain}
\bibliography{the_ref,dnn,mypub}

\end{document}